\pdfoutput=1

\RequirePackage{fix-cm}

\documentclass[leqno, fleqn, centertags, 12pt]{article}

\usepackage{latexsym}
\usepackage{amsmath}
\usepackage{amssymb}
\usepackage{verbatim}

\usepackage[T1]{fontenc}

\usepackage[upright]{fourier}

\usepackage{fullpage}

\newskip\nineskipamount \nineskipamount=9pt plus 0pt minus 0pt
\newskip\zeroskipamount \zeroskipamount=0pt plus 0pt minus 0pt
\usepackage[noorphans,vskip=\nineskipamount]{quoting}

\makeatletter
\renewcommand{\@makefntext}[1]{\vspace*{0.5ex}\parindent=0em
\hspace*{-0.4em}
\hbox to 0.4em{\hss\@makefnmark}\hspace*{0.4em}{#1}
}
\makeatother

\newcounter{mysectionnumber}
\newcommand{\mysection}[2]{\setcounter{footnote}{0}
\setcounter{myparnum}{0}
\refstepcounter{mysectionnumber}
\vspace{21pt}{\Large {\themysectionnumber.} {#1}}\label{#2}\vspace*{15pt}}

\newcommand{\mynonumbersection}[1]{\setcounter{footnote}{0}
\setcounter{myparnum}{0}
\vspace{21pt}{\Large \dnsp$\phantom{1.}$\dnsp\ {#1}}\vspace*{15pt}}

\newcommand{\myit}[1]{\textbf{\textit{#1}}\hspace{0.0em}}

\newcounter{myparnum}[mysectionnumber]
\renewcommand{\themyparnum}{\themysectionnumber.\arabic{myparnum}}
\newcommand{\mypar}[2]{\refstepcounter{myparnum}{\vspace{\medskipamount}\textbf{{\themyparnum. #1}\label{#2}}\hspace{0.5em}}}

\newcounter{mylemmanum}[myparnum]
\newcommand{\myuppar}[1]{\vspace{\medskipamount}\textbf{#1}\hspace*{0.5em}}

\newcounter{myappendnumber}
\newcounter{myaparnum}[myappendnumber]
\newcommand{\myappend}[2]{\setcounter{footnote}{0}
\setcounter{myaparnum}{0}
\setcounter{myparnum}{0}
\refstepcounter{myappendnumber}
\vspace{21pt}{\Large A\dff.{\themyappendnumber.}\oss {#1}}\label{#2}\vspace*{15pt}}

\newcommand{\myapar}[2]{\refstepcounter{myaparnum}{\vspace{\medskipamount}\textbf{{\themyaparnum. #1}\label{#2}}\hspace{0.5em}}}

\renewcommand{\themyaparnum}{A\halfff\fff.\fff\themyappendnumber.\arabic{myaparnum}}

\newcommand{\proof}{\vspace{\medskipamount}{\textbf{{\emph{Proof}.}}\hspace*{1em}}}
\newcommand{\prooftitle}[1]{\vspace{\medskipamount}{\textbf{{\emph{#1}.}}\hspace*{1em}}}

\newcommand{\eproof}{ $\blacksquare$}

\newcommand{\dis}{\displaystyle}

\def\sss{\hspace{0.05em}\ }
\def\dss{\hspace{0.1em}\ }
\def\trs{\hspace{0.15em}\ }
\def\qss{\hspace{0.2em}\ }
\def\pss{\hspace{0.3em}\ }
\def\oss{\hspace{0.4em}\ }

\def\halfff{\hspace*{0.025em}}
\def\fff{\hspace*{0.05em}}
\def\dff{\hspace*{0.1em}}
\def\trf{\hspace*{0.15em}}
\def\qff{\hspace*{0.2em}}
\def\pff{\hspace*{0.3em}}
\def\off{\hspace*{0.4em}}

\newcommand{\hnsp}{\hspace*{-0.05em}}
\newcommand{\nsp}{\hspace*{-0.1em}}
\newcommand{\nnsp}{\hspace*{-0.15em}}
\newcommand{\snsp}{\hspace*{-0.175em}}
\newcommand{\dnsp}{\hspace*{-0.2em}}

\renewcommand{\leq}{\leqslant}
\renewcommand{\geq}{\geqslant}

\newcommand{\zzz}{\mathbb{Z}}

\newcommand{\rrr}{\mathbb{R}}

\newcommand{\ftwo}{\mathbb{F}_{\dff 2}}

\newcommand{\num}[1]{|\qff #1 \qff|}

\newcommand{\fclass}[1]{[\snsp [\dff #1\dff]\snsp ]}
\newcommand{\norm}[1]{\|\qff #1 \qff\|}

\newcommand{\ttoo}{\hspace*{0.2em}\longrightarrow\hspace*{0.2em}}

\DeclareMathOperator{\pr}{pr}
\newcommand{\bd}{\operatorname{bd}}

\begin{document}

\setlength{\baselineskip}{12pt plus 0pt minus 0pt}
\setlength{\parskip}{12pt plus 0pt minus 0pt}
\setlength{\abovedisplayskip}{12pt plus 0pt minus 0pt}
\setlength{\belowdisplayskip}{12pt plus 0pt minus 0pt}

\newskip\smallskipamount \smallskipamount=3pt plus 0pt minus 0pt
\newskip\medskipamount   \medskipamount  =6pt plus 0pt minus 0pt
\newskip\bigskipamount   \bigskipamount =12pt plus 0pt minus 0pt

\author{Nikolai\qss V.\qss Ivanov}
\title{Cubes\qss and\qss cubical\qss chains\qss and\qss cochains\\ in\qss combinatorial\pss topology}
\date{}

\footnotetext{\hspace*{-0.65em}\copyright\oss 
Nikolai\qss V.\qss Ivanov,\oss 2019.\trs 
Neither\dss the work reported\dss in\dss the present paper\halfff,\qss
nor\dss its preparation were supported\dss by\dss any corporate entity.}

\maketitle

\vspace*{16ex}

{\renewcommand{\baselinestretch}{1}
\selectfont

\myit{\hspace*{0em}\large Contents}\vspace*{1.5ex} \\ 
\hbox to 0.8\textwidth{\myit{Preface} \hfil 2}\hspace*{0.5em} \vspace*{1.5ex}\\
\hbox to 0.8\textwidth{\myit{\phantom{1}1.}\hspace*{0.5em} Du\sss c\^{o}t\'{e}\sss de\sss chez\qss 
Lebesgue \hfil 5}\hspace*{0.5em} \vspace*{0.25ex}\\
\hbox to 0.8\textwidth{\myit{\phantom{1}2.}\hspace*{0.5em} Lebesgue\trs tilings \hfil 18}\hspace*{0.5em} \vspace*{0.25ex}\\
\hbox to 0.8\textwidth{\myit{\phantom{1}3.}\hspace*{0.5em} Modification of\qss Lebesgue\dss methods\sss by\qss W.\dss Hurewicz  \hfil 26}\hspace*{0.5em}\vspace*{0.25ex}\\
\hbox to 0.8\textwidth{\myit{\phantom{1}4.}\hspace*{0.5em} Nerves of\dss coverings\sss and\trs
Hurewicz's\dss theorems \hfil 31}\hspace*{0.5em} \vspace*{0.25ex}\\
\hbox to 0.8\textwidth{\myit{\phantom{1}5.}\hspace*{0.5em} Lusternik--Schnirelmann\qss theorems \hfil 35}\hspace*{0.5em} \vspace*{0.25ex}\\
\hbox to 0.8\textwidth{\myit{\phantom{1}6.}\hspace*{0.5em} First\sss applications of\qss Lusternik--Schnirelmann\qss theorems \hfil 41}\hspace*{0.5em} \vspace*{0.25ex}\\
\hbox to 0.8\textwidth{\myit{\phantom{1}7.}\hspace*{0.5em} The discrete cube 
and\dss products of\dss cubical\sss cochains \hfil 44}\hspace*{0.5em} \vspace*{0.25ex}\\
\hbox to 0.8\textwidth{\myit{\phantom{1}8.}\hspace*{0.5em} The discrete sphere 
and\dss products of\dss cubical\sss cochains  \hfil 57}\hspace*{0.5em} \vspace*{0.25ex}\\
\hbox to 0.8\textwidth{\myit{\phantom{1}9.}\hspace*{0.5em} Kuhn's\qss cubical\qss Sperner\trs lemmas \hfil 64}\hspace*{0.5em} \vspace*{0.25ex}\\
\hbox to 0.8\textwidth{\myit{10.}\hspace*{0.5em} Ky\qss Fan's\qss cubical\qss Sperner\qss lemma \hfil 69}\hspace*{0.5em} \vspace*{1.5ex}\\
\myit{Appendices}\hspace*{0.5em}  \hspace*{0.5em} \vspace*{1.5ex}\\
\hbox to 0.8\textwidth{\myit{\phantom{1}A.1.}\hspace*{0.5em} Freudenthal's\dss triangulations
of\dss cubes and simplices\hfil 75}\hspace*{0.5em} \vspace*{0.25ex}\\
\hbox to 0.8\textwidth{\myit{\phantom{1}A.2.}\hspace*{0.5em} The solid cube and\dss the discrete cube \hfil 79}\hspace*{0.5em}  \vspace*{1.5ex}\\
\hbox to 0.8\textwidth{\myit{References}\hspace*{0.5em}\hfil 84}\hspace*{0.5em}  \vspace*{0.25ex}  

}

\renewcommand{\baselinestretch}{1}
\selectfont

\newpage
\mynonumbersection{Preface}

\myuppar{Simplices\dss and\dss cubes.}
The present\dss paper\sss can\dss be considered as a continuation of\dss
author's paper\qss \cite{i2}\qss devoted\dss to\sss the lemmas of\qss
Alexander\dss and\dss Sperner\dss and\dss their applications,\oss
but\trs is\dss independent\dss from\qss \cite{i2}.\oss
As\sss a\dss motivation,\pss 
it\trs is\dss worth\dss to\sss quote\sss 
the first\dss few phrases of\qss \cite{s}.\oss\vspace{-8.25pt}

\begin{quoting}
Lebesgue\dss has shown\dss that\dss the invariance of\trs the dimension 
and also\sss the invariance of\dss domains can be reduced\dss in 
a simple way to a fundamental theorem of\trs topology\qss (see Section 1).\oss 
This theorem\dss was also already proved\dss by\qss Brouwer\trs and\dss Lebesgue.\oss 
But\dss both\dss these proofs present\sss significant\sss difficulties\sss 
to\sss the readers\sss who are not\dss familiar\dss with\dss this circle of\dss questions.\oss 
In\dss the present\dss work,\oss first of\dss all,\oss a new,\oss 
simple proof\dss of\trs this important\dss theorem\dss 
is\dss given\dss in\dss Sections\qss 1 -- 3.\oss 
It\dss turns out\dss to be appropriate\sss to place\sss the 
$n$\dnsp-dimensional\sss simplex at\dss the center of\dss consideration\dss 
instead of\dss the $n$\dnsp-dimensional\sss cube.\oss 
\end{quoting}

\vspace{-8.25pt}
In\dss the present\dss paper\dss $n$\dnsp-dimensional\sss cubes
are\dss put\dss back at\dss the center of\dss consideration.\oss
It\dss
begins\sss by\sss a\sss step back\dss from\dss
Alexander\dss and\trs Sperner\trs to\dss
Lebesgue\dss papers\qss \cite{l-one},\pss \cite{l-two},\oss
one of\trs the main sources of\trs
inspiration\dss for\dss Sperner\halfff.\oss
One could\dss say\dss that\dss the present\dss paper\dss is\dss
devoted\dss to cubical\sss analogues of\trs lemmas of\qss Alexander\trs
and\qss Sperner\halfff,\oss
if\qss not\dss the fact\dss that\dss the first\sss such analogues
were proved\dss by\qss Lebesgue\dss several\dss years\sss earlier\dss
than\dss these\sss lemmas.\oss

\myuppar{A\dss fragment\sss of\trs history\halfff.}
The\qss ``fundamental\dss theorem of\trs topology'',\pss
referred\dss to\sss by\trs Sperner\halfff,\oss
is\dss the\qss Lebesgue\trs covering\dss theorem\qss
or\halfff,\oss rather\halfff,\oss two\sss covering\trs theorems.\oss
See\dss Section\qss \ref{lebesgue}.\oss
Lebesgue\dss results were originally\dss published\qss \cite{l-one}\qss as
an excerpt\dss from a\sss letter\dss to\dss O.\dss Blumenthal,\pss one of\dss the editors of\qss
\emph{Mathematische\dss Annalen},\pss complementing\qss
Brouwer's\trs paper\qss \cite{b1}\qss
about\dss the\sss topological\dss invariance of\trs the dimension.\oss 
In\dss his\sss letter\trs Lebesgue\dss suggested an approach\dss to\sss
the\sss topological\dss invariance based\sss on a\sss completely\dss new\sss idea.\pss
Namely\halfff,\pss 
Lebesgue\dss suggested\dss to prove\sss the\sss topological\dss
invariance of\trs the dimension\dss by\dss using coverings of\dss domains\sss 
by\sss small\sss sets\sss and\dss investigating\dss the pattern of\dss 
intersections of\dss these small\sss sets.\oss

This\dss idea\dss had\dss very\dss far-reaching\dss influence.\oss
One can\dss trace\sss to\sss it\dss not\sss only\dss the notion of\dss the covering
dimension,\pss but\sss also\sss the notion of\trs the nerve of\dss a covering and,\pss
eventually\halfff,\pss the notion of\qss a\dss sheaf\halfff.\oss
In contrast\halfff,\pss the influence of\qss Lebesgue\dss proofs\sss was\sss rather\dss limited.\oss
It\dss looks\sss like almost\dss nobody\dss read\dss his detailed\dss exposition\qss
\cite{l-two},\pss in contrast\dss with\dss the outline\qss \cite{l-one}.\oss
The only\sss exceptions known\dss to\sss the author\sss are\dss W.\dss Hurewicz\qss \cite{h2},\oss
who modified\qss Lebesgue\dss proofs,\oss
and\qss L.\dss Lusternik\dss and\dss L.\dss Schnirelmann\qss \cite{ls},\oss
who applied\qss Lebesgue\dss methods\sss to a different\dss problem.\oss
Lusternik\dss and\dss Schnirelmann\dss did\dss not\dss refer\dss to\dss Lebesgue,\oss
but\dss the\dss influence of\qss \cite{l-two}\qss is\dss quite obvious;\oss
they\sss even use some of\qss Lebesgue notations.\oss

Lebesgue\dss main\dss technical\dss tools were\sss
the standard\dss
partitions of\dss a given cube in\dss $\rrr^{\fff n}$\dss into small\sss cubes
of\dss equal\sss sizes
and\dss an approximation of\dss closed sets by\dss unions of\trs these small\sss cubes.\oss
The use of\dss cubes by\qss Lebesgue\dss strongly\sss contrasted\dss with\dss
the\sss domination of\dss simplicial\dss methods in\dss topology\halfff.\oss
But\sss somewhat\dss later\dss it\dss became a desirable feature at\dss least\dss
in\dss some quarters.\oss
This\sss lead\dss to searching\dss for\dss cubical\sss analogues of\pss
Sperner's\dss results\sss and\sss methods.\oss
See,\pss for example,\oss the works of\pss A.W.\qss Tucker\qss \cite{t1},\pss \cite{t2},\oss
H.W.\qss Kuhn\qss \cite{ku},\pss and\qss Ky\dss Fan\qss \cite{kf}.\oss
Apparently\halfff,\oss none of\trs these mathematicians\dss 
realized\dss that\sss such analogues existed\dss before\dss 
Sperner's\dss work.\oss\vspace{-0.125pt}

\myuppar{Lebesgue,\oss Brouwer\halfff,\oss and\qss Sperner\halfff.}
Brouwer\dss deemed\qss Lebesgue\dss outline\qss \cite{l-one}\qss faulty\dss
and\dss two years\sss later\dss published\qss \cite{b2}\qss
a\dss fairly\sss cumbersome\sss simplicial\dss version of\qss Lebesgue\sss theorem.\oss
See\qss \cite{d},\oss Section\qss II.5\qss for a modern\dss rendering\dss of\qss
Brouwer's\dss version.\oss
Brouwer\dss also\sss included\dss in\qss \cite{b2}\qss
an almost\dss page-long\dss footnote devoted\dss to a counterexample
to\dss Lebesgue\dss proofs\qss
({\fff}but\dss not\dss to\sss the\sss theorems).\oss
The validity of\trs this counterexample depends on\dss interpreting\qss
French\qss \emph{un ensemble}\qss 
(a set\halfff)\trs 
as\dss German\qss
\emph{ein\dss Teilkontinuum}\pss (a connected closed subset\fff).\oss
But\qss Lebesgue\dss accepted\qss Brouwer's\trs critique\qss
(see\qss \cite{l-two},\oss pp.\qss 264--265)\qss
and\dss was much\dss more careful\dss in\qss \cite{l-two}.

Brouwer's\dss proof\dss of\qss his simplicial\dss version\dss is\dss
compressed\dss into 13\trs lines.\oss
J.\dss Dieudonne\dss writes\sss that\oss
``the proof\dss uses\sss the properties of\trs the degree of\dss a map,\pss
and,\pss as usual,\pss 
is\dss very sketchy and\dss has\sss to be\sss interpreted\dss to make sense''\qss
(see\qss \cite{d},\pss Section\qss II.5).\oss
In contrast\halfff,\oss Lebesgue\trs paper\qss \cite{l-two}\qss
is\dss self-contained and\sss detailed.\oss
If\trs it\dss appears\sss to be demanding,\pss then not\dss
because\sss of\trs prerequisites,\pss
but\dss precisely\dss because\sss of\trs its\dss elementary\dss language.\oss 
At\dss the same\sss time\dss Lebesgue\dss proofs\sss are\dss direct\sss and\dss natural.\oss
Sperner\qss \cite{s}\qss replaced\dss his\sss 
direct\dss methods 
by\dss a special\sss case of\pss
Brouwer's\qss degree\sss theory\sss admitting\dss an\dss independent\dss proof\halfff.\oss
See\qss \cite{i2}\qss for\dss the details.\oss\vspace{-0.125pt}

\myuppar{Lebesgue\dss methods and\dss their\sss extensions and\sss applications.}
The arguments of\qss
Lebesgue\trs and\trs Lusternik--Schnirelmann\trs can\dss be substantially\sss
clarified\dss by\sss using\dss the language of\dss cubical\sss chains.\oss
This language\dss is\dss introduced\dss in\dss Section\qss \ref{lebesgue},\pss
devoted\dss to an exposition of\qss Lebesgue\dss results,\pss
and\dss is\dss used\dss throughout\dss the paper\halfff.\oss
Following\dss the norms of\trs the modern\dss mathematical\dss exposition,\oss
the continuous narrative of\qss Lebesgue\dss is\dss split\dss into several\dss
lemmas and\dss theorems.\oss
A similar\dss treatment\sss of\qss
the\qss Lusternik--Schnirelmann covering\dss theorems\qss \cite{ls}\qss
is\dss the\sss topic of\qss Section\qss \ref{ls}.\oss
In\dss Section\qss \ref{ls-applications}\qss these\sss theorems are complemented\dss
by several\dss popular\sss and\sss elementary\sss applications,\oss
usually\sss going\dss under\dss the name\qss ``Borsuk--Ulam\qss theorems''.\oss

Section\qss \ref{tilings}\qss is\dss devoted\dss to\sss the easier\dss
part\sss of\qss Lebesgue\dss results\fff:\oss 
Lebesgue\dss construction of\dss coverings\qss
(also known as\sss tilings)\qss 
of\trs $\rrr^{\fff n}$\dss by\sss small\sss
closed sets having\dss the property\dss that\dss no more\sss than\dss $n\qff +\qff 1$\dss
of\trs these sets may\dss have a common\dss point\halfff.\oss
The usual\sss presentations of\trs Lebesgue\dss construction are\sss
limited\dss by\sss a\sss picture illustrating\dss the case\qss $n\off =\off 2$\nnsp.\oss 
In\dss fact\halfff,\oss it\dss is\dss much\sss more interesting\dss in\dss higher
dimensions.\oss
We complement\dss Lebesgue\dss construction\dss by a full\sss description of\trs the intersection\dss
pattern of\trs these sets,\oss which,\pss
strangely\sss enough,\oss appears\sss to\sss be new.

Section\qss \ref{wh}\qss is\dss devoted\dss to\qss
Hurewicz's\dss modification\dss \cite{h2}\dss
of\qss Lebesgue\dss theory.\oss
Hurewicz's\dss proofs\dss bear a substantial\sss similarity\dss to\dss Sperner's\dss
ones,\pss but\sss neither simplices,\pss nor\dss the double counting\dss are\sss used.\oss 
Section\qss \ref{nerve}\qss is\dss
devoted\dss to\sss combinatorics of\qss Hurewicz's\dss results and\sss an\sss
application of\qss
Hurewicz's\trs theorems
to\trs Lebesgue\dss coverings of\trs $\rrr^{\fff n}$\dnsp.\oss
One would\dss think\dss that\dss this\dss was done around\qss 1930,\oss 
but\halfff,\oss
apparently\halfff,\oss this\dss is\dss not\dss the case.\oss
Surprisingly\halfff,\oss this application\dss leads\sss to\sss the strong\qss
Kuhn's\dss lemma\qss \cite{ku}\qss from\qss 1960,\oss
which appears\sss in\dss Section\qss \ref{nerve}\qss as\qss Theorem\qss \ref{hurewicz-lemma-nerve}.\oss

\myuppar{Hurewicz\dss and\dss path-following\sss algorithms.}
Like\dss Sperner\halfff,\pss Hurewicz\dss used\dss the\sss principle\qss 
``an\sss odd\dss number\dss cannot\dss be equal\dss to $0$\nnsp''\qss
in order\dss to prove\sss the existence of\dss desired\dss intersection\dss points.\oss
But\dss along\dss the way\dss he constructed\sss a path\dss leading\dss
to\sss these points,\oss
and\dss his proofs almost\sss explicitly\sss contain a\qss
\emph{path-following\dss algorithm}\qss of\trs the variety\dss which\dss
became\sss popular\dss after\qss H.\dss Scarf\qss \cite{sc}.\oss
Hurewicz's\dss path
consists of\dss several\sss segments\qss
(see\qss Lemma\qss \ref{intersection-path}),\oss
and\dss passing\dss from a segment\dss to\sss the next\sss one\dss
is\dss nothing else but\sss a\qss \emph{pivot\dss step}\qss of\trs
this\dss algorithm.\oss
After\dss reading\dss Section\qss \ref{wh}\qss
all\dss this should\dss be obvious\sss to anyone\sss 
familiar\dss
with\dss such\sss algorithms.\oss

\myuppar{Cubical\sss cochains\sss and\dss products of\dss cubical\sss cochains.}
The most\sss conceptual\dss 
way\dss to prove covering\dss theorems of\qss Lusternik--Schnirelmann\dss
is\dss to use\sss the multiplication\dss
in\dss the cohomology\dss rings of\trs the real\dss projective spaces.\oss
The multiplication\dss of\trs the relative cohomology\dss classes
allows\sss to prove\dss Lebesgue\dss covering\dss theorems\sss in\sss
a similar\dss manner\halfff.\oss

In\dss Sections\qss \ref{cubical-cochains}\qss and\qss \ref{sphere-cochains}\qss
we adapt\qss Serre's\qss \cite{se}\qss definition of\dss products of\dss singular\sss
cubical\sss cochains\sss to\sss discrete setting\dss 
and combinatorial\sss cubical\sss cochains.\pss
By\dss the reasons partially\sss outlined\dss in\dss Appendix\qss \ref{two-cubes}\qss
(morally\sss speaking\halfff,\oss in\sss order\dss to be in agreement\dss
with\qss Poincar\'{e}\dss duality),\oss
for\dss the discussion of\trs products\sss the
solid cubes of\qss Lebesgue\dss are\sss replaced\dss
by\dss discrete\sss cubes consisting in dimension $n$ of\dss
$2^{\fff n}$ points\sss
corresponding\dss to\sss the centers of\trs adjacent\qss Lebesgue\dss cubes.\oss

Products of\dss cubical\sss cochains\sss lead\dss to a new\sss approach\dss
to\trs Lebes\-gue\dss and\qss Lusternik--Schnirel\-mann\qss covering\dss theorems\dss
which\dss is\dss both conceptual\sss and el\-e\-men\-tary\halfff.\oss
The main\dss results are\dss new\dss purely\sss combinatorial\qss 
``cubical\dss lemmas''.\oss
Their\sss cochain\dss versions are\qss 
Theorems\qss \ref{products-induction}\qss and\qss \ref{products-faces}\qss
in\qss Lebesgue\dss context\sss and\qss Theorem\qss \ref{power-of-generator}\qss
in\dss Lusternik--Schnirel\-mann\dss context\halfff,\oss 
and\sss combinatorial\dss versions are\qss
Theorems\qss \ref{separation-products}\qss
and\qss \ref{subsets-of-sphere-products}\qss respectively\halfff.\oss
The\sss classical\dss topological\sss covering\dss theorems can\dss be deduced\dss from\qss
Theorems\qss \ref{separation-products}\qss
and\qss \ref{subsets-of-sphere-products}\qss
by\dss standard\dss elementary\dss topological\dss tools.\oss
See\sss the proof\dss of\qss Lebesgue\dss first\sss covering\dss theorem\dss
in\dss Section\qss \ref{cubical-cochains}\qss
and\qss Theorem\qss \ref{ls-product}.\oss

The little\sss theory\sss developed\dss in\dss Section\qss \ref{cubical-cochains}\qss
was partially\dss motivated\dss by\sss the desire\dss to elucidate\qss
Kuhn's\qss cubical\dss Sperner\dss lemma\qss \cite{ku}\qss
with\dss its\dss rather\sss strange,\pss at\dss least\dss to a\sss topologist\halfff,\pss
state\-ment\halfff.\oss
See\dss Section\qss \ref{kuhn}.\oss 
It\dss is\dss not\dss hard\dss to\sss translate\qss Kuhn's\qss lemma\dss into\sss
the language of\dss subsets,\oss and\dss then\dss it\dss turns into 
an\dss immediate corollary\sss of\qss Theorem\qss \ref{separation-products}.\oss
At\dss the same\dss time\dss a weakened\dss version of\qss Kuhn's\dss lemma,\oss
sufficient\dss for all\sss applications,\pss easily\dss follows from\qss
Lebesgue\dss results.\oss
A similar\dss translation\dss turns\qss
Kuhn's\qss strong\dss lemma into
an\dss immediate\sss corollary\sss of\qss Hurewicz's\dss theorems\sss and combinatorics
of\qss Lebesgue\dss tilings.\oss
We present\sss also\qss Kuhn's\dss own\dss proofs,\oss
based\sss on\dss Freudenthal's\dss triangulations\qss
(see\dss Appendix\qss \ref{freudenthal})\qss 
of\dss cubes and\dss a\qss Sperner-like\dss argument\halfff.\oss

Ky\dss Fan's\qss cubical\qss Sperner\qss lemma\qss \cite{kf}\qss
appears\sss to be even\dss more mysterious\dss than\qss Kuhn's\dss one,\oss
despite a fairly simple proof\halfff.\oss
At\dss the first\sss sight\dss it\dss seems\sss that\qss Ky\dss Fan's\qss
\emph{cubical\dss vertex\dss maps}\qss
are cubical\sss analogues of\dss simplicial\dss maps,\oss
but\dss then\dss it\dss turns out\dss that\dss they are not\halfff.\oss
By\dss this reason\dss the author\dss took\dss the liberty\dss to rename\sss
them as\qss \emph{adjacency-preserving\dss maps}.\oss 
The property\sss of\trs being adjacency-preserving\dss can\dss be interpreted
as a sort\sss of\trs transversality\halfff,\pss
and\dss the\sss theory of\qss Section\qss \ref{cubical-cochains}\qss
allows\sss to understand\qss Ky\dss Fan's\qss lemma\sss as\sss a\sss 
natural\sss strengthening\sss
of\qss Lebesgue\dss or\qss Kuhn's\qss results\sss under\dss this\sss transversality\sss
assumption.\oss
See\dss Section\qss \ref{kyfan}.\oss

\newpage
\mysection{Du\qss c\^{o}t\'{e}\qss de\qss chez\qss Lebesgue}{lebesgue}

\myuppar{The big cube.}
Let\dss us fix a natural\dss number\dss $n$\nnsp,\pss
and\dss let\pss
$I\off =\off \{\qff 1\fff,\pff 2\fff,\pff \ldots\fff,\pff n\qff \}$\nnsp.\oss
Let\dss us also fix another natural\dss number\dss $l$\nnsp.\oss
The\qss \emph{big\dss cube}\pss of\dss size\sss $l$\qss is\dss\vspace{3pt} 
\[
\quad
Q
\off =\off
[\trf 0\fff,\pff l\trf]^{\fff n}\qff \subset\qff \rrr^{\fff n}
\qff.
\]

\vspace{-9pt}
A\qss \emph{face}\qss of\dss $Q$\dss is a product\sss of\dss the form\vspace{-0.5pt}
\[
\quad
\prod_{\dff i\qff =\qff 1}^{\dff n}\qff
J_{\fff i}
\qff,
\]

\vspace{-12pt}
where each\dss $J_{\fff i}$\dss is either\dss the interval\dss
$[\trf 0\fff,\pff l\qff]$\nnsp,\oss
or\dss the one-point\sss set\dss $\{\dff 0\dff\}$\nnsp,\oss
or\dss $\{\dff l\qff\}$\nnsp.\oss
The\qss \emph{dimension}\qss of\dss a face of\dss $Q$\dss is\sss defined as\sss
the number of\dss intervals in\dss the corresponding\sss product\halfff.\oss
A face of\dss dimension\dss $m$\dss is also called an\qss \emph{$m$\dnsp-face}.\oss

For\qss $i\qff \in\qff I$\qss let\qss $\pr_i\dff \colon\dff Q\qff \ttoo\qff [\trf 0\fff,\pff l\qff]$\qss
be\sss the projection\qss 
$(\dff x_{\dff 1}\fff,\pff x_{\dff 2}\fff,\pff \ldots\fff,\pff x_{\dff n}\dff)
\off \longmapsto\off
x_{\dff i}$\qss
and\dss let\vspace{4.5pt}
\[
\quad
A_{\dff i}\off =\off \pr_i^{\dff -\dff 1}\dff (\dff 0\dff)
\qff\quad
\mbox{and}\qff\quad
B_{\dff i}\off =\off \pr_i^{\dff -\dff 1}\dff (\dff l\qff)
\qff.
\]

\vspace{-7.5pt}
Clearly,\oss $A_{\dff i}$\dss and\dss $B_{\dff i}$\dss are\dss
$(\fff n\dff -\dff 1\fff)$\dnsp-faces of\dss $Q$\nnsp,\oss and every\dss
$(\fff n\dff -\dff 1\fff)$\dnsp-face of\dss $Q$\dss is equal\dss to
either\dss $A_{\dff i}$\dss or\dss $B_{\dff i}$\dss for some\dss $i$\nnsp.\oss
For\sss every\qss $i\qff \in\qff I$\qss the faces\qss $A_{\dff i}\fff,\pff B_{\dff i}$\dss are said\dss 
to be\qss \emph{opposite}.\oss
The usual\dss boundary\dss
$\bd Q$\dss of\dss $Q$\dss is equal\dss to\sss the union of\trs
all\sss $(\fff n\dff -\dff 1\fff)$\dnsp-faces\qss $A_{\dff i}\fff,\pff B_{\dff i}$\nsp.\oss

\myuppar{The cubes.}
An \emph{$n$\dnsp-cube}\pss
is a product of\dss the form\qss\vspace{-3pt}
\[
\quad
\prod_{\dff i\qff =\qff 1}^{\dff n}\qff
[\trf a_{\dff i}\fff,\pff a_{\dff i}\qff +\qff 1\qff ]
\qff,
\]

\vspace{-12pt}
where\qss $a_{\dff 1}\fff,\pff a_{\dff 2}\fff,\pff \ldots\fff,\off a_{\dff n}$\qss
are non-negative integers\qss $\leq\qff l\qff -\qff 1$\nnsp.\oss
The\sss $n$\dnsp-cubes form a partition of\dss $Q$\dss in an obvious sense.\oss
A\qss \emph{cubical\dss set}\pss is a subset\sss of\dss $Q$\dss 
equal\dss to\sss a union
of\dss $n$\dnsp-cubes.\oss
If\dss $e$\dss is a cubical\sss set\halfff,\oss then\dss $\overline{e}$\dss
denotes\sss the union of\dss all $n$-cubes not\sss contained\dss in\dss $e$\nnsp.\oss

A\qss \emph{cube}\pss ({\fff}which is not\dss necessarily\dss an $n$\dnsp-cube)\qss 
is a product\sss of\dss the form\vspace{-3pt}
\[
\quad
c
\off =\off
\prod_{\dff i\qff =\qff 1}^{\dff n}\qff
J_{\fff i}
\qff,
\]

\vspace{-12pt}
where\dss for every\qss $i\qff \in\qff I$\qss
either\qss
$J_{\fff i}
\off =\off
[\trf a_{\dff i}\fff,\pff a_{\dff i}\qff +\qff 1\trf]$\qss  
for a non-negative integer\qss $a_{\dff i}\qff \leq\qff l\qff -\qff 1$\nnsp,\oss
or\dss the one-point\sss set\qss 
$J_{\fff i}
\off =\off
\{\dff a_{\dff i}\trf\}$\qss
for a non-negative integer\qss $a_{\dff i}\qff \leq\qff l$\nnsp.\oss
The\qss \emph{dimension}\qss of\dss $c$\dss is\sss defined\sss as\sss the number of\dss intervals in\dss
the above\dss product\halfff.\oss
A cube of\dss dimension\dss $m$\dss is also called a\dss
\emph{$m$\dnsp-cube}.\oss
A\qss \emph{$k$\dnsp-face}\qss of\dss $c$\dss
is a\dss $k$\dnsp-cube contained\dss in\dss $c$\nnsp.\oss
Every\dss $k$\dnsp-face of\dss an $m$\dnsp-cube\dss $c$\dss is obtained\dss by\dss replacing\sss
in\dss the above product\qss $m\qff -\qff k$\qss intervals\qss
$[\trf a_{\dff i}\fff,\pff a_{\dff i}\qff +\qff 1\trf]$\qss by\sss one of\dss their endpoints.

\myuppar{Cubical chains.}
A\qss \emph{cubical\sss chain of\dss dimension\dss $m$}\nnsp,\oss
or simply an\dss
\emph{$m$\dnsp-chain},\oss is\sss a formal sum
of\dss several\dss $m$\dnsp-cubes with coefficients
in\dss $\ftwo$\nnsp.\oss
An\dss $m$\dnsp-chain\dss $\gamma$\dss is essentially\dss a set of\dss $m$\dnsp-cubes,\oss
namely\dss the set\sss of\dss all\dss $m$\dnsp-cubes entering\dss $\gamma$\dss
with\dss the coefficient\dss $1$\nnsp.\oss
The union of\dss these\dss $m$\dnsp-cubes is called\dss the\qss
\emph{support}\qss of\dss $\gamma$\dss and\dss is denoted\dss by\dss $\num{\gamma}$\nnsp.\oss
Clearly,\oss an\dss $m$\dnsp-chain\dss $\gamma$\dss is equal\dss to\sss the sum of\dss
all\dss $m$\dnsp-cubes\dss $c$\dss such\dss that\qss
$c\qff \subset\qff \num{\gamma}$\nnsp.\oss
We will call\dss such\dss $m$\dnsp-cubes\sss $c$\sss the\dss
\emph{$m$\dnsp-cubes of}\qss $\gamma$\nnsp,\oss
or\dss simply\dss the\qss \emph{cubes}\qss of\dss $\gamma$\nnsp.\oss
Clearly,\pss $\gamma$\dss is equal\dss to\sss the sum of\dss the cubes of\dss $\gamma$\nnsp.\oss

The boundary\dss $\partial\fff c$\dss of\dss an\sss $m$\dnsp-cube\dss $c$\dss
is defined as\sss the sum of\dss all\dss its\dss $(\fff m\dff -\dff 1\fff)$\dnsp-faces.\oss
The\qss \emph{boundary\dss operator}\qss $\partial$\dss is\sss the extension of\dss
the boundary of\dss $m$\dnsp-cubes\sss to all\sss $m$\dnsp-chains by\dss the linearity.\oss
As always,\oss the main\sss property of\trs the boundary operator $\partial$\dss is\vspace{3pt} 
\begin{equation}
\label{cubes-dd}
\quad
\partial\dff \circ\dff \partial\off =\off 0
\qff.
\end{equation}

\vspace{-9pt}
It\dss is sufficient\dss to check\dss this property\dss for\dss $m$\dnsp-cubes.\oss
In\dss this case it\dss follows from\dss the fact\dss that\dss
an\sss $(\fff m\dff -\dff 2\fff)$\dnsp-face of\dss an\dss $m$\dnsp-cube\dss $c$\dss
is contained\dss in exactly\dss two\dss $(\fff m\dff -\dff 1\fff)$\dnsp-faces of\dss $c$\nnsp.\oss

\myuppar{Projections.}
We\sss need a limited class of\dss maps similar\dss to simplicial\dss maps.\oss
Let\qss $k\qff \leq\qff n$\qss and\dss let\vspace{3pt}
\[
\quad
p\dff \colon\dff
[\trf 0\fff,\pff l\qff]^{\fff n}
\qff \ttoo\qff
[\trf 0\fff,\pff l\qff]^{\fff k}
\]

\vspace{-9pt}
be\sss the projection defined\dss by\sss a set\sss of\sss $k$\sss coordinates.\oss
If\sss $c$\sss is an\sss $m$\dnsp-cube,\oss
then\dss the image\dss $p\dff(\dff c\trf)$\dss is an\sss $m'$\dnsp-cube for some\dss
$m'\qff \leq\qff m$\nnsp.\qff\oss
For an $m$\dnsp-cube\dss $c$\dss let\vspace{3pt}
\[
\quad
p_{\fff *}\dff(\dff c\trf)
\off =\off
c
\hspace*{1.2em}\mbox{if}\hspace*{0.6em}
p\dff(\dff c\trf)\hspace*{0.6em}
\mbox{is\dss an\dss $m$-cube\dss and}
\]

\vspace{-36pt}
\[
\quad
p_{\fff *}\dff(\dff c\trf)
\off =\off
0
\hspace*{1.2em}\mbox{if}\hspace*{0.6em}
p\dff(\dff c\trf)\hspace*{0.6em}
\mbox{is\dss an\dss $m'$-cube\dss with}\hspace*{0.6em}
m'\qff <\qff m\qff.
\]

\vspace{-9pt}
The map\dss $p_{\fff *}$\dss extends\dss to all\dss $m$\dnsp-chains by\dss linearity.\oss
The extended\dss map is also denoted\dss by\dss $p_{\fff *}$\dss
and is called\dss the map\qss \emph{induced}\pss by\dss $p$\nnsp.\oss
Taking\dss the induced\dss maps respects\sss the compositions in\dss the sense\sss that\qss
if\vspace{3pt}
\[
\quad
r\dff \colon\dff
[\trf 0\fff,\pff l\qff]^{\fff k}
\qff \ttoo\qff
[\trf 0\fff,\pff l\qff]^{\fff l}
\]

\vspace{-9pt}
is a similar\dss projection,\oss
then\qss
$(\dff r\dff \circ\dff p\trf)_*
\off =\off
r_{\dff *}\dff \circ\qff p_{\fff *}$\nnsp,\oss
as a\sss trivial\sss verification\sss shows.\oss

\myuppar{Lemma.} 
\emph{The map\dss $p_{\fff *}$\dss
commutes with\qss $\partial$\qss in\dss the sense\sss that\dss for every\dss chain\dss $\gamma$}\vspace{3pt}
\begin{equation}
\label{cube-d-f}
\quad
\partial\dff \circ\dff p_{\fff *}\dff(\dff \gamma\trf)
\off =\off
p_{\fff *}\dff \circ\qff \partial\qff (\dff \gamma\trf)
\qff.
\end{equation}

\vspace{-9pt}
\proof
It\dss is sufficient\dss to prove\qss (\ref{cube-d-f})\qss when\dss
$\gamma$\dss is equal\dss to an\dss $m$\dnsp-cube\dss $c$\nnsp.\oss
If\dss $p\dff(\dff c\trf)$\dss is also an\sss $m$\dnsp-cube,\oss (\ref{cube-d-f})\qss is obvious.\oss
Also,\oss if\dss $p\dff(\dff c\trf)$\dss is an\sss $m'$\dnsp-cube with\qss
$m'\qff \leq\qff m\qff -\qff 2$\nnsp,\oss
then\sss both sides of\qss (\ref{cube-d-f})\qss are equal\dss to\dss $0$\dss and\dss hence\qss
(\ref{cube-d-f})\qss holds.\oss
It\dss remains\sss to consider\dss the case when\dss 
$p\dff(\dff c\trf)$\dss is an\sss $(\fff m\dff -\dff 1\fff)$\dnsp-cube.\oss
In\dss this case\dss 
$p\dff(\dff d\trf)$\dss is an\sss $(\fff m\dff -\dff 1\fff)$\dnsp-cube
for exactly\dss two\dss $(\fff m\dff -\dff 1\fff)$\dnsp-faces\qss 
$d\off =\off d_{\dff 1}\fff,\pff d_{\dff 2}$\qss
of\qss $c$\dss and\dss is\sss an\sss $(\fff m\dff -\dff 2\fff)$\dnsp-cube 
for all other\sss $(\fff m\dff -\dff 1\fff)$\dnsp-faces\dss $d$\dss
of\dss $c$\nnsp.\oss
Clearly,\oss
$p\dff(\dff d_{\dff 1}\trf)\off =\off p\dff(\dff d_{\dff 2}\trf)$\qss
and\dss
$p\dff(\dff d\trf)\off =\off 0$\dss
if\dss $d\off \neq\off d_{\dff 1}\fff,\off d_{\dff 2}$\nsp.\oss
Hence\dss this case\qss\vspace{3pt}
\[
\quad
p_{\fff *}\dff(\dff \partial\dff c\trf)
\off =\off
p_{\fff *}\dff(\dff d_{\dff 1}\trf)
\qff +\qff
p_{\fff *}\dff(\dff d_{\dff 2}\trf)
\off =\off
0
\qff.
\]

\vspace{-9pt}
Also,\oss in\dss this case\qss
$\partial\dff p_{\fff *}\dff(\dff c\trf)\off =\off \partial\dff 0\off =\off 0$\nnsp.\oss
Hence\dss the lemma\sss
holds in\dss this case also.\oss  \eproof

\myuppar{Lebesgue\dss construction.}
Suppose\sss that\sss 
$\gamma$ is an $m$\dnsp-chain
and $e$ is a\sss cubical\sss set\halfff.\pss 
Let\dss $\gamma_{\fff e}$\dss be\sss the sum of\dss 
all\dss $m$\dnsp-cubes
of\dss $\gamma$\dss
contained\dss in\dss $e$\dss and\trs 
$\gamma'$\dss be\sss the sum of\dss all\sss other $m$\dnsp-cubes
of\dss $\gamma$\nnsp.\oss
Then\vspace{3pt}
\[
\quad
\gamma\off =\off \gamma_{\fff e}\qff +\qff \gamma'
\]

\vspace{-9pt}
and\dss hence\qss 
$\partial\dff \gamma
\off =\off 
\partial\dff \gamma_{\fff e}\qff +\qff \partial\dff \gamma'$\nnsp.\oss
Clearly,\oss
$\num{\gamma_{\fff e}}\off \subset\off \num{\gamma}\qff \cap\qff e$\qss
and\qss
$\num{\gamma'}\off \subset\off \num{\gamma}\qff \cap\qff \overline{e}$\nnsp.\oss
Let\dss $\delta$\dss be\sss the sum of\dss all
common\ $(\fff m\dff -\dff 1\fff)$\dnsp-cubes of\dss the chains\dss
$\partial\dff \gamma_{\fff e}$\dss and\dss $\partial\dff \gamma'$\nnsp.\oss
Then\dss
$\num{\delta}$\dss is contained\dss in\dss both\qss $\num{\gamma}\qff \cap\qff e$\qss
and\qss
$\num{\gamma}\qff \cap\qff \overline{e}$\nnsp.\oss
Let\dss 
$\delta_{\dff e}$\dss be the sum of\dss all\dss
$(\fff m\dff -\dff 1\fff)$\dnsp-cubes of\dss $\partial\dff \gamma_{\fff e}$\dss
which are not\dss $(\fff m\dff -\dff 1\fff)$\dnsp-cubes of\dss $\delta$\nnsp.\oss
Clearly,\vspace{3pt}
\begin{equation}
\label{gamma-e-decomposition-d}
\quad
\partial\dff \gamma_{\fff e}
\off =\off
\delta_{\dff e}
\qff +\qff
\delta
\qff.
\end{equation}

\vspace{-9pt}
A $m$\dnsp-chain $\gamma$ is\dss said\dss to\sss be\qss
\emph{proper}\oss if\pss 
$\num{\partial\dff \gamma}\off =\off \num{\gamma}\qff \cap\qff \bd\fff Q$\qss
and\qss
\emph{special}\oss if\qss it\trs is\dss proper\sss and\dss
$\num{\gamma}$\dss
is disjoint\dss from\qss
$A_{\dff i}\qff \cup\qff B_{\dff i}$\qss
for\qss
$i\qff \leq\qff n\qff -\qff m$\nnsp.\oss
The following\dss two\sss lemmas play a key role.

\mypar{Lemma.}{proper}
\emph{If\pss $\gamma$\dss is\dss proper\halfff,\oss
then\dss $\delta$\dss
is\dss also proper\halfff.}

\proof
Let\dss us prove\sss first\dss 
that\qss
$\num{\partial\dff \delta}\off \subset\off \bd\fff Q$\nnsp.\oss
The\dss $(\fff m\dff -\dff 1\fff)$\dnsp-cubes of\dss $\delta$\dss
are exactly\dss the cubes entering\dss the sum\qss
$\partial\dff \gamma_{\fff e}\qff +\qff \partial\dff \gamma'$\qss
twice,\oss i.e.\qss cancelling\dss in\dss this sum.\oss
The other\dss $(\fff m\dff -\dff 1\fff)$\dnsp-cubes in\dss this sum do not\sss
cancel and\sss are\dss $(\fff m\dff -\dff 1\fff)$\dnsp-cubes of\dss
$\partial\dff \gamma$\nnsp.\oss
Since\dss $\gamma$\dss is proper\halfff,\oss
these\dss $(\fff m\dff -\dff 1\fff)$\dnsp-cubes are contained\dss in\dss $\bd\fff Q$\nnsp.\oss
It\dss follows\sss that\halfff,\pss
in\dss particular\halfff,\oss
$\num{\delta_{\dff e}}\qff \subset\pff \bd\fff Q$\nnsp.\oss
Since\qss
$\partial\dff \circ\dff \partial\off =\off 0$\nnsp,\oss
the equality\qss 
(\ref{gamma-e-decomposition-d})\qss implies\dss that\vspace{3pt}
\[
\quad
\partial\dff \delta_{\dff e}
\qff +\qff
\partial\dff \delta
\off =\off
\partial\dff \circ\dff \partial\dff \left(\dff \gamma_{\fff e}\dff\right)
\off =\off
0
\]

\vspace{-9pt}
and\dss hence\qss
$\partial\dff \delta
\off =\off
-\qff
\partial\dff \delta_{\dff e}
\off =\off
\partial\dff \delta_{\dff e}$\nnsp.\oss
Therefore\qss
$\num{\partial\dff \delta}
\off =\off 
\num{\partial\dff \delta_{\dff e}}
\off \subset\off
\num{\delta_{\dff e}}
\qff \subset\pff \bd\fff Q
$\qss
and\dss hence\vspace{3pt}
\[
\quad
\num{\partial\dff \delta}
\off \subset\off
\bd\fff Q \qff.
\]

\vspace{-9pt}
Every\dss $(\fff m\dff -\dff 1\fff)$\dnsp-cube of\dss $\delta$\dss
is\sss an\dss $(\fff m\dff -\dff 1\fff)$\dnsp-face of\dss 
one $m$\dnsp-cube of\dss $\gamma_{\dff E}$\dss
and\dss one of\dss $\gamma'$\nnsp.\oss
Since\dss $\gamma$\dss is a proper\halfff,\oss 
these\dss $m$\dnsp-cubes are not\sss contained\dss in\dss $\bd\fff Q$\nnsp.\oss
It\dss follows\dss that\dss no\dss $(\fff m\dff -\dff 1\fff)$\dnsp-cube of\dss $\delta$\dss
is contained\dss in\dss $\bd\fff Q$\nnsp.\oss
Together\sss with\qss 
$\num{\partial\dff \delta}
\off \subset\off
\bd\fff Q$\nnsp,\oss
this implies\sss that\dss $\delta$\dss is\dss proper\halfff.\oss  \eproof

\mypar{Lemma.}{special}
\emph{Suppose\sss that\dss $\gamma$\dss is
a special\dss $m$\dnsp-chain.\qff\oss
If\halfff,\oss moreover\halfff,\oss}\qss\vspace{4.5pt} 
\[
\quad
e\qff \cap\qff B_{\dff n\dff -\dff m\dff +\dff 1}
\off =\off
\varnothing
\hspace*{1.25em}
\mbox{\emph{and}}
\hspace*{1.5em}
\num{\gamma\fff}\qff \cap\qff A_{\dff n\dff -\dff m\dff +\dff 1}\off \subset\off e
\qff,
\]
 
\vspace{-7.5pt}
\emph{then\dss
$\num{\gamma'}$\dss is disjoint\trs from\dss $A_{\dff n\dff -\dff m\dff +\dff 1}$\dss
and\dss the\qss
$(\fff m\dff -\dff 1\fff)$\dnsp-chain\dss $\delta$\dss is special.\oss}

\proof
Let\dss us\sss prove\sss first\dss that\dss
$\num{\gamma'}$\dss is disjoint\dss from\dss $A_{\dff n\dff -\dff m\dff +\dff 1}$\nnsp.\oss
Suppose\sss that\dss $c$\dss is an $m$\dnsp-cube of\dss $\gamma$\dss 
intersecting\dss $A_{\dff n\dff -\dff m\dff +\dff 1}$\nnsp.\oss
Since\dss $A_{\dff n\dff -\dff m\dff +\dff 1}$\dss is an\dss
$(\fff n\dff -\dff 1\fff)$\dnsp-face of\dss $Q$\nnsp,\oss
in\dss this case either\dss 
$c$\dss is contained\dss in\dss $A_{\dff n\dff -\dff m\dff +\dff 1}$\nnsp,\oss
or\qss
$c\qff \cap\qff A_{\dff n\dff -\dff m\dff +\dff 1}$\qss 
is an\dss $(\fff n\dff -\dff 1\fff)$\dnsp-face of\dss $c$\nnsp.\oss
In\dss the first\sss case\dss\vspace{2.675pt} 
\[
\quad
c
\off \subset\off
\num{\gamma\fff}\qff \cap\qff A_{\dff n\dff -\dff m\dff +\dff 1}
\off \subset\off
e
\qff,
\] 

\vspace{-9.375pt}
and\dss hence\dss $c$\dss is a cube of\dss $\gamma_{\fff e}$\nnsp.\oss
In\dss the second case\qss 
$c\qff \cap\qff A_{\dff n\dff -\dff m\dff +\dff 1}
\off \subset\off
e$\qss
and
since\dss $e$\dss is a cubical set\halfff,\oss
the intersection\qss
$c\qff \cap\qff A_{\dff n\dff -\dff m\dff +\dff 1}$\qss
is contained\dss in\dss an $(\fff n\dff -\dff 1\fff)$\dnsp-face of\qss an\dss $n$\dnsp-cube\qss
$d\qff \subset\qff e$\nnsp.\oss
It\dss follows\dss that\qss 
$c\qff \subset\qff d\qff \subset\qff e$\qss
and\dss hence\dss $c$\dss is a cube of\dss $\gamma_{\fff e}$\dss in\dss this case also.\oss
It\dss follows\sss that\dss no\dss $m$\dnsp-cube of\dss $\gamma'$\dss intersects\dss
$A_{\dff n\dff -\dff m\dff +\dff 1}$\dss
and\dss hence\dss
$\num{\gamma'}$\dss is disjoint\dss from\dss $A_{\dff n\dff -\dff m\dff +\dff 1}$\nnsp.\oss

Since\qss
$\num{\delta}
\off \subset\off 
\num{\gamma}$\nnsp,\oss
the support\dss $\num{\delta}$\dss
is disjoint\dss from\qss
$A_{\dff i}\qff \cup\qff B_{\dff i}$\qss
for\qss
$i\qff \leq\qff n\qff -\qff m$\nnsp.\oss
Since\dss
$\num{\delta}$\dss is contained\dss in\qss 
$\num{\partial\dff \gamma_{\fff e}}
\off \subset\off
e$\qss
and\dss $e$\dss is disjoint\dss from\qss
$B_{\dff n\dff -\dff m\dff +\dff 1}$\nnsp,\oss
the support\dss $\num{\delta}$\dss
is disjoint\dss from\qss
$B_{\dff n\dff -\dff m\dff +\dff 1}$\nnsp.\oss
By\dss the construction,\pss
$\num{\delta} 
\qff \subset\qff 
\num{\partial\dff \gamma'}
\qff \subset\qff 
\num{\gamma'}$\nnsp.\oss
By\dss the previous paragraph\dss
$\num{\gamma'}$\dss is disjoint\dss from\dss
$A_{\dff n\dff -\dff m\dff +\dff 1}$\nnsp,\oss
and\dss hence\dss $\num{\delta}$\dss is disjoint\dss from\dss
$A_{\dff n\dff -\dff m\dff +\dff 1}$\dss also.\oss  \eproof

\mypar{Corollary.}{sides-and-chains}
\emph{Suppose\sss that\dss the assumptions of\dss the previous lemma hold.\oss
Let\pss
$\beta$\dss
be\dss the sum of\dss all $(\fff m\dff -\dff 1\fff)$\dnsp-cubes\dss
of\oss $\partial\dff \gamma_{\fff e}$\trs contained\dss in\qss
$A_{\dff n\dff -\dff m\dff +\dff 1}$\nnsp.\qff\oss
Then}\qss\vspace{2.675pt}
\begin{equation}
\label{sides-decomposition}
\quad
\partial\dff \gamma_{\fff e}
\off =\off
\beta
\qff +\pff 
\delta
\qff +\qff
\delta_{\dff 0}
\qff,
\end{equation}

\vspace{-9.375pt}
\emph{for some\dss $(\fff m\dff -\dff 1\fff)$\dnsp-chain\dss 
$\delta_{\dff 0}$\dss
such\dss that\dss
$\num{\delta_{\dff 0}}$\dss
is\dss contained\trs in\oss  
$\dis
\bigcup\nolimits_{\dff k\qff >\qff n\qff -\qff m\qff +\qff 1}\qff
A_{\dff k}\qff \cup\qff B_{\dff k}$.}

\proof
By\qss Lemma\qss \ref{special}\qss the support\dss
$\num{\gamma'}$\dss is disjoint\trs from\dss $A_{\dff n\dff -\dff m\dff +\dff 1}$\dss
and\dss hence
an\dss $(\fff m\dff -\dff 1\fff)$\dnsp-cube 
contained\trs in\dss $A_{\dff n\dff -\dff m\dff +\dff 1}$\dss
is\sss a\sss cube\sss of\dss $\partial\dff \gamma$\dss
if\dss and\sss only\trs if\trs it\dss is a cube of\trs
$\partial\dff \gamma_{\dff E}$\nnsp.\oss
Therefore\dss
$\beta$\dss
is\sss equal\dss to\sss the sum of\trs all $(\fff m\dff -\dff 1\fff)$\dnsp-cubes
of\dss $\partial\dff \gamma_{\fff e}$\sss contained\sss in\sss
$B_{\dff n\dff -\dff m\dff +\dff 1}$\nnsp.\oss

The support\dss
$\num{\gamma_{\fff e}}$\dss
is contained\dss in\sss
$\num{\gamma}$\sss
and\dss hence is disjoint\dss from\qss
$A_{\dff i}\qff \cup\qff B_{\dff i}$\qss
for\qss
$i\qff \leq\qff n\qff -\qff m$\nnsp.\oss
Also,\oss
it\dss is contained\dss in\dss $e$\dss
and\dss hence is disjoint\dss from\dss $B_{\dff n\dff -\dff m\dff +\dff 1}$\nnsp.\oss
It\dss follows\dss that\sss
every\dss $(\fff m\dff -\dff 1\fff)$\dnsp-cube
of\qss $\partial\dff \gamma_{\fff e}$\dss
is either\sss an\dss $(\fff m\dff -\dff 1\fff)$\dnsp-cube\dss
of\qss $\delta$\nnsp,\oss
or\dss is\dss contained\dss in\dss $A_{\dff n\dff -\dff m\dff +\dff 1}$\dss
and\dss hence is\dss an $(\fff m\dff -\dff 1\fff)$\dnsp-cube of\dss $\beta$\nnsp,\oss
or\dss is\dss contained\trs in\qss $A_{\dff k}\qff \cup\qff B_{\dff k}$\qss
for\dss some\qss $k\qff >\qff n\qff -\qff m\qff +\qff 1$\nnsp.\oss
Therefore one can\dss take as\dss $\delta_{\dff 0}$\dss the\sss sum 
of\trs the\sss latter\dss $(\fff m\dff -\dff 1\fff)$\dnsp-cubes.\oss  \eproof

\myuppar{Essential\sss chains.}
The\qss \emph{essential\dss chains},\oss to be defined\dss in a moment\halfff,\oss
are a\sss tool\dss for\dss proving\dss that\dss some intersections of\dss cubical\sss sets
are non-empty.\oss
Let\dss us begin\dss with\qss Lebesgue\dss definition\qss
(who did\dss not\dss use\sss this\sss term).\oss
Let\dss us choose arbitrary\sss non-integer\sss numbers\qss
$a_{\dff 1}\fff,\pff a_{\dff 2}\fff,\pff \ldots\fff,\pff a_{\dff n}
\qff \in\qff
[\trf 0\fff,\pff l\trf]$\nnsp.\oss
For an\dss integer\dss $u$\dss 
between $0$ and $n$\dss let\dss
$L^u$\dss be\sss
the plane in\dss $\rrr^{\fff n}$\dss defined\dss by\dss the equations\vspace{1.5pt}
\[
\quad
x_{\dff u\dff +\dff 1}\off =\off a_{\dff u\dff +\dff 1}\qff,\hspace*{1.2em}
x_{\dff u\dff +\dff 2}\off =\off a_{\dff u\dff +\dff 2}\qff,\hspace*{1.2em}
\ldots\fff,\hspace*{1.2em}
x_{\dff n}\off =\off a_{\dff n}
\qff,
\]

\vspace{-10.5pt}
where\qss
$(\dff x_{\dff 1}\fff,\pff x_{\dff 2}\fff,\pff 
\ldots\fff,\pff x_{\dff n}\dff)$\qss
are\sss the coordinates\sss in\dss $\rrr^{\fff n}$\dnsp.\oss
In\dss particular\halfff,\pss
$L^u$\dss is\dss a $u$\dnsp-dimensional\dss plane,\pss
$L^0$\dss consists of\dss a single point\halfff,\oss namely\qss
$(\dff a_{\dff 1}\fff,\pff a_{\dff 2}\fff,\pff 
\ldots\fff,\pff a_{\dff n}\dff)$\nnsp,\oss
and\qss
$L^n\off =\off \rrr^{\fff n}$\dnsp.\oss
The numbers\dss $a_{\dff i}$\dss are chosen\dss to be non-integers
in order\dss to ensure\sss the following\dss property\fff:\oss
if\dss $c$\sss is an $(\fff n\dff -\dff u\fff)$\dnsp-cube,\oss
then\dss the intersection\qss
$c\qff \cap\qff L^u$\qss is either empty,\oss 
or consists of\dss a single point\halfff.\oss
More generally,\oss if\dss $c$\sss is an $(\fff n\dff -\dff k\fff)$\dnsp-cube,\oss
then\qss
$c\qff \cap\qff L^u$\qss is either empty,\oss
or a cube of\dss dimension\qss $u\qff -\qff k$\qss in\dss the usual\dss geometric sense\qss
(of\dss course,\oss this intersection is not\sss a cube according\dss to our\sss definition).\oss
Using\dss the modern\dss language,\oss
one can say\dss that\sss all cubes from our\dss partition of\dss $Q$\dss are in\qss
\emph{general\dss position}\qss with respect\dss to\sss the planes\dss $L^u$\nnsp.\oss

Now we are ready\dss for\dss the definition.\oss
An $m$\dnsp-chain $\gamma$ is said\dss to be\qss 
\emph{essential}\oss if\trs $\gamma$\dss is\dss special\sss and\vspace{1.5pt}
\[
\quad
\num{\gamma}\qff \cap\qff L^{n\dff -\dff m}
\]

\vspace{-10.5pt}
consists of\dss an\sss odd\dss number of\dss points.\oss
In\dss particular\halfff,\oss if\dss $\gamma$\sss is essential,\oss
then\dss $\num{\gamma}$\dss is non-empty.\oss
In\qss Lebesgue\dss terminology\dss
a\sss set\sss of\qss the\dss form\sss $\num{\gamma}$\nnsp,\oss
where\sss $\gamma$\dss is\dss an\sss essential\dss
$(\fff n\dff -\dff p\fff)$\dnsp-chain,\oss
is\sss called\sss a\sss set\qss (of\dss the\sss type)\qss $J_{\fff p}$\nsp.\oss
This\sss is\sss a\sss central\dss notion of\dss his paper\qss \cite{l-two}.\oss

An obvious drawback of\dss this definition\sss is\dss the dependence
on\dss the choice of\dss the numbers\dss $a_{\dff i}$\nnsp.\oss
In\dss fact\halfff,\pss the property of\dss being an essential\sss chain
does not\sss depend on\dss this choice 
and\qss Lebesgue\dss arguments work  
for any\dss
fixed choice.\oss
In\dss the context\sss of\qss Lebesgue\dss papers\qss \cite{l-one},\pss \cite{l-two}\qss
there is another drawback\halfff,\oss
the dependence on\dss geometric\dss ideas.\oss
After\dss reliance on\sss geometric\sss language rendered\dss to\sss be unconvincing\qss
Lebesgue\dss outline\qss \cite{l-one}\qss of\dss his ideas about\dss
the invariance of\dss the dimension,\oss
he decided\dss to\qss ``d'arithm\'{e}tiser\dss la d\'{e}monstration\halfff''.\dff\oss
While\trs Lebesgue\dss arguments in\qss \cite{l-two}\qss
are complete and\dss never\dss met\sss any\sss objection,\pss
the following definition is essentially a combinatorial one and\dss allows\dss to
complete\qss Lebesgue\dss plan of\qss
``arithmetization\halfff'' of\dss his arguments.\oss

For\sss every\qss $m\qff \in\qff I$\qss let\qss
$P_{\dff m}
\off =\off
A_{\dff 1}\qff \cap\qff A_{\dff 2}\qff \cap\qff
\ldots\qff \cap\qff
A_{\dff n\dff -\dff m}$\nnsp.\oss
Clearly,\pss $P_{\dff m}$\dss is\dss an $m$\dnsp-face of\dss $Q$\nnsp.\qff\oss
Let\qss\vspace{1.5pt} 
\[
\quad
p_{\fff m}\dff \colon\dff
Q\qff \ttoo\qff P_{\dff m} 
\]

\vspace{-10.5pt}
be\sss the canonical\dss projection.\oss
An\dss $m$\dnsp-chain\dss $\gamma$\dss is said\dss to\sss be a\qss
\emph{essential}\oss  
if\qss it\dss is\dss special\sss and\vspace{1.5pt}
\begin{equation*}
\quad
p_{\fff m\dff *}\dff(\trf \gamma\trf)
\off =\dff\off
\fclass{\fff P_{\dff m}\fff}
\qff,
\end{equation*}
 
\vspace{-10.5pt}
where\dss $\fclass{\fff P_{\dff m}\fff}$\dss is\sss the sum of\dss all\sss
$m$\dnsp-cubes contained\dss in\dss $P_{\dff m}$\nsp.\oss
Clearly,\oss if\dss $\gamma$\dss is essential\halfff,\oss
then\qss $\gamma\off \neq\off 0$\qss
and\dss hence\qss $\num{\gamma\fff}\off \neq\off \varnothing$\nnsp.\oss
The main\dss result\sss about\sss essential\sss chains is\qss
Lemma\qss \ref{essential}\qss below,\oss
in which one can understand\dss the\sss term\qss ``essential''\qss
in either of\dss the above\sss two ways.\oss
We will\dss give\sss two proofs of\dss this lemma,\oss
the first\sss one is\sss following\qss Lebesgue\qss \cite{l-two},\oss
while\dss the second one is based on\dss the second definition.\oss
We will\dss not\dss use\sss the equivalence of\trs these\sss two definitions of\dss
essential\sss chains and\dss leave\sss the proof\dss of\trs their equivalence
as an exercise\sss to\sss the reader\halfff.\oss

\myuppar{Intersections\sss of\dss chains\dss with\dss the\sss planes\dss $L^u$\dnsp.}
The first\dss proof\trs requires some additional\dss preparations.\oss
For an\dss integer\dss $u$\dss 
between $0$ and $n$\dss let\qss
$Q^{\fff u}
\off =\off
Q\qff \cap\qff L^u$\nnsp.\oss
Clearly,\pss\vspace{1.5pt}
\[
\quad
Q^{\fff u}
\off =\off
[\trf 0\fff,\pff l\trf]^{u}\dff \times\qff
(\dff 
a_{\dff u\dff +\dff 1}\dff,\pff
a_{\dff u\dff +\dff 2}\dff,\pff
\ldots\dff,\pff
a_{\dff n}
\dff)
\qff.
\]

\vspace{-10.5pt}
Equivalently,\pss
$Q^{\fff u}$\dss
is\dss the $u$\dnsp-dimensional cube\qss
$[\trf 0\fff,\pff l\trf]^{u}
\qff \subset\qff \rrr^{\fff u}
\qff \subset\qff \rrr^{\fff n}$\qss
shifted\dss by\dss the vector\vspace{1.5pt}
\[
\quad
(\dff 
0\dff,\pff 
0\dff,\pff 
\ldots\fff,\pff 
0\dff,\pff
a_{\dff u\dff +\dff 1}\dff,\pff
a_{\dff u\dff +\dff 2}\dff,\pff
\ldots\dff,\pff
a_{\dff n}
\dff)
\qff \in\qff \rrr^{\fff n}
\qff.
\]

\vspace{-10.5pt}
This allows\sss to\sss transfer\dss to\dss $Q^{\fff u}$\sss the notions
of\dss cubes,\oss chains,\oss etc.\oss
Moreover\halfff,\oss we can\dss take\sss the intersection of\dss chains in\dss $Q$\dss
with\dss $Q^{\fff u}$\dnsp.\oss
Namely,\oss let\sss $\gamma$\sss be an $(\fff n\dff -\dff k\fff)$\dnsp-chain\dss in\dss
$Q$\nnsp,\oss i.e.\qss let\vspace{3pt}
\[
\quad
\gamma
\off =\dff\off
\sum\nolimits_{\dff i}\dff c_{\dff i}
\qff,
\] 

\vspace{-9pt}
where\dss $c_{\dff i}$\dss are\dss $(\fff n\dff -\dff k\fff)$\dnsp-cubes in\dss $Q$\nnsp.\oss
The fact\dss that\sss all\sss cubes in\sss $Q$\sss are in general\dss position\dss with\dss
respect\dss to\dss $L^u$\dss allows\dss to define\sss 
the\qss \emph{intersection}\qss
$\gamma\dff \cap\dff L^u$\qss 
as\dss the $(\fff u\dff -\dff k\fff)$\dnsp-chain\vspace{3pt}
\[
\quad
\gamma\qff \cap\qff L^u
\off =\qff\off
\sum\nolimits_{\dff i}\dff c_{\dff i}\qff \cap\qff L^u
\qff,
\] 

\vspace{-9pt}
where\sss the empty\dss intersections are interpreted as\dss $0$\nnsp'{\halfff}s.\oss
The operation of\dss taking\dss the intersection of\dss chains with\dss $L^u$\dss
commutes with\dss $\partial$\nnsp,\oss i.e.\vspace{1.5pt}
\[
\quad
\left(\dff \partial\dff \gamma\dff\right)\qff \cap\qff L^u
\off =\off
\partial\dff \left(\dff \gamma \qff \cap\qff L^u\dff\right)
\]

\vspace{-10.5pt}
for any chain\dss $\gamma$\dss in\dss $Q$\nnsp.\oss
This is\dss trivial\dss for cubes,\oss and extends\sss to\sss the general case
by\dss linearity.\oss

\mypar{Lemma.}{essential}
\emph{Under\dss the assumptions of\pss Lemma\qss \ref{special},\oss  
if\pss the chain\dss $\gamma$\dss is\sss essential,\oss
then\trs the chain\trs $\delta$\trs is\sss also essential.\oss}

\prooftitle{Proof\qss based\sss on\dss intersections}
In\dss this\sss proof\dss we use\qss Lebesgue\dss 
definition of\dss essential\sss chains.\oss
The idea is\sss to consider\dss the intersections of\dss
chains\sss with\dss $L^{n\dff -\dff m\dff +\dff 1}$\nnsp.\oss
For a chain $\varepsilon$ in\dss $Q$\dss
let\dss\vspace{1.5pt}
\[
\quad
\cap\trf{\varepsilon}
\off =\off
\varepsilon\qff \cap\qff L^{n\dff -\dff m\dff +\dff 1}
\qff.
\]

\vspace{-10.5pt}
Then\trs $\cap\dff{\gamma}$\trs and\trs $\cap\dff{\gamma}_{\fff e}$\trs are 
$1$\dnsp-chains\dss  
and\trs $\cap\dff{\delta}$\trs
is\dss a $0$\dnsp-chain\dss in\dss $Q^{\fff n\dff -\dff m\dff +\dff 1}$\nnsp.\oss
By\dss the definition,\pss $\delta$\dss is essential\dss if\dss and\dss only\trs if\qss
$\num{\nsp\cap\dff{\delta}}$\trs consists of\dss an\sss odd\dss number\sss of\dss points.\oss
Since\dss $\gamma$\dss and\dss $\delta$\dss are proper
and\dss $\partial$\dss commutes with\dss taking\sss intersections,\pss
$\cap\dff{\gamma}$\trs and\trs $\cap\dff{\delta}$\trs are
proper as chains in 
$Q^{\fff n\dff -\dff m\dff +\dff 1}$\dnsp.\qff\oss
Since\dss $\gamma$\dss and\dss $\delta$\dss are
special\halfff,\oss
this implies\dss that\trs $\cap\dff{\gamma}$\trs and\trs $\cap\dff{\delta}$\trs are
special\dss as chains\dss in\dss 
$Q^{\fff n\dff -\dff m\dff +\dff 1}$\qss
and\dss hence\dss
$\num{\nsp\cap\dff{\gamma}}$\dss
is disjoint\dss from all $(\fff n\dff -\dff m\fff)$\dnsp-faces of\dss $Q^{\fff n\dff -\dff m\dff +\dff 1}$\dss
except\qss\vspace{1.5pt}
\[
\quad 
A_{\dff n\dff -\dff m\dff +\dff 1}\qff \cap\qff Q^{\fff n\dff -\dff m\dff +\dff 1}
\hspace*{1em}
\mbox{and}
\hspace*{1.2em} 
B_{\dff n\dff -\dff m\dff +\dff 1}\qff \cap\qff Q^{\fff n\dff -\dff m\dff +\dff 1}
\qff,
\]

\vspace{-10.5pt}
and\dss $\num{\cap\dff{\delta}}$\dss is disjoint\dss from\dss the whole
boundary\dss $\bd\fff Q^{\fff n\dff -\dff m\dff +\dff 1}$\dnsp.\oss 

Let\dss us\sss consider\dss now\dss the chains $\beta$
and $\delta_{\dff 0}$
from\dss Corollary\qss \ref{sides-and-chains}.\oss
The intersection\qss
$\delta_{\dff 0}\qff \cap\qff L^{n\dff -\dff m\dff +\dff 1}$\qss
is empty\sss and\dss hence\qss
$\cap\dff{\delta}_{\dff 0}\off =\off 0$\nnsp.\oss
By\dss intersecting\dss both sides of\qss (\ref{sides-decomposition})\qss
with\dss $L^{n\dff -\dff m\dff +\dff 1}$\qss we see\sss that\vspace{1.5pt}
\[
\quad 
\partial\dff (\dff \cap\dff{\gamma}_{\fff e}\dff)
\off =\off
\cap\qff\partial\dff \gamma_{\fff e}
\off =\dff\off
\cap\qff{\beta}
\pff +\pff 
\cap\dff{\delta}
\pff +\pff
\cap\dff{\delta}_{\dff 0}
\off =\off
\cap\qff{\beta}
\pff +\pff 
\cap\dff{\delta}
\qff.
\]

\vspace{-10.5pt}
The boundary of\dss every $1$\dnsp-chain\dss is a sum of\dss an even\dss
numbers of\dss points.\oss
It\dss follows\sss that\trs $\cap\dff{\delta}$\trs
is a sum of\dss an odd\dss number\sss of\dss points if\dss and\sss only\trs if\qss
$\cap\qff{\beta}$\trs is.\oss

Since\dss $\gamma$\dss is essential\halfff,\pss
$\num{\gamma}\qff \cap\qff L^{n\dff -\dff m}$\qss
consists of\dss an odd\dss number of\dss points.\oss
But\dss 
$\num{\gamma}\qff \cap\qff L^{n\dff -\dff m}$\dss is equal\dss to\qss
$\num{\nsp\cap\dff{\gamma}}\qff \cap\qff Q^{\fff n\dff -\dff m}$\dnsp.\oss
It\dss follows\dss that\dss the number of\sss $1$\dnsp-cubes of\sss $\cap\dff{\gamma}$\dss
intersecting\dss $Q^{\fff n\dff -\dff m}$\dss is odd.\oss
Let\dss $\eta$\dss be\sss the sum of\dss these $1$\dnsp-cubes and\dss $\theta$\dss
be\sss the sum of\dss $1$\dnsp-cubes 
of\dss $\cap\dff{\gamma}$\dss contained\dss in\vspace{3pt}
\[
\quad
Q^{\fff n\dff -\dff m}
\off =\off
[\trf 0\fff,\pff l\trf]^{n\dff -\dff m}\dff \times\qff
[\trf 0\dff,\pff a_{\dff n\dff -\dff m\dff +\dff 1}\trf]\qff \times\qff
(\dff 
a_{\dff n\dff -\dff m\dff +\dff 2}\dff,\pff
\ldots\dff,\pff
a_{\dff n}
\dff)
\qff.
\]

\vspace{-9pt}
By\dss the construction,\oss 
out\sss of\dss two endpoints of\dss a\sss
$1$\dnsp-cube of\dss $\eta$\dss
one is\sss in\dss $Q^{\fff n\dff -\dff m}$\dss
and\dss the other\sss is\sss not\halfff.\oss
Let\dss $\alpha$\dss be\sss the sum of\dss the endpoints which are not\dss in\dss
$Q^{\fff n\dff -\dff m}$\dnsp.\oss
The number of\dss them is
equal\dss to\sss the number of\dss $1$\dnsp-cubes of\dss $\eta$\dss
and\dss hence is odd\halfff.\oss
Since\dss $\gamma$\dss is a proper chain,\vspace{1.5pt}
\[
\quad
\partial\dff(\dff \eta\qff +\qff \theta\trf)
\off =\off
\alpha\qff +\qff \beta
\qff.
\]

\vspace{-10.5pt}
The chain\dss
$\partial\dff(\dff \eta\qff +\qff \theta\trf)$\dss
is a boundary\sss and\dss hence is a sum of\dss an even\dss number of\dss points.\oss
Since\dss $\alpha$\dss is a sum of\dss an odd\dss number of\dss points,\oss
it\dss follows\dss that\dss $\beta$\dss is a sum of\dss an odd\dss number of\dss
points.\oss
By\dss the previous paragraph,\oss this implies\dss that\trs $\cap\dff{\delta}$\trs
is a sum of\dss an odd\dss number of\dss points.\oss
Some of\dss these points cancel\dss under\dss the summation,\oss 
but\dss this does not\sss affect\dss
the parity of\dss their\dss number\halfff.\oss
Therefore\dss $\num{\nsp\cap\dff{\delta}}$\dss consists of\dss an odd\dss
number of\dss points and\dss $\delta$\dss is\dss
essential\halfff.\oss  \eproof

\prooftitle{Proof\qss based\sss on\dss projections}
Let\qss $p\off =\off p_{\fff m}$\nnsp.\oss
Since $\gamma$ is  essential,\pss 
$p_{\fff *}\dff(\trf \gamma\trf)
\off =\dff\off
\fclass{\fff P_{\dff m}\fff}$\qss
and\vspace{2.5pt}
\[
\quad
p_{\fff *}\dff(\trf \partial\dff \gamma\trf)
\off =\dff\off
\partial\dff p_{\fff *}\dff(\trf \gamma\trf)
\off =\dff\off
\partial\trf \fclass{\fff P_{\dff m}\fff}
\qff,
\]

\vspace{-9.5pt}
as follows\sss from\qss (\ref{cube-d-f}).\oss
Since\qss
$A_{\dff n\dff -\dff m\dff +\dff 1}
\off =\off
p^{\dff -\dff 1}\dff(\trf P_{\dff m\dff -\dff 1} \trf)$\nnsp,\oss
this implies\sss that\vspace{2.5pt}
\begin{equation*}
\quad
p_{\fff *}\dff(\trf \beta \trf)
\off =\dff\off
\fclass{\fff P_{\dff m\dff -\dff 1}\fff}
\qff,
\end{equation*}

\vspace{-9.5pt}
where
$\beta$
is\dss the chain\dss from\trs Corollary\qss \ref{sides-and-chains}.\oss 
Next\halfff,\oss let\qss $q\off =\off p_{\fff m\dff -\dff 1}$\nnsp.\oss
Then\qss $q\off =\off q\dff \circ\dff p$\qss
and\dss hence\vspace{2.5pt}
\[
\quad
q_{\fff *}\dff(\trf \beta \trf)
\off =\off
q_{\fff *}\dff(\trf p_{\fff *}\dff(\trf \beta\trf)\trf)
\off =\off
q_{\fff *}\dff\bigl(\trf \fclass{P_{\dff m\dff -\dff 1}}\trf\bigr)
\qff.
\] 

\vspace{-9.5pt}
Since\dss $q$\dss is equal\dss to\sss the identity on\dss $P_{\dff m\dff -\dff 1}$\nnsp,\oss
this implies\dss that\vspace{2.5pt}
\[
\quad
q_{\fff *}\dff(\trf \beta \trf)
\off =\off
\fclass{P_{\dff m\dff -\dff 1}}
\qff.
\]

\vspace{-9.5pt}
Therefore,\oss in order\dss to prove\sss that\dss $\delta$\dss is essential\dss
it\dss is sufficient\dss to prove\sss that\qss
$q_{\fff *}\dff(\trf \delta\trf)
\off =\off
q_{\fff *}\dff(\trf \beta\trf)$\nnsp.\vspace{2.25pt}

By\dss applying\dss $q_{\fff *}$\dss to\qss (\ref{sides-decomposition})\qss
we see\sss that\vspace{3pt}
\[
\quad
q_{\fff *}\dff(\trf \beta\trf)
\qff +\qff 
q_{\fff *}\dff(\trf \delta\trf)
\qff +\pff
q_{\fff *}\dff(\trf \delta_{\dff 0}\trf)
\off =\dff\off
q_{\fff *}\dff(\trf \partial\dff \gamma_{\dff E}\trf)
\off =\dff\off
\partial\dff p_{\fff *}\dff(\trf \gamma_{\dff E}\trf)
\qff.
\]

\vspace{-9pt}
But\dss $q$\dss takes every\dss $m$\dnsp-cube\sss to an\dss $m'$\dnsp-cube
with\qss $m'\qff \leq\qff m\qff -\qff 1$\nnsp.\oss
Since\dss $\gamma_{\dff E}$\dss is an\dss $m$\dnsp-chain,\oss
it\dss follows\dss that\qss
$q_{\fff *}\dff(\trf \gamma_{\dff E}\trf)
\off =\off
0$\qss
and\dss hence\qss
$\partial\dff q_{\fff *}\dff(\trf \gamma_{\dff E}\trf)
\off =\off
0$\nnsp.\oss
Therefore\vspace{3pt}
\[
\quad
q_{\fff *}\dff(\trf \beta\trf)
\qff +\qff 
q_{\fff *}\dff(\trf \delta\trf)
\qff +\qff
q_{\fff *}\dff(\trf \delta_{\dff 0}\trf)
\off =\off
0
\qff.
\]

\vspace{-9pt}
Also,\oss if\qss an\sss $(\fff m\dff -\dff 1\fff)$\dnsp-cube\sss $c$\dss
is contained\dss in\dss $A_{\dff k}$\qss or\qss $B_{\dff k}$\qss
with\qss $k\qff >\qff n\qff -\qff m\qff +\qff 1$\nnsp,\pss
then\dss $q\dff(\dff c\dff)$\dss
is\sss an\dss $m'$\dnsp-cube\dss
with\qss $m'\qff \leq\qff m\qff -\qff 2$\nnsp.\oss
It\dss follows\dss that\qss
$q_{\fff *}\dff(\trf \delta_{\dff 0}\trf)
\off =\off
0$\qss
and\dss therefore\vspace{3pt}
\[
\quad
q_{\fff *}\dff(\trf \beta\trf)
\qff +\qff
q_{\fff *}\dff(\trf \delta\trf)
\off =\off
0
\qff,
\]

\vspace{-9pt}
or\halfff,\oss equivalently,\oss 
$q_{\fff *}\dff(\trf \delta\trf)
\off =\off
q_{\fff *}\dff(\trf \beta\trf)$\nnsp.\oss
As we saw,\oss
this implies\sss that\dss $\delta$\dss is essential.\oss  \eproof\vspace{3pt}

\mypar{Theorem.}{e-coverings}
\emph{Suppose\dss that\pss
$e_{\dff 1}\fff,\pff e_{\dff 2}\dff,\pff \ldots\fff,\pff e_{\dff n\dff +\dff 1}$\qss
are cubical\sss sets covering\pss $Q$\dss and such\dss that}
\vspace{3.5pt}

\hspace*{1em}({\fff}i{\fff})\phantom{i{\dff}}\oss 
\emph{$e_{\dff 1}\qff \cup\qff e_{\dff 2}\qff \cup\qff \ldots\qff \cup\qff e_{\dff i}$\qss
contains\qss 
$A_{\dff i}$\qss
and\qss $e_{\dff i}$\dss 
is\dss disjoint\dss from\qss
$B_{\dff i}$\oss for\qss $i\qff \leq\qff n$\nnsp,\oss and}
\vspace{3.0pt}

\hspace*{1em}({\fff}i{\fff}i{\fff})\oss
\emph{$e_{\dff i}$\dss is\dss disjoint\qss from\qss
$A_{\fff j}$\pss if\oss $i\qff >\qff j$\nnsp.\oss}
\vspace{3.5pt}

\emph{Then\dss the\sss intersection\sss of\qss the\sss sets\pss
$e_{\dff 1}\fff,\pff e_{\dff 2}\dff,\pff \ldots\fff,\pff e_{\dff n\dff +\dff 1}$\qss
is\sss non-empty.\oss}\vspace{3pt}

\proof
We\sss will\dss recursively\dss construct\sss 
essential\dss $(\fff n\dff -\dff m\fff)$\dnsp-chains\dss $\gamma_{\fff m}$\dss in\dss $Q$\dss
such\dss that\qss\vspace{4.5pt}
\[
\quad 
\num{\gamma_{\fff m}}
\off \subset\off 
e_{\dff 1}\qff \cap\qff 
e_{\dff 2}\qff \cap\qff
\ldots
\qff \cap\qff
e_{\dff m}
\]

\vspace{-13.5pt}
and\vspace{0pt}
\begin{equation}
\label{other-e}
\quad
\num{\gamma_{\fff m}}
\off \subset\off 
e_{\dff m\dff +\dff 1}\qff \cup\qff 
e_{\dff m\dff +\dff 2}\qff \cup\qff
\ldots
\qff \cup\qff
e_{\dff n\dff +\dff 1}
\qff.
\end{equation}

\vspace{-7.5pt}
Let\qss $\gamma_{\dff 0}\off =\off \fclass{Q}$\nnsp,\oss
where\dss $\fclass{Q}$\dss is\dss the sum of\dss all\sss $n$\dnsp-cubes of\dss $Q$\nnsp.\oss
Clearly,\oss the $n$\dnsp-chain\dss $\gamma_{\dff 0}$\dss is essential\dss
and\sss the\sss two other conditions are vacuous 
for\qss $m\off =\off 0$\nnsp.\oss
Suppose\sss that\dss $\gamma_{\fff m}$\dss  
is already constructed\halfff.\oss
Let\dss us apply\qss Lebesgue\dss construction\dss to\qss
$\gamma\off =\off \gamma_{\fff m}$\qss
and\qss
$e\off =\off e_{\dff m\dff +\dff 1}$\nsp,\oss
and\dss let\dss $\delta$\dss be\sss the resulting $(\fff n\dff -\dff m\dff -\dff 1\fff)$\dnsp-chain.\oss
By\trs Lemma\qss \ref{proper}\qss the chain\dss $\delta$\dss is proper\halfff.\oss

Let\dss us check\dss the assumptions of\qss Lemma\qss \ref{special}.\oss
The assumption\dss ({\fff}i{\fff})\dss of\dss the\dss theorem\dss implies\sss that\qss
$e_{\dff m\dff +\dff 1}\qff \cap\qff B_{\dff m\dff +\dff 1}
\off =\off
\varnothing$\nnsp.\oss
The assumption\dss ({\fff}i{\fff}i{\fff})\dss of\dss the\dss theorem\dss implies\sss that\qss
$e_{\dff m\dff +\dff 2}\qff \cup\qff
\ldots
\qff \cup\qff
e_{\dff n\dff +\dff 1}$\qss
is disjoint\dss from\dss $A_{\dff m\dff +\dff 1}$\qss
and\dss hence\qss (\ref{other-e})\qss implies\sss that\qss
$\num{\gamma_{\fff m}}
\qff \cap\qff
A_{\dff m\dff +\dff 1}$\qss
is contained\dss in\dss
$e_{\dff m\dff +\dff 1}$\nsp.\oss
Therefore\sss the assumptions of\qss Lemma\qss \ref{special}\qss
hold\sss and\dss this lemma implies\sss the chain\dss $\delta$\dss is special\halfff.\oss 
Now we can apply\qss Lemma\qss \ref{essential}\qss and conclude\sss that\dss
the chain\dss $\delta$\dss is essential\halfff.\oss

By\dss the construction\dss $\num{\delta}$\dss is contained\dss in\qss
$\num{\gamma_{\fff m}}\qff \cap\qff e_{\dff m\dff +\dff 1}$\qss
and\dss hence\qss\vspace{2.25pt}
\[
\quad 
\num{\delta}
\off \subset\off 
e_{\dff 1}\qff \cap\qff 
e_{\dff 2}\qff \cap\qff
\ldots
\qff \cap\qff
e_{\dff m\dff +\dff 1}
\qff.
\]

\vspace{-9.75pt}
Also by\dss the construction,\oss every cube of\dss $\delta$\dss
is a common face of\dss two cubes of\dss $\gamma_{\dff m}$\nnsp,\oss
one of\dss which is contained\dss in\dss $e_{\dff m\dff +\dff 1}$\dss
and\dss the other does not\halfff.\oss
By\qss (\ref{other-e})\qss that\sss other cube is contained\dss in\qss
$e_{\dff m\dff +\dff 2}\qff \cup\qff
\ldots
\qff \cup\qff
e_{\dff n\dff +\dff 1}$\nsp.\oss
It\dss follows\sss that\vspace{2.25pt}
\[
\quad 
\num{\delta}
\off \subset\off 
e_{\dff m\dff +\dff 2}\qff \cup\qff
\ldots
\qff \cup\qff
e_{\dff n\dff +\dff 1}
\qff.
\]

\vspace{-10pt}
Therefore we can set\qss
$\gamma_{\fff m\dff +\dff 1}
\off =\off
\delta$\nnsp.\oss
This completes\sss the construction of\dss the chains\dss $\gamma_{\dff m}$\nsp.\oss

Let\dss us consider\dss the last\sss of\dss these chains,\oss
namely\dss $\gamma_{\fff n}$\nsp.\oss
The support\dss $\num{\gamma_{\fff n}}$\dss 
is contained\dss in\dss the intersection\qss
$e_{\dff 1}\qff \cap\qff 
e_{\dff 2}\qff \cap\qff
\ldots
\qff \cap\qff
e_{\dff n}$\qss
and\dss by\qss (\ref{other-e})\qss also\dss in $e_{\dff n\dff +\dff 1}$\nsp.\oss
Since\dss $\gamma_{\fff n}$\dss is essential\sss and\dss hence\qss
$\num{\gamma_{\fff n}}\off \neq\off \varnothing$\nsp,\oss
this implies\sss that\qss
$e_{\dff 1}\qff \cap\qff 
e_{\dff 2}\qff \cap\qff
\ldots
\qff \cap\qff
e_{\dff n\dff +\dff 1}
\off \neq\off
\varnothing$\nnsp.\oss  \eproof

\myuppar{Lebesgue fusion of\dss sets.}
Let\qss
$d_{\dff 1}\fff,\pff d_{\dff 2}\dff,\pff \ldots\fff,\pff d_{\dff r}$\qss
be a covering\dss of\pss $Q$\dss
by several\dss sets\qss
(in\dss the applications\sss these sets are closed or even cubical\fff).\oss
Suppose\dss that\dss none of\qss the sets\dss $d_{\dff i}$\dss intersects\sss
two opposite\dss $(\fff n\dff -\dff 1\fff)$\dnsp-faces of\pss $Q$\nnsp.\oss
An elementary\sss construction of\qss Lebesgue\dss allows\sss to
construct\sss sets\qss
$e_{\dff 1}\fff,\pff e_{\dff 2}\fff,\pff \ldots\fff,\pff e_{\dff n\dff +\dff 1}$\qss
covering\dss $Q$\dss and\sss
satisfying\dss the conditions\qss ({\fff}i{\fff})\qss and\qss ({\fff}i{\fff}i{\fff})\qss
of\qss Theorem\qss \ref{e-coverings}\qss
by\dss taking\sss unions of\dss disjoint\sss
collections of\dss sets\dss $d_{\dff i}$\nsp.\oss  

Let\dss $e_{\dff 1}$\dss be\sss the union of\dss all sets\dss $d_{\fff j}$\dss
having non-empty intersection with\dss $A_{\dff 1}$\nnsp.\oss
Then\dss $e_{\dff 1}$\dss is disjoint\dss from\dss $B_{\dff 1}$\nnsp.\oss
Let\dss $e_{\dff 2}$\dss be\sss the union of\dss all sets\dss $d_{\fff j}$\dss
which have non-empty intersection with\dss $A_{\dff 2}$\dss
and\dss were\sss not\dss used\dss to form\dss $e_{\dff 1}$\nnsp.\oss
Since\sss the sets\dss $d_{\dff i}$\dss form a covering of\dss $Q$\nnsp,\oss
the set\dss $A_{\dff 2}$\dss
is contained\dss in\qss
$e_{\dff 1}\qff \cup\qff e_{\dff 2}$\nsp.\oss
Clearly,\oss
$e_{\dff 2}$\dss is disjoint\dss from\dss $A_{\dff 1}$\dss and\trs $B_{\dff 2}$\nsp.\oss

Suppose\sss that\dss $k\qff \leq\qff n\qff -\qff 1$\nnsp,\oss 
the sets\qss
$e_{\dff 1}\fff,\pff e_{\dff 2}\fff,\pff \ldots\fff,\pff e_{\dff k}$\qss
are already\sss defined and\dss
satisfy\dss the conditions\qss ({\fff}i{\fff})\qss and\qss ({\fff}i{\fff}i{\fff})\qss
of\qss Theorem\qss \ref{e-coverings}\qss for\qss $i\qff \leq\qff k$\nnsp.\oss 
Let\dss $e_{\dff k\dff +\dff 1}$\dss be\sss the union of\dss all sets\dss $d_{\fff j}$\dss
which have non-empty intersection with\dss $A_{\dff k\dff +\dff 1}$\dss
and\dss were\sss not\dss used\dss to form sets\qss 
$e_{\dff 1}\fff,\pff e_{\dff 2}\fff,\pff \ldots\fff,\pff e_{\dff k}$\nsp.\oss
By\dss the construction,\oss $e_{\dff k\dff +\dff 1}$\dss
satisfies\dss the conditions\qss ({\fff}i{\fff})\qss and\qss ({\fff}i{\fff}i{\fff})\qss
of\qss Theorem\qss \ref{e-coverings}\qss for\qss $i\off =\off k\qff +\qff 1$\nnsp.\oss 

Finally,\oss let\dss $e_{\dff n\dff +\dff 1}$\sss be\sss the union
of\dss the sets\dss $d_{\fff j}$\dss not\dss used\dss to form\qss
$e_{\dff 1}\fff,\pff e_{\dff 2}\fff,\pff \ldots\fff,\pff e_{\dff n}$\nnsp.\oss
The condition\qss ({\fff}i{\fff})\qss is vacuous for\qss $i\off =\off n\qff +\qff 1$\nnsp,\oss
and\trs $e_{\dff n\dff +\dff 1}$\dss obviously satisfies\sss the condition\qss ({\fff}i{\fff}i{\fff})\qss
for\qss $i\off =\off n\qff +\qff 1$\nnsp.\oss

\mypar{Theorem.}{collecting-sets}
\emph{Let\qss
$d_{\dff 1}\fff,\pff d_{\dff 2}\dff,\pff \ldots\fff,\pff d_{\dff r}$\qss
be a covering\dss of\pss $Q$\dss
by cubical\dss sets.\oss
Suppose\dss that\dss none of\qss the sets\dss $d_{\dff i}$\dss intersects\sss
two opposite\dss $(\fff n\dff -\dff 1\fff)$\dnsp-faces of\pss $Q$\nnsp.\oss
Then\dss among\dss the sets\qss $d_{\dff i}$\qss there are\qss 
$n\qff +\qff 1$\qss sets with non-empty\dss intersection.\oss}

\proof
Let\qss
$e_{\dff 1}\fff,\pff e_{\dff 2}\fff,\pff \ldots\fff,\pff e_{\dff n\dff +\dff 1}$\qss
be\sss the sets obtained\dss from\qss
$d_{\dff 1}\fff,\pff d_{\dff 2}\dff,\pff \ldots\fff,\pff d_{\dff r}$\qss
by\qss Lebesgue\dss fusion\sss construction.\oss
Then\dss the assumptions of\qss Theorem\qss \ref{e-coverings}\qss hold
and\dss hence\qss
$e_{\dff 1}\qff \cap\qff 
e_{\dff 2}\qff \cap\qff
\ldots
\qff \cap\qff
e_{\dff n\dff +\dff 1}$\qss
is non-empty\halfff.\oss
Since\sss the sets\dss $e_{\dff i}$\dss
are unions of\dss disjoint\sss collections of\dss sets\dss $d_{\fff j}$\nsp,\oss
this\sss implies\sss that\dss 
among\dss the sets\qss $d_{\dff i}$\qss there are\qss 
$n\qff +\qff 1$\qss sets with non-empty\dss intersection.\oss  \eproof

\myuppar{Lebesgue\dss first\dss covering\dss theorem.} 
\emph{Let\qss
$D_{\dff 1}\fff,\pff D_{\dff 2}\dff,\pff \ldots\fff,\pff D_{\dff r}$\qss
be a covering\dss of\pss the unit\sss cube\qss
$[\trf 0\fff,\qff 1\trf]^{n}$\qss
by closed\dss sets.\oss
Suppose\dss that\dss none of\qss the sets\dss $D_{\dff i}$\sss intersects\sss
two opposite\sss $(\fff n\dff -\dff 1\fff)$\dnsp-faces of\qss $[\trf 0\fff,\qff 1\trf]^{n}$\dnsp.\oss
Then\dss among\dss the sets\qss $D_{\dff i}$\qss there are\qss 
$n\qff +\qff 1$\qss sets with non-empty\dss intersection.\oss}

\proof
Let\qss $\varepsilon\qff >\qff 0$\qss be a\dss Lebesgue\dss number of\dss the covering\sss
of\qss the unit\sss cube\dss $[\trf 0\fff,\qff 1\trf]^{n}$\dss by\dss the sets\dss
$D_{\dff 1}\fff,\pff D_{\dff 2}\dff,\pff \ldots\fff,\pff D_{\dff r}$\dss
and $(\fff n\dff -\dff 1\fff)$\dnsp-faces of\qss $[\trf 0\fff,\qff 1\trf]^{n}$\dnsp.\oss
Suppose\sss that\dss the number\dss $l$\dss is chosen\dss to\sss be so large\sss that\dss
an $n$\dnsp-dimensional\sss cube in\dss $\rrr^{\fff n}$\dss with\dss
the sides of\dss the length\dss $1/l$\dss has diameter\qss
$<\qff \varepsilon$\nnsp.\oss
Let\qss 
$p\dff \colon\dff
Q\qff \ttoo\qff [\trf 0\fff,\qff 1\trf]^{n}$\qss
be\sss the map defined\dss by\vspace{0pt}
\[
\quad
p\dff(\dff x_{\dff 1}\fff,\pff x_{\dff 2}\fff,\pff \ldots\fff,\pff x_{\dff n}\dff)
\off =\off
(\dff x_{\dff 1}/l\fff,\off x_{\dff 2}/l\fff,\off \ldots\fff,\off x_{\dff n}/l\qff)
\qff.
\]

\vspace{-12pt}
For each\qss $i\qff \leq\qff r$\qss let\dss $d_{\dff i}$\dss
be\sss the union of\dss all\sss $n$\dnsp-cubes of\dss $Q$\dss intersecting\dss
the preimage\dss
$p^{\dff -\dff 1}\dff (\trf D_{\dff i}\trf)$\nnsp.\oss
The sets\qss
$d_{\dff 1}\fff,\pff d_{\dff 2}\dff,\pff \ldots\fff,\pff d_{\dff r}$\qss
are cubical sets.\oss
Obviously,\oss they cover\dss $Q$\nnsp.\oss
If\dss a cube $d_{\dff i}$\dss intersects an $(\fff n\dff -\dff 1\fff)$\dnsp-face of\dss $Q$\nnsp,\oss
then\dss $D_{\dff i}$\dss intersects\sss the image of\dss this face in\dss
$[\trf 0\fff,\qff 1\trf]^{n}$\dss
({\fff}in\dss fact\halfff,\oss the $(\fff n\dff -\dff 1\fff)$\dnsp-faces of\qss
$[\trf 0\fff,\qff 1\trf]^{n}$\dss were included\dss in\dss the covering 
defining\dss $\varepsilon$\dss precisely\dss to ensure\sss this property\fff).\oss
It\dss follows\dss that\dss
none of\dss the sets\dss $d_{\dff i}$\dss  intersects\sss
two opposite\dss $(\fff n\dff -\dff 1\fff)$\dnsp-faces\dss of\dss $Q$\nnsp.\oss
By\trs Theorem\qss \ref{collecting-sets}\qss among\dss these cubical sets\sss there 
are\dss $n\qff +\qff 1$\trs sets with\dss non-empty\sss intersection.\oss
If\halfff,\oss say\halfff,\pss
$d_{\dff 1}\qff \cap\qff 
d_{\dff 2}\qff \cap\qff
\ldots
\qff \cap\qff
d_{\dff n\dff +\dff 1}$\qss
is\dss non-empty\dss
and\pss
$x\qff \in\qff
d_{\dff 1}\qff \cap\qff 
d_{\dff 2}\qff \cap\qff
\ldots
\qff \cap\qff
d_{\dff n\dff +\dff 1}$\nnsp,\oss
then\dss the distance of\trs the point\qss 
$p\dff(\dff x\trf)\qff \in\qff [\trf 0\fff,\qff 1\trf]^{n}$\qss
from each of\dss the sets\qss
$D_{\dff 1}\fff,\pff D_{\dff 2}\dff,\pff \ldots\fff,\pff D_{\dff n\dff +\dff 1}$
is\qss $<\qff \varepsilon$\nnsp.\oss
Now\trs Lebesgue\dss lemma for closed sets implies\sss that\dss the intersection\qss
$D_{\dff 1}\qff \cap\qff 
D_{\dff 2}\qff \cap\qff
\ldots
\qff \cap\qff
D_{\dff n\dff +\dff 1}$\qss
is\dss non-empty\halfff.\oss  \eproof \vspace{0pt}

\myuppar{Lebesgue\dss second\dss covering\dss theorem.} 
\emph{If\oss $\varepsilon\trf >\trf 0$\qss is\dss sufficiently\sss small\halfff,\oss
then every\dss closed\dss $\varepsilon$\dnsp-covering\dss of\oss $[\trf 0\fff,\qff 1\trf]^{n}$\dss
has\dss order\qss $\geq\qff n\qff +\qff 1$\nnsp.} \vspace{1pt}

\proof
It\dss is\dss sufficient\dss to note\sss that\dss 
no set\sss of\dss diameter\qss 
$<\qff 1$\qss can\dss intersect\dss two opposite
\dss $(\fff n\dff -\dff 1\fff)$\dnsp-faces of\pss $[\trf 0\fff,\qff 1\trf]^{n}$\qss
and\dss hence one can simply\dss take\qss $\varepsilon\off =\off 1$\nnsp.\oss  
\eproof \vspace{0pt}

\myuppar{Further\dss applications.}
The above\dss results\sss are contained\dss in\qss Lebesgue\dss
papers\qss \cite{l-one},\oss \cite{l-two}\qss
either explicitly\sss or\dss implicitly.\oss
Now\dss we will\dss present\sss some further\sss applications of\qss Theorem\qss \ref{e-coverings}.\oss

\mypar{Theorem.}{cubes-partitions}
\emph{Let\pss
$d_{\dff 1}\fff,\pff d_{\dff 2}\fff,\pff \ldots\fff,\pff d_{\dff n\dff +\dff 1}$\qss
be cubical\sss sets covering\dss $Q$\dss and\dss such\dss that\dss the union\qss
$d_{\dff 1}\qff \cup\qff d_{\dff 2}\qff \cup\qff \ldots\qff \cup\qff d_{\dff i}$\qss
contains\qss 
$A_{\dff i}$\qss
and\qss $d_{\dff i}$\dss 
is\dss disjoint\dss from\qss
$B_{\dff i}$\oss for every\qss $i\qff \leq\qff n$\nnsp.\oss 
Then\dss the intersection\qss
$d_{\dff 1}\qff \cap\qff 
d_{\dff 2}\qff \cap\qff
\ldots
\qff \cap\qff
d_{\dff n\dff +\dff 1}$\qss
is\dss non-empty\halfff.\oss}

\proof
For\qss $i\qff \leq\qff n\qff +\qff 1$\qss
let\dss $e_{\dff i}$\dss be\sss the union of\dss all\sss $n$\dnsp-cubes contained\dss in\dss
$d_{\dff i}$\nnsp,\oss
but\dss not\dss 
in any\dss $d_{\fff j}$\dss with\qss
$j\qff <\qff i$\nnsp.\oss
Then\dss $e_{\dff i}$\dss is disjoint\dss from\dss $A_{\dff j}$\qss if\oss
$i\qff >\qff j$\nnsp,\oss and\qss\vspace{0pt} 
\[
\quad
e_{\dff 1}\qff \cup\qff e_{\dff 2}\qff \cup\qff \ldots\qff \cup\qff e_{\dff i}
\off\qff =\off\qff
d_{\dff 1}\qff \cup\qff d_{\dff 2}\qff \cup\qff \ldots\qff \cup\qff d_{\dff i}
\]

\vspace{-12pt}
for every\qss $i\qff \leq\qff n$\nnsp.\oss
It\dss follows\dss that\qss
$e_{\dff 1}\fff,\pff e_{\dff 2}\dff,\pff \ldots\fff,\pff e_{\dff n\dff +\dff 1}$\qss
satisfy\dss the assumptions of\qss Theorem\qss \ref{e-coverings}\qss
and\dss hence\sss the intersection\qss 
$e_{\dff 1}\qff \cap\qff 
e_{\dff 2}\qff \cap\qff
\ldots
\qff \cap\qff
e_{\dff n\dff +\dff 1}$\qss
is\dss non-empty\halfff.\oss
Since\qss $e_{\dff i}\qff \subset\qff d_{\dff i}$\qss for every\dss $i$\nnsp,\oss
it\dss follows\dss that\dss the intersection\qss
$d_{\dff 1}\qff \cap\qff 
d_{\dff 2}\qff \cap\qff
\ldots
\qff \cap\qff
d_{\dff n\dff +\dff 1}$\qss
is\dss non-empty\halfff.\oss \eproof \vspace{1pt}

\mypar{Theorem.}{closed-partitions}
\emph{Let\pss
$D_{\dff 1}\fff,\pff D_{\dff 2}\fff,\pff \ldots\fff,\pff D_{\dff n\dff +\dff 1}$\qss
be closed\sss sets covering\dss $[\trf 0\fff,\qff 1\trf]^{n}$\dss and\dss such\dss that\dss the union\qss
$D_{\dff 1}\qff \cup\qff D_{\dff 2}\qff \cup\qff \ldots\qff \cup\qff D_{\dff i}$\qss
contains\qss 
$A_{\dff i}$\qss
and\qss $D_{\dff i}$\dss 
is\dss disjoint\dss from\qss
$B_{\dff i}$\oss for every\qss $i\qff \leq\qff n$\nnsp.\oss 
Then\dss the intersection\qss
$D_{\dff 1}\qff \cap\qff 
D_{\dff 2}\qff \cap\qff
\ldots
\qff \cap\qff
D_{\dff n\dff +\dff 1}$\qss
is\dss non-empty\halfff.\oss} \vspace{1pt}

\proof
The proof\trs is\dss completely\dss similar\dss to\sss the proof\dss of\qss Lebesgue\dss
first\dss covering\dss theorem,\oss
with\qss Theorem\qss \ref{cubes-partitions}\qss playing\dss the role of\qss
Theorem\qss \ref{collecting-sets}.\oss  \eproof \vspace{0pt}

\mypar{Theorem.}{early-h-th}
\emph{Let\pss
$D_{\dff 1}\fff,\pff D_{\dff 2}\fff,\pff \ldots\fff,\pff D_{\dff n\dff +\dff 1}$\qss
be closed\sss sets covering\dss $[\trf 0\fff,\qff 1\trf]^{n}$\dss and\dss such\dss that\dss
the set\dss 
$D_{\dff k}$\dss 
is\dss disjoint\dss from\dss
$A_{\dff i}$\dss with\qss $i\qff <\qff k$\qss
for every\qss $k\qff \leq\qff n\qff +\qff 1$\qss 
and\dss is\dss disjoint\dss from\dss
$B_{\dff k}$\dss for every\qss $k\qff \leq\qff n$\nnsp.\oss
Then\dss the intersection\qss
$D_{\dff 1}\qff \cap\qff 
D_{\dff 2}\qff \cap\qff
\ldots
\qff \cap\qff
D_{\dff n\dff +\dff 1}$\qss
is\dss non-empty\halfff.\oss}

\proof
For every\dss $i\qff \leq\qff n$\nnsp,\pss the union\dss
$D_{\dff 1}\qff \cup\qff
D_{\dff 2}\qff \cup\qff
\ldots\qff \cup\qff
D_{\dff n\dff +\dff 1}$\dss
contains\sss $A_{\dff i}$\sss and\dss hence\sss
the assumptions of\dss the\sss theorem\dss imply\dss that\dss
$D_{\dff 1}\qff \cup\qff D_{\dff 2}\qff \cup\qff \ldots\qff \cup\qff D_{\dff i}$\dss
contains\dss 
$A_{\dff i}$\nsp.\pss
Hence\dss Theorem\dss \ref{closed-partitions}\qss applies\sss to\sss the sets\dss
$D_{\dff 1}\fff,\qff
D_{\dff 2}\fff,\qff
\ldots\fff,\qff
D_{\dff n\dff +\dff 1}$\sss
and\sss implies\sss that\dss their\dss intersection\sss is\sss non-empty\halfff.\oss \eproof

\myuppar{Separation of\dss faces.}
Let\qss $A\fff,\pff B$\qss be\sss two disjoint\sss subsets of\dss
a closed subset\qss $X\qff \subset\qff \rrr^{\fff n}$\qss
(or\sss a\sss topological\sss space\dss $X$\nsp).\oss
A closed subset\dss $C\qff \subset\qff X$\qss 
is\dss said\dss to be a\qss \emph{partition\dss between\qss $A$\qss and\dss $B$\dss
in}\qss $X$\qss if\qss $A$\dss and\dss $B$\dss
are contained\dss in\dss two different\sss components
of\trs the complement\qss
$X\qff \smallsetminus\qff C$\nnsp.\oss
The next\dss theorem\dss 
is\dss a\sss basic result\sss of\trs dimension\dss theory.\oss
See,\oss for example,\oss \cite{hw},\oss Section\qss IV.{\dff}1.\oss
It\dss is\dss concerned\dss with\dss  
sets separating\dss $A_{\dff i}$\dss from\dss $B_{\dff i}$\qss
in\qss $Q$\dss 
in\dss the case\qss $l\off =\off 1$\nnsp,\oss
i.e.\qss 
$Q\off =\off [\trf 0\fff,\qff 1\trf]^{n}$\dnsp.\oss \vspace{1pt}

\myuppar{Theorem about\dss partitions.}
\emph{Let\pss
$C_{\dff 1}\fff,\pff C_{\dff 2}\fff,\pff \ldots\fff,\pff C_{\dff n}$\qss
be closed\sss subsets of\pss $[\trf 0\fff,\qff 1\trf]^{n}$\dss such\dss that\dss 
for every\dss $i$\dss the set\dss $C_{\dff i}$\dss
is\dss a\dss partition\dss between\dss $A_{\dff i}$\qss and\trs $B_{\dff i}$\nsp.\qff\oss 
Then\oss
$C_{\dff 1}\qff \cap\qff 
C_{\dff 2}\qff \cap\qff
\ldots
\qff \cap\qff
C_{\dff n}
\off \neq\off
\varnothing$\nnsp.\oss} \vspace{1pt}

\proof
Let\dss $U_{\fff i}$\dss be\sss the component\sss of\qss
$[\trf 0\fff,\qff 1\trf]^{n}\pff \smallsetminus\qss C_{\dff i}$\qss
containing\dss $A_{\dff i}$\nsp,\oss
and\dss let\dss $D_{\dff i}$\dss be its closure.\oss
Let\dss $E_{\dff i}$\dss be\sss the closure of\qss
$[\trf 0\fff,\qff 1\trf]^{n}\pff \smallsetminus\qss D_{\dff i}$\nsp.\oss
Then\dss $B_{\dff i}$\dss is contained\dss in\dss $E_{\dff i}$\nsp.\oss
If\qss
$x\qff \in\qff
[\trf 0\fff,\qff 1\trf]^{n}$\qss  
is\sss not\dss in\dss $C_{\dff i}$\nsp,\oss
then\sss either\qss $x\qff \in\qff U_{\fff i}$\nsp,\oss
or\dss $x$\dss belongs\sss to some other component\sss of\qss 
$[\trf 0\fff,\qff 1\trf]^{n}\pff \smallsetminus\qss C_{\dff i}$\nsp.\oss
In\dss the first\sss case\qss
$x\qff \not\in\qff E_{\dff i}$\nsp,\oss
in\dss the second case\qss
$x\qff \not\in\qff D_{\dff i}$\nsp.\oss
It\dss follows\dss that\qss 
$D_{\dff i}\qff \cap\qff E_{\dff i}
\off \subset\off
C_{\dff i}$\nsp.\oss 
Let\oss
$D_{\dff n\dff +\dff 1}
\off =\off
E_{\dff 1}\qff \cap\qff 
E_{\dff 2}\qff \cap\qff
\ldots
\qff \cap\qff
E_{\dff n}$\nsp.\oss
The sets\qss
$D_{\dff 1}\fff,\pff D_{\dff 2}\fff,\pff \ldots\fff,\pff D_{\dff n\dff +\dff 1}$\qss
form a closed\sss covering of\dss
$[\trf 0\fff,\qff 1\trf]^{n}$\dss
satisfying\dss the assumptions of\qss Theorem\qss \ref{closed-partitions}\qss
and\dss hence\sss 
their\sss intersection\dss is\sss non-empty.\oss
Let\dss $x$\dss be a point\dss in\dss the intersection of\dss these sets.\oss
Then\oss\vspace{4pt}
\[
\quad
x\qff \in\qff
D_{\dff i}\qff \cap\qff D_{\dff n\dff +\dff 1}
\off \subset\off
D_{\dff i}\qff \cap\qff E_{\dff i}
\off \subset\off
C_{\dff i} 
\]

\vspace{-8pt}
for every\dss $i$\nnsp.\oss
It\dss follows\dss that\qss
$C_{\dff 1}\qff \cap\qff 
C_{\dff 2}\qff \cap\qff
\ldots
\qff \cap\qff
C_{\dff n}
\off \neq\off
\varnothing$\nnsp.\oss  \eproof \vspace{0pt}

\myuppar{Two cartesian\sss versions of\pss KKM\dss argument\halfff.}
Suppose\sss that\qss
$f\dff \colon\dff
[\trf 0\fff,\qff 1\trf]^{n}\qff \ttoo\qff [\trf 0\fff,\qff 1\trf]^{n}$\qss
is\sss a continuous map without\dss fixed\dss points.\oss
By\sss slightly\dss perturbing\dss $f$\dnsp,\oss if\dss necessary\halfff,\oss
we may assume\sss that\dss the image of\dss $f$\dss is
disjoint\dss from\dss 
$\bd\dff [\trf 0\fff,\qff 1\trf]^{n}$\nnsp.\oss
For\qss 
$x\qff \in\qff [\trf 0\fff,\qff 1\trf]^{n}$\pss 
let\qss 
$(\dff x_{\dff 1}\fff,\pff x_{\dff 2}\fff,\pff \ldots\fff,\pff x_{\dff n}\dff)$\qss
and\qss
$(\dff y_{\dff 1}\fff,\pff y_{\dff 2}\fff,\pff \ldots\fff,\pff y_{\dff n}\dff)$\qss
be\sss the usual coordinates of\dss  
$x$\dss and\qss
$y\off =\off f\dff(\dff x\trf)$\qss respectively.

\prooftitle{Deducing\dss Brouwer\dss fixed\dss point\dss theorem\dss from\qss
Theorem\qss \ref{early-h-th}}
For each\qss
$i\off =\off 1\fff,\pff 2\fff,\pff \ldots\fff,\pff n$\qss
let\qss $D_{\dff i}\off \subset\off [\trf 0\fff,\qff 1\trf]^{n}$\qss
be\sss the set\sss of\dss points\dss $x$\dss such\dss that\qss
$y_{\dff i}\qff \geq\qff x_{\dff i}$\nnsp.\oss
Let\qss $D_{\dff n\dff +\dff 1}\off \subset\off [\trf 0\fff,\qff 1\trf]^{n}$\qss
be\sss the set\sss of\dss points\dss $x$\dss such\dss that\qss
$y_{\dff i}\qff \leq\qff x_{\dff i}$\qss
for every\qss
$i\off =\off 1\fff,\pff 2\fff,\pff \ldots\fff,\pff n$\nnsp.\oss

If\qss $x\qff \in\qff A_{\dff i}$\nnsp,\oss
then\qss $x_{\dff i}\off =\off 0$\qss
and\qss
$y_{\dff i}\qff >\qff 0\off =\off x_{\dff i}$\qss
because\sss the image of\dss $f$\dss is disjoint\dss from\dss $A_{\dff i}$\dss
by\dss the above assumption.\oss
It\dss follows\dss that\qss
$A_{\dff i}$\dss is\dss disjoint\dss from\dss $D_{\dff n\dff +\dff 1}$\dss
for\sss every\qss $i\qff \leq\qff n$\nnsp.\oss
If\qss $x\qff \in\qff B_{\dff i}$\nnsp,\oss
then\qss $x_{\dff i}\off =\off 1$\qss
and\qss 
$y_{\dff i}\off <\off 1\off =\off x_{\dff i}$\qss
because\dss $B_{\dff i}$\dss
is disjoint\dss from\dss the image of\dss $f$\sss by\dss the above assumption.\oss
It\dss follows\dss that\dss $D_{\dff i}$\dss
is disjoint\dss from\dss $B_{\dff i}$\nnsp.\oss

Therefore\qss Theorem\qss \ref{early-h-th}\qss applies\sss
and\dss hence\sss there exists a point\qss
$x\qff \in\qff 
D_{\dff 1}\qff \cap\qff 
D_{\dff 2}\qff \cap\qff
\ldots
\qff \cap\qff
D_{\dff n\dff +\dff 1}$\nnsp.\oss
Now\qss $x\qff \in\qff D_{\dff i}$\qss
implies\sss that\qss
$y_{\dff i}\qff \geq\qff x_{\dff i}$\qss
and\qss $x\qff \in\qff D_{\dff n\dff +\dff 1}$\qss
implies\sss that\qss
$y_{\dff i}\qff \geq\qff x_{\dff i}$\qss
for every\dss
$i$\nnsp.\oss
It\dss follows\dss that\qss 
$y\off =\off x$\nnsp,\oss
i.e.\qss $x$\dss is a fixed\dss point\sss of\dss $f$\nsp\dnsp.\oss  \eproof

\prooftitle{Deducing\dss Brouwer\dss fixed\dss point\dss theorem\dss from\qss
Theorem\sss about\dss partitions}
Let\dss $C_{\dff i}$\dss be\sss the set\sss of\dss points\dss $x$\dss such\dss that\qss
$y_{\dff i}\off =\off x_{\dff i}$\nnsp.\oss
The complement\qss $[\trf 0\fff,\qff 1\trf]^{n}\qff \smallsetminus\qff C_{\dff i}$\qss
is equal\dss to\sss the union of\dss two disjoint\sss open sets,\oss
one defined\dss by\dss the inequality\qss 
$y_{\dff i}\qff >\qff x_{\dff i}$\qss
and\dss another\dss by\dss the inequality\qss
$y_{\dff i}\qff <\qff x_{\dff i}$\nnsp.\oss
Since\sss the image of\dss $f$\dss is\sss disjoint\dss
from\dss the boundary\dss $\bd\dff [\trf 0\fff,\qff 1\trf]^{n}$\dnsp,\oss
the face\dss $A_{\dff i}$\dss is contained\dss in\dss the first\sss of\dss
these open sets and\dss the face\dss $B_{\dff i}$\dss in\dss the second.\oss
Therefore\dss $C_{\dff i}$\dss is\dss a\sss partition\dss between\dss
$A_{\dff i}$\dss and\dss $B_{\dff i}$\dss and\dss hence\qss
$C_{\dff 1}\qff \cap\qff 
C_{\dff 2}\qff \cap\qff
\ldots
\qff \cap\qff
C_{\dff n}$\qss
is non-empty.\oss
If\dss $x$\dss belongs\sss to\sss this intersection,\oss
then\qss $y_{\dff i}\off =\off x_{\dff i}$\qss for all\dss $i$\dss
and\dss hence\qss $y\off =\off x$\nnsp,\oss
i.e.\qss $f\dff(\dff x\dff)\off =\off x$\nnsp.\oss  \eproof

\myuppar{The unit\dss cube and\dss the standard simplex.}
Let\dss $\delta$\dss be\sss the standard $n$\dnsp-simplex in\dss $\rrr^{\fff n\dff +\dff 1}$
defined\dss in\dss the standard coordinates\qss
$y_{\dff 0}\fff,\pff
y_{\dff 1}\fff,\pff
\ldots\fff,\pff
y_{\dff n}$\qss
by\dss the conditions\vspace{4.5pt}
\[
\quad
y_{\dff 0}\qff +\qff
y_{\dff 1}\qff +\qff
\ldots\qff +\qff
y_{\dff n}
\off =\off
1
\hspace{1.2em}\mbox{and}\hspace{1.2em}
y_{\dff 0}\fff,\pff
x_{\dff 1}\fff,\pff
\ldots\fff,\pff
y_{\dff n}
\off \geq\off
0
\qff,
\]

\vspace{-7.5pt}
and\dss let\dss $\delta_{\dff i}$\dss be\sss the face of\dss $\delta$\dss
defined\dss by\dss the equation\qss $y_{\dff i}\off =\off 0$\nnsp,\oss
where\qss $i\off =\off 0\fff,\pff 1\fff,\pff \ldots\fff,\pff n$\nnsp.\oss
Let\qss\vspace{3pt} 
\[
\quad
\lambda\dff \colon\dff 
[\trf 0\fff,\qff 1\trf]^{n}
\qff \ttoo\qff
\delta
\]

\vspace{-9pt}
be\sss the map defined\dss by\qss
$\lambda\dff
(\dff
x_{\dff 1}\fff,\pff
x_{\dff 2}\fff,\pff
\ldots\fff,\pff
x_{\dff n}
\dff)
\off =\off 
(\dff
y_{\dff 0}\fff,\pff
y_{\dff 1}\fff,\pff
\ldots\fff,\pff
y_{\dff n}
\dff)$\nnsp,\oss 
where\vspace{3pt}
\[
\quad
y_{\dff 0}
\off =\off
1\qff -\qff x_{\dff 1}
\qff,
\]

\vspace{-36pt}
\[
\quad
y_{\dff i}
\off =\off
(\qff 1\qff -\qff x_{\dff i\dff +\dff 1}\qff)\qff
x_{\dff 1}\qff x_{\dff 2}\pff \ldots\pff x_{\dff i}
\hspace*{1.2em}\mbox{for}\hspace*{1.2em}
0\off <\off i\off <\off n
\qff,
\]

\vspace{-36pt}
\[
\quad
y_{\dff n}
\off =\off
x_{\dff 1}\qff x_{\dff 2}\pff \ldots\pff x_{\dff n}
\qff.
\]

\vspace{-12pt}
Clearly\halfff,\qss\vspace{-2.25pt} 
\[
\quad
y_{\dff 0}\qff +\qff
y_{\dff 1}\qff +\qff
\ldots\qff +\qff
y_{\dff i}
\off =\off
1
\pff -\pff
x_{\dff 1}\qff x_{\dff 2}\pff \ldots\pff x_{\dff i\dff +\dff 1}
\]

\vspace{-10.5pt}
for\qss $i\qff <\qff n$\qss
and\dss hence\qss
$y_{\dff 0}\qff +\qff
y_{\dff 1}\qff +\qff
\ldots\qff +\qff
y_{\dff n}
\off =\off
1$\nnsp.\oss

The preimages under\dss the map\dss $\lambda$\dss of\trs the faces\dss $\delta_{\dff i}$\dss are\sss
the following\fff:\vspace{3pt}
\[
\quad
\lambda^{\dff -\dff 1}\dff(\trf \delta_{\dff 0}\trf)
\off =\off
B_{\dff 1}
\qff,
\]

\vspace{-36pt}
\begin{equation}
\label{preimages}
\quad
\lambda^{\dff -\dff 1}\dff(\trf \delta_{\dff i}\trf)
\off =\off
B_{\dff i\dff +\dff 1}\qff \cup\qff
A_{\dff 1}\qff \cup\qff
A_{\dff 2}\qff \cup\qff
\ldots\qff \cup\qff
A_{\dff i}
\hspace*{1.2em}\mbox{for}\hspace*{1.2em}
0\off <\off i\off <\off n
\qff,
\hspace*{1.2em}\mbox{and}
\end{equation}

\vspace{-36pt}
\[
\quad
\lambda^{\dff -\dff 1}\dff(\trf \delta_{\dff n}\trf)
\off =\off
A_{\dff 1}\qff \cup\qff
A_{\dff 2}\qff \cup\qff
\ldots\qff \cup\qff
A_{\dff n}
\qff.
\]

\vspace{-9pt}
\myuppar{Sperner's\dss version of\qss Lebesgue\dss first\dss covering\dss theorem.}
\emph{Suppose\sss that\qss
$F_{\dff 0}\fff,\pff
F_{\dff 1}\fff,\pff
\ldots\fff,\pff
F_{\dff n}$\qss 
is\dss a\sss closed\sss covering\sss of\dss $\delta$\dss such\dss that\dss $F_{\dff i}$\dss
is\dss disjoint\dss from\dss $\delta_{\dff i}$\dss for all\dss $i$\nnsp.\oss
Then\dss the intersection of\dss the sets\qss 
$F_{\dff 0}\fff,\pff
F_{\dff 1}\fff,\pff
\ldots\fff,\pff
F_{\dff n}$\qss 
is\dss non-empty\halfff.\oss}

\proof
Let\qss
$D_{\dff i}
\off =\off
\lambda^{\dff -\dff 1}\dff(\dff F_{\dff i\dff -\dff 1}\dff)$\qss
for\sss every\qss
$i\off =\off 1\fff,\pff 2\fff,\pff \ldots\fff,\pff n\qff +\qff 1$\nnsp.\oss
Then\qss
$D_{\dff 1}\fff,\pff
D_{\dff 2}\fff,\pff
\ldots\fff,\pff
D_{\dff n\dff +\dff 1}$\qss
is\dss a\sss closed covering\sss of\dss $[\trf 0\fff,\qff 1\trf]^{n}$\dnsp.\qff\oss 
Since\dss $F_{\dff i}$\dss
is\dss disjoint\dss from\dss $\delta_{\dff i}$\dss
for every\dss $i$\nnsp,\oss
(\ref{preimages})\pss 
implies\dss that\dss\vspace{3pt}
\[
\quad
D_{\dff n\dff +\dff 1}
\hspace*{0.8em}\mbox{is\qss disjoint\qss from}\hspace*{0.8em}
A_{\dff 1}\qff \cup\qff
A_{\dff 2}\qff \cup\qff
\ldots\qff \cup\qff
A_{\dff n}\qff,
\hspace*{1.2em}\mbox{and}
\]

\vspace{-36pt}
\[
\quad
D_{\dff k}
\hspace*{0.8em}\mbox{is\qss disjoint\qss from}\hspace*{0.8em}
B_{\dff k}
\hspace*{0.8em}\mbox{and}\hspace*{0.8em}
A_{\dff i}
\hspace*{0.8em}\mbox{with}\hspace*{0.8em}
i\off <\off k
\]

\vspace{-9pt}
for\qss
$k\qff \leq\qff n$\nnsp.\oss
Therefore\sss the sets\qss
$D_{\dff 1}\fff,\pff
D_{\dff 2}\fff,\pff
\ldots\fff,\pff
D_{\dff n\dff +\dff 1}$\qss
satisfy\dss the assumptions of\qss Theorem\qss \ref{early-h-th}\qss
and\dss hence\qss 
$D_{\dff 1}\qff \cap\qff 
D_{\dff 2}\qff \cap\qff
\ldots
\qff \cap\qff
D_{\dff n\dff +\dff 1}
\off \neq\off
\varnothing$\nnsp.\oss
But\dss if\qss
$x\qff \in\qff
D_{\dff 1}\qff \cap\qff 
D_{\dff 2}\qff \cap\qff
\ldots
\qff \cap\qff
D_{\dff n\dff +\dff 1}$\nsp,\oss
then\qss 
$\lambda\qff(\dff x\trf)\qff \in\pff
F_{\dff i}$\qss
for\sss every\dss $i$\dss
and\dss hence\sss
the intersection of\dss the sets\qss 
$F_{\dff 0}\fff,\pff
F_{\dff 1}\fff,\pff
\ldots\fff,\pff
F_{\dff n}$\qss 
is\dss non-empty\halfff.\oss \eproof

\myuppar{Knaster--Kuratowski--Mazurkiewicz\qss and\qss Lebesgue\dss paper\halfff.}
Knaster\halfff,\pss Kuratowski,\qss and\qss Mazur\-kiewicz\dss
were aware of\trs Lebesgue\dss paper\qss \cite{l-two},\oss
but\dss in\qss \cite{kkm}\qss this paper\dss 
is\sss only\dss mentioned.\oss
But\dss they included\qss Theorem\qss \ref{early-h-th}\qss in\qss 
\cite{kkm},\oss 
referring\dss to\qss W.\qss Hurewicz\pss \cite{h1},\pss \cite{h2}.\oss
W.\qss Hurewicz\dss stated\qss Theorem\qss \ref{early-h-th}\qss
in a footnote in\qss \cite{h1}\qss as a\sss theorem of\qss
Lebesgue\dss and\dss Brouwer\halfff,\oss and\dss referred\dss to\qss \cite{h2}\qss
for\dss this\sss particular\dss form of\dss their\dss results.\oss
As we saw,\oss
the main\dss idea of\qss Knaster--Kuratowski--Mazurkiewicz\dss argument\dss
can\dss be easily\dss adapted\dss to use\dss 
Theorem\dss \ref{early-h-th}\dss instead\dss 
of\qss KKM\qss theorem.\oss
But\qss Knaster\halfff,\pss Kuratowski,\pss and\qss Mazurkiewicz\dss 
missed\dss this opportunity\dss to
use\sss more familiar 
cartesian coordinates instead of\dss bary\-cen\-tric ones.\oss
They\dss went\dss in\dss the opposite direction
and deduced\qss 
Theorem\qss \ref{early-h-th}\qss 
from\dss Sperner's\dss version of\qss Lebesgue\trs theorems.\oss
It\dss seems\sss that\qss H.W.\qss Kuhn\qss \cite{ku}\qss was\dss the first\dss to
adapt\dss the\qss KKM\qss argument\dss to cartesian coordinates.

Here\dss is\dss their\sss deduction.\oss
If\trs the sets\qss
$D_{\dff 1}\fff,\pff
D_{\dff 2}\fff,\pff
\ldots\fff,\pff
D_{\dff n\dff +\dff 1}$\qss
satisfy\dss the assumptions of\qss Theorem\qss \ref{early-h-th},\oss
then\qss 
$F_{\dff i}
\off =\off
\lambda\trf(\dff D_{\dff i\dff +\dff 1}\dff)$\qss 
is\dss disjoint\dss from\dss $\delta_{\dff i}$\dss for every\qss
$i\off =\off 0\fff,\pff 1\fff,\pff \ldots\fff,\pff n$\qss
and\qss
$F_{\dff 0}\fff,\pff
F_{\dff 1}\fff,\pff
\ldots\fff,\pff
F_{\dff n}$\qss
cover\dss $\delta$\nnsp.\oss
By\qss Sperner's\dss version of\qss Lebesgue\dss first\dss covering\dss theorem\qss
$F_{\dff 0}\qff \cap\qff 
F_{\dff 1}\qff \cap\qff
\ldots
\qff \cap\qff
F_{\dff n}
\off \neq\off
\varnothing$\nnsp.\oss
If\dss $x$\dss belongs\sss to\sss this\sss intersection,\oss
then\dss $x$\dss cannot\dss belong\dss to any\dss face\dss $\delta_{\dff i}$\dss
and\dss hence belongs\sss to\sss the interior of\dss $\delta$\nnsp.\oss
But\dss $\lambda$\dss induces a bijection\dss between\dss the interiors of\dss
$[\trf 0\fff,\qff 1\trf]^{n}$\dss and\dss $\delta$\nnsp.\oss
It\dss follows\dss that\dss $\lambda^{\dff -\dff 1}\dff(\trf x\trf)$\dss
belongs\sss to\sss the intersection\qss
$D_{\dff 1}\qff \cap\qff 
D_{\dff 2}\qff \cap\qff
\ldots
\qff \cap\qff
D_{\dff n\dff +\dff 1}$\nsp,\oss
which\dss is\dss therefore\sss non-empty\halfff.\oss
Probably\halfff,\oss this\dss deduction\qss \cite{kkm}\qss
is\dss the earliest\sss appearance of\trs
the map\dss $\lambda$\nnsp.\oss
About\dss two decades later\dss the map\dss $\lambda$\dss was\dss rediscovered\dss by\qss
Serre\qss \cite{se}.\oss

\newpage
\mysection{Lebesgue\qss tilings}{tilings}

\myuppar{Lebesgue\qss coverings of\dss order\qss $n\qff +\qff 1$\nnsp.}
In order\dss to prove\sss the\sss topological\dss invariance of\dss dimension,\oss
the second\dss covering\dss theorem of\qss Lebesgue\dss should\dss be complemented\dss
by a construction of\dss $\varepsilon$\dnsp-coverings of\dss the cube\dss 
$[\trf 0\fff,\qff 1\trf]^{n}$\dss
of\dss order\qss $n\qff +\qff 1$\qss for every\qss $\varepsilon\qff >\qff 0$\nnsp.\oss
Actually\halfff,\oss Lebesgue constructed\sss coverings of\dss order\dss
$n\qff +\qff 1$\dss of\trs the whole\dss
$\rrr^{\fff n}$\dss by cubes with\dss the side $\varepsilon'$ for every\dss $\varepsilon'\qff >\qff 0$\nnsp.\oss
Moreover\halfff,\oss his coverings are\qss \emph{partitions},\oss
i.e.\qss consist\sss of\dss cubes with disjoint\dss interiors.\oss
We will\sss call\dss them\qss \emph{Lebesgue\dss tilings}.\oss
The second\dss Lebesgue\dss covering\dss theorem\dss is\dss often
called\qss \emph{Lebesgue\dss tiling\dss theorem}\qss
({\fff}it\dss claims\sss that\sss arbitrary closed $\varepsilon$\dnsp-coverings\sss 
behave\sss like\sss these\sss tilings).\oss

We will\sss construct\qss Lebesgue\trs tilings only\dss for\qss
$\varepsilon'\off =\off 1$\nnsp.\oss
Obviously\halfff,\oss they can\dss be scaled\dss to make\sss the cubes of\dss
scaled\dss partitions arbitrarily\dss small.\oss

Let\dss us consider\dss for every\dss $n$\dss real\dss numbers\qss
$a_{\dff 1}\fff,\pff a_{\dff 2}\fff,\pff \ldots\fff,\off a_{\dff n}
\qff \in\qff \rrr$\qss
the cube\vspace{-3pt}
\[
\quad
c\qff(\dff a_{\dff 1}\fff,\pff a_{\dff 2}\fff,\pff \ldots\fff,\off a_{\dff n} \dff)
\off =\off
\prod_{\dff i\qff =\qff 1}^{\dff n}\qff
[\trf a_{\dff i}\fff,\pff a_{\dff i}\qff +\qff 1\qff ]
\qff.
\]

\vspace{-13.5pt}
Clearly,\oss the cubes\dss
$c\qff(\dff a_{\dff 1}\fff,\pff a_{\dff 2}\fff,\pff \ldots\fff,\off a_{\dff n} \dff)$\dss
with integer\dss $a_{\dff 1}\fff,\pff a_{\dff 2}\fff,\pff \ldots\fff,\off a_{\dff n}$\dss
form a partition of\qss $\rrr^{\fff n}$\dnsp,\oss
but\dss its order\dss $2^n$\dss is\dss too big.\oss
Lebesgue\dss idea\sss is\dss to\sss translate\sss these cubes in\dss the directions
of\dss coordinate axes.\oss
Let\dss us fix\dss $n\qff -\qff 1$\dss real\dss numbers\qss
$\varepsilon_{\dff 1}\fff,\pff \varepsilon_{\dff 2}\fff,\pff 
\ldots\fff,\off \varepsilon_{\dff n\dff -\dff 1}
\qff \in\qff \rrr$\dss 
and consider\dss the cubes\vspace{1.5pt}
\[
\quad
e\qff(\dff a_{\dff 1}\fff,\pff a_{\dff 2}\fff,\pff \ldots\fff,\off a_{\dff n} \dff)
\off =\off
c\qff(\dff 
u_{\dff 1}\fff,\pff 
u_{\dff 2}\fff,\pff  
\ldots\fff,\off 
u_{\dff n}
\dff)
\qff,
\]

\vspace{-13.5pt}
where\vspace{-6pt} 
\[
\quad
u_{\dff 1}
\off =\off
a_{\dff 1}\qff +\qff 
a_{\dff 2}\qff \varepsilon_{\dff 1}\qff +\qff 
a_{\dff 3}\qff \varepsilon_{\dff 2}\qff +\qff 
\ldots\qff +\qff 
a_{\dff n}\qff \varepsilon_{\dff n\dff -\dff 1}
\qff,
\]

\vspace{-39pt}
\[
\quad
u_{\dff 2}
\off =\off
a_{\dff 2}\qff +\qff 
a_{\dff 3}\qff \varepsilon_{\dff 2}\qff +\qff 
\ldots\qff +\qff 
a_{\dff n}\qff \varepsilon_{\dff n\dff -\dff 1}
\qff,
\]

\vspace{-42pt}
\begin{equation}
\label{a-to-u}
\quad
\hspace{1.5em}
\ldots\ldots\ldots
\end{equation}

\vspace{-39pt}
\[
\quad
u_{\dff n\dff -\dff 1}
\off =\off
a_{\dff n\dff -\dff 1}\qff +\qff 
a_{\dff n}\qff \varepsilon_{\dff n\dff -\dff 1}
\qff,
\]

\vspace{-39pt}
\[
\quad
u_{\dff n}
\off =\off
a_{\dff n}
\qff,
\]

\vspace{-12pt}
and\dss $a_{\dff 1}\fff,\pff a_{\dff 2}\fff,\pff \ldots\fff,\off a_{\dff n}$\dss
are integers.\oss

The collection of\dss cubes\dss
$e\qff(\dff a_{\dff 1}\fff,\pff a_{\dff 2}\fff,\pff \ldots\fff,\off a_{\dff n} \dff)$\dss
can\dss be constructed\dss recursively\sss starting\dss with\dss
the partition of\trs  
$\rrr$\dss into intervals\dss
$[\dff a\fff,\pff a\qff +\qff 1\trf]$\nnsp,\oss where\dss $a$\dss is\dss an\dss integer\halfff.\oss
Let\dss us\sss consider $\rrr^{\fff n}$ as\sss the union of\trs the\qss
\emph{layers}\qss
$\mathcal{L}\dff(\dff a\trf)
\off =\off
\rrr^{\fff n\dff -\dff 1} \times\dff [\dff a\fff,\pff a\qff +\qff 1\trf]$\nnsp,\oss
where\dss $a$\dss is\sss an\dss integer\halfff.\oss
Then\dss
the cube\dss 
$e\qff(\dff a_{\dff 1}\fff,\pff a_{\dff 2}\fff,\pff \ldots\fff,\off a_{\dff n} \dff)$
is\dss contained\dss in\dss the\sss layer\dss
$\mathcal{L}\dff(\dff a_{\dff n}\dff)$\nnsp.\oss
The cubes contained\dss in\dss the\sss layer\sss
$\mathcal{L}\dff(\dff 0\dff)$\sss
are\sss the products of\dss similar cubes in\dss $\rrr^{\fff n\dff -\dff 1}$\dss
(constructed\dss using\dss 
$\varepsilon_{\dff 1}\fff,\pff \varepsilon_{\dff 2}\fff,\pff 
\ldots\fff,\off \varepsilon_{\dff n\dff -\dff 2}$\nsp)\qss
with\dss the interval\qss $[\dff 0\fff,\pff 1\trf]$\nnsp.\oss
The cubes\dss in\dss the\sss layer\sss
$\mathcal{L}\dff(\dff a_{\dff n}\dff)$\sss
are obtained\dss from\dss the cubes\sss in\dss the layer\sss 
$\mathcal{L}\dff(\dff 0\dff)$\sss
by\sss a\sss translation\dss   
along\dss the $n${\nnsp}th coordinate axis\dss
moving\dss this\sss layer\dss to\sss the\sss layer\sss
$\mathcal{L}\dff(\dff a_{\dff n}\dff)$\sss
followed\dss by\sss a\sss translation\dss inside\sss the\sss layer\sss 
$\mathcal{L}\dff(\dff a_{\dff n}\dff)$\sss by\dss the vector\vspace{1.5pt}
\[
\quad
(\qff 
a_{\dff n}\qff \varepsilon_{\dff n\dff -\dff 1}\fff,\pff 
a_{\dff n}\qff \varepsilon_{\dff n\dff -\dff 1}\fff,\pff  
\ldots\fff,\off
a_{\dff n}\qff \varepsilon_{\dff n\dff -\dff 1}\fff,\pff 
0
\qff)
\qff.
\]

\vspace{-10.5pt}
In\dss particular\halfff,\oss
we see\sss that\dss the cubes\dss
$e\qff(\dff a_{\dff 1}\fff,\pff a_{\dff 2}\fff,\pff \ldots\fff,\off a_{\dff n} \dff)$\dss
indeed\dss form a partition of\dss $\rrr^{\fff n}$\dnsp.\oss
It\dss turns out\dss that\dss the order of\trs
this partition\dss is\dss equal\dss to\qss $n\qff +\qff 1$\qss
under some mild assumptions about\dss the numbers\qss 
$\varepsilon_{\dff 1}\fff,\pff \varepsilon_{\dff 2}\fff,\pff 
\ldots\fff,\off \varepsilon_{\dff n\dff -\dff 1}$\nsp.\oss
Let\dss us say\dss that\dss the numbers\qss
$\varepsilon_{\dff 1}\fff,\pff \varepsilon_{\dff 2}\fff,\pff 
\ldots\fff,\off \varepsilon_{\dff n\dff -\dff 1}$\qss
are\qss \emph{generic}\oss if\qss
for\dss every\qss $i\qff \leq\qff n\qff -\qff 1$\qss
the number\dss 
$\varepsilon_{\dff i}$\dss
is\dss not\sss equal\dss to any\dss linear combination\oss\vspace{1.25pt}
\[
\quad
m_{\dff 0}\qff +\qff
m_{\dff 1}\dff \varepsilon_{\dff 1}\qff +\qff
m_{\dff 2}\dff \varepsilon_{\dff 2}\qff +\qff
\ldots\qff +\qff
m_{\dff i\dff -\dff 1}\dff \varepsilon_{\dff i\dff -\dff 1} 
\]

\vspace{-10.75pt}
of\trs the numbers\dss
$1\fff,\pff 
\varepsilon_{\dff 1}\fff,\pff \varepsilon_{\dff 2}\fff,\pff 
\ldots\fff,\off \varepsilon_{\dff i\dff -\dff 1}$\dss
with\dss integer coefficients\qss
$m_{\dff 0}\fff,\pff m_{\dff 1}\fff,\pff m_{\dff 2}\fff,\pff 
\ldots\fff,\pff m_{\dff i\dff -\dff 1}$\nnsp.\oss
For example,\oss
Lebesgue\dss choice\qss
$\varepsilon_{\dff i}\off =\off 2^{\dff -\dff i}$\qss is generic.\oss

\mypar{Theorem.}{tiling-intersections}
\emph{If\qss the choice of\qss numbers\qss
$\varepsilon_{\dff 1}\fff,\pff \varepsilon_{\dff 2}\fff,\pff 
\ldots\fff,\off \varepsilon_{\dff n\dff -\dff 1}$\qss
is\dss generic,\oss then\dss
the intersection of\qss every\dss $r$ different\sss cubes of\qss the form\dss
$e\qff(\dff a_{\dff 1}\fff,\pff a_{\dff 2}\fff,\pff \ldots\fff,\off a_{\dff n} \dff)$\nnsp,\oss
where\qss $a_{\dff 1}\fff,\pff a_{\dff 2}\fff,\pff \ldots\fff,\off a_{\dff n}$\qss
are\dss integers,\oss
is\dss contained\dss a\sss plane\sss of\qss dimension\qss $n\qff +\qff 1\qff -\qff r$\qss
defined\qss by\qss $r\qff -\qff 1$\qss 
equations\dss of\qss the\dss form}\vspace{0pt}
\[
\quad
x_{\dff i}
\off =\off
A_{\dff i}
\qff,
\]

\vspace{-12pt}
\emph{where\dss  
$A_{\dff i}$\sss
is\dss a\dss linear combination\dss of\pss the numbers\pss
$1\fff,\pff \varepsilon_{\dff 1}\fff,\pff \varepsilon_{\dff 2}\fff,\pff 
\ldots\fff,\off \varepsilon_{\dff n\dff -\dff 1}$\dss
with\dss integer\dss coefficients.\oss
In\dss par\-tic\-u\-lar\halfff,\oss every\qss $n\qff +\qff 2$\qss different\dss cubes
have empty\dss intersection and\dss the\dss order\sss of\pss the covering\dss of\oss $\rrr^{\fff n}$\dss 
by\dss these cubes\dss
is\dss equal\dss to\qss $n\qff +\qff 1$\nnsp.\oss}

\proof
The\sss theorem\sss is\sss obvious\sss for\qss $n\off =\off 1$\nnsp.\oss
Suppose\sss that\dss the\sss theorem\sss is\dss true with\qss $n\qff -\qff 1$\qss 
in\dss the role\sss of\dss $n$\nnsp.\oss
Then\dss the intersections of\dss cubes contained\dss in\sss a
single\sss layer\sss
$\mathcal{L}\dff(\dff a_{\dff n}\trf)$\sss
have\sss the required\dss property\halfff.\oss
In\dss more details,\oss the intersections of\dss any\dss $r$\dss 
different\sss cubes\dss in\dss the\sss layer\dss
$\mathcal{L}\dff(\dff a_{\dff n}\trf)$\dss
is\dss contained\dss in\sss a plane defined\dss by equations of\trs the form\vspace{3pt}
\begin{equation}
\label{tiling-equations}
\quad
x_{\dff i}
\off =\off
A_{\dff i}
\qff +\qff
a_{\dff n}\qff \varepsilon_{\dff n\dff -\dff 1}
\qff,
\end{equation}

\vspace{-9pt}
where\dss $A_{\dff i}$\dss is\dss a\dss linear combination\dss of\oss 
$1\fff,\pff \varepsilon_{\dff 1}\fff,\pff \varepsilon_{\dff 2}\fff,\pff 
\ldots\fff,\off \varepsilon_{\dff n\dff -\dff 2}$\dss
with\dss integer\dss coefficients.\oss
Suppose\sss that\sss our\dss $r$\dss cubes contained\dss in\dss two layers one
next\dss to\sss the other\qss
(otherwise\sss their\dss intersection\dss is\sss empty\fff),\oss
say,\oss in\dss the\sss layers\dss
$\mathcal{L}\dff(\dff a_{\dff n}\qff -\qff 1\dff)$\dss
and\dss
$\mathcal{L}\dff(\dff a_{\dff n}\trf)$\nnsp,\oss
but\dss not\dss in\sss one of\trs them.\oss
Suppose\sss that\dss $r_{\dff 0}$\dss of\trs these cubes are contained\dss in\dss
$\mathcal{L}\dff(\dff a_{\dff n}\qff -\qff 1\dff)$\dss
and\dss $r_{\dff 1}$\dss
in\dss $\mathcal{L}\dff(\dff a_{\dff n}\trf)$\nnsp.\oss
Then\qss
$r\off =\off r_{\dff 0}\qff +\qff r_{\dff 1}$\qss
and\qss $r_{\dff 0}\dff,\off r_{\dff 1}\qff >\pff 0$\nnsp.\oss
Clearly,\oss the intersection of\trs these\dss $r$\dss cubes\sss
is\sss contained\dss in\dss the hyperplane\dss
$\rrr^{\fff n\dff -\dff 1} \times\dff a_{\dff n}$\nsp.\oss
By\dss the previous paragraph,\oss 
the intersection of\dss $r_{\dff 1}$\dss cubes\dss in\dss
$\mathcal{L}\dff(\dff a_{\dff n}\trf)$\dss
is contained\dss in\dss plane defined\dss by\dss the equations of\trs the
form\qss (\ref{tiling-equations})\qss
with\qss
$i\off \neq\off n$\nnsp,\oss
and\dss the intersection of\dss $r_{\dff 0}$\dss cubes\dss in\dss
$\mathcal{L}\dff(\dff a_{\dff n}\qff -\qff 1\dff)$\dss
is\dss contained\dss in\sss a plane with\dss the equations of\trs the form\vspace{2.25pt}
\[
\quad
x_{\fff j}
\off =\off
A_{\fff j}'
\qff +\qff
(\dff a_{\dff n}\qff -\qff 1\dff)\qff \varepsilon_{\dff n\dff -\dff 1}
\qff,
\]

\vspace{-9.75pt}
where\dss $A_{\fff j}'$\dss is\dss also a\dss linear combination\dss of\oss 
$1\fff,\pff \varepsilon_{\dff 1}\fff,\pff \varepsilon_{\dff 2}\fff,\pff 
\ldots\fff,\off \varepsilon_{\dff n\dff -\dff 2}$\dss
with\dss integer\dss coefficients,\oss
and\qss
$j\off \neq\off n$\qss
also.\oss
If\qss the same coordinate\dss $x_{\dff i}$\dss occurs\dss in each of\trs these\sss
two systems of\dss equations,\oss
then either\dss they are not\sss compatible and\dss the intersection of\trs
all\sss $r$\sss cubes\sss is\sss empty\halfff,\oss
or\vspace{3pt}
\[
\quad
A_{\dff i}
\qff +\qff
a_{\dff n}\qff \varepsilon_{\dff n\dff -\dff 1}
\off =\off
A_{\dff i}'
\qff +\qff
(\dff a_{\dff n}\qff -\qff 1\dff)\qff \varepsilon_{\dff n\dff -\dff 1}
\qff,
\]

\vspace{-9pt}
and\dss hence\qss
$\varepsilon_{\dff n\dff -\dff 1}
\off =\off
A_{\dff i}'
\pff -\pff
A_{\dff i}$\qss
is\dss a\dss linear combination\dss of\oss 
$1\fff,\pff \varepsilon_{\dff 1}\fff,\pff \varepsilon_{\dff 2}\fff,\pff 
\ldots\fff,\off \varepsilon_{\dff n\dff -\dff 2}$\dss
with\dss integer\dss coefficients,\oss
contrary\dss to\sss the assumption.\oss
The contradiction shows\sss that\dss if\trs the intersection\dss is\dss non-empty\halfff,\oss
then it\dss is contained\dss in\dss the plane defined\dss
by\qss
$r_{\dff 0}\qff -\qff 1\qff +\qff r_{\dff 1}\qff -\qff 1
\off =\off
r\qff -\qff 2$\qss
linearly\dss independent\sss equations.\oss
Together\dss with\qss
$x_{\dff n}\off =\off a_{\dff n}$\qss
we get\qss $r\qff -\qff 1$\qss equations,\oss
as claimed.\oss  \eproof

\myuppar{Lebesgue\dss tiling\dss theorem.}
\emph{For every\qss $\varepsilon\qff >\qff 0$\qss
there\dss exists\dss a\dss closed\dss $\varepsilon$\dnsp-covering\sss
of\pss $[\trf 0\fff,\qff 1\trf]^{n}$\dss of\dss order\dss $n\qff +\qff 1$\dnsp.}

\proof
A scaled\dss version of\trs the above covering of\qss $\rrr^{\fff n}$\dss
consists of\dss cubes with sides of\dss any\sss given\dss length.\oss
The cube\qss $[\trf 0\fff,\qff 1\trf]^{n}$\qss is\dss contained\dss in\dss
the union of\dss finitely\dss many\sss cubes of\trs the scaled covering\halfff,\oss
and\dss it\dss is\sss sufficient\dss to\sss take\sss their\dss intersections
with\qss $[\trf 0\fff,\qff 1\trf]^{n}$\dnsp.\oss  \eproof

\myuppar{Remark.}
A converse\sss to\dss Theorem\dss \ref{tiling-intersections}\dss is\dss also\sss true.\oss
Namely\halfff,\oss if\qss
for every $r$ the\sss intersection
of\dss $r$ different\sss cubes of\qss the form\dss
$e\qff(\dff a_{\dff 1}\fff,\pff a_{\dff 2}\fff,\pff \ldots\fff,\off a_{\dff n} \dff)$\qss
({\fff}with\dss integer\qss $a_{\dff 1}\fff,\pff a_{\dff 2}\fff,\pff \ldots\fff,\off a_{\dff n}$\nsp)\qss
is\dss contained\dss in\sss a\sss plane\sss of\qss dimension\qss $n\qff +\qff 1\qff -\qff r$\dnsp,\oss
then\qss
$\varepsilon_{\dff 1}\fff,\pff \varepsilon_{\dff 2}\fff,\pff 
\ldots\fff,\off \varepsilon_{\dff n\dff -\dff 1}$\qss
are generic.\oss
This can\dss be proved\dss by\dss inverting\dss the arguments in\dss the proof\trs of\qss
Theorem\qss \ref{tiling-intersections}.\oss
We leave details\sss to\sss reader\halfff.\oss

\myuppar{Intersections of\dss cubes in\qss Lebesgue\dss tilings.}
Now\sss we\sss turn\dss to\sss
the pattern of\dss intersections of\dss cubes\dss of\dss
a\dss Lebesgue\dss tiling\halfff.\oss
Of\dss course,\oss  
the choice of\qss
$\varepsilon_{\dff 1}\fff,\pff \varepsilon_{\dff 2}\fff,\pff 
\ldots\fff,\off \varepsilon_{\dff n\dff -\dff 1}$\qss
is\dss assumed\dss to be generic.\oss
In\dss particular\halfff,\oss none of\dss numbers\sss $\varepsilon_{\dff i}$\sss is\dss an\dss integer\halfff.\oss
Clearly\halfff,\oss replacing\dss $\varepsilon_{\dff i}$\dss
by\qss $\varepsilon_{\dff i}\qff +\qff k_{\dff i}$\qss
with\dss integer\dss $k_{\dff i}$\dss does not\sss change\sss the partition.\oss
We prefer\dss $\varepsilon_{\dff i}$\dss to be negative
in order\dss to match\dss better some further constructions
and\dss will\sss assume\sss that\qss
${}-\qff 1\qff <\qff \varepsilon_{\dff i}\qff <\qff 0$\qss
for every\qss $i\qff \leq\qff n\qff -\qff 1$\nnsp.\oss
We will\sss also assume\sss that\oss 
$\num{\varepsilon_{\dff 1}}\qff +\qff 
\num{\varepsilon_{\dff 2}}\qff +\qff
\ldots\qff +\qff 
\num{\varepsilon_{\dff n\dff -\dff 1}} 
\off <\off 1$\nnsp.\oss

\mypar{Lemma.}{intersection-inequality}
\emph{If\pss the cubes\qss 
$e\qff(\dff a_{\dff 1}\fff,\pff a_{\dff 2}\fff,\pff \ldots\fff,\off a_{\dff n} \dff)$\dss
and\qss
$e\qff(\dff b_{\dff 1}\fff,\pff b_{\dff 2}\fff,\pff \ldots\fff,\off b_{\dff n} \dff)$\qss
intersect\halfff,\oss
then\dss}\vspace{3pt}
\begin{equation}
\label{int-ineq}
\quad
b_{\dff i}
\off\qff \leq\off\qff
a_{\dff i}\qff +\qff
1\qff +\qff
(\dff a_{\dff i\dff +\dff 1}\qff -\qff b_{\dff i\dff +\dff 1} \dff)\qff \varepsilon_{\dff i}
\qff +\qff 
\ldots\qff +\qff 
(\dff a_{\dff n}\qff -\qff b_{\dff n} \dff)\qff \varepsilon_{\dff n\dff -\dff 1}
\end{equation}

\vspace{-9pt}
\emph{for\dss every\qss $i\qff \leq\qff n$\qss
(without\sss any assumptions about\dss
the choice of\pss
$\varepsilon_{\dff 1}\fff,\pff \varepsilon_{\dff 2}\fff,\pff 
\ldots\fff,\off \varepsilon_{\dff n\dff -\dff 1}$\nsp).\oss}

\proof
If\qss these cubes intersect\halfff,\oss
then\dss the intervals\vspace{4.5pt}
\[
\quad
[\qff
a_{\dff i}\qff +\qff 
a_{\dff i\dff +\dff 1}\qff \varepsilon_{\dff i}\qff +\qff 
\ldots\qff +\qff 
a_{\dff n}\qff \varepsilon_{\dff n\dff -\dff 1}
\dff,\qff\off
a_{\dff i}\qff +\qff 
1\qff +\qff 
a_{\dff i\dff +\dff 1}\qff \varepsilon_{\dff i}\qff +\qff 
\ldots\qff +\qff 
a_{\dff n}\qff \varepsilon_{\dff n\dff -\dff 1}
\qff]
\]

\vspace{-15.25pt}
and\vspace{-3.25pt}
\[
\quad
[\qff
b_{\dff i}\qff +\qff 
b_{\dff i\dff +\dff 1}\qff \varepsilon_{\dff i}\qff +\qff 
\ldots\qff +\qff 
b_{\dff n}\qff \varepsilon_{\dff n\dff -\dff 1}
\dff,\qff\off
b_{\dff i}\qff +\qff 
1\qff +\qff 
b_{\dff i\dff +\dff 1}\qff \varepsilon_{\dff i}\qff +\qff 
\ldots\qff +\qff 
b_{\dff n}\qff \varepsilon_{\dff n\dff -\dff 1}
\qff]
\]

\vspace{-9pt}
overlap\dss for every\qss $i\qff \leq\qff n$\nnsp.\oss
It\dss follows\dss that\vspace{3pt}
\[
\quad
b_{\dff i}\qff +\qff 
b_{\dff i\dff +\dff 1}\qff \varepsilon_{\dff i}\qff +\qff 
\ldots\qff +\qff 
b_{\dff n}\qff \varepsilon_{\dff n\dff -\dff 1}
\off\qff \leq\off\qff
a_{\dff i}\qff +\qff 
1\qff +\qff 
a_{\dff i\dff +\dff 1}\qff \varepsilon_{\dff i}\qff +\qff 
\ldots\qff +\qff 
a_{\dff n}\qff \varepsilon_{\dff n\dff -\dff 1}
\qff.
\]

\vspace{-9pt}
The last\dss inequality\sss obviously\dss implies\qss (\ref{int-ineq}).\oss  \eproof

\mypar{Lemma.}{intersecting-cubes-close}
\emph{Under\dss the above assumptions,\oss if\qss 
$e\qff(\dff a_{\dff 1}\fff,\pff a_{\dff 2}\fff,\pff \ldots\fff,\off a_{\dff n} \dff)$\dss
and\qss
$e\qff(\dff b_{\dff 1}\fff,\pff b_{\dff 2}\fff,\pff \ldots\fff,\off b_{\dff n} \dff)$\qss
intersect\halfff,\oss
then\qss 
$\num{a_{\dff i}\qff -\qff b_{\dff i}}
\off \leq\off
1$\qss
for every\qss $i\qff \leq\qff n$\nnsp.\oss}

\proof
Suppose\sss that\qss $i\qff \leq\qff n$\qss
and\qss
$\num{a_{\dff j}\qff -\qff b_{\dff j}}\qff \leq\qff 1$\pss
for\oss $j\qff \geq\qff i\qff +\qff 1$\nnsp.\oss
Lemma\qss \ref{intersection-inequality}\qss implies\sss that\vspace{3pt}
\[
\quad
b_{\dff i}
\off\qff \leq\off\qff
a_{\dff i}\qff +\qff
1\qff +\qff
\num{a_{\dff i\dff +\dff 1}\qff -\qff b_{\dff i\dff +\dff 1}}\dff \cdot\dff \num{\varepsilon_{\dff i}}
\qff +\qff 
\ldots\qff +\qff 
\num{a_{\dff n}\qff -\qff b_{\dff n}}\dff \cdot\dff \num{\varepsilon_{\dff n\dff -\dff 1}}
\qff.
\]

\vspace{-36pt}
\[
\quad
\phantom{b_{\dff i}
\off\qff }
\leq\off\qff
a_{\dff i}\qff +\qff
1\qff +\qff
\num{\varepsilon_{\dff i}}
\qff +\qff 
\ldots\qff +\qff 
\num{\varepsilon_{\dff n\dff -\dff 1}}
\off\qff 
\leq\off\qff
a_{\dff i}\qff +\qff
1\qff +\qff
\num{\varepsilon_{\dff 1}}
\qff +\qff 
\ldots\qff +\qff 
\num{\varepsilon_{\dff n\dff -\dff 1}}
\qff.
\]

\vspace{-9pt}
But\qss
$\num{\varepsilon_{\dff 1}}
\qff +\qff 
\ldots\qff +\qff 
\num{\varepsilon_{\dff n\dff -\dff 1}}
\off <\off 1$\qss
and\dss hence\qss
$b_{\dff i}\qff <\qff a_{\dff i}\qff +\qff 2$\nnsp.\oss
Since\dss $a_{\dff i}$\dss and\dss $b_{\dff i}$\dss
are integers,\oss
in\dss fact\qss
$b_{\dff i}\qff \leq\qff a_{\dff i}\qff +\qff 1$\nnsp.\oss
By\dss interchanging\dss the roles of\dss $a_{\dff i}$\dss and\dss $b_{\dff i}$\nsp,\oss
we conclude\sss that\sss also\qss
$a_{\dff i}\qff \leq\qff b_{\dff i}\qff +\qff 1$\nnsp.\oss
It\dss follows\dss that\qss
$\num{a_{\dff i}\qff -\qff b_{\dff i}}
\off \leq\off
1$\nnsp.\oss
It\dss remains\sss to use\sss the descending\dss induction\dss by\dss $i$\nnsp.\oss  \eproof

\mypar{Lemma.}{intersecting-cubes-comparable}
\emph{Under\dss the above assumptions,\oss if\qss 
$e\qff(\dff a_{\dff 1}\fff,\pff a_{\dff 2}\fff,\pff \ldots\fff,\off a_{\dff n} \dff)$\dss
and\qss
$e\qff(\dff b_{\dff 1}\fff,\pff b_{\dff 2}\fff,\pff \ldots\fff,\off b_{\dff n} \dff)$\qss
intersect\halfff,\oss
then\sss either\pss 
$a_{\dff i}\qff \leq\qff b_{\dff i}$\pss
for\dss every\qss $i\qff \leq\qff n$\nnsp,\oss
or\pss
$b_{\dff i}\qff \leq\qff a_{\dff i}$\pss
for\dss every\qss $i\qff \leq\qff n$\nnsp.\oss}

\proof
Suppose\sss that\dss the\sss lemma holds with\qss $n\qff -\qff 1$\qss in\dss the role 
of\dss $n$\nnsp.\oss
Then\dss it\dss is\dss true in\dss the case when\dss both cubes are contained\dss in\dss the same\sss layer\sss
$\mathcal{L}\dff(\dff a\qff)$\nnsp.\oss
The cubes\qss
$e\qff(\dff a_{\dff 1}\fff,\pff a_{\dff 2}\fff,\pff \ldots\fff,\off a_{\dff n} \dff)$\dss
and\qss
$e\qff(\dff b_{\dff 1}\fff,\pff b_{\dff 2}\fff,\pff \ldots\fff,\off b_{\dff n} \dff)$\qss
are contained\dss in\dss the layers\dss $\mathcal{L}\dff(\dff a_{\dff n}\dff)$\dss and\qss
$\mathcal{L}\dff(\trf b_{\dff n}\dff)$\dss respectively\halfff.\oss
If\qss these layers are different\halfff,\oss
then\qss
$\num{a_{\dff n}\qff -\qff b_{\dff n}}
\off =\off
1$\nnsp.\oss
We may assume\sss that\qss $a_{\dff n}\off =\off b_{\dff n}\qff +\qff 1$\nnsp.\oss
Suppose\sss that\qss $i\qff <\qff n$\qss and\qss
$a_{\fff j}\qff \geq\qff b_{\fff j}$\pss
for\oss $j\qff \geq\qff i\qff +\qff 1$\nnsp.\oss
By\qss Lemma\qss \ref{intersection-inequality}\vspace{3pt}
\[
\quad
b_{\dff i}
\off\qff \leq\off\qff
a_{\dff i}\qff +\qff
1\qff +\qff
(\dff a_{\dff i\dff +\dff 1}\qff -\qff b_{\dff i\dff +\dff 1} \dff)\qff \varepsilon_{\dff i}
\qff +\qff 
\ldots\qff +\qff 
(\dff a_{\dff n}\qff -\qff b_{\dff n} \dff)\qff \varepsilon_{\dff n\dff -\dff 1}
\qff.
\]

\vspace{-9pt}
Since\sss the numbers\qss
$\varepsilon_{\dff i}\fff,\pff \varepsilon_{\dff i\dff +\dff 1}\fff,\pff 
\ldots\fff,\off \varepsilon_{\dff n\dff -\dff 1}$\qss
are negative,\pss 
$a_{\dff n}\qff -\qff b_{\dff n}\off =\off 1$\nnsp,\oss
and\qss
$a_{\fff j}\qff -\qff b_{\fff j}\qff \geq\qff 0$\qss
for\qss $j\qff \geq\qff i\qff +\qff 1$\nnsp,\oss
this inequality\dss implies\sss that\qss
$b_{\dff i}\qff \leq\qff
a_{\dff i}\qff +\qff
1\qff +\qff
\varepsilon_{\dff n\dff -\dff 1}
\qff <\qff a_{\dff i}\qff +\qff
1$\nnsp.\oss
Since\dss $a_{\dff i}$\dss and\dss $b_{\dff i}$\dss
are integers,\oss
it\dss follows\dss that\qss
$b_{\dff i}\qff \leq\qff a_{\dff i}$\nsp.\oss
It\dss remains\sss to use\sss the descending\dss induction\dss by\dss $i$\nnsp.\oss  \eproof

\myuppar{Discrete cubes and a partial\dss order on\dss $\rrr^{\fff n}$\dnsp.}
The last\dss two lemmas suggest\dss the following\dss two definitions.\oss
A\qss \emph{discrete $n$\dnsp-cube}\qss is\dss defined\sss as\sss a\sss subset\sss of\qss $\zzz^{\fff n}$\dss 
of\trs the form\vspace{1.5pt}
\[
\quad
\prod_{\dff i\qff =\qff 1}^{\dff n}\qff
\{\trf z_{\dff i}\fff,\pff z_{\dff i}\qff +\qff 1\qff \}
\qff,
\]

\vspace{-12.25pt}
where\qss
$z_{\dff 1}\fff,\pff z_{\dff 2}\fff,\pff \ldots\fff,\off z_{\dff n}
\off \in\off \zzz$\nnsp.\oss
By\dss the definition\dss 
$(\dff a_{\dff 1}\fff,\pff a_{\dff 2}\fff,\pff \ldots\fff,\off a_{\dff n} \dff)
\qff\dff \leq\qff\dff
(\dff b_{\dff 1}\fff,\pff b_{\dff 2}\fff,\pff \ldots\fff,\off b_{\dff n} \dff)$\dss
means\sss that\dss
$a_{\dff i}\qff \leq\qff b_{\dff i}$\dss
for\dss all\dss $i\qff \leq\qff n$\nnsp.\oss
As usual,\pss $u\qff <\qff v$
means\sss that\sss $u\qff \leq\qff v$\dss and\dss $u\off \neq\off v$\nnsp.\oss

\myuppar{Pivot\dss (sub)sequences.}
Let\qss
$v\dff(\dff 0\dff)\fff,\off
v\dff(\dff 1\dff)\fff,\off
\ldots\fff,\off
v\dff(\dff r\dff)
\qff \in\qff
\zzz^{\fff n}$\qss
be pairwise\qss $\leq$\dnsp-comparable distinct\dss $n$\dnsp-tuples of\dss integers.\oss
After\dss reordering\dss them,\oss if\trs necessary\halfff,\oss
we may assume\sss that\vspace{3pt}
\[
\quad
v\dff(\dff 0\dff)\off < \off
v\dff(\dff 1\dff)\off < \off
\ldots\off < \off
v\dff(\dff r\dff)
\qff.
\]

\vspace{-9pt}
Since\dss $v\dff(\dff i\qff -\qff 1\dff)\qff <\qff v\dff(\dff i\dff)$\nnsp,\oss
the $n$\dnsp-tuple $v\dff(\dff i\trf)$
can\dss be\sss obtained\dss from $v\dff(\dff i\qff -\qff 1\dff)$
by\dss increasing\sss several\sss coordinates.\oss
Suppose\sss that\halfff,\oss in addition,\oss
the $n$\dnsp-tuples\qss
$v\dff(\dff 0\dff)\fff,\off
v\dff(\dff 1\dff)\fff,\off
\ldots\fff,\off
v\dff(\dff r\dff)$\qss
are contained\dss in\dss a\sss discrete $n$\dnsp-cube.\oss
Then each coordinate may increase only\dss by\sss $1$\sss
and\sss only\sss once along\dss the sequence\qss
$v\dff(\dff 0\dff),\off
v\dff(\dff 1\dff),\off
\ldots,\off
v\dff(\dff r\dff)$\nsp.\oss
This\dss implies,\oss in\dss particular\halfff,\pss
that\dss $r\qff \leq\qff n$\nnsp.\oss
Clearly\halfff,\oss the directions of\dss vectors from $v\dff(\dff i\qff -\qff 1\dff)$
to $v\dff(\dff i\trf)$ are all\sss different\halfff.\oss
Let\dss us\sss define a\qss \emph{pivot\dss subsequence}\qss 
as\dss an\dss $<$\dnsp-increasing sequence of\dss 
elements of\dss a discrete $n$\dnsp-cube.\oss
A\qss \emph{pivot\sss sequence}\qss is\dss a pivot\sss subsequence 
consisting of\qss $n\qff +\qff 1$\qss terms.\oss

In a pivot\sss sequence each coordinate increases by $1$ as we go along\dss the sequence,\oss
and only once.\oss
Increasing $m$ coordinates at\dss one step in a pivot\sss subsequence
can\dss be replaced\dss by\dss increasing\dss them one by one\sss in $m$ steps.\oss
It\dss follows\dss that\sss every\dss pivot\sss subsequence\dss is\dss
a\sss subsequence of\dss a\sss pivot\sss sequence,\oss
justifying\dss the\sss term.\oss

\myuppar{Sequences of\dss cubes associated\dss with\dss pivot\dss subsequences.}
Let\qss
$v\dff(\dff 0\dff)\qff < \qff
v\dff(\dff 1\dff)\qff < \qff
\ldots\qff < \qff
v\dff(\dff r\dff)$\qss
be a pivot\sss subsequence.\oss
For\qss $i\off =\off 0\fff,\pff 1\fff,\pff \ldots\fff,\pff r$\qss
let\vspace{3pt}
\[
\quad
v\dff(\dff i\trf)
\off =\off
(\trf
a_{\dff 1}\dff(\dff i\trf),\off
a_{\dff 2}\dff(\dff i\trf),\off
\ldots,\off
a_{\dff n}\dff(\dff i\trf)
\qff)
\qff,
\]

\vspace{-9pt}
and\trs let\dss the real\dss numbers\dss\vspace{3pt}
\[
\quad
u_{\dff 1}\dff(\dff i\trf),\off
u_{\dff 2}\dff(\dff i\trf),\off
\ldots,\off
u_{\dff n}\dff(\dff i\trf)
\]

\vspace{-9pt}
be\sss related\dss to\sss the integers\dss \vspace{3pt}
\[
\quad
a_{\dff 1}\dff(\dff i\trf),\off
a_{\dff 2}\dff(\dff i\trf),\off
\ldots,\off
a_{\dff n}\dff(\dff i\trf)
\]

\vspace{-9pt}
by\dss an obvious modification of\trs the formulas\qss (\ref{a-to-u}).\oss
Let\vspace{3pt}
\[
\quad
e\qff(\dff i\trf)
\off =\off
e\qff(\trf
a_{\dff 1}\dff(\dff i\trf),\off
a_{\dff 2}\dff(\dff i\trf),\off
\ldots,\off
a_{\dff n}\dff(\dff i\trf)
\qff)
\off =\off
c\qff(\trf
u_{\dff 1}\dff(\dff i\trf),\off
u_{\dff 2}\dff(\dff i\trf),\off
\ldots,\off
u_{\dff n}\dff(\dff i\trf)
\qff)
\qff.
\]

\vspace{-6pt}
\mypar{Theorem.}{rotating-intersections}
\emph{Under\dss the above assumptions}\qss\vspace{2.25pt}
\[
\quad
e\qff(\dff 0\trf)\qff \cap\qff
e\qff(\dff 1\trf)\qff \cap\qff
\ldots\qff \cap\qff
e\qff(\dff r\trf)
\qff\off =\dff\off
\prod_{\dff i\qff =\qff 1}^{\dff n}\qff
J_{\fff i}
\off,
\]

\vspace{-9.875pt}
\emph{where\qss $n\qff -\qff r$\qss factors\qss $J_{\dff i}$\dss 
are non-degenerate closed\dss intervals in\dss $\rrr$\dss
and\dss $r$\dss factors\qss $J_{\dff i}$\dss are one-point\sss subsets of\pss $\rrr$\nnsp.\oss
In\dss particular\halfff,\oss
$e\qff(\dff 0\trf)\qff \cap\qff
e\qff(\dff 1\trf)\qff \cap\qff
\ldots\qff \cap\qff
e\qff(\dff r\trf)$\qss
is\dss non-empty\halfff.\oss}

\proof
{\qff}By\dss the definitions\vspace{0pt}
\[
\quad
e\qff(\dff 0\trf)\qff \cap\qff
e\qff(\dff 1\trf)\qff \cap\qff
\ldots\qff \cap\qff
e\qff(\dff r\trf)
\off =\off
\bigcap_{i\qff =\off 0}^r\qff
c\qff(\trf
u_{\dff 1}\dff(\dff i\trf),\off
u_{\dff 2}\dff(\dff i\trf),\off
\ldots,\off
u_{\dff n}\dff(\dff i\trf)
\qff)
\qff.
\]

\vspace{-10.5pt}
Clearly\halfff,\oss the\dss latter\dss intersection\dss is\dss a\sss product\sss of\sss $n$ intervals,\oss
namely\halfff,\oss of\trs the intervals\vspace{1.5pt}
\[
\quad
U_{\fff k}
\off =\off
\bigcap_{i\qff =\qff 0}^r\off
[\trf 
u_{\dff k}\dff(\dff i\trf)\dff,\off 
u_{\dff k}\dff(\dff i\trf)
\qff +\qff
1
\qff]
\]

\vspace{-12pt}
with\qss $k\off =\off 1\fff,\pff 2\fff,\pff \ldots\fff,\pff n$\nnsp.\oss

Let\dss us\dss temporarily\dss fix\dss $k$\dss and write\dss $u\trf(\dff i\trf)$\dss
for\dss $u_{\dff k}\dff(\dff i\trf)$\nsp.\oss
For an\dss integer\dss $i$\dss between $1$ and\dss $r$\dss let\dss
$\sigma\dff(\dff i\trf)$\dss be\sss the sum of\trs the numbers
$\varepsilon_{\fff j}$
such\dss that\qss $j\qff \geq\qff k\qff +\qff 1$\qss
and\qss\vspace{1.5pt}
\[
\quad
a_{\dff j}\dff(\dff i\trf)
\off =\off 
a_{\dff j}\dff(\dff i\qff -\qff 1\trf)\qff +\qff 1
\qff.
\]

\vspace{-10.5pt}
The sums\dss $\sigma\dff(\dff i\trf)$\dss are negative and\dss 
since\dss $v\dff(\dff i\dff)$\dss is\dss a\sss pivot subsequence,\oss
different\sss sums\dss $\sigma\dff(\dff i\trf)$\dss 
involve different\dss numbers\dss $\varepsilon_{\fff j}$\nsp.\oss
Together with our assumptions about\dss $\varepsilon_{\fff j}$\dss
this implies\sss that\vspace{3pt}
\begin{equation*}
\quad
\sigma\dff(\dff s\qff +\qff 1\trf)\qff +\qff \sigma\dff(\dff s\qff +\qff 2\trf)\qff +\qff
\ldots\qff +\qff
\sigma\dff(\dff t\trf)\qff +\qff 1
\off >\off 
0
\qff.
\end{equation*}

\vspace{-9pt}
if\pss $0\qff \leq\qff s\qff <\qff t\qff \leq\qff r$\nnsp.\oss
Suppose\sss now\dss that\qss\vspace{1.5pt} 
\[
\quad
a_{\dff k}\dff(\dff s\trf)
\off =\off 
a_{\dff k}\dff(\dff s\qff +\qff 1\trf)
\off =\off 
\ldots
\off =\off
a_{\dff k}\dff(\dff t\trf)
\qff.
\]

\vspace{-10.5pt}
Then\dss the interval\vspace{1.5pt}
\[
\quad
[\trf 
u\trf(\dff i\trf)\dff,\off 
u\trf(\dff i\trf)
\qff +\qff
1
\qff]
\]

\vspace{-10.5pt}
with\qss $s\qff +\qff 1\qff \leq\qff i\qff \leq\qff t$\qss
is obtained\dss by\sss shifting\dss the interval\qss\vspace{1.5pt}
\[
\quad
[\trf 
u\trf(\dff i\qff -\qff 1\trf)\dff,\off 
u\trf(\dff i\qff -\qff 1\trf)
\qff +\qff
1
\qff]
\]

\vspace{-10.5pt}
by\dss the amount\dss $\num{\sigma\dff(\dff i\trf)}$\dss to\sss the left\halfff.\oss
The\sss total\sss amount\sss of\qss shifting\qss\vspace{4pt}
\[
\quad
\mbox{from}\hspace*{1.2em}
[\trf 
u\trf(\dff s\trf)\dff,\off 
u\trf(\dff s\trf)
\qff +\qff
1
\qff]
\hspace*{1.2em}\mbox{to}\hspace*{1.2em}
[\trf 
u\trf(\dff t\trf)\dff,\off 
u\trf(\dff t\trf)
\qff +\qff
1
\qff]
\]

\vspace{-8pt}
is\dss equal\dss to\qss
$\num{\sigma\dff(\dff s\qff +\qff 1\trf)\qff +\qff \sigma\dff(\dff s\qff +\qff 2\trf)\qff +\qff
\ldots\qff +\qff
\sigma\dff(\dff t\trf)}
\off <\off 1$\nnsp.\qff\oss
It\dss follows\dss that\vspace{1.5pt}
\begin{equation}
\label{no-change-intersection}
\quad
\bigcap_{i\qff =\qff s}^t\off
[\trf 
u\trf(\dff i\trf)\dff,\off 
u\trf(\dff i\trf)
\qff +\qff
1
\qff]
\off =\off
[\qff
u\trf(\dff s\trf)\dff,\off 
u\trf(\dff t\trf)
\qff +\qff
1
\qff]
\qff
\end{equation}

\vspace{-9pt}
and\dss this interval\dss is\dss non-empty\dss and,\pss
moreover\halfff,\pss proper\halfff.\oss

In\dss particular\halfff,\oss
if\pss\vspace{3pt}
\[
\quad
a_{\dff k}\dff(\dff 0\trf)
\off =\off 
a_{\dff k}\dff(\dff 1\trf)
\off =\off 
\ldots
\off =\off
a_{\dff k}\dff(\trf r\trf)
\qff,
\]

\vspace{-9pt}
then\dss 
$U_{\fff k}$\dss is equal\dss to\sss the proper\dss interval\qss
$[\qff
u\trf(\trf r\trf)\dff,\off 
u\trf(\dff 0\trf)
\qff +\qff
1
\qff]$\nnsp.\oss

Suppose now\dss that\dss not\sss all\dss values\dss $a_{\dff k}\dff(\dff i\trf)$\dss
are equal.\oss
Then\dss there\dss is\sss a unique $m$ such\dss that\vspace{3pt}
\[
\quad
a_{\dff k}\dff(\dff 0\trf)
\off =\off 
a_{\dff k}\dff(\dff 1\trf)
\off =\off 
\ldots
\off =\off
a_{\dff k}\dff(\dff m\qff -\qff 1\trf)
\qff,
\]

\vspace{-36pt}
\[
\quad
a_{\dff k}\dff(\dff m\trf)
\off =\off 
a_{\dff k}\dff(\dff m\qff -\qff 1\trf)\qff +\qff 1
\qff,
\hspace*{1.2em}\mbox{and}
\]

\vspace{-36pt}
\[
\quad
a_{\dff k}\dff(\dff m\trf)
\off =\off 
a_{\dff k}\dff(\dff m\qff +\qff 1\trf)
\off =\off 
\ldots
\off =\off
a_{\dff k}\dff(\dff r\trf)
\qff.
\]

\vspace{-12pt}
Clearly\halfff,\vspace{-1.5pt}
\[
\quad
u\trf(\dff m\trf)
\off =\off
u\trf(\dff m\qff -\qff 1\trf)
\qff +\qff
\sigma\dff(\dff m\dff )
\qff +\qff 1
\]

\vspace{-37.5pt}
\[
\quad
\phantom{u\trf(\dff m\trf)
\off }
=\off
u\trf(\dff 0\dff)
\qff +\qff
\sigma\dff(\dff 1\dff)
\qff +\qff
\ldots
\qff +\qff
\sigma\dff(\dff m\qff -\qff 1\dff )
\qff +\qff
\sigma\dff(\dff m\dff )
\qff +\qff 1
\qff.
\]

\vspace{-15.25pt}
It\dss follows\dss that\qss\vspace{1.5pt}
\[
\quad
u\trf(\dff 0\dff)
\off <\off
u\trf(\dff m\trf)
\off \leq\off
u\trf(\dff m\qff -\qff 1\dff)
\qff +\qff
1
\]

\vspace{-9pt}
and\qss
$u\trf(\dff m\trf)
\off =\off
u\trf(\dff m\qff -\qff 1\dff)
\qff +\qff
1$\qss
if\qss and\dss only\qss if\pss
$\sigma\dff(\dff m\dff )
\off =\off
0$\nnsp.\oss
Also,\vspace{3pt}
\[
\quad
u\trf(\dff r\dff)
\off =\off
u\trf(\dff m\trf)
\qff +\qff
\sigma\dff(\dff m\qff +\qff 1\trf)\qff +\qff 
\ldots\qff +\qff
\sigma\dff(\dff r\trf)
\]

\vspace{-36pt}
\[
\quad
\phantom{u\trf(\dff r\dff)
\off }
=\off
u\trf(\dff m\qff -\qff 1\trf)
\qff +\qff
1
\qff +\qff
\sigma\dff(\dff m\trf)
\qff +\qff
\sigma\dff(\dff m\qff +\qff 1\trf)\qff +\qff 
\ldots\qff +\qff
\sigma\dff(\dff r\trf)
\off >\off
u\trf(\dff m\qff -\qff 1\trf)
\qff.
\]

\vspace{-9pt}
By applying\qss (\ref{no-change-intersection})\qss with\qss $s\off =\off 0$\qss
and\qss $t\off =\off m\qff -\qff 1$\qss we see\sss that\vspace{1.5pt}
\[
\quad
\bigcap_{i\qff =\qff 0}^{m\qff -\qff 1}\off
[\trf 
u\trf(\dff i\trf)\dff,\off 
u\trf(\dff i\trf)
\qff +\qff
1
\qff]
\off =\off
[\qff
u\trf(\dff 0\trf)\dff,\off 
u\trf(\dff m\qff -\qff 1\trf)
\qff +\qff
1
\qff]
\qff.
\]

\vspace{-9pt}
Similarly\halfff,\oss by applying\qss (\ref{no-change-intersection})\qss 
with\qss $s\off =\off m$\qss
and\qss $t\off =\off r$\qss we see\sss that\vspace{3pt}
\[
\quad
\bigcap_{i\qff =\qff m}^r\off
[\trf 
u\trf(\dff i\trf)\dff,\off 
u\trf(\dff i\trf)
\qff +\qff
1
\qff]
\off =\off
[\qff
u\trf(\dff m\trf)\dff,\off 
u\trf(\dff r\trf)
\qff +\qff
1
\qff]
\qff
\]

\vspace{-9pt}
By\dss combining\dss the last\dss two displayed\dss 
formulas conclude\sss that\vspace{3pt}
\[
\quad
U_{\fff k}
\off =\off
\bigcap_{i\qff =\qff 0}^r\off
[\trf 
u\trf(\dff i\trf)\dff,\off 
u\trf(\dff i\trf)
\qff +\qff
1
\qff]
\off =\off
[\qff
u\trf(\dff 0\trf)\dff,\off 
u\trf(\dff m\qff -\qff 1\trf)
\qff +\qff
1
\qff]
\qff \cap\qff
[\qff
u\trf(\dff m\trf)\dff,\off 
u\trf(\dff r\trf)
\qff +\qff
1
\qff]
\qff.
\]

\vspace{-9pt}
The inequalities\qss
$u\trf(\dff 0\dff)
\qff <\qff
u\trf(\dff m\trf)
\qff \leq\qff
u\trf(\dff m\qff -\qff 1\dff)
\qff +\qff
1
\qff <\qff
u\trf(\dff r\dff)
\qff +\qff
1$\nnsp,\oss
proved\sss above,\oss imply\dss that\vspace{3pt}
\[
\quad
[\qff
u\trf(\dff 0\trf)\dff,\off 
u\trf(\dff m\qff -\qff 1\trf)
\qff +\qff
1
\qff]
\qff \cap\qff
[\qff
u\trf(\dff m\trf)\dff,\off 
u\trf(\dff r\trf)
\qff +\qff
1
\qff]
\off =\off
[\qff
u\trf(\dff m\trf)\dff,\off
u\trf(\dff m\qff -\qff 1\trf)
\qff +\qff
1
\qff]
\qff
\]

\vspace{-9pt}
and\dss hence\qss
$U_{\fff k}
\off =\off
[\qff
u\trf(\dff m\trf)\dff,\off
u\trf(\dff m\qff -\qff 1\trf)
\qff +\qff
1
\qff]$\nnsp.\oss
Since\qss
$u\trf(\dff m\trf)
\qff \leq\qff
u\trf(\dff m\qff -\qff 1\dff)$\nnsp,\oss
this interval\dss is\dss non-empty\halfff.\oss
It\dss is\dss proper\dss if\qss and\dss only\trs if\pss
$u\trf(\dff m\trf)
\qff <\qff
u\trf(\dff m\qff -\qff 1\dff)
\qff +\qff
1$\nnsp,\oss
or\halfff,\oss
equivalently\halfff,\oss
$\sigma\dff(\dff m\dff )
\off <\off
0$\nnsp.\oss
The\sss latter\sss condition is equivalent\dss to\qss
$a_{\dff j}\dff(\dff m\qff +\qff 1\trf)
\off =\off 
a_{\dff j}\dff(\dff m\trf)\qff +\qff 1$\qss
for some\qss $j\qff \geq\qff k\qff +\qff 1$\nnsp.\oss
Equiv\-a\-lent\-ly\halfff,\oss the interval\dss $U_{\dff k}$\dss
is\dss proper\dss if\qss and\dss only\trs if\pss
$a_{\dff k}\dff(\dff i\dff)$\dss increases at\dss the same step as some\dss
$a_{\fff j}\dff(\dff i\dff)$\dss with\qss $j\qff >\qff k$\nnsp.\oss

Since\dss $k$\dss was arbitrary\halfff,\oss this proves\sss that\dss the intersection\qss
$e\qff(\dff 0\trf)\qff \cap\qff
e\qff(\dff 1\trf)\qff \cap\qff
\ldots\qff \cap\qff
e\qff(\dff r\trf)$\qss
is\dss non-empty\sss and\dss is\dss a product\sss of\dss closed\dss intervals.\oss
It\dss remains\sss to determine\sss the number of\dss proper\sss intervals among\dss them.\oss
Let\dss us\dss return\dss to\sss treating\dss $k$\dss as an arbitrary\dss integer\dss
between $1$ and $n$\nnsp.\oss
As we saw,\oss 
if\dss
$a_{\dff k}\dff(\dff i\dff)$\dss is\dss independent\sss of\dss $i$\nnsp,\oss
then\dss the interval\dss $U_{\dff k}$\dss is\dss proper\halfff.\oss
Also,\oss 
let\vspace{1.5pt}
\[
\quad
a_{\dff k_{\dff 1}}\dff(\dff i\dff)\dff,\qff\off
a_{\dff k_{\dff 2}}\dff(\dff i\dff)\dff,\qff\off
\ldots\dff,\qff\off
a_{\dff k_{\trf N}}\dff(\dff i\dff)
\]

\vspace{-10.5pt}
be\sss the coordinates increasing\dss when\dss $i$\dss changes from\qss $m\qff -\qff 1$\qss to\dss $m$\nnsp.\oss
If\qss $N\qff >\qff 1$\nnsp,\oss
then\dss\vspace{1.5pt} 
\[
\quad
U_{\dff k_{\dff 1}}
\dff,\qff\off
U_{\dff k_{\dff 2}}\dff(\dff i\dff)\dff,\qff\off
\ldots\dff,\qff\off
U_{\dff k_{\trf N\dff -\dff 1}}\dff(\dff i\dff)
\]

\vspace{-10.5pt}
are proper\dss intervals.\oss 
There are no other\dss proper\dss intervals.\oss
Let\dss us\sss look\dss how\dss the sum\vspace{3pt}
\[
\quad
\Sigma\dff(\dff i\trf)
\off =\off
a_{\dff 1}\dff(\dff i\trf)\qff +\qff
a_{\dff 2}\dff(\dff i\trf)\qff +\qff
\ldots\qff +\qff
a_{\dff n}\dff(\dff i\trf)
\qff
\]

\vspace{-9pt}
changes when $i$ changes from $0$ to $r$\dnsp.\oss
Each $a_{\dff k}\dff(\dff i\trf)$ may\sss stay\sss constant\sss 
or\dss increase\dss by $1$ and\dss hence\qss
$\Sigma\dff(\dff r\trf)
\off =\off
\Sigma\dff(\dff 0\dff)\qff +\qff n\qff -\qff c$\nnsp,\oss
where $c$ is\dss the number of\dss coordinates\dss 
staying\sss constant\halfff.\oss 
At\dss the same\sss time\qss
$\Sigma\dff(\dff m\trf)
\off =\off
\Sigma\dff(\dff m\qff -\qff 1\trf)
\qff +\pff N_{\dff m}$\nsp,\oss
where\dss $N_{\dff m}$\dss is\dss the number of\dss coordinates increasing\dss
from\qss $m\qff -\qff 1$\qss to\dss $m$\nnsp.\oss
It\dss follows\dss that\qss 
$n\qff -\qff c
\off =\off
N_{\dff 1}\qff +\qff
N_{\dff 2}\qff +\qff
\ldots\qff +\qff
N_{\dff r}$\qss
and\dss hence\qss\vspace{1.5pt}
\[
\quad
n\qff -\qff r
\off =\off
(\qff N_{\dff 1}\qff -\qff 1\qff)\qff +\qff
(\qff N_{\dff 2}\qff -\qff 1\qff)\qff +\qff
\ldots\qff +\qff
(\qff N_{\dff r}\qff -\qff 1\qff)\qff +\qff
c
\qff.
\]

\vspace{-10.5pt}
But\dss the right\dss hand of\trs this equality\dss is\dss equal\dss to\sss the number of\dss proper\dss
intervals\sss in\dss the product\halfff.\oss
It\dss follows\sss that\dss this number\dss is\sss equal\dss to\qss $n\qff -\qff r$\dnsp.\oss
This completes\sss the proof\halfff.\oss  \eproof

\mypar{Corollary\halfff.}{tiling-combinatorics}
\emph{Under\dss the above assumptions,\oss the intersection of\qss several\sss cubes of\qss the\dss form\dss
$e\qff(\dff a_{\dff 1}\fff,\pff a_{\dff 2}\fff,\pff \ldots\fff,\off a_{\dff n} \dff)$\dss
is\dss non-empty\dss if\qss and\dss only\trs if\pss
the corresponding\dss $n$\dnsp-tuples\qss
$(\dff a_{\dff 1}\fff,\pff a_{\dff 2}\fff,\pff \ldots\fff,\off a_{\dff n} \dff)$\qss
are\dss the\dss terms of\qss a\dss pivot\dss subsequence.\oss
In\dss particular\halfff,\oss
the intersection of\qss
$n\qff +\qff 2$\qss distinct\sss cubes of\qss the\dss form\dss
$e\qff(\dff a_{\dff 1}\fff,\pff a_{\dff 2}\fff,\pff \ldots\fff,\off a_{\dff n} \dff)$\dss
is\dss empty\halfff.\oss}  \eproof

\myuppar{Remark.}
Theorem\qss \ref{rotating-intersections}\qss and\dss
Corollary\qss \ref{tiling-combinatorics}\qss
provide another\dss proof\dss of\qss Lebesgue\dss tiling\dss theorem.\oss

\newpage
\mysection{Modification\qss of\pss Lebesgue\qss methods\qss by\pss W.\qss Hurewicz}{wh}

\myuppar{Tilings.}
The intersections of\dss $n$\dnsp-cubes of\dss $Q$ are rather excessive,\oss
if\dss compared\dss with\dss Lebesgue\dss covering\dss theorems\fff:\oss
vertices of\dss $n$\dnsp-cubes,\oss except\sss of\dss ones belonging\dss
to\sss the boundary\dss $\bd\fff Q$\nnsp,\oss are contained\dss in\dss $2^n$\sss
different\sss $n$\dnsp-cubes.\oss
At\dss the same\sss time,\oss Lebesgue\trs tilings consist\sss of\dss similar cubes
with no more\sss than\dss $n\qff +\qff 1$\dss  
intersecting\halfff.\oss
W.\dss Hurewicz\qss \cite{h2}\qss suggested\dss to start\dss with a\sss partition
of\dss a cube or\sss of\dss a product\sss of\dss intervals 
already\dss having\dss this property\halfff.\oss

Following\qss W.\dss Hurewicz\halfff,\oss
let\dss us call\dss \emph{$n$\dnsp-intervals}\oss the products of\trs the form\vspace{-1.5pt}
\[
\quad
Z
\off\qff =\off
\prod_{\dff i\qff =\qff 1}^{\dff n}\qff
J_{\fff i}
\off \subset\off
\rrr^{\fff n}
\qff,
\]

\vspace{-12pt}
where\qss 
$J_{\fff i}
\off =\off
[\trf a_{\dff i}\fff,\pff b_{\dff i}\trf]$\qss 
for some\qss
$a_{\dff i}\qff <\qff b_{\dff i}$\qss
for every\qss 
$i\off =\off 1\fff,\pff 2\fff,\pff \ldots\fff,\pff n$\nnsp.\oss
More generally\halfff,\oss 
an\dss \emph{$m$\dnsp-interval}\oss is\dss a product\sss of\trs the same form,\oss
but\dss with\sss $m$ factors\dss 
$J_{\fff i}$\dss
being intervals\dss $[\trf a_{\dff i}\fff,\pff b_{\dff i}\trf]$\dss 
with\qss $a_{\dff i}\qff <\qff b_{\dff i}$\qss 
as before
and\dss the other\qss $n\qff -\qff m$\qss factors being\dss one-point\dss
subsets of\trs $\rrr$\nnsp.\oss
By\dss an\qss \emph{interval}\pss we will\dss understand an $m$\dnsp-interval\dss
for some $m$\nnsp.\oss
The usual\dss intervals will\dss be called\qss \emph{intervals\dss in}\qss $\rrr$\nnsp.\oss
The\qss \emph{faces}\qss and\qss the\qss \emph{boundary}\dss $\bd\fff Z$\dss 
of\dss an interval\dss $Z$\dss
are defined\dss in\dss the obvious\sss manner\halfff.

For\dss the rest\sss of\trs this section we will\sss assume\dss that\sss
an $n$\dnsp-interval\dss $Z$\dss as above is fixed.\oss
We will\sss denote by\dss $A_{\dff i}$\dss and\dss  $B_{\dff i}$\dss
the $(\fff n\dff -\dff 1\fff)$\dnsp-faces of\trs $Z$\dss obtained\dss by\dss
replacing\dss in\dss the above product\dss the factor\qss
$[\trf a_{\dff i}\fff,\pff b_{\dff i}\trf]$\qss by\dss $\{\trf a_{\dff i}\trf\}$\dss
and\dss $\{\trf b_{\dff i}\trf\}$\dss respectively\halfff.\oss
A\qss \emph{partition}\qss of\trs $Z$\dss is\dss a\sss finite collection of\dss
$n$\dnsp-intervals such\dss that\dss $Z$\dss is\dss equal\dss to\sss their union
and\dss their\dss interiors are pair-wise disjoint\halfff.\oss
A partition of\trs $Z$\dss is\dss called\sss a\qss \emph{tiling}\pss
if\qss its\dss intervals satisfy\dss the following\dss two conditions\fff:\vspace{3.75pt}

\hspace*{1em}({\fff}a{\fff})\oss 
No point\dss belongs\sss 
to more\sss than\qss $n\qff +\qff 1$\qss
of\qss intervals.\oss\vspace{3.75pt}

\hspace*{1em}({\fff}b{\fff})\oss
The intersection of\dss $k$\sss intervals\dss is\dss
either\sss empty\halfff,\oss
or\trs is\dss an $(\fff n\dff +\qff 1\dff -\dff k\fff)$\dnsp-interval.\oss\vspace{3.75pt}

From\dss now\sss on\sss we will denote\sss by\qss
$Z_{\dff 1}\fff,\off Z_{\dff 2}\dff,\off \ldots\fff,\off Z_{\dff m}$\qss
a\dss tiling\sss of\dss $Z$\dss
and\dss call\dss its\dss elements\dss $Z_{\dff i}$\dss \emph{tiles}.

\myuppar{Tilings\dss into\sss small\sss $n$\dnsp-intervals.}
Of\dss course,\oss one needs\sss to know\dss that\dss tilings consisting of\dss 
arbitrarily\sss small $n$\dnsp-intervals exist\halfff.\oss
Theorem\qss \ref{rotating-intersections}\qss 
and\qss Corollary\qss \ref{tiling-combinatorics}\qss
imply\sss that\qss Lebesgue\dss tilings satisfy\sss conditions\qss
({\fff}a{\fff})\qss and\qss ({\fff}b{\fff}).\oss
By\dss scaling\dss and\dss intersecting\dss the scaled cubes with\dss $Z$\dss
one can get\sss a\sss tiling of\dss $Z$\dss by arbitrary small\dss tiles.\oss
It\dss seems\sss that\trs Hurewicz\dss wasn't\dss interested\dss in\dss
the detailed\sss description of\dss intersections\sss in\dss Lebesgue\dss tilings
provided\dss by\qss Theorem\qss \ref{rotating-intersections}\qss
(and\dss its proof\dff),\oss but\sss at\dss his disposal\dss was\qss
Theorem\qss \ref{tiling-intersections},\oss
due\sss to\dss Lebesgue\qss \cite{l-two}.\oss
Hurewicz\qss \cite{h2}\qss mentioned\sss in\sss a\sss footnote\sss
that\dss
if\dss no interval\sss of\dss a partition\dss intersects\sss two
opposite $(\fff n\dff -\dff 1\fff)$\dnsp-faces of\trs $Z$\nnsp,\oss
then\qss ({\fff}a{\fff})\qss implies\qss ({\fff}b{\fff}).\oss
While\sss this fact\dss together\dss with\dss 
Lebesgue\dss results\qss
({\fff}i.e.\qss together\dss with\qss
Theorem\qss \ref{tiling-intersections})\qss
also implies\sss that\qss Lebesgue\dss tilings satisfy\sss conditions\qss
({\fff}a{\fff})\qss and\qss ({\fff}b{\fff}),\oss
Hurewicz\dss did\dss not\dss wrote down a proof\halfff.\oss
Instead,\oss he outlined\qss (in another\dss footnote)\qss 
a\sss more simple and\dss general\dss way\dss to construct\dss
the needed\dss tilings.\oss
See\sss the next\dss theorem.\oss

\mypar{Theorem.}{hurewicz-tilings}
\emph{There exist\qss tilings\dss of\oss $Z$\qss
consisting\dss of\dss arbitrarily\dss small\dss $n$\dnsp-intervals.\oss}

\prooftitle{Hurewicz's\qss proof\halfff}
Let\qss
$(\dff x_{\dff 1}\fff,\pff x_{\dff 2}\fff,\pff \ldots\fff,\off x_{\dff n} \dff)$\qss
be\sss the standard coordinates in\dss $\rrr^{\fff n}$\dnsp.\oss
To begin\dss with,\pss the partition consisting only\sss of\qss $Z$\dss
is\dss trivially a\sss tiling.\oss
Let\qss
$Z_{\dff 1}\fff,\off Z_{\dff 2}\dff,\off \ldots\fff,\off Z_{\dff m}$\qss
be\sss a\sss tiling\sss of\dss $Z$\nnsp.\oss
Let\dss $Z_{\dff k}$\dss be an $n$\dnsp-interval\sss of\trs this\sss tiling\halfff,\oss
and\dss let\dss us divide\dss $Z_{\dff k}$\dss into\sss two $n$\dnsp-intervals\qss
\[
Z_{\dff k}'\qff,\off\off Z_{\dff k}''
\]
by\sss a\sss hyperplane\dss 
$H\pff \subset\off \rrr^{\fff n}$\dss defined\dss by\sss an equation of\trs
the form\qss $x_{\dff i}\off =\off a$\nnsp,\oss where\qss $a\qff \in\qff \rrr$\nnsp.\oss 
If\trs the hyperplane\dss $H$\dss does not\sss contain any $(\fff n\dff -\dff 1\fff)$\dnsp-face
of\dss any\dss tile\qss 
$Z_{\dff 1}\dff,\off Z_{\dff 2}\dff,\off \ldots\fff,\off Z_{\dff m}$\nsp,\oss 
then\dss replacing\dss $Z_{\dff k}$\dss by\dss the\sss pair\qss
$Z_{\dff k}'\dff,\qff\off Z_{\dff k}''$\qss results\dss in\sss
a\dss tiling\halfff.\oss
By\dss repeating\dss this procedure\sss with\dss hyperplanes of\dss various directions
one can get\dss tilings consisting of\dss arbitrarily small\dss intervals.\oss  \eproof

\mypar{Lemma.}{tiling-face}
\emph{Let\qss $B$\dss be an\dss $(\fff n\dff -\dff 1\fff)$\dnsp-face\dss of\pss $Z$\nnsp.\oss
The non-empty\dss sets among\dss the intersections\qss $Z_{\dff i}\qff \cap\qff B$\nnsp,\oss
where\qss
$i\off =\off 1\fff,\pff 2\fff,\pff \ldots\fff,\pff m$\nnsp,\qff\oss
form\dss a\dss tiling\dss of\pss $B$\nnsp.\oss}

\proof
If\qss $x\qff \in\pff Z_{\dff i}\qff \cap\qff B$\nnsp,\oss
then\dss $Z_{\dff i}$\dss contains a $1$\dnsp-interval\dss
having\dss $x$\dss as one of\dss its endpoints and orthogonal\dss to\dss $B$\nnsp.\oss
It\dss follows\dss that\dss taking\dss the intersection with\dss $B$\dss
decreases\sss the dimension of\dss any\dss interval\sss of\trs the form\qss
$Z_{\dff i_{\dff 1}}\off \cap\off
Z_{\dff i_{\dff 2}}\off \cap\off
\ldots\off \cap\off
Z_{\dff i_{\dff s}}$\qss
by $1$\nnsp.\oss
In\dss turn,\oss this implies\dss the lemma.\oss  \eproof

\mypar{Lemma.}{one-intervals}
\emph{Let\sss $z$\sss be\sss the intersection of\qss $n$ tiles.\oss
If\qss 
$z\off \neq\off \varnothing$\nnsp,\oss then\sss $z$\sss is\dss a\dss $1$\dnsp-interval\dss
such\dss that\dss each\sss of\pss its\dss endpoints\sss  
either\dss belongs\dss to\sss an $(\fff n\dff -\dff 1\fff)$\dnsp-face\dss of\pss $Z$\nnsp,\oss
or\dss is\dss contained\dss in\dss $n\qff +\qff 1$\dss
tiles.\oss
Other\dss points\sss of\dss $z$\sss are contained\dss in\sss only $n$
tiles.\oss}

\proof
Without\sss any\dss loss of\dss generality\sss we can assume\sss that\qss
$z\off =\off
Z_{\dff 1}\pff \cap\pff
Z_{\dff 2}\pff \cap\pff
\ldots\pff \cap\pff
Z_{\dff n}$\nsp.\oss
Suppose\sss that\dss $p$\dss is\dss an endpoint\sss of\dss $z$\dss not\sss contained\dss in\dss an
$(\fff n\dff -\dff 1\fff)$\dnsp-face\dss of\pss $Z$\nnsp.\oss
Let\dss $H$\dss be a\sss hyperplane\sss in\dss $\rrr^{\fff n}$\dss
containing\sss $p$\dss and orthogonal\trs to\dss $z$\nnsp,\oss
and\dss let\dss us extend $z$\sss to a $1$\dnsp-interval\sss $z'$\sss
containing\dss $p$\dss in\dss its\sss interior\halfff.\oss
Every\dss point\sss of\trs the difference\qss $z'\qff \smallsetminus\qff z$\qss
is not\dss contained\dss in\dss tiles\dss $Z_{\dff i}$\dss with\qss $i\qff \leq\qff n$\nnsp.\oss
It\dss follows\dss that\dss for some\dss $j\qff \leq\qff n$\qss 
all\dss points\sss of\qss $z'\qff \smallsetminus\qff z$\qss
sufficiently\sss close\sss to\dss $p$\sss do not\dss belong\dss to\dss $Z_{\fff j}$\nsp.\oss
Since\qss $p\qff \in\qff Z_{\fff j}$\nsp,\oss
the intersection\qss $Z_{\fff j}\pff \cap\pff H$\qss
is\dss an $(\fff n\dff -\dff 1\fff)$\dnsp-face of\pss $Z_{\fff j}$\nsp.\oss

Since\sss the\sss tiles form a partition,\oss
there is a\sss tile\dss $Z_{\dff k}$\dss 
different\dss from\dss $Z_{\dff j}$\dss and\sss
such\dss that\qss
$Z_{\dff k}\pff \cap\pff H$\qss
is\dss an $(\fff n\dff -\dff 1\fff)$\dnsp-face of\pss $Z_{\dff k}$\nsp.\oss
Clearly\halfff,\oss $Z_{\dff k}$\dss and\dss $Z_{\fff j}$\dss
cannot\dss be at\dss the same side of\qss $H$\nnsp.\oss 
It\dss follows\dss that\dss $Z_{\dff k}$\dss is different\dss from all\dss 
$Z_{\dff i}$\dss with\qss $i\qff \leq\qff n$\qss
and\dss hence\sss $p$\sss is\dss contained\dss in\qss $n\qff +\qff 1$\qss tiles.\oss

Suppose now\dss that\qss
$p\qff \in\qff Z_{\dff k}$\qss for some\qss $k\qff >\qff n$\nnsp.\oss
Then\dss $p$\dss is\dss the only\dss point\sss of\trs the intersection\qss
$Z_{\dff 1}\pff \cap\pff
Z_{\dff 2}\pff \cap\pff
\ldots\pff \cap\pff
Z_{\dff n}\pff \cap\pff
Z_{\dff k}
\off =\off
z\pff \cap\pff Z_{\dff k}$\nsp.\oss
But\dss the intersection of\dss a $1$\dnsp-interval\sss and\sss an $n$\dnsp-interval\dss
may consist\sss of\dss only\sss one point\sss only\dss if\trs this point\dss is\dss an
endpoint\sss of\trs the $1$\dnsp-interval.\oss  \eproof

\myuppar{Tiled\dss subsets\sss of\trs $Z$\nnsp.}
A\qss \emph{tiled}\pss subset\sss of\trs $Z$\dss is\sss defined as\sss the union
of\dss several\dss tiles.\oss
Tiled subsets of\trs $Z$\dss are\sss the analogues of\dss cubical\sss subsets
of\dss $Q$\dss from\dss Section\qss \ref{lebesgue}.\oss
Two\sss tiled subsets are said\dss to be\qss \emph{essentially\dss disjoint}\qss
if\dss they are unions of\trs two disjoint\sss collections of\trs tiles.\oss

\mypar{Lemma.}{intersection-path}
\emph{Suppose\dss that\pss
$e_{\dff 1}\fff,\off e_{\dff 2}\dff,\off \ldots\fff,\off e_{\dff n\dff +\dff 1}$\qss
is a covering of\pss $Z$\dss by\qss $n\qff +\qff 1$\qss tiled subsets,\oss
and\dss that\dss these\sss tiled subsets are
pairwise essentially\dss disjoint\halfff.\qff\oss
Let}\qss\vspace{0pt}
\[
\quad
P
\off =\off
e_{\dff 1}\pff \cap\pff
e_{\dff 2}\pff \cap\pff
\ldots\pff \cap\pff
e_{\dff n}
\qff.
\]

\vspace{-12pt}
\emph{Then\dss $P$\dss is\dss  the union of\dss several\pss $1$\dnsp-intervals\sss
$z_{\dff 1}\fff,\off z_{\dff 2}\dff,\off \ldots\fff,\off z_{\dff f}$\qss
intersecting\dss only\sss at\dss their endpoints and\sss
such\dss that\dss
each endpoint\sss of\dss each\dss interval\dss $z_{\dff i}$\dss
is\dss the end\-point\sss of\oss $1$\dss or\qss $2$\dss of\qss these intervals.\oss
If\qss $p$\dss is\dss the endpoint\sss of\dss only\sss one  interval\trs
and\qss $p\pff \not\in\off \bd\dff Z$\nnsp,\qff\oss
then\oss
$p
\pff \in\off
e_{\dff 1}\pff \cap\pff
e_{\dff 2}\pff \cap\pff
\ldots\pff \cap\pff
e_{\dff n\dff +\dff 1}$\nsp.\oss}

\proof
Let\dss us consider\dss intersections of\trs the form\qss
$Z_{\dff i_{\dff 1}}\off \cap\off
Z_{\dff i_{\dff 2}}\off \cap\off
\ldots\off \cap\off
Z_{\dff i_{\dff n}}$\qss
with\qss
$Z_{\dff i_{\dff k}}\pff \subset\pff e_{\dff k}$\qss
for every\qss $k\qff \leq\qff n$\nnsp.\oss
Let\qss
$z_{\dff 1}\fff,\off z_{\dff 2}\dff,\off \ldots\fff,\off z_{\dff f}$\qss
be\dss the non-empty\sss sets among\dss these intersections.\oss
Since\dss $E_{\dff i}$\dss are\sss tiled sets,\oss
$P$\dss is\dss equal\dss to\sss the union\qss
$z_{\dff 1}\off \cup\off
z_{\dff 2}\off \cup\off
\ldots\off \cup\off
z_{\dff f}$\nsp.\oss
Since\qss
$e_{\dff 1}\fff,\off e_{\dff 2}\dff,\off \ldots\fff,\off e_{\dff n\dff +\dff 1}$\qss
are pairwise essentially\dss disjoint\halfff,\oss
the property\qss ({\fff}b{\fff})\qss of\trs tilings implies\sss that\sss
each\dss $z_{\dff i}$\dss is a $1$\dnsp-interval.\oss
Let\dss $p$\dss be an endpoint\sss of\dss some such $1$\dnsp-interval\oss\vspace{1.5pt}
\[
\quad
z_{\dff t}
\off =\off
Z_{\dff i_{\dff 1}}\off \cap\off
Z_{\dff i_{\dff 2}}\off \cap\off
\ldots\off \cap\off
Z_{\dff i_{\dff n}}
\qff.
\]

\vspace{-10.5pt}
Suppose\sss that\qss $p\qff \in\pff \bd\dff Z$\nnsp.\oss
Lemma\qss \ref{tiling-face}\qss implies\sss that\dss $\bd\dff Z$\dss does not\sss
contain any of\dss the intervals\dss $z_{\dff i}$\nsp.\oss
It\dss follows\dss that\dss in\dss this case\dss $p$\dss is\dss the endpoint\sss
only\sss of\dss $z_{\dff t}$\nsp.\oss

Suppose\sss now\dss that\qss $p\qff \not\in\pff \bd\dff Z$\nnsp.\oss
Then\dss Lemma\qss \ref{one-intervals}\qss implies\sss that\dss there\dss is\dss a\sss tile\qss
$Z_{\dff i_{\dff n\dff +\dff 1}}$\dss containing\dss $p$\dss and\dss
not\sss equal\dss to any\qss $Z_{\dff i_{\dff k}}$\dss with\qss $k\qff \leq\qff n$\nnsp.\oss
By\dss the property\qss ({\fff}b{\fff})\qss of\trs tilings,\oss 
the intersection\vspace{1.75pt}
\[
\quad
Z_{\dff i_{\dff 1}}\off \cap\off
Z_{\dff i_{\dff 2}}\off \cap\off
\ldots\off \cap\off
Z_{\dff i_{\dff n}}\off \cap\off
Z_{\dff i_{\dff n\dff +\dff 1}}
\qff
\]

\vspace{-10.25pt}
is\dss a $0$\dnsp-interval\dss and\dss hence contains only one point\halfff,\oss
namely\dss $p$\nnsp.\oss
The property\qss ({\fff}a{\fff})\qss of\trs tilings implies\sss that
$p$\dss does not\dss belong\dss to any\dss tiles other\dss than\dss 
$Z_{\dff i_{\dff k}}$\dss with\qss $k\qff \leq\qff n\qff +\qff 1$\nnsp.\oss

There are\sss two cases\sss to consider\halfff.\oss
First\halfff,\oss suppose\dss that\qss $p\pff \in\off e_{\dff n\dff +\dff 1}$\nsp.\oss
Since\dss $e_{\dff n\dff +\dff 1}$\dss is\sss essentially\dss disjoint\dss from
each\dss $e_{\dff i}$\dss with\qss $i\qff \leq\qff n$\nnsp,\oss
in\dss this case\qss
$Z_{\dff i_{\dff n\dff +\dff 1}}\pff \subset\off e_{\dff n\dff +\dff 1}$\nsp.\oss
In view of\trs the observation at\dss the end of\trs the previous paragraph,\oss
this\dss implies\sss that\dss $p$\dss is\dss not\dss an endpoint\sss
of\dss any $1$\dnsp-interval\dss $z_{\dff i}$\dss different\dss from\dss $z_{\dff t}$\nsp.\oss
In other\dss words,\oss in\dss this case\dss $p$\dss is\dss also\sss the endpoint\sss
only\sss of\dss $z_{\dff t}$\nsp.\oss

If\qss $p\qff \not\in\qff e_{\dff n\dff +\dff 1}$\nsp,\oss
then\qss
$Z_{\dff i_{\dff n\dff +\dff 1}}\pff \subset\off e_{\dff k}$\qss
for some\qss $k\qff \leq\qff n$\qss
and\dss among\dss the\qss $n\qff +\qff 1$\qss tiles\vspace{1.75pt}
\[
\quad
Z_{\dff i_{\dff 1}}\dff,\qff\off
Z_{\dff i_{\dff 2}}\dff,\qff\off
\ldots\dff,\qff\off
Z_{\dff i_{\dff n}}\dff,\qff\off
Z_{\dff i_{\dff n\dff +\dff 1}}
\qff
\]

\vspace{-10.25pt}
there are exactly\dss $2$\dss sets of\dss
$n$\dss tiles including one\sss tile from each of\qss
$e_{\dff 1}\fff,\off e_{\dff 2}\dff,\off \ldots\fff,\off e_{\dff n}$\nsp,\oss
namely\dss\vspace{1.75pt}
\[
\quad
Z_{\dff i_{\dff 1}}\dff,\qff\off
Z_{\dff i_{\dff 2}}\dff,\qff\off
\ldots\dff,\qff\off
Z_{\dff i_{\dff n}}
\hspace*{1.5em}\mbox{and}\hspace*{1.5em}
Z_{\dff i_{\dff 1}}\dff,\qff\off
\ldots\dff,\qff\off
Z_{\dff i_{\dff k\dff -\dff 1}}\dff,\qff\off
Z_{\dff i_{\dff k\dff +\dff 1}}\dff,\qff\off
\ldots\dff,\qff\off
Z_{\dff i_{\dff n\dff +\dff 1}}
\qff.
\]

\vspace{-10.25pt}
The intersection of\dss tiles from\trs the first\sss set\dss is\dss $z_{\dff t}$\dss
and\dss the intersection of\dss tiles from\trs the second set\dss is\sss some other\dss 
$z_{\dff u}$\nsp.\oss
Clearly,\pss $p$\dss is\dss the endpoint\sss of\dss $z_{\dff t}$\dss
and\dss $z_{\dff u}$\nnsp,\oss but\sss of\trs no other\dss interval\dss $z_{\dff i}$\nsp.\oss
In\dss this case\dss $p$\dss is\dss the endpoint\sss of\dss exactly $2$ intervals\dss
$z_{\dff i}$\nsp.\oss
The lemma follows.\oss  \eproof

\mypar{Theorem.}{e-coverings-hurewicz}
\emph{Let\pss
$e_{\dff 1}\fff,\pff e_{\dff 2}\dff,\pff \ldots\fff,\pff e_{\dff n\dff +\dff 1}$\qss
be as\sss in\dss Lemma\qss \ref{intersection-path}.\oss
Suppose\sss that\dss $e_{\dff i}$\dss is\dss disjoint\dss from\dss $B_{\dff i}$\dss
for\dss every\qss $i\qff \leq\qff n$\qss
and\dss is\dss disjoint\dss from\dss $A_{\dff i\dff -\dff 1}$\dss for\dss
every\qss $i\qff \geq\qff 2$\nnsp.\qff\oss
Then\dss the\sss intersection\sss of\qss the\sss sets\pss
$e_{\dff 1}\fff,\pff e_{\dff 2}\dff,\pff \ldots\fff,\pff e_{\dff n\dff +\dff 1}$\qss
consists of\dss an odd\dss number of\qss points.\oss}

\proof
Let\dss us\sss use\sss the induction\dss by\dss $n$\nnsp.\oss
Suppose\sss that\qss $n\off =\off 1$\nnsp.\oss
Then\dss $Z$\dss is\sss a $1$\dnsp-interval\dss partitioned\dss into several\sss
$1$\dnsp-intervals with disjoint\dss interiors.\oss
There are\sss two sets,\oss $e_{\dff 1}$\dss and\dss $e_{\dff 2}$\dss
which are unions of\trs two disjoint\sss collections of\trs these intervals
and such\dss that\qss $Z\off =\off e_{\dff 1}\qff \cup\qff e_{\dff 2}$\nsp.\oss
It\dss is\dss assumed\dss that\sss one endpoint\sss of\trs $Z$\dss
belongs\sss to\dss $e_{\dff 1}$\dss and\dss the other\dss to\dss $e_{\dff 2}$\nsp.\oss
Therefore,\oss if\dss one moves along\dss $Z$\dss from one endpoint\dss to\sss
the other\halfff,\oss one has\sss to pass from\dss $e_{\dff 1}$\dss
to\dss $e_{\dff 2}$\dss or\dss vice versa an odd\dss number of\trs times.\oss
This means\sss that\qss $e_{\dff 1}\qff \cap\qff e_{\dff 2}$\qss
consists of\dss an odd\dss number of\dss points.\oss

Suppose now\dss that\qss $n\qff >\qff 1$\qss and\dss that\dss the\sss
theorem\dss is\dss already\dss proved\dss for\qss $n\qff -\qff 1$\qss
in\dss the role of\dss $n$\nnsp.\oss
By\dss Lemma\qss \ref{tiling-face}\qss the non-empty sets
among\dss the intersections\qss
$Z_{\dff i}\qff \cap\qff A_{\dff n}$\qss
form a\sss tiling of\dss the $(\fff n\dff -\dff 1\fff)$\dnsp-face\dss $A_{\dff n}$\dss
of\trs $Z$\nnsp.\oss
Since\dss $e_{\dff n\dff +\dff 1}$\dss is\dss disjoint\dss from\dss $A_{\dff n}$\nsp,\oss 
the sets\qss\vspace{2pt}
\[
\quad
e_{\dff 1}\qff \cap\qff A_{\dff n}\dff,\qff\off
e_{\dff 2}\qff \cap\qff A_{\dff n}\dff,\qff\off
\ldots\dff,\qff\off
e_{\dff n}\qff \cap\qff A_{\dff n}
\]

\vspace{-10.25pt}
cover\dss $A_{\dff n}$\nsp.\oss
Obviously\halfff,\oss they are\sss tiled subsets of\dss $A_{\dff n}$\dss
and are essentially\sss disjoint\halfff.\oss
By\dss the\dss inductive assumption\dss the intersection\qss 
$P\qff \cap\qff A_{\dff n}$\qss 
of\dss these sets
consists of\dss an odd\dss number of\dss points.\oss

At\dss the same\sss time\sss
the assumptions\qss of\trs the\sss theorem\dss imply\dss that\dss $P$\dss
is\dss disjoint\dss from\dss $B_{\dff i}$\dss for all\qss $i\qff \leq\qff n$\qss
and\dss is\dss disjoint\dss from\dss $A_{\dff i}$\dss for all\qss $i\qff \leq\qff n\qff -\qff 1$\nnsp.\oss
Therefore\qss
$P\qff \cap\qff \bd\dff Z
\off =\off
P\qff \cap\qff A_{\dff n}$\nsp.\oss

Let\qss
$z_{\dff 1}\fff,\off z_{\dff 2}\dff,\off \ldots\fff,\off z_{\dff f}$\qss
be\sss the intervals from\dss Lemma\qss \ref{intersection-path}.\oss
Let\dss $e$\dss and\dss $h$\dss be\sss the numbers of\dss endpoints of\trs these intervals
belonging\dss to\qss
$e_{\dff 1}\pff \cap\pff
e_{\dff 2}\pff \cap\pff
\ldots\pff \cap\pff
e_{\dff n\dff +\dff 1}$\qss
and\dss $A_{\dff n}$\dss respectively\halfff.\oss
By\dss the previous paragraph\dss $h$\dss is\sss also equal\dss to\sss the number
of\dss endpoints of\trs these intervals belonging\dss to\sss the boundary\dss $\bd\dff Z$\nnsp.\oss
Let\dss $r$\dss be\sss the number of\dss endpoints in\dss
the interior\qss 
$Z\qff \smallsetminus\qff \bd\dff Z$\qss of\trs $Z$\dss which are\sss endpoints of\dss
exactly\dss two of\trs the intervals.\oss
Now\qss Lemma\qss \ref{intersection-path}\qss implies\sss that\vspace{1.32pt}
\[
\quad
e\qff +\qff h\qff +\qff  2\dff r
\off =\off
2\dff f
\qff.
\]

\vspace{-10.68pt}
Since\dss $h$\dss is\dss odd\dss by\dss the inductive assumption,\oss
it\dss follows\dss that\dss $e$\dss is\dss also odd.\oss
Since\sss the sets\dss $e_{\dff i}$\dss are essentially\dss disjoint\halfff,\oss
every\dss point\sss of\trs the intersection\qss
$e_{\dff 1}\pff \cap\pff
e_{\dff 2}\pff \cap\pff
\ldots\pff \cap\pff
e_{\dff n\dff +\dff 1}$\qss
is\sss contained\dss in\qss $n\qff +\qff 1$\qss different\dss tiles,\oss
one from each\dss $e_{\dff i}$\nnsp.\oss
In view of\qss Lemma\qss \ref{one-intervals}\qss this implies\sss that\sss
every\dss point\sss of\trs this intersection\dss is\dss an\sss endpoint\sss
of\dss one of\dss the $1$\dnsp-intervals\dss $z_{\dff i}$\nsp.\oss
It\dss follows\dss that\dss $e$\dss is\dss the number of\dss points of\trs
this intersection,\oss and\dss hence\qss
$e_{\dff 1}\pff \cap\pff
e_{\dff 2}\pff \cap\pff
\ldots\pff \cap\pff
e_{\dff n\dff +\dff 1}$\qss 
consists of\dss an odd\dss number of\dss points.\oss
This completes\sss the step of\trs the induction and\dss hence\sss
the proof\halfff.\oss  \eproof

\mypar{Theorem.}{collecting-tiled-sets}
\emph{Let\qss
$d_{\dff 1}\fff,\pff d_{\dff 2}\dff,\pff \ldots\fff,\pff d_{\dff r}$\qss
be a covering\dss of\pss $Z$\dss
by\dss tiled\dss sets.\oss
Suppose\dss that\dss 
none of\qss the sets\dss 
$d_{\dff i}$\dss intersects\sss
two opposite\dss $(\fff n\dff -\dff 1\fff)$\dnsp-faces of\pss $Z$\nnsp.\oss
Then\dss among\dss the sets\qss 
$d_{\dff i}$\qss there are\qss 
$n\qff +\qff 1$\qss sets with non-empty\dss intersection.\oss}

\proof
It\dss is\dss sufficient\dss to combine\dss Theorem\qss \ref{e-coverings-hurewicz}\qss 
with\qss Lebesgue\dss fusion\sss construction.\oss  \eproof

\myuppar{Remark.}
At\dss the corresponding place of\trs his paper\qss \cite{h2}\qss
Hurewicz\dss wrote\sss that\qss\vspace{-10.5pt}

\begin{quoting}
One should\dss point\dss out\sss empathically\dss
the combinatorial\dss nature of\trs the preceding\sss considerations.\oss
The apparent\sss use of\trs the geometry\sss of\dss continuous figures served
only\dss the goal\sss of\dss simplifying\dss the exposition.
\end{quoting}

\vspace{-10.5pt}
See\qss \cite{h2},\oss footnote\dss 28\dss on\dss p.\qss 217.\oss
The combinatorial\dss nature of\qss Hurewicz's\dss proofs manifests itself\qss best\trs
if\trs they are\sss translated\dss into\sss the simplicial\dss language.\oss 
See\qss Section\qss \ref{nerve}.\oss

\myuppar{Lebesgue\dss first\dss covering\dss theorem.} 
\emph{Let\qss
$D_{\dff 1}\fff,\pff D_{\dff 2}\dff,\pff \ldots\fff,\pff D_{\dff r}$\qss
be a covering\dss of\pss the unit\sss cube\qss
$[\trf 0\fff,\qff 1\trf]^{n}$\qss
by closed\dss sets.\oss
Suppose\dss that\dss none of\qss the sets\dss $D_{\dff i}$\sss intersects\sss
two opposite\sss $(\fff n\dff -\dff 1\fff)$\dnsp-faces of\qss $[\trf 0\fff,\qff 1\trf]^{n}$\dnsp.\oss
Then\dss among\dss the sets\qss $D_{\dff i}$\qss there are\qss 
$n\qff +\qff 1$\qss sets with non-empty\dss intersection.\oss}

\prooftitle{Hurewicz's\qss proof\halfff}
The main\dss part\sss of\dss work\dss is\dss already\sss done.\oss
Given\qss Theorem\qss \ref{collecting-tiled-sets},\oss
the rest\sss of\dss
proof\dss is\sss almost\dss the same as\sss the proof\dss in\dss
Section\qss \ref{lebesgue}.\oss
Let\qss
$Z\off =\off [\trf 0\fff,\qff 1\trf]^{n}$\qss
and\dss let\dss us choose a\sss tiling of\trs $Z$\dss consisting of\trs tiles\sss
of\dss diameter\qss $<\qff \varepsilon$\nnsp,\oss
where $\varepsilon$ is\sss the same as\dss in\dss Section\qss \ref{lebesgue}.\oss
For each\qss $i\qff \leq\qff r$\qss let\dss $d_{\dff i}$\dss
be\sss the union of\dss all\dss tiles intersecting\dss $D_{\dff i}$\nsp.\oss
Now one can use\dss  Theorem\qss \ref{collecting-tiled-sets}\qss
instead of\qss Theorem\qss \ref{collecting-sets}.\oss  
We leave details\sss to\sss the reader\halfff.\oss  \eproof

\mypar{Theorem.}{hurewicz-lemma}
\emph{Let\pss
$d_{\dff 1}\fff,\off d_{\dff 2}\dff,\off \ldots\fff,\off d_{\dff n\dff +\dff 1}$\qss
be\sss a covering\dss of\pss $Z$\dss by\qss $n\qff +\qff 1$\qss tiled subsets.\oss
Suppose\sss that\dss $d_{\dff i}$\dss is\dss disjoint\dss from\dss $B_{\dff i}$\dss
for\dss every\qss $i\qff \leq\qff n$\qss
and\dss is\dss disjoint\dss from\dss $A_{\dff i\dff -\dff 1}$\dss for\dss
every\qss $i\qff \geq\qff 2$\nnsp.\qff\oss
Then\dss the\sss intersection\sss of\qss the\sss sets\pss
$d_{\dff 1}\fff,\pff d_{\dff 2}\dff,\pff \ldots\fff,\pff d_{\dff n\dff +\dff 1}$\qss
is\dss non-empty\halfff.\oss}

\proof
It\dss is\dss based\sss on\dss the same idea as\sss the proof\dss of\qss
Theorem\qss \ref{cubes-partitions}.\oss
For\qss $i\qff \leq\qff n\qff +\qff 1$\qss
let\dss $e_{\dff i}$\dss be\sss the union of\dss all\dss tiles contained\dss in\dss
$d_{\dff i}$\nnsp,\oss
but\dss not\dss 
in any\dss $d_{\fff j}$\dss with\qss
$j\qff <\qff i$\nnsp.\oss
Then\qss\vspace{1.125pt} 
\[
\quad
e_{\dff 1}\qff \cup\qff e_{\dff 2}\qff \cup\qff \ldots\qff \cup\qff e_{\dff i}
\off\qff =\off\qff
d_{\dff 1}\qff \cup\qff d_{\dff 2}\qff \cup\qff \ldots\qff \cup\qff d_{\dff i}
\]

\vspace{-10.875pt}
for every\dss $i$\nnsp.\oss
In\dss particular\halfff,\oss
$e_{\dff 1}\fff,\off e_{\dff 2}\dff,\off \ldots\fff,\off e_{\dff n\dff +\dff 1}$\qss
is a covering of\pss $Z$\dss by\qss $n\qff +\qff 1$\qss tiled subsets.\oss
By\dss the construction,\oss these sets are essentially\sss disjoint\halfff.\oss
On\dss the other\dss hand,\pss $e_{\dff i}\qff \subset\qff d_{\dff i}$\qss
and\dss hence
$e_{\dff i}$ is\dss disjoint\dss from $B_{\dff i}$
for\dss every\qss $i\qff \leq\qff n$\qss
and\dss is\dss disjoint\dss from $A_{\dff i\dff -\dff 1}$ for\sss
every\dss $i\qff \geq\qff 2$\nnsp.\oss
It\dss follows\dss that\dss the sets\qss
$e_{\dff 1}\fff,\pff e_{\dff 2}\dff,\pff \ldots\fff,\pff e_{\dff n\dff +\dff 1}$\qss
satisfy\dss the assumptions of\qss Theorem\qss \ref{e-coverings-hurewicz}\qss
and\dss hence\qss 
$e_{\dff 1}\qff \cap\qff 
e_{\dff 2}\qff \cap\qff
\ldots
\qff \cap\qff
e_{\dff n\dff +\dff 1}
\off \subset\off
d_{\dff 1}\qff \cap\qff 
d_{\dff 2}\qff \cap\qff
\ldots
\qff \cap\qff
d_{\dff n\dff +\dff 1}$\qss
is\dss non-empty\halfff.\oss  \eproof

\mypar{Theorem.}{hurewicz-theorem}
\emph{Let\pss
$D_{\dff 1}\fff,\off D_{\dff 2}\dff,\off \ldots\fff,\off D_{\dff n\dff +\dff 1}$\qss
be a covering of\pss the unit\sss cube\qss
$[\trf 0\fff,\qff 1\trf]^{n}$\qss
by closed\dss sets.\oss
Suppose\sss that\dss $D_{\dff i}$\dss is\dss disjoint\dss from\dss $B_{\dff i}$\dss
for\dss every\qss $i\qff \leq\qff n$\qss
and\dss is\dss disjoint\dss from\dss $A_{\dff i\dff -\dff 1}$\dss for\dss
every\qss $i\qff \geq\qff 2$\nnsp.\qff\oss
Then\dss the\sss intersection\sss of\qss the\sss sets\pss
$D_{\dff 1}\fff,\pff D_{\dff 2}\dff,\pff \ldots\fff,\pff D_{\dff n\dff +\dff 1}$\qss
is\dss non-empty\halfff.\oss}

\proof
It\dss is\dss completely\dss similar\dss to\dss Hurewicz's\dss proof\dss
of\qss Lebesgue\dss first\dss covering\dss theorem,\oss
with\qss Theorem\qss \ref{hurewicz-lemma}\qss
playing\dss the role of\qss Theorem\qss \ref{collecting-tiled-sets}.\oss  \eproof

\myuppar{Remark\halfff.}
This\dss is\dss a\sss stronger\dss version of\qss
Theorem\qss \ref{early-h-th}.\oss

\mysection{Nerves\qss of\pss coverings\qss and\qss Hurewicz's\qss theorems}{nerve}

\myuppar{The nerve of\dss a\sss finite system\sss of\dss sets.}
Let\vspace{1.5pt}
\[
\quad
M_{\dff 1}\dff,\off
M_{\dff 2}\dff,\off
\ldots\dff,\off
M_{\dff r}
\]

\vspace{-10.5pt}
be a system of\dss sets.\oss
In\dss the present\sss context\dss such system\dss usually\sss arises
as a covering of\dss a set\dss under consideration.\oss
The\qss \emph{nerve}\qss of\trs this system\dss is\dss an\sss abstract\sss
simplicial\sss complex constructed as follows.\oss
Its vertices\qss
$v_{\dff 1}\fff,\pff
v_{\dff 2}\fff,\pff
\ldots\fff,\pff
v_{\dff r}$\qss
are in a one-to-one correspondence with\dss the sets of\trs the system.\oss
A set\qss
$\{\qff
v_{\dff i_{\dff 1}}\fff,\pff
v_{\dff i_{\dff 2}}\fff,\pff
\ldots\fff,\pff
v_{\dff i_{\dff s}}
\qff\}$\qss
of\dss vertices\dss is\dss a\sss simplex\dss if\dss and\dss only\trs if\vspace{1.5pt}
\[
\quad
M_{\dff i_{\dff 1}}\qff \cap\qff
M_{\dff i_{\dff 2}}\qff \cap\qff
\ldots\qff \cap\qff
M_{\dff i_{\dff s}}
\off \neq\off
\varnothing
\qff.
\]

\vspace{-10.5pt}
This notion was introduced\dss in\qss 1926\qss by\qss P.\trs Alexandroff\pss \cite{a-nerve}.\oss
He had\dss in\dss mind\dss
the goal\sss of\dss relating\sss general\sss topological\sss spaces
with simplicial\sss complexes,\oss either abstract\sss ones,\oss
which are pure combinatorial\dss finite objects,\oss
or\dss geometric ones,\oss
which are geometric figures of\dss finite nature.\oss
The notion of\dss a nerve\sss turned out\dss to be extremely\dss fruitful\dss
far\dss beyond\dss this initial\qss (already\dss very\sss ambitious)\qss goal,\oss
but\dss we will\dss use it\sss only\sss as a convenient\dss language.\oss

Probably\halfff,\oss
the simplest\sss example\dss is\dss provided\dss by\dss the collection of\dss
$(\fff n\dff -\dff 1\fff)$\dnsp-faces of\dss a geometric $n$\dnsp-simplex\dss $\Delta$\nnsp.\oss
It\dss is\sss a covering\dss of\trs the boundary\dss $\bd\fff \Delta$\nnsp.\oss
The intersection of\dss all $(\fff n\dff -\dff 1\fff)$\dnsp-faces is empty\halfff,\oss
but\dss the intersection of\dss any\dss proper subset\sss of\dss the set\sss of\dss
$(\fff n\dff -\dff 1\fff)$\dnsp-faces\dss is\dss not\halfff.\oss
Therefore,\oss the nerve of\dss this covering\sss of\trs the boundary\dss  $\bd\fff \Delta$\dss
is\dss the boundary\sss of\dss an abstract\sss $n$\dnsp-simplex\qss
(considered as a pseudo-manifold\fff).\oss
The reader\dss may\dss try\dss to identify\dss the nerve of\trs the covering\sss
of\trs the boundary\sss of\dss an $n$\dnsp-cube by\dss its
$(\fff n\dff -\dff 1\fff)$\dnsp-faces.\oss
As a more general\sss example,\oss
let\dss us consider a geometric simplicial\sss complex\dss $S$\dss
and\dss the covering of\trs its polyhedron\dss $\norm{S}$\dss by\dss the closed\dss
barycentric stars.\oss
The discussion\dss in\dss Section\qss 3\qss of\pss \cite{i2}\qss
implies\sss that\dss the nerve of\dss this covering is isomorphic\sss to\sss the
corresponding abstract\sss simplicial\sss complex.\oss

\myuppar{Infinite systems of\dss sets.}
The notion of\dss the nerve easily\sss extends\sss to infinite systems of\dss sets,\oss
once\sss the notion of\dss an\sss abstract\sss simplicial\sss complex is extended\dss
to allow\dss infinite sets of\dss vertices.\oss
This\sss is\sss also easy\fff:\oss the sets of\dss vertices and simplices are allowed\dss
to be infinite,\oss but\dss the simplices are still\dss required\dss to be finite.\oss
So,\oss a potentially\dss infinite abstract\sss simplicial\sss complex\dss $K$\dss
is defined as a collection\sss $K$\sss of\oss \emph{finite}\qss 
subsets of\dss a set\sss $v\dff(\trf K\trf)$\sss
such\dss that\dss if\qss $\sigma\qff \in\qff K$\qss and\qss
$\sigma'\qff \subset\qff \sigma$\nnsp,\oss then\qss
$\sigma'\qff \in\qff K$\nnsp,\oss 
and\sss $v\dff(\trf K\trf)$\sss is equal\dss to\sss the union of\dss all\sss subsets\sss in\dss $K$\nnsp.\oss
The elements of\dss $v\dff(\trf K\trf)$\dss are called\qss \emph{vertices},\oss
and\dss the elements of\trs $K$\dss the\qss \emph{simplices}\qss of\trs $K$\nnsp.\oss

Now\dss the nerve of\dss an arbitrary\sss system\dss $\mathfrak{M}$\dss 
of\dss sets\sss is\dss defined\dss
in\dss the same way as before.\oss
Its vertices\dss $v_{\dff m}$\dss are\sss in\sss a one-to-one correspondence with\dss the sets of\trs
the system,\oss i.e.\qss with\dss the elements\qss $m\qff \in\qff \mathfrak{M}$\nnsp.\qff\oss 
A\dss set\qss
$\{\qff
v_{\dff m_{\dff 1}}\fff,\pff
v_{\dff m_{\dff 2}}\fff,\pff
\ldots\fff,\pff
v_{\dff m_{\dff s}}
\qff\}$\qss
of\dss vertices\dss is\dss a\sss simplex\dss if\dss and\dss only\trs if\vspace{1.5pt}
\[
\quad
m_{\dff 1}\qff \cap\qff
m_{\dff 2}\qff \cap\qff
\ldots\qff \cap\qff
m_{\dff s}
\off \neq\off
\varnothing
\qff.
\]

\vspace{-10.5pt}
\myuppar{The nerve of\qss Lebesgue\dss tilings.}
We are going\dss to discuss essentially\sss only\sss one infinite system of\dss sets,\oss
namely\halfff,\oss Lebesgue\dss tilings of\dss $\rrr^{\fff n}$\dnsp.\oss
We\sss keep\sss the notations and assumptions of\qss Section\qss \ref{tilings}.\oss
In\dss particular\halfff,\oss
we will\sss assume\sss that\trs the numbers\qss
$\varepsilon_{\dff 1}\fff,\pff \varepsilon_{\dff 2}\fff,\pff 
\ldots\fff,\off \varepsilon_{\dff n\dff -\dff 1}$\qss
are\sss generic,\oss such\dss that\qss
${}-\qff 1\qff <\qff \varepsilon_{\dff i}\qff <\qff 0$\qss
for every\qss $i\qff \leq\qff n\qff -\qff 1$\nnsp,\oss
and such\dss that\oss 
$\num{\varepsilon_{\dff 1}}\qff +\qff 
\num{\varepsilon_{\dff 2}}\qff +\qff
\ldots\qff +\qff 
\num{\varepsilon_{\dff n\dff -\dff 1}} 
\off <\off 1$\nnsp.\oss

The cubes\qss
$e\qff(\dff a_{\dff 1}\fff,\pff a_{\dff 2}\fff,\pff \ldots\fff,\off a_{\dff n} \dff)$\qss 
of\dss a\dss Lebesgue\sss tiling\dss are in\dss a canonical\sss one-to-one
correspondence with elements of\dss $\zzz^{\fff n}$\dnsp.\oss
Therefore,\oss the nerve of\dss a\dss Lebesgue\dss tiling\sss may\dss be considered
as an\dss infinite abstract\sss simplicial\sss complex with\dss $\zzz^{\fff n}$\dss
being\dss the set\sss of\dss vertices.\oss
Corollary\qss \ref{tiling-combinatorics}\qss 
provides a\sss complete combinatorial\sss description of\dss this nerve.\oss
Namely\halfff,\pss
a subset\sss of\trs $\zzz^{\fff n}$\dss is a simplex\dss
if\trs and\dss only\trs if\qss it\dss is\dss the set\sss of\trs terms of\dss
a pivot subsequence\qss (as defined\dss in\dss Section\qss \ref{tilings}).\oss
In\dss particular\halfff,\oss
under\dss the above assumptions about\dss the parameters\dss $\varepsilon_{\dff i}$\dss
the nerve is\sss independent\sss on\dss the choice of\dss these parameters.\oss

\myuppar{The nerve of\dss a\dss tiling\halfff.}
As in\dss Section\qss \ref{wh},\oss
let\dss $Z$\dss be an $n$\dnsp-interval\sss and\dss let\qss
$Z_{\dff 1}\fff,\off Z_{\dff 2}\dff,\off \ldots\fff,\off Z_{\dff m}$\qss
be a\dss tiling\sss of\dss $Z$\nnsp.\oss
We will\sss
assume\sss that\dss no\sss tile\dss $Z_{\dff i}$\dss
intersects\sss two opposite faces of\trs $Z$\nnsp.\oss
Let\dss $\mathcal{Z}$\dss be\sss the nerve of\trs the covering of\trs $Z$\dss
by\dss the\sss tiles,\oss
and\qss let\qss
$\mathcal{V}
\off =\off
\{\qff
v_{\dff 1}\fff,\pff
v_{\dff 2}\fff,\pff
\ldots\fff,\pff
v_{\dff m}
\qff\}$\qss
be\sss the set\sss of\trs its vertices,\pss
the vertex\dss $v_{\dff i}$\dss corresponding\dss to\sss the\sss tile\dss $Z_{\dff i}$\nsp.\oss
By\dss the definition of\trs tilings,\oss the dimension of\dss $\mathcal{Z}$\dss
is\qss $\leq\qff n$\qss
({\dff}Lebesgue covering\dss theorems\dss imply\dss that\dss it\dss is\sss equal\dss to\dss $n$\nsp).\oss 

For\dss each\sss integer $i$\trs between $1$ and $n$\dss let\dss
$\mathcal{A}_{\dff i}$\sss 
be\sss the nerve of\dss the system of\dss sets consisting\sss 
of\trs tiles intersecting\dss $A_{\dff i}$\nnsp.\oss 
Let\dss us define\trs $\mathcal{B}_{\dff i}$\dss in\dss the same manner\halfff,\oss
and\dss let\qss \emph{boundary}\qss $\bd\dff \mathcal{Z}$\dss be\sss 
the union of\dss simplicial\sss complexes\qss 
$\mathcal{A}_{\dff i}$\nsp,\qss $\mathcal{B}_{\dff i}$\qss 
over all\qss $1\qff \leq\qff i\qff \leq\qff n$\nnsp.\oss
The assumption\dss that\dss no\sss tile intersects\sss two opposite faces of\trs $Z$\dss
implies\sss that\dss $\mathcal{A}_{\dff i}$\dss is\dss disjoint\dss from\dss
$\mathcal{B}_{\dff i}$\dss for every\dss $i$\nnsp.\oss
Lemma\qss \ref{tiling-face}\qss  
implies\sss that\dss the dimension of\dss
$\bd\dff \mathcal{Z}$\dss is\qss $\leq\qff n\qff -\qff 1$\qss
({\fff}in\dss fact\halfff,\oss it\dss is\sss equal\dss to\qss
$n\qff -\qff 1$\nsp).\oss

By\dss the definition of\trs the simplicial\sss complex\dss $\mathcal{Z}$\nnsp,\qff\oss 
its $n$\dnsp-simplices 
correspond\dss to\sss non-empty\dss intersections of\qss $n\qff +\qff 1$\qss tiles,\oss
and\dss its $(\fff n\dff -\dff 1\fff)$\dnsp-\dnsp-simplices 
correspond\dss to non-empty\dss intersections of\dss $n$\dss tiles.\oss
The latter are exactly\dss the $1$\dnsp-intervals considered\dss in\qss
Lemma\qss \ref{one-intervals}.\oss 
Since we assumed\dss that\dss no\sss tiles
intersect\dss both\dss
$A_{\dff i}$\sss and\trs $B_{\dff i}$\nsp,\oss
these $1$\dnsp-intervals\sss cannot\dss have endpoints in\dss both\dss
$A_{\dff i}$\sss and\trs $B_{\dff i}$\nsp.\oss
If\dss such an $1$\dnsp-interval\dss has an endpoint\dss in\dss 
$A_{\dff i}$\sss or\trs $B_{\dff i}$\nsp,\oss
then\dss the\sss the corresponding $(\fff n\dff -\dff 1\fff)$\dnsp-\dnsp-simplex\dss
is\dss a\sss simplex of\qss
$\mathcal{A}_{\dff i}$\dss or\dss $\mathcal{B}_{\dff i}$\dss
respectively\halfff.\oss

Therefore,\pss
in\dss the language of\dss $\mathcal{Z}$\nnsp,\oss 
Lemma\qss \ref{one-intervals}\qss  
means\sss that\dss an
$(\fff n\dff -\dff 1\fff)$\dnsp-simplex 
of\dss $\mathcal{Z}$\dss
is\dss either contained\dss in\dss the boundary\dss $\bd\dff \mathcal{Z}$\dss
and\dss then\dss it\dss is\sss a face of\dss one $n$\dnsp-simplex,\oss
or\dss it\dss is\sss a face of\trs two $n$\dnsp-simplices.\oss
In other\dss words,\oss this lemma establishes an analogue of\dss
the non-branching\dss property\sss of\trs triangulations of\dss a geometric 
simplex\qss (see\qss \cite{i2},\oss Section\qss 2).\oss

\myuppar{The\sss combinatorics of\qss Hurewicz's\trs method.}
Let\qss
$e_{\dff 1}\fff,\off e_{\dff 2}\dff,\off \ldots\fff,\off e_{\dff n\dff +\dff 1}$\qss
be\dss a\sss covering of\trs the set\dss $\mathcal{V}$\dss
of\dss vertices of\qss $\mathcal{Z}$\dss 
by\qss $n\qff +\qff 1$\qss pairwise\sss
disjoint\sss subsets.\oss
Such a covering can\sss be interpreted as a\sss labeling of\dss $\mathcal{V}$\dnsp\dnsp,\oss
where\sss the label\sss of\trs the vertex\dss $v_{\dff i}$\dss is\dss the unique number\dss $k$\dss
such\dss that\qss $v_{\dff i}\qff \in\qff e_{\dff k}$\nsp.\oss
Even\dss better\halfff,\oss this covering\sss can\sss be interpreted as a\sss
simplicial\dss map\qss\vspace{3pt}
\[
\quad
\varphi\dff \colon\dff
\mathcal{Z}\qff \ttoo\qff \Delta
\qff,
\]

\vspace{-9pt}
where\dss $\Delta$\dss is\dss the simplicial\sss complex having\dss 
$\{\qff 1\fff,\pff 2\fff,\pff \ldots\fff,\pff n\qff +\qff 1\qff\}$\dss
as\sss the set\sss of\dss vertices and
all\sss sets of\trs
vertices as simplices\qss
({\fff}i.e.\qss the simplex\sss 
$\{\qff 1\fff,\pff 2\fff,\pff \ldots\fff,\pff n\qff +\qff 1\qff\}$\sss
considered as a simplicial\sss complex).\oss
It\dss is\dss worth\dss to point\sss out\dss that\dss these interpretations
depend on\dss the disjointness assumption.\oss

Let\dss 
$\Delta_{\dff n\dff +\dff 1}$\dss
be\sss the $(\fff n\dff -\dff 1\fff)$\dnsp-face\qss
$\{\qff 1\fff,\pff 2\fff,\pff \ldots\fff,\pff n\qff\}$\qss
of\dss $\Delta$\nnsp.\oss
Then\dss the $1$\dnsp-intervals\sss $z_{\dff k}$\sss of\qss Lemma\qss \ref{intersection-path}\qss
correspond\dss to $(\fff n\dff -\dff 1\fff)$\dnsp-simplices\dss $\tau$\dss
of\qss $\mathcal{Z}$\dss such\dss that\qss
$\varphi\dff(\dff \tau\dff)\off =\off \Delta_{\dff n\dff +\dff 1}$\nsp.\oss
In\dss the language of\dss $\mathcal{Z}$\dss Lemma\qss \ref{intersection-path}\qss means\sss
that\dss if\trs an $n$\dnsp-simplex\dss $\sigma$\dss has an
$(\fff n\dff -\dff 1\fff)$\dnsp-face\dss $\tau$\dss such\dss that\qss
$\varphi\dff(\dff \tau\dff)\off =\off \Delta_{\dff n\dff +\dff 1}$\nsp,\oss
then\sss either\sss $\sigma$\sss has\sss two such\dss faces,\oss
or\qss
$\varphi\dff(\dff \sigma\dff)\off =\off \Delta$\nnsp.\oss

Now\sss everything\dss is\dss ready\dss to carry\sss out\sss a counting\sss
argument\sss in\dss the spirit\sss of\qss Sperner's\dss one,\oss
and\dss this\sss is\sss done in\dss the proof\dss of\qss Theorem\qss
\ref{e-coverings-hurewicz}.\oss 
Let\qss $e\fff,\off h\fff,\off r$\qss and\dss $f$\dss be\sss the numbers
from\dss this proof\halfff,\oss
and\dss let\dss $g$\dss be\sss the number of\dss endpoints of\trs the intervals\dss
$z_{\dff i}$\dss contained\dss in\dss the interior of\trs $Z$\nnsp.\oss
Then\qss $r\off =\off g\qff -\qff e$\qss
and\dss the equality\qss\vspace{1.5pt}
\[
\quad
e\qff +\qff h\qff +\qff  2\dff r
\off =\off
2\dff f
\]

\vspace{-10.5pt}
from\dss this proof\dss
can\dss be rewritten as\qss
$e\qff +\qff h\qff +\qff  2\qff (\dff g\qff -\qff e\dff)
\off =\off
2\dff f$\nnsp,\oss
or\halfff,\oss
what\dss is\dss the same,\oss as\vspace{1.5pt}
\[
\quad
h\qff +\qff 2\dff g
\off =\off 
e\qff +\qff 2\dff f
\qff.
\]

\vspace{-10.5pt}
This\dss is\dss the same equality\sss as\sss in\trs Sperner's\trs proof\dss of\trs his\dss lemma.\oss
The\dss Hurewicz's\dss and\dss Sperner's\dss proofs\dss differ\dss in\dss the context\sss
and\dss the\sss language,\oss but\sss on\sss a\sss deeper\dss level\dss they are nearly\dss
the same.

\myuppar{The\sss methods\sss of\pss Hurewicz\qss and\qss Sperner\qss in\qss 1928.}
W.\dss Hurewicz\dss himself\dss was\sss well\sss aware of\trs the analogy\dss between\dss his
methods and\dss Sperner's\dss ones.\oss
He wrote\vspace{-9pt}

\begin{quoting}
During\dss the preparation of\dss the page proofs of\trs this paper\halfff,\oss
an article of\qss E.\dss Sperner\dss appeared\qss
(Abhandlungen des Hamburgischen Math.\qss Sem.\qss 6,\oss pp.\qss 265-472,\oss
submitted\dss in\dss June\qss 1928),\oss
in\dss which one\sss finds a proof\dss of\qss Lebesgue--Brouwer\dss theorems
exhibiting\sss a\sss far-reaching analogy with\dss the proof\dss presented\dss here.\oss
\end{quoting}

\vspace{-9pt}
See\qss \cite{h2},\oss footnote\dss 13a\dss on\dss p.\qss 211.\oss
Hurewicz\dss paper was submitted in\dss January\qss 1928,\oss
several\dss months earlier\dss than\dss Sperner's\dss paper\halfff,\oss
but\dss was published at\dss least\sss several\dss months\sss later\halfff.\oss
In\dss the meantime\dss
K.\dss Menger\dss included\dss both\dss proofs in\dss his\dss monograph\qss \cite{me}.\oss
K.\dss Menger\dss wrote\vspace{-9pt}

\begin{quoting}
The first\dss proof\dss is\dss based on\dss Sperner's\dss work\qss
(already\dss appeared\dss in\dss Hamburger Abhandl.)\qss
and\dss deals with \dnsp$n$\dnsp-dimensional\sss simplices.\oss
It\dss simplifies\dss Brouwer's\dss proof\trs by\sss using\sss
Brouwer's\sss notion of\sss simplicial\sss decompositions\qss
(Math.\qss Ann.\qss 71,\pss 1911,\pss p.\qss 161)\dss
and\dss suppressing\sss the notion
of\dss degrees of\dss maps.\oss
The second\dss proof\dss is\dss Hure\-wicz's\dss simplification\qss 
({\fff}to appear shortly\sss in\dss Math.\qss Ann.)\qss
of\qss Lebesgue\dss proof\halfff,\oss
dealing with \dnsp$n$\dnsp-dimensional\sss cubes and using\dss Lebesgue\dss
canonical\sss cubical\sss decompositions.
\end{quoting}

\vspace{-9pt}
In\qss 1928\qss the\sss two methods were standing on equal\dss footing,\oss
but\dss their\dss further\dss fates are different\halfff.\oss
Hurewicz's\qss beautiful\dss proof\trs wasn't\trs included even in his
book\qss \cite{hw}\qss with\dss H.\dss Wallman.\oss

\myuppar{Lebesgue\dss tilings\dss and\dss Hurewicz's\dss theorems.}
It\dss is\dss only\dss natural\dss to apply\dss Hurewicz's\dss theorems\sss
to\dss Lebesgue\dss tilings.\oss
This\sss leads\sss to purely\sss combinatorial\qss 
Theorems\qss \ref{wh-lebesgue-tilings}\trs and\qss \ref{hurewicz-lemma-nerve}\qss
below.\oss
A\-par\-ent\-ly\halfff,\pss
those\sss who\sss knew\dss 
Hurewicz's\dss work\dss 
were not\dss interested\dss in\sss such combinatorial\dss results,\oss
and\dss those\sss who were interested\dss in\dss them,\pss
such as\qss H.W.\qss Kuhn\qss \cite{ku},\pss
were not\sss aware of\qss his\dss work.

Let\dss $Z$\dss be an $n$\dnsp-interval.\oss
The intersections of\dss cubes of\dss
a\dss Lebesgue\dss tiling\dss with\dss $Z$\dss are $m$\dnsp-intervals with\qss $m\qff \leq\qff n$\nnsp,\oss
and\dss the intersections which are $n$\dnsp-intervals form a\sss tiling of\trs $Z$\nnsp.\oss
Let\dss $k$\dss be a natural\dss number\sss and\qss
$K\off =\off
\{\qff 0\fff,\pff 1\fff,\pff \ldots\fff,\pff k\qff\}^{\fff n}$\dnsp.\oss
We would\dss like\sss to choose an $n$\dnsp-interval\dss $Z$\dss and a\dss Lebesgue\dss tiling\dss in such a way\dss
that\dss the cubes\dss
$e\qff(\dff a_{\dff 1}\fff,\pff a_{\dff 2}\fff,\pff \ldots\fff,\off a_{\dff n} \dff)$\dss
with\qss 
$(\dff a_{\dff 1}\fff,\pff a_{\dff 2}\fff,\pff \ldots\fff,\off a_{\dff n} \dff)
\qff \in \qff K$\qss
would\dss be\sss the ones which are intersecting\dss $Z$\dss in $n$\dnsp-intervals.\oss
A natural\dss way\dss to do\sss this\dss is\dss the following.\oss
Let\dss $\varepsilon$\dss be a small\dss positive number\halfff.\oss
Assuming\dss that\qss $\varepsilon\qff <\qff 1$\qss is\dss sufficient\halfff.\oss
Let\vspace{4.75pt}
\[
\quad
Z
\off =\off
[\qff \varepsilon\dff,\off k\qff +\qff \varepsilon\qff]^{\fff n}
\qff.
\]

\vspace{-7.75pt}
Let\qss
$\varepsilon_{\dff 1}\fff,\pff \varepsilon_{\dff 2}\fff,\pff 
\ldots\fff,\off \varepsilon_{\dff n\dff -\dff 1}
\qff \in\qff
\rrr$\qss 
be\sss such\dss that\qss 
${}-\qff 1\qff <\qff \varepsilon_{\dff i}\qff <\qff 0$\qss
for every\dss $i$\dss and\vspace{4.5pt}
\[
\quad
k\qff \num{\varepsilon_{\dff 1}}\qff +\qff 
k\qff \num{\varepsilon_{\dff 2}}\qff +\qff 
\ldots\qff +\qff 
k\qff \num{\varepsilon_{\dff n\dff -\dff 1}}
\off <\off
1\qff -\qff \varepsilon
\qff.
\]

\vspace{-7.5pt}
In\dss particular\halfff,\pss
$\num{\varepsilon_{\dff 1}}\qff +\qff 
\num{\varepsilon_{\dff 2}}\qff +\qff
\ldots\qff +\qff 
\num{\varepsilon_{\dff n\dff -\dff 1}} 
\off <\off 1$\qss
and\dss hence\dss $\varepsilon_{\dff i}$\dss satisfy\dss all\sss assumptions
of\qss Section\qss \ref{tilings}.\oss
Let\dss us\dss consider\dss the\trs
Lebesgue\sss tiling with\dss the parameters\qss
$\varepsilon_{\dff 1}\fff,\pff \varepsilon_{\dff 2}\fff,\pff 
\ldots\fff,\off \varepsilon_{\dff n\dff -\dff 1}$\qss
and\dss its cubes\dss
$e\qff(\dff a_{\dff 1}\fff,\pff a_{\dff 2}\fff,\pff \ldots\fff,\off a_{\dff n} \dff)$\nnsp.\oss
A\sss trivial\sss check shows\sss that\dss $Z$\dss
is contained\dss in\dss the union of\dss cubes\dss
$e\qff(\dff a_{\dff 1}\fff,\pff a_{\dff 2}\fff,\pff \ldots\fff,\off a_{\dff n} \dff)$\dss
with\qss
$(\dff a_{\dff 1}\fff,\pff a_{\dff 2}\fff,\pff \ldots\fff,\off a_{\dff n} \dff)
\qff \in \qff K$\qss
and\dss such cubes intersect\dss $Z$\dss in $n$\dnsp-intervals.\oss 

By\qss Theorem\qss \ref{rotating-intersections}\qss and\dss Corollary\qss 
\ref{tiling-combinatorics}\qss
these intersections form a\sss tiling of\trs $Z$\nnsp.\oss
This\sss tiling\dss is\dss a part\sss of\trs the\dss Lebesgue\sss tiling of\trs $\rrr^{\fff n}$\dss
with\dss the same parameters.\oss
This allows\sss to identify\dss the nerve\dss $\mathcal{K}$\dss of\trs this\sss tiling\halfff.\oss
Namely\halfff,\oss a subset\sss of\trs $K$\dss is\dss a\sss simplex of\dss
$\mathcal{K}$\dss
if\trs and\dss only\trs if\qss it\dss is\dss the set\sss of\trs terms of\dss
a pivot subsequence.\oss
The complexes\dss
$\mathcal{A}_{\dff i}$\sss 
and\dss
$\mathcal{B}_{\dff i}$\sss 
are also easy\dss to identify\halfff.\oss
Namely\halfff,\pss
$(\dff a_{\dff 1}\fff,\pff a_{\dff 2}\fff,\pff \ldots\fff,\off a_{\dff n} \dff)
\qff \in\qff K$\qss
is\dss a\sss vertex of\dss $\mathcal{A}_{\dff i}$\qss if\trs and\dss only\trs if\qss 
$a_{\dff i}\off =\off 0$\nnsp,\oss
and\dss is\dss a\sss vertex of\dss $\mathcal{B}_{\dff i}$\dss 
if\trs and\dss only\trs if\qss $a_{\dff i}\off =\off k$\nnsp.\oss
The simplices are\sss sets of\trs terms of\dss pivot\sss subsequences.\oss

\mypar{Theorem.}{wh-lebesgue-tilings}
\emph{Suppose\dss that\pss
$e_{\dff 1}\fff,\pff e_{\dff 2}\dff,\pff \ldots\fff,\pff e_{\dff n\dff +\dff 1}$\qss
are subsets of\pss $K$\dss such\dss that\dss their union\dss is\dss equal\dss to\dss $K$\dss
and\dss they\sss are\dss pairwise disjoint\halfff.\oss
Suppose\sss that\dss $e_{\dff i}$\dss is\dss disjoint\dss from\qss
({\fff}the set\sss of\qss vertices\sss of{\trf})\dss $\mathcal{B}_{\dff i}$\dss
for\dss every\qss $i\qff \leq\qff n$\qss
and\dss  
from\dss $\mathcal{A}_{\dff i\dff -\dff 1}$\dss for\dss
every\qss $i\qff \geq\qff 2$\nnsp.\qff\oss
Then\dss the number of\dss pivot sequences having\sss a\sss term\dss belonging\dss to\dss
$e_{\dff i}$\dss for every\qss
$i\off =\off 1\fff,\pff 2\fff,\pff \ldots\fff,\pff n\qff +\qff 1$\qss
is\dss odd.\oss}

\mypar{Theorem.}{hurewicz-lemma-nerve}
\emph{Let\pss
$d_{\dff 1}\fff,\off d_{\dff 2}\dff,\off \ldots\fff,\off d_{\dff n\dff +\dff 1}$\qss
be\dss a\sss covering of\oss $K$\nnsp.\oss
Suppose\sss that\dss $d_{\dff i}$\dss is\dss disjoint\dss from\dss $\mathcal{B}_{\dff i}$\dss
for\dss every\qss $i\qff \leq\qff n$\qss
and\dss is\dss disjoint\dss from\dss $\mathcal{A}_{\dff i\dff -\dff 1}$\dss for\dss
every\qss $i\qff \geq\qff 2$\nnsp.\qff\oss
Then\dss there\dss is\dss a\sss pivot sequence having\sss a\sss term\dss belonging\dss to\dss
$d_{\dff i}$\dss for every\qss
$i\off =\off 1\fff,\pff 2\fff,\pff \ldots\fff,\pff n\qff +\qff 1$\nnsp.\oss}

\prooftitle{Proofs}
It\dss is\dss sufficient\dss to
apply\qss Theorems\qss \ref{e-coverings-hurewicz}\qss and\qss \ref{hurewicz-lemma}.\oss
The resulting\dss proofs\sss  
in\-volve geom\-e\-try\halfff,\oss
but\dss mostly\dss to simplify\dss the exposition.\oss
The main input\sss of\trs geometry\dss is\dss the non-branching\trs
property\dss of\trs $\mathcal{K}$\dnsp.\oss
See\qss Lemma\qss \ref{cube-non-branching}\qss  
for a combinatorial\dss 
proof\halfff.\oss  \eproof

\newpage
\mysection{Lusternik--Schnirelmann\qss theorems}{ls}

\myuppar{The\dss big\dss sphere\dss $S^n$\dnsp.}
Let\dss us fix a natural\dss number\dss $n$\nnsp,\pss
and\dss let\oss
$I\off =\off \{\qff 1\fff,\pff 2\fff,\pff \ldots\fff,\pff n\qff +\qff 1\qff \}$\nnsp.\oss
In\dss this section we will\sss work in\dss $\rrr^{\fff n\dff +\dff 1}$\dnsp.\oss
Let\dss us also fix another natural\dss number\dss $l$\nnsp,\oss
and\dss let\dss\vspace{3pt} 
\[
\quad
Q
\off =\off
[\trf -\qff l\qff -\qff 1/2\fff,\pff l\qff +\qff 1/2\trf]^{\fff n\dff +\dff 1}
\qff \subset\qff 
\rrr^{\fff n\dff +\dff 1}
\qff.
\]

\vspace{-9pt}
Obviously,\pss $Q$\dss is invariant\dss under\dss the\qss \emph{antipodal\dss involution}\qss
$\iota\dff \colon\dff
\rrr^{\dff n\dff +\dff 1}\qff \ttoo\qff \rrr^{\dff n\dff +\dff 1}$\qss defined\dss by
\[
\quad
\iota\dff \colon\dff
x\off \longmapsto\off -\qff x
\qff.
\]
For\qss $X\qff \subset\qff \rrr^{\fff n\dff +\dff 1}$\qss we will abbreviate\dss
$\iota\dff (\dff X\dff)$\dss to\dss $\overline{X}$\nnsp.\oss
The reason for using\dss a half-integer\qss $l\qff +\qff 1/2$\qss instead of\dss an integer\dss
will\dss be clear\dss later\halfff.\oss
Let\qss 
$\pr_i\dff \colon\dff Q
\qff \ttoo\qff 
[\trf -\qff l\qff -\qff 1/2\fff,\pff l\qff +\qff 1/2\trf]$\qss
be\sss the projection\qss
\[
\quad
(\dff x_{\dff 1}\fff,\pff x_{\dff 2}\fff,\pff \ldots\fff,\pff x_{\dff n\dff +\dff 1}\dff)
\off \longmapsto\off
x_{\dff i}
\qff.
\] 
A\qss \emph{face}\qss of\dss $Q$\dss is a product\sss of\dss the form\vspace{-3pt}
\[
\quad
\prod_{\dff i\qff =\qff 1}^{n\qff +\qff 1}\qff
J_{\fff i}
\qff,
\]

\vspace{-12pt}
where each\dss $J_{\fff i}$\dss is either\dss the interval\dss
$[\trf -\qff l\qff -\qff 1/2\fff,\pff l\qff +\qff 1/2\trf]$\nnsp,\oss
or\dss the one-point\sss set\dss $\{\dff -\qff l\qff -\qff 1/2\dff\}$\nnsp,\oss
or\dss $\{\dff l\qff +\qff 1/2\qff\}$\nnsp.\oss
The\qss \emph{dimension}\qss of\dss a face of\dss $Q$\dss is\sss 
the number of\dss intervals in\dss the corresponding\sss product\halfff.\oss
A face of\dss dimension\dss $m$\dss is also called an\qss \emph{$m$\dnsp-face}.\oss
For\qss $i\qff \in\qff I$\oss let\dss\vspace{3pt}
\[
\quad
A_{\dff i}\off =\off \pr_i^{\dff -\dff 1}\dff (\dff -\qff l\qff -\qff 1/2\dff)
\qff\quad
\mbox{and}\qff\quad
B_{\dff i}\off =\off \pr_i^{\dff -\dff 1}\dff (\dff l\qff +\qff 1/2\qff)
\qff.
\]

\vspace{-9pt}
Clearly,\oss $A_{\dff i}$\dss and\dss $B_{\dff i}$\dss are\dss
$n$\dnsp-faces of\dss $Q$\nnsp,\oss and\dss there are no other\dss
$n$\dnsp-faces of\dss $Q$\nnsp.\oss
The usual\dss boundary\dss
$\bd Q$\dss of\dss $Q$\dss is equal\dss to\sss the union of\dss
all\dss faces\qss $A_{\dff i}\fff,\pff B_{\dff i}$\nsp.\oss
The\qss \emph{big\dss sphere}\qss is\vspace{-1.5pt}
\[
\quad
S^n
\off =\off
\bigcup_{i\off =\off 1\vphantom{K^K}}^{n\qff +\qff 1\vphantom{j_j}}
\qff
A_{\dff i}\qff \cup\qff B_{\dff i}
\qff.
\]

\vspace{-12pt}
Obviously,\oss the big sphere\dss 
$S^n$\dss is 
homeomorphic\sss to\sss the standard\dss $n$\dnsp-dimensional sphere
and\dss is invariant\dss under\dss the antipodal\dss involution\dss $\iota$\nnsp.\oss

\myuppar{Cubes.}
A\dss \emph{$n$\dnsp-cube}\qss of\dss $S^n$\dss
is a product of\dss the form\qss\vspace{-1.5pt}
\[
\quad
\prod_{\dff i\qff =\qff 1}^{n\qff +\qff 1}\qff
J_{\dff i} 
\qff,
\]

\vspace{-13.5pt}
where each\dss $J_{\dff i}$\dss is either an\dss interval of\dss the form\qss
$[\trf a\qff -\qff 1/2\fff,\pff a\qff +\qff 1/2\trf]$\nnsp,\oss
where\dss $a$\dss is an\sss integer\dss 
between\dss $-\qff l$\dss and\dss $l$\nnsp,\oss
or\sss a boundary\dss point\sss of\dss the interval\dss
$[\trf -\qff l\qff -\qff 1/2\fff,\pff l\qff +\qff 1/2\trf]$\nnsp,\oss
and exactly\sss one of\dss the sets\dss $J_{\dff i}$\dss is such a boundary\dss point\halfff.\oss
Clearly,\oss every\dss $n$\dnsp-cube is contained\dss
in\dss $S^n$\dss
and\sss $n$\dnsp-cubes form a partition of\dss $S^n$\dss in an obvious sense.\oss
The\qss \emph{cubical\dss sets}\qss are defined as\sss the unions of\dss $n$\dnsp-cubes.\oss
A\qss \emph{cube}\qss
is a product\dss of\dss the same form\dss in which some of\dss the sets\dss $J_{\dff i}$\dss
are allowed\dss to be one-point\sss sets\dss $\{\dff b\qff -\qff 1/2\trf\}$\nnsp,\oss
where\dss $b$\dss is an\sss integer\dss 
between\dss $-\qff l$\dss and\dss $l\qff +\qff 1$\nnsp,\oss
and at\dss least\sss one of\dss these numbers\dss $b$\dss is equal\dss to\dss $-\qff l$\dss
or\qss $l\qff +\qff 1$\nnsp.\oss
The\qss \emph{dimension}\qss of\dss a cube and\dss
\emph{$m$\dnsp-cubes}\qss are defined as before.\oss 
The \emph{$m$\dnsp-chains}\qss are formal sums of\dss $m$\dnsp-cubes with coefficients in\dss
$\ftwo$\nnsp.\oss
The\qss \emph{supports}\qss and\dss the\qss \emph{boundary\sss operator}\dss 
$\partial$\dss are defined as before and\qss
$\partial\dff \circ\dff \partial
\off =\off
0$\qss
as before.\oss

Clearly,\pss the antipodal\dss involution\dss $\iota$\dss maps\dss 
$m$\dnsp-cubes\sss to $m$\dnsp-cubes,\qss $m$\dnsp-chains\sss to\dss $m$\dnsp-chains
etc.\oss
Two objects related\dss by\dss $\iota$\dss are said\dss to be\qss
\emph{antipodal}\pss to each other\halfff,\oss
and\dss an object\dss is called\qss \emph{symmetric}\pss if\dss $\iota$\dss
leaves it\dss invariant\halfff.\oss
The image of\dss a chain $c$ under\dss $\iota$\dss
is\sss denoted\dss by\dss $\iota_{\dff *}\dff(\dff c\dff)$\dnsp.

\myuppar{Intersections with spheres\sss of\qss lower\dss dimensions.}
For\qss $k\qff <\qff n$\qss let\dss $S^k$\dss be\sss the intersection of\dss $S^n$\dss
with\dss the subspace\dss $\rrr^{\fff k\dff +\dff 1}\times\dff 0$\dss of\dss
$\rrr^{\fff n\dff +\dff 1}$\dss
defined\dss by\dss the equations\qss $x_{\dff i}\off =\off 0$\qss
with\qss
$i\qff \geq\qff
k\qff +\qff 2$\nnsp.\oss 
Clearly,\oss the intersection of\dss an $n$\dnsp-cube of\dss $S^n$\dss
with\dss $S^k$\dss is either empty\sss or\sss a $k$\dnsp-cube of\dss $S^k$\nnsp.\oss
More generally,\oss the intersection of\dss an\dss $(\fff n\dff -\dff m\fff)$\dnsp-cube
in\dss $S^n$\dss with\dss the sphere\dss $S^k$\dss 
is either empty\sss or\dss is\sss a\dss $(\fff k\dff -\dff m\fff)$\dnsp-cube
in\dss $S^k$\dss
({\fff}for\qss $k\qff <\qff m$\qss it\dss is always empty\fff).\oss
Using\dss the modern\dss language,\oss
one can say\dss that\sss all cubes in\dss $S^n$\dss are in\qss
\emph{general\dss position}\qss with respect\dss to\sss the sphere\dss $S^k$\dnsp.\oss
This wouldn't\dss be\sss the case if\dss we would\dss use an integer instead of\qss
$l\qff +\qff 1/2$\qss in\dss
the definition of\dss $Q$\nnsp.\oss

Let\dss $\gamma$\dss be\sss an\dss $(\fff n\dff -\dff m\fff)$\dnsp-chain in\dss $S^n$\dnsp,\oss
i.e.\qss let\vspace{3pt}
\[
\quad
\gamma
\off =\dff\off
\sum\nolimits_{\dff i}\dff c_{\dff i}
\qff,
\] 

\vspace{-9pt}
where\dss $c_{\dff i}$\dss are\dss $(\fff n\dff -\dff m\fff)$\dnsp-cubes in\dss $S^n$\dnsp.\oss
The fact\dss that\sss all\sss cubes in\dss $S^n$\dss are in general\dss position\dss with\dss
respect\dss to\dss $S^k$\dss allows\dss to define\sss 
the\qss \emph{intersection}\qss
$\gamma\dff \cap\dff S^k$\qss of\dss $\gamma$\dss with\dss $S^k$\dnsp.\oss
Namely,\vspace{3.5pt}
\[
\quad
\gamma\qff \cap\qff S^k
\off =\qff\off
\sum\nolimits_{\dff i}\dff c_{\dff i}\qff \cap\qff S^k
\qff,
\] 

\vspace{-8.5pt}
where\sss the empty\dss intersections are interpreted as\dss $0$\nnsp'{\halfff}s,\oss
is an $(\fff k\dff -\dff m\fff)$\dnsp-chain\sss
in\dss $S^k$\dnsp.\oss
The operation of\dss taking\dss the intersection of\dss chains with\dss $S^k$\dss
commutes with\dss $\partial$\nnsp,\oss i.e.\vspace{3.5pt}
\[
\quad
\left(\dff \partial\dff \gamma\dff\right)\qff \cap\qff S^k
\off =\off
\partial\dff \left(\dff \gamma \qff \cap\qff S^k\dff\right)
\]

\vspace{-8.5pt}
for any chain\dss $\gamma$\dss in\dss $S^n$\dnsp.\oss
This is\dss trivial\dss for cubes,\oss and extends\sss to\sss the general case
by\dss linearity.\oss

\mypar{Lemma.}{symmetric-chains}
\emph{Every\dss symmetric\dss $m$\dnsp-chain\dss $\gamma$\dss in\dss $S^n$\dss
can\dss be represented\dss in\dss the form}\vspace{3pt}
\begin{equation}
\label{symmetric-split}
\quad
\gamma
\off =\off
\delta\qff +\fff\qff \iota_{\dff *}\dff(\dff \delta\dff)
\qff,
\end{equation}

\vspace{-9pt}
\emph{where\sss the\dss $m$\dnsp-chain\dss $\delta$\dss is such\dss that\dss $\delta$\dss
and\qss $\iota_{\dff *}\dff(\dff \delta\dff)$\dss
have no common\dss $m$\dnsp-cubes.\oss
If\qss $\partial\dff \gamma\off =\off 0$\nnsp,\oss
then\dss $\partial\dff \delta$\dss is a symmetric\dss
$(\fff m\dff -\dff 1\fff)$\dnsp-chain\qss
(which\dss may\dss be equal\dss to\dss zero).\oss}

\proof
Since\dss $\gamma$\dss is symmetric,\oss the set\sss of\dss $m$\dnsp-cubes of\dss $\gamma$\dss
is\dss invariant\dss under\dss the antipodal\dss involution\dss $\iota$\nnsp.\oss
Since\dss $\iota$\dss leaves no $m$\dnsp-cube invariant\halfff,\oss
this set\dss is\sss
equal\dss to\sss the union
of\dss a collection of\dss disjoint\dss pairs of\dss the form\dss
$\{\dff c\fff,\pff \iota\dff(\dff c\dff)\qff\}$\nnsp.\oss
If\dss we select\sss one\dss $m$\dnsp-cube\sss from each such pair\sss
and\dss define\dss $\delta$\dss as\sss the sum of\dss these selected cubes,\oss
then\qss (\ref{symmetric-split})\qss obviously\dss holds.\oss
If\qss (\ref{symmetric-split})\qss holds,\oss then\vspace{3pt}
\[
\quad 
\partial\dff \gamma
\off =\off
\partial\dff \delta\qff +\fff\qff \partial\qff \iota_{\dff *}\dff(\dff \delta\dff)
\off =\off
\partial\dff \delta\qff +\fff\qff \iota_{\dff *}\dff(\dff \partial\dff \delta\dff)
\qff.
\]

\vspace{-9pt}
If\halfff,\oss moreover\halfff,\pss
$\partial\dff \gamma\off =\off 0$\nnsp,\oss
then\qss
$\partial\dff \delta\qff +\fff\qff \iota_{\dff *}\dff(\dff \partial\dff \delta\dff)
\off =\off
0$\nnsp.\oss
Since we are working over\dss $\ftwo$\nsp,\oss
this implies\sss that\qss
$\partial\dff \delta
\off =\fff\off 
\iota_{\dff *}\dff(\dff \partial\dff \delta\dff)$\nnsp,\oss
i.e.\dss $\partial\dff \delta$\dss is symmetric.\oss  \eproof

\myuppar{The\dss odd $0$\dnsp-chains.}
In\dss what\dss follows $0$\dnsp-chains will appear as\sss the intersections\qss
$\omega
\off =\off 
\gamma\qff \cap\qff S^{m}$\qss
of\dss $(\fff n\dff -\dff m\fff)$\dnsp-chains\dss $\gamma$\dss in\dss $S^n$\dss
with\dss $S^m$\nnsp.\oss
Clearly,\oss such a chain\dss $\omega$\dss is\sss symmetric\dss if\dss $\gamma$\dss is.\oss
If\dss $\omega$\dss is\sss a\sss symmetric $0$\dnsp-chain,\oss
then\dss $\omega$\dss is\sss the sum of\dss $0$\dnsp-cubes in a set\sss of\dss $0$\dnsp-cubes invariant\dss
under\dss the action of\dss $\iota$\nnsp.\oss
Clearly,\oss such a set\dss is equal\dss to\sss the union of\dss a collection
of\dss pairs of\dss antipodal\dss points\qss (\nsp$0$\dnsp-cubes).\oss
If\dss the number of\dss these pairs is odd,\oss
then\dss $\omega$\dss is\dss said\dss to be an\qss \emph{odd}\qss $0$\dnsp-chain.\oss

\mypar{Lemma.}{ls-step}
\emph{Let\qss $m\qff \leq\qff n$\qss and\dss
let\dss $\gamma$\dss be a\sss symmetric\sss
$(\fff n\dff -\dff m\fff)$\dnsp-chain\dss in\dss $S^n$\dss
such\dss that\qss $\partial\dff \gamma\off =\off 0$\qss and\qss
$\gamma\qff \cap\qff S^{m}$\qss
is\dss an odd\qss $0$\dnsp-chain.\oss 
Suppose\sss that\qss $\gamma$\dss is represented\dss in\dss the form 
}\vspace{3pt}
\[
\quad
\gamma
\off =\off
\delta\qff +\fff\qff \iota_{\dff *}\dff(\dff \delta\dff)
\qff
\]

\vspace{-9pt}
\emph{and\dss the\dss $(\fff n\dff -\dff m\fff)$\dnsp-chains\dss $\delta$\dss 
and\qss $\iota_{\dff *}\dff(\dff \delta\dff)$\dss
have no common\dss $(\fff n\dff -\dff m\fff)$\dnsp-cubes.\oss
Then}\qss
\vspace{3pt}
\[
\quad
\gamma\fff'\off =\dff\off \partial\dff \delta
\]

\vspace{-9pt}
\emph{is a symmetric\dss $(\fff n\dff -\dff m\dff -\dff 1\fff)$\dnsp-chain\dss and\qss
$\gamma\fff'\qff \cap\qff S^{\fff m\dff +\dff 1}$\qss
is\dss an odd\qss $0$\dnsp-chain.\oss 
In\dss particular\halfff,\qss
$\gamma\fff'
\off \neq\off
0$\nnsp.}

\proof
Lemma\qss \ref{symmetric-chains}\qss implies\sss that\dss $\gamma\fff'$\dss is a symmetric
$(\fff n\dff -\dff m\dff -\dff 1\fff)$\dnsp-chain.\oss
Since an odd\dss $0$\dnsp-chain\dss is obviously\dss non-zero,\oss
it\dss remains\sss to prove\dss that\dss $\gamma\fff'\qff \cap\qff S^{\fff m\dff +\dff 1}$\qss
is an odd\dss $0$\dnsp-chain.\oss

Let\dss us\dss think about\dss the sphere\dss $S^m$\dss
as\dss the\qss \emph{equator}\pss dividing\dss $S^{m\dff +\dff 1}$\dss
into\dss two hemispheres,\oss
the\qss \emph{northern}\qss one,\oss defined\dss by\dss the inequality\qss 
$x_{\dff m\dff +\dff 2}\qff \geq\qff 0$\qss
and\sss denoted\dss by\dss $N^{\fff m\dff +\dff 1}$\nsp,\oss
and\dss the\qss \emph{southern}\qss one.\oss
The chain\dss $\gamma\fff'$\dss is an\dss $(\fff n\dff -\dff m\dff -\dff 1\fff)$\dnsp-chain
and\dss hence\dss
$\num{\gamma\fff'}$\dss is disjoint\dss from\dss $S^m$\dnsp.\oss
Therefore\dss it\dss is sufficient\dss to prove\sss that\qss
$\num{\gamma\fff'}\qff \cap\qff N^{\fff m\dff +\dff 1}$\qss
consists of\dss an odd number of\dss points.\oss

The main idea is\sss to consider\dss the intersection of\dss everything\dss in sight\dss
with\dss $S^{m\dff +\dff 1}$\dnsp.\oss
Let\dss\vspace{3pt} 
\[
\quad
\gamma_{\dff 1}
\off =\off
\gamma\qff \cap\qff S^{m\dff +\dff 1}
\hspace*{1.2em}\mbox{and}\hspace*{1.2em}
\delta_{\dff 1}
\off =\off
\delta\qff \cap\qff S^{m\dff +\dff 1} 
\qff.
\]

\vspace{-9pt}
Then\dss $\gamma_{\dff 1}$\dss and\dss $\delta_{\dff 1}$\dss are 
$1$\dnsp-chains such\dss that\vspace{3pt}
\[
\quad
\gamma_{\dff 1}
\off =\off
\delta_{\dff 1}\qff +\dff\qff \iota_{\dff *}\dff(\dff \delta_{\dff 1}\dff)
\dff,
\]

\vspace{-9pt}
the chains\dss 
$\delta_{\dff 1}$\dss and\qss $\iota_{\dff *}\dff(\dff \delta_{\dff 1}\dff)$\dss
have no common\dss $1$\dnsp-cubes,\oss 
$\gamma\qff \cap\qff S^{m}
\off =\off
\gamma_{\dff 1}\qff \cap\qff S^{m}$\nnsp,\oss 
and\vspace{3pt}
\[
\quad
\gamma\fff'\qff \cap\qff S^{\fff m\dff +\dff 1}
\off =\off
\partial\dff \delta\qff \cap\qff S^{\fff m\dff +\dff 1}
\off =\off
\partial\dff \delta_{\dff 1}\qff \cap\qff S^{\fff m\dff +\dff 1}
\dff.
\]

\vspace{-9pt}
In\dss particular\halfff,\oss
it\dss is sufficient\dss to prove\sss that\dss the number of\dss points in\qss
$\num{\partial\dff \delta_{\dff 1}}\qff \cap\qff N^{\fff m\dff +\dff 1}$\qss
is odd.\oss

Let\dss $\gamma_{\dff E}$\dss  
be\sss the sums of\dss all\dss $1$\dnsp-cubes of\dss
$\gamma_{\dff 1}$\dss  
intersecting\dss 
$S^m$\nnsp.\oss
Then\dss $\gamma_{\dff E}$\dss is a symmetric $1$\dnsp-chain
and\dss hence
is\sss the sum of\dss $1$\dnsp-cubes\sss in\dss an\dss $\iota$\dnsp-invariant\sss set\halfff.\oss 
Clearly,\oss 
such a set\dss is 
equal\dss to\sss the union of\dss a set\sss
of\dss pairs of\dss antipodal\dss $1$\dnsp-cubes.\oss
Let\dss $p$\dss be\sss the number of\dss these pairs.\oss
Then\dss the number of\dss $1$\dnsp-cubes of\dss 
$\gamma_{\dff E}$\dss is\sss 
$2\fff p$\dss
and\dss $p$\dss
is equal\dss to\sss
the number of\dss pairs of\dss antipodal\dss points in\vspace{3pt}
\[
\quad
\gamma_{\dff E}\qff \cap\qff S^{m}
\off =\dff\off
\gamma_{\dff 1}\qff \cap\qff S^{m}
\off =\dff\off
\gamma\qff \cap\qff S^{m}
\qff.
\]

\vspace{-9pt}
Since\qss
$\gamma_{\dff 1}\qff \cap\qff S^{m}$\qss
is an odd\dss $0$\dnsp-chain,\pss
$p$\dss is odd.\oss

Let\dss $\delta_{\dff E}$\dss
be\sss the sums of\dss all\sss $1$\dnsp-cubes\sss of\dss
$\delta_{\dff 1}$\dss
intersecting\dss $S^m$\dnsp.\qff\oss 
Then\vspace{3pt}
\[
\quad
\gamma_{\dff E}
\off =\off
\delta_{\dff E}\qff +\fff\qff \iota_{\dff *}\dff(\dff \delta_{\dff E}\trf)
\qff,
\]

\vspace{-9pt}
and\dss the chains\dss 
$\delta_{\dff E}$\dss and\qss $\iota_{\dff *}\dff(\dff \delta_{\dff E}\dff)$\dss
have no common\dss $1$\dnsp-cubes.\oss
It\dss follows\dss that\dss the number of\dss $1$\dnsp-cubes of\dss $\delta_{\dff E}$\dss
is equal\dss to\sss the half\trs of\trs the number of\dss 
$1$\dnsp-cubes of\dss $\gamma_{\dff E}$\nsp,\oss
i.e.\qss to\dss $p$\nnsp.\oss

Let\dss $\delta_{\dff N}$\dss
be\sss the sums of\dss all\dss $1$\dnsp-cubes of\dss
$\delta_{\dff 1}$\dss
contained\dss in\dss the northern\dss hemisphere\dss $N^{\fff m\dff +\dff 1}$\nsp\dnsp.\dff\oss
Then\qss
$\delta_{\dff E}\qff +\fff\qff \delta_{\dff N}$\qss
is\dss the sum of\dss all\dss $1$\dnsp-cubes of\trs
$\delta_{\dff 1}$\dss
intersecting\dss $N^{\fff m\dff +\dff 1}$\nnsp,\oss
and\dss hence\qss\vspace{3pt}
\begin{equation}
\label{northern}
\quad
\num{\partial\dff \delta_{\dff 1}}\qff \cap\qff N^{\fff m\dff +\dff 1}
\off =\off
\num{\partial\dff (\dff \delta_{\dff E}\qff +\qff \delta_{\dff N}\dff)}
\qff \cap\qff N^{\fff m\dff +\dff 1}
\qff
\end{equation}

\vspace{-9pt}
The\dss $1$\dnsp-chain\qss
$\partial\dff (\trf \delta_{\dff E}\qff +\qff \delta_{\dff N}\trf)$\qss
is a boundary\dss and\dss hence is a sum of\dss an even\dss number of\dss points.\oss
Some of\dss these points cancel\dss under\dss the summation,\oss 
but\dss this does not\sss affect\dss the parity\sss of\dss the number of\dss points.\oss
Since\dss $\num{\delta_{\dff N}}$\dss is contained\dss inside of\dss the northern\dss
hemisphere,\oss
the endpoints of\dss $1$\dnsp-cubes of\dss $\delta_{\dff E}$\dss
contained\dss in\dss the southern\dss hemisphere do not\dss cancel.\oss
There are\dss $p$\dss such endpoints.\oss
Since\dss $p$\dss is odd,\oss
it\dss follows\dss that\dss the\sss number of\dss points of\qss
$\num{\partial\dff (\dff \delta_{\dff E}\qff +\qff \delta_{\dff N}\dff)}$\qss
contained\dss in\dss the northern\dss hemisphere is odd.\oss
In view of\qss (\ref{northern})\qss this implies\sss that\dss the number of\dss points of\qss
$\num{\partial\dff \delta_{\dff 1}}\qff \cap\qff N^{\fff m\dff +\dff 1}$\qss
is odd.\oss
As was noted above,\oss this is sufficient\dss to complete\sss the proof\halfff.\oss  \eproof

\myuppar{Distances.}
Let\dss us use as\dss the distance in\dss $\rrr^{\fff n\dff +\dff 1}$\dss
the\dss \emph{$l_{\dff \infty}$\dnsp-distance}\pss defined as\vspace{3pt}
\[
\quad
\operatorname{dist}\qff(\dff x\fff,\pff y\dff)
\off =\off
\max\nolimits_{\dff i}\qff \num{x_{\dff i}\qff -\qff y_{\dff i}}
\qff,
\]

\vspace{-9pt}
where\qss
$x
\off =\off 
(\dff x_{\dff 1}\fff,\pff x_{\dff 2}\fff,\pff \ldots\fff,\pff x_{\dff n\dff +\dff 1} \dff)$\qss
and\qss
$y
\off =\off 
(\dff y_{\dff 1}\fff,\pff y_{\dff 2}\fff,\pff \ldots\fff,\pff y_{\dff n\dff +\dff 1} \dff)$\nnsp.\oss

Obviously,\oss if\dss the points\qss $x\fff,\pff y$\qss belong\dss to\sss
two disjoint\sss cubes of\dss $S^n$\dnsp,\oss
then\qss
$\operatorname{dist}\qff(\dff x\fff,\pff y\dff)
\qff \geq\qff
1$\nnsp.\oss
It\dss follows\dss that\dss the distance between\dss two
disjoint\sss cubical subsets of\dss $S^n$\dss is\qss $\geq\qff 1$\nnsp.\oss

\mypar{Theorem.}{borsuk-lebesgue-cubes}
\emph{Let\pss
$e_{\dff 1}\fff,\pff e_{\dff 2}\fff,\pff \ldots\fff,\pff e_{\dff n}$\qss
be cubical\sss subsets of\dss the sphere\qss $S^n$\nnsp.\oss
Suppose\dss that\qss none of\qss them contains a pair of\trs antipodal\dss points.\oss
Then\dss the sets\qss
$e_{\dff i}\qff \cup\qff \overline{e}_{\dff i}$\qss
do not\dss cover\dss $S^n$\dnsp.\oss}

\proof
The set\dss 
$e_{\dff i}$\dss
is disjoint\dss from\dss $\overline{e}_{\dff i}$\dss for every\dss $i$\nnsp,\oss
and\dss hence\sss the distance between\dss $e_{\dff i}$\dss
and\dss $\overline{e}_{\dff i}$\dss is\qss $\geq\qff 1$\qss for every\dss $i$\nnsp.\oss
Suppose\sss that\dss the sets\qss
$e_{\dff i}\qff \cup\qff \overline{e}_{\dff i}$\qss
cover\dss $S^n$\dnsp.\oss

Let\dss us consider\dss the images of\dss $S^n$\dss and\dss the sets\dss $e_{\dff i}$\dss
under\dss the map\qss
$\mu\dff \colon\dff
\rrr^{\fff n\dff +\dff 1}\qff \ttoo\qff \rrr^{\fff n\dff +\dff 1}$\qss
multiplying every coordinate by\dss $3$\nnsp,\oss
i.e.\qss under\dss the map\vspace{3pt}
\[
\quad
\mu\dff \colon\dff
(\dff x_{\dff 1}\fff,\pff x_{\dff 2}\fff,\pff \ldots\fff,\pff x_{\dff n\dff +\dff 1} \dff)
\off \longmapsto\off
(\trf 3\fff x_{\dff 1}\fff,\pff 3\fff x_{\dff 2}\fff,\pff \ldots\fff,\pff 3\fff x_{\dff n\dff +\dff 1} \dff)
\qff.
\]

\vspace{-9pt}
Since\qss
$3\dff (\dff l\qff +\qff 1/2\dff)\off =\off (\dff 3\dff l\dff +\dff 1\dff)\qff +\qff 1/2$\nnsp,\oss
the image\dss $\mu\dff(\dff S^n\dff)$\dss is\dss the sphere constructed\dss
with\dss the number\qss $3\dff l\qff +\qff 1$\qss in\dss the role of\dss $l$\nnsp.\oss
Obviously\halfff,\oss the images\dss $\mu\dff(\dff e_{\dff i}\dff)$\dss of\dss the sets\dss $e_{\dff i}$\dss
are disjoint\dss from\dss their antipodal sets,\oss
and\dss it\dss is sufficient\dss to prove\sss that\dss the sets\qss
$\mu\dff(\dff e_{\dff i}\dff)\qff \cup\qff \mu\dff(\dff \overline{e}_{\dff i}\dff)$\qss
do not\sss cover\dss $\mu\dff(\dff S^n\dff)$\nnsp.\oss
Hence\sss we can replace\dss $l$\dss by\qss $3\dff l\dff +\dff 1$\qss
in\dss the definition of\dss $S^n$\dss and\dss the sets\dss $e_{\dff i}$\dss
by\dss their images\dss $\mu\dff(\dff e_{\dff i}\dff)$\nnsp.\oss
We will use for\dss the new sets\dss $\mu\dff(\dff e_{\dff i}\dff)$\dss
the old\dss notation\dss $e_{\dff i}$\nsp.\oss

Now\dss the distance between\dss $e_{\dff i}$\dss
and\qss $\overline{e}_{\dff i}$\dss is\qss $\geq\qff 3$\qss for every\dss $i$\nnsp,\oss
while\dss the cubes have\sss the same size as before.\oss
For every\dss $i$\dss let\dss $E_{\dff i}$\dss be\sss the union of\dss
$e_{\dff i}$\sss and\dss all $n$\dnsp-cubes\sss intersecting\dss $e_{\dff i}$\nnsp.\oss
Then\dss $e_{\dff i}$\dss is contained\dss in\dss
the interior of\dss $E_{\dff i}$\nnsp.\oss
Clearly,\oss the distance of\dss every\dss point\sss of\dss $E_{\dff i}$\dss from\dss $e_{\dff i}$\dss
is\qss $\leq\qff 1$\nnsp,\oss
and\dss the distance of\dss every\dss point\sss of\dss $\overline{E}_{\dff i}$\dss 
from\dss $\overline{e}_{\dff i}$\dss
is also\qss $\leq\qff 1$\nnsp.\oss
Since\sss the distance between\dss $e_{\dff i}$\dss
and\qss $\overline{e}_{\dff i}$\dss is\qss $\geq\qff 3$\nnsp,\oss
it\dss follows\sss that\dss
$E_{\dff i}$\dss
is disjoint\dss from\dss $\overline{E}_{\dff i}$\dss for every\dss $i$\nnsp.\oss

Let\dss us recursively\dss construct\sss 
for each\qss
$m\off =\off 0\fff,\pff 1\fff,\pff \ldots\fff,\pff n$\qss
a symmetric $(\fff n\dff -\dff m\fff)$\dnsp-chain\dss $\gamma_{\fff m}$\dss in\dss $S^n$\dss
such\dss that\qss $\partial\dff \gamma_{\fff m}\off =\off 0$\nnsp,\oss
the $0$\dnsp-chain\qss
$\gamma_{\fff m}\qff \cap\qff S^{m}$\qss
is odd\halfff,\oss
and\dss the support\dss $\num{\gamma_{\fff m}}$\dss
is disjoint\dss from\qss
$e_{\dff i}\qff \cup\qff \overline{e}_{\dff i}$\qss
for\qss $1\qff \leq\qff i\qff \leq\qff m$\qss 
after a possible renumbering of\dss the sets\dss $e_{\dff i}$\nnsp.\oss

Let\qss 
$\gamma_{\dff 0}
\off =\off
\fclass{S^n}$\nnsp,\oss
where\dss $\fclass{S^n}$\dss is\dss
the sum of\dss all\dss $n$\dnsp-cubes in\dss $S^n$\dnsp.\oss
Clearly,\pss $\partial\dff \gamma_{\fff 0}\off =\off 0$\nnsp,\oss
the $0$\dnsp-chain\qss $\gamma_{\fff 0}\qff \cap\qff S^{0}$\qss
is equal\dss to\dss $\fclass{S^0}$\dss and\dss hence is odd\halfff,\oss
and\dss the condition\dss concerned\sss with\dss $\num{\gamma_{\fff 0}}$\dss
is vacuous.\oss
Therefore\dss $\gamma_{\fff 0}$\dss has\sss the required\dss properties.\oss

Suppose\sss that\dss $\gamma_{\fff m}$\dss is already constructed.\oss
Since\sss the sets\qss
$e_{\dff i}\qff \cup\qff \overline{e}_{\dff i}$\qss
are covering\dss $S^n$\dss
and are disjoint\dss from\dss $\num{\gamma_{\fff m}}$\dss
for\qss $i\qff \leq\qff m$\nnsp,\oss
there exists\qss $j\qff \geq\qff m\qff +\qff 1$\qss
such\dss that\qss
$e_{\fff j}\qff \cup\qff \overline{e}_{\fff j}$\qss
has non-empty\dss intersection with\dss $\num{\gamma_{\fff m}}$\nnsp.\oss
After\dss renumbering\dss the sets\sss
$e_{\fff j}$\dss
with\qss $j\qff \geq\qff m\qff +\qff 1$\qss
we may assume\sss that\qss
$e_{\dff m\dff +\dff 1}\qff \cup\qff \overline{e}_{\dff m\dff +\dff 1}$\qss
has non-empty\dss intersection with\dss $\num{\gamma_{\fff m}}$\nnsp.\oss
Since\dss $\gamma_{\fff m}$\dss is symmetric,\oss
this implies\sss that\qss
$e_{\dff m\dff +\dff 1}\qff \cap\qff \num{\gamma_{\fff m}}
\off \neq\off
\varnothing$\qss
and\dss hence\qss
$E_{\dff m\dff +\dff 1}\qff \cap\qff \num{\gamma_{\fff m}}
\off \neq\off
\varnothing$\nnsp.\oss
Let\dss $D_{\dff m}$\dss be\dss the closure of\vspace{4pt}
\[
\quad
\num{\gamma_{\fff m}}
\off \smallsetminus\off
\left(\dff E_{\dff m\dff +\dff 1}\qff \cup\qff \overline{E}_{\dff m\dff +\dff 1}\dff\right)
\qff.
\]

\vspace{-8pt}
Clearly,\pss $D_{\dff m}$\dss is\sss the union of\dss a symmetric set\sss
of\dss $(\fff n\dff -\dff m\fff)$\dnsp-cubes.\oss
Let\dss us\sss select\sss one $(\fff n\dff -\dff m\fff)$\dnsp-cube
from each pair of\dss the form\qss $c\fff,\off \overline{c}$\qss contained\dss
in\dss this set\sss of\dss cubes
and consider\dss the union\dss $F_{\dff m\dff +\dff 1}$\dss of\trs 
$E_{\dff m\dff +\dff 1}\qff \cap\qff \num{\gamma_{\fff m}}$\dss 
with\dss these selected\dss $(\fff n\dff -\dff m\fff)$\dnsp-cubes.\oss
Then\dss 
$F_{\dff m\dff +\dff 1}\qff \cup\pff \overline{F}_{\dff m\dff +\dff 1}
\pff =\off 
\num{\gamma_{\fff m}}$\nnsp,\pss
and\sss since\dss $E_{\dff m\dff +\dff 1}$\dss 
is\dss disjoint\dss from\dss $\overline{E}_{\dff m\dff +\dff 1}$\nsp,\oss 
the set\dss $F_{\dff m\dff +\dff 1}$\dss has no common $(\fff n\dff -\dff m\fff)$\dnsp-cubes\dss 
with\dss $\overline{F}_{\dff m\dff +\dff 1}$\nsp.\oss
Let\dss $\delta_{\dff m}$\dss be\sss the sum of\dss all\dss
$(\fff n\dff -\dff m\fff)$\dnsp-cubes contained\dss in\dss $F_{\dff m\dff +\dff 1}$\nnsp.\oss
Then\vspace{3pt} 
\[
\quad
\gamma_{\fff m}
\off =\off
\delta_{\dff m}\qff +\qff \iota_{\dff *}\dff(\dff \delta_{\dff m}\dff)
\qff
\]

\vspace{-9pt}
and\dss the\dss $(\fff n\dff -\dff m\fff)$\dnsp-chains\dss $\delta_{\dff m}$\dss 
and\qss $\iota_{\dff *}\dff(\dff \delta_{\dff m}\dff)$\dss
have no common\dss $(\fff n\dff -\dff m\fff)$\dnsp-cubes.\oss
Let\qss\vspace{3pt} 
\[
\quad
\gamma_{\fff m\dff +\dff 1}
\off =\dff\off 
\partial\dff \delta_{\dff m}
\qff.
\]

\vspace{-9pt}
By\dss Lemma\dss \ref{ls-step}\qss the {\nsp}$(\fff n\dff -\dff m\dff -\dff 1\fff)$\dnsp-chain
$\gamma_{\fff m\dff +\dff 1}$ is\sss symmetric\sss
and\sss
$\gamma_{\fff m\dff +\dff 1}\qff \cap\qff S^{\fff m\dff +\dff 1}$
is\dss an odd $0$\dnsp-chain.\oss 
The identity\qss
$\partial\dff \circ\dff \partial
\off =\off
0$\qss
implies\sss that\qss
$\partial\dff \gamma_{\fff m\dff +\dff 1}
\off =\off
0$\nnsp.\oss
Clearly,\pss
$\num{\gamma_{\fff m\dff +\dff 1}}
\qff \subset\qff
\num{\gamma_{\fff m}}$\qss
and\dss hence\dss
$\num{\gamma_{\fff m\dff +\dff 1}}$\dss
is disjoint\dss from\qss
$e_{\dff i}\qff \cup\qff \overline{e}_{\dff i}$\qss
for\qss $i\qff \leq\qff m$\nnsp.\oss
By\dss the construction,\pss
$\num{\gamma_{\fff m\dff +\dff 1}}$\dss
may\dss intersect\dss the set\qss
$E_{\dff m\dff +\dff 1}\qff \cup\qff \overline{E}_{\dff m\dff +\dff 1}$\qss
only\sss along\dss its boundary.\oss 
Since\qss
$e_{\dff m\dff +\dff 1}\qff \cup\qff \overline{e}_{\dff m\dff +\dff 1}$\qss
is contained\dss in\dss its interior\halfff,\oss
$\num{\gamma_{\fff m\dff +\dff 1}}$\dss
is disjoint\dss from\qss
$e_{\dff m\dff +\dff 1}\qff \cup\qff \overline{e}_{\dff m\dff +\dff 1}$\qss
also.\oss
Therefore\dss $\gamma_{\fff m\dff +\dff 1}$\dss has\sss the required\dss properties.\oss

This completes\sss the construction of\dss chains\dss $\gamma_{\fff m}$\nsp.\oss
The last\sss of\dss them\dss is\dss the $0$\dnsp-chain\dss $\gamma_{\fff n}$\nsp.\oss
Its support\dss $\num{\gamma_{\fff n}}$\dss
is disjoint\dss from all sets\qss
$e_{\dff i}\qff \cup\qff \overline{e}_{\dff i}$\nnsp.\oss
At\dss the same\sss time\dss $\num{\gamma_{\fff n}}$\dss
is non-empty\dss because\qss
$\gamma_{\fff n}
\off =\off
\gamma_{\fff n}\qff \cap\qff S^{n}$\qss
is an odd\dss $0$\dnsp-chain.\oss
It\dss follows\dss that\dss the sets\qss
$e_{\dff i}\qff \cup\qff \overline{e}_{\dff i}$\qss
do not\sss cover\dss $S^n$\dnsp.\oss  \eproof\vspace{1.25pt}

\myuppar{The\sss first\qss Lusternik--Schnirelmann\dss theorem.} 
\emph{Let\qss $\mathbb{S}^n$\dss be\sss the standard\dss unit\dss sphere in\dss 
$\rrr^{\fff n\dff +\dff 1}$\dnsp.\oss
Let\pss
$F_{\dff 1}\fff,\pff F_{\dff 2}\fff,\pff \ldots\fff,\pff F_{\dff n}$\qss
be closed\sss subsets of\qss $\mathbb{S}^n$\nnsp.\oss
Suppose\dss that\qss none of\qss them contains a pair of\trs antipodal\dss points.\oss
Then\dss the sets\qss
$F_{\dff i}\qff \cup\qff \overline{F}_{\dff i}$\qss
do not\dss cover\dss $\mathbb{S}^n$\dnsp.\oss}\vspace{1.25pt}

\proof
For every\dss $i$\dss the sets\dss $F_{\dff i}$\dss and\dss $\overline{F}_{\dff i}$\dss
are disjoint\halfff.\oss
Since\sss they are closed subsets\dss $\mathbb{S}^n$\dss 
and\dss $\mathbb{S}^n$\dss is compact\halfff,\oss 
the distance between\dss them\dss 
is\qss $>\qff 0$\nnsp.\oss
Let\qss $\varepsilon\qff >\qff 0$\qss be a real\dss number smaller\dss than\dss
the distance between\dss $F_{\dff i}$\dss and\dss $\overline{F}_{\dff i}$\dss
for every\dss $i$\nnsp.\oss

Let\qss
$r\dff\colon\dff
S^n\qff \ttoo \mathbb{S}^n$\qss
be\sss the radial\dss projection.\oss
The images\dss $r\dff(\dff c\trf)$\dss of\dss $n$\dnsp-cubes of\dss $S^n$\dss
form a partition of\dss 
the sphere\dss $\mathbb{S}^n$\dnsp.\oss
Let\dss $E_{\dff i}$\dss
be\dss the union of\dss all such images intersecting\dss $F_{\dff i}$\nnsp.\oss
Clearly,\oss $F_{\dff i}\qff \subset\qff E_{\dff i}$\nnsp.\oss
If\dss $l$\dss is sufficiently\dss big,\oss
then\dss the diameter of\dss $r\dff(\dff c\dff)$\dss
is\qss $<\qff \varepsilon/3$\qss for every $c$\nnsp.\oss
In\dss this case\sss the distance of\dss every point\sss of\dss $E_{\dff i}$\dss 
from\dss $F_{\dff i}$\dss is\qss $<\qff \varepsilon/3$\nnsp,\oss
and\dss hence\dss $E_{\dff i}$\dss
and\dss $\overline{E}_{\dff i}$\dss are disjoint\halfff.\oss

For every\dss $i$\dss let\qss 
$e_{\dff i}
\off =\off
r^{\fff -\dff 1}\dff(\dff E_{\dff i}\dff)$\nnsp.\oss
Then\dss the sets $e_{\dff i}$ are cubical sets,\oss
and\dss for every\sss $i$\sss the sets\sss $e_{\dff i}$\sss
and\dss $\overline{e}_{\dff i}$\dss are disjoint\halfff.\oss 
By\trs Theorem\qss \ref{borsuk-lebesgue-cubes}\qss
the sets\qss
$e_{\dff i}\qff \cup\qff \overline{e}_{\dff i}$\qss
do not\dss cover\dss $S^n$\dnsp.\oss
It\dss follows\dss that\dss
the sets\qss
$E_{\dff i}\qff \cup\qff \overline{E}_{\dff i}$\qss
do not\dss cover\dss $\mathbb{S}^n$\dnsp,\oss
and\dss hence\sss their subsets\qss
$F_{\dff i}\qff \cup\qff \overline{F}_{\dff i}$\qss
also do not\dss cover\dss $\mathbb{S}^n$\dnsp.\oss  \eproof\vspace{1.25pt}

\myuppar{The second\qss Lusternik--Schnirelmann\dss theorem.} 
\emph{Let\pss
$F_{\dff 1}\fff,\pff F_{\dff 2}\fff,\pff \ldots\fff,\pff F_{\dff n\dff +\dff 1}$\qss
be closed\sss subsets of\qss $\mathbb{S}^n$\nnsp.\oss
If\qss none of\qss them contains a pair of\trs antipodal\dss points,\oss
then\dss these sets
do not\dss cover\dss $\mathbb{S}^n$\dnsp.\oss}\vspace{1.25pt}

\proof
By\trs the previous\dss theorem\qss 
the sets\qss
$F_{\dff i}\qff \cup\qff \overline{F}_{\dff i}$\qss
with\qss $i\qff \leq\qff n$\dss do not\sss cover\dss $\mathbb{S}^n$\dnsp,\oss
i.e.\qss the complement\dss $C$\dss of\dss their\dss union\dss is non-empty.\oss
If\dss the sets\qss
$F_{\dff 1}\fff,\pff F_{\dff 2}\fff,\pff \ldots\fff,\pff F_{\dff n\dff +\dff 1}$\qss
cover\dss $\mathbb{S}^n$\dnsp,\oss
then\dss $C$\dss is contained\dss in\dss $F_{\dff n\dff +\dff 1}$\nsp.\oss
But\dss $C$\dss is obviously\dss invariant\dss under\dss the antipodal\dss
involution,\oss and\dss hence in\dss this case\dss $F_{\dff n\dff +\dff 1}$\dss
contains a pair of\dss antipodal\dss points,\oss contrary\dss to\sss
the assumption.\oss  \eproof

\mysection{First\qss applications\qss of\pss Lusternik--Schirelmann\qss theorems}{ls-applications}

\mypar{Theorem.}{antipodal-maps}
\emph{If\qss $m\qff <\qff n$\nnsp,\oss
then\dss there exists no continuous map\qss
$f\dff \colon\dff
\mathbb{S}^n\qff \ttoo\qff \mathbb{S}^{m}$\qss
taking pairs of\dss antipodal\dss points\sss to pairs of\dss antipodal\dss points,\oss
i.e.\qss such\dss that\qss
$f\dff(\dff -\qff x\trf)\off =\off -\qff f\dff(\dff x\trf)$\qss
for all\qss $x$\nnsp.}

\proof
Let\qss 
$f\dff \colon\dff
\mathbb{S}^n\qff \ttoo\qff \mathbb{S}^{m}$\qss
be a continuous map.\oss
Let\qss $\varepsilon\qff >\qff 0$\qss be a small\dss number\halfff.\oss
For every\qss
$i\off =\off 1\fff,\pff 2\fff,\pff \ldots\fff,\pff m\qff +\qff 1$\qss
let\dss
$H_{\dff i}$\dss
be\dss the subset\sss of\qss $\mathbb{S}^{m}$\dss defined\dss by\dss the inequality\qss
$x_{\dff i}\qff \geq\qff \varepsilon$\nnsp,\oss
where\qss
$x_{\dff 1}\dff,\pff x_{\dff 2}\dff,\pff \ldots\dff,\pff x_{\dff m\dff +\dff 1}$\qss
are\sss the standard coordinates in\dss $\rrr^{\fff m\dff +\dff 1}$\nnsp.\oss
Clearly,\pss $H_{\dff i}$\dss is disjoint\dss from\dss 
$\overline{H}_{\dff i}$\nnsp.\oss
If\dss $\varepsilon$\dss is small enough,\oss
then\dss the sets\qss
$H_{\dff i}\qff \cup\qff \overline{H}_{\dff i}$\qss
cover\dss $\mathbb{S}^{m}$\dss
and\dss hence\sss
their\dss preimages\dss\vspace{0pt}
\[
\quad
f^{\dff -\dff 1}\dff\left(\trf H_{\dff i}\qff \cup\qff \overline{H}_{\dff i} \trf\right)
\off =\off
f^{\dff -\dff 1}\dff(\trf H_{\dff i} \trf)
\qff \cup\qff 
f^{\dff -\dff 1}\dff\bigl(\trf \overline{H}_{\dff i} \trf\bigr)
\]

\vspace{-12pt}
cover\dss  $\mathbb{S}^n$\dnsp.\oss
Also,\oss the preimages of\dss $H_{\dff i}$\dss
and\dss $\overline{H}_{\dff i}$\dss are disjoint\dss
because\dss $H_{\dff i}$\dss
and\dss $\overline{H}_{\dff i}$\dss are.\oss

Let\qss
$F_{\dff i}
\off =\off
f^{\dff -\dff 1}\dff(\trf H_{\dff i} \trf)$\nnsp.\oss
If\qss
$f$\dss takes pairs of\dss antipodal\dss points\sss to pairs of\dss antipodal\dss points,\oss
then\vspace{0pt}
\[
\quad
f^{\dff -\dff 1}\dff\bigl(\trf \overline{H}_{\dff i} \trf\bigr)
\off =\qff\off
\overline{F}_{\dff i}
\]

\vspace{-12pt}
for every\dss $i$\nnsp.\oss
It\dss follows\dss that\dss $F_{\dff i}$\dss
is a closed set\sss disjoint\dss from\dss $\overline{F}_{\dff i}$\dss 
for every\dss $i$\dss
and\dss the sets\qss
$F_{\dff i}\qff \cup\qff \overline{F}_{\dff i}$\qss
cover\dss $\mathbb{S}^n$\dnsp.\oss
But\dss there are\qss
$m\qff +\qff 1\qff \leq\qff n$\qss sets\dss $F_{\dff i}$\dss
and\dss hence we reached\sss a contradiction\dss with\trs
the first\qss Lusternik--Schnirelmann\dss theorem\qss
(see\dss Section\qss \ref{ls}).\oss 
Therefore no such map\dss $f$\dss exists.\oss \eproof

\myuppar{Borsuk--Ulam\dss theorem.} 
\emph{Let\qss
$f\dff \colon\dff
\mathbb{S}^n\qff \ttoo\qff \rrr^{\fff n}$\qss
be a continuous map.\oss
Then\trs there\dss is\dss a\sss point\qss $x\qff \in\qff \mathbb{S}^n$\qss
such\dss that\qss
$f\dff(\dff -\qff x\trf)\off =\off f\dff(\dff x\trf)$\qss
for some\dss $x$\nnsp.}

\proof
Suppose\sss that\qss
$f\dff(\dff -\qff x\trf)\off \neq\off f\dff(\dff x\trf)$\qss
for all\dss $x\qff \in\qff \mathbb{S}^n$\dnsp.\oss
Then\vspace{2.25pt}
\[
\quad
g\dff \colon\dff
x\off \longmapsto\off
\frac{f\dff(\dff x\trf)\qff -\qff f\dff(\dff -\qff x\trf)}
{\norm{f\dff(\dff x\trf)\qff -\qff f\dff(\dff -\qff x\trf)}}
\off,
\]

\vspace{-9.75pt}
where\dss $\norm{\bullet}$\dss is\sss the usual euclidean\dss norm,\oss
is a well defined continuous map\qss
$\mathbb{S}^n\qff \ttoo\qff \mathbb{S}^{n\dff -\dff 1}$\dnsp.\oss
Obviously,\oss
$g\dff(\dff -\qff x\trf)\off =\off -\qff g\dff(\dff x\trf)$\qss
for all\qss $x\qff \in\qff \mathbb{S}^n$\dnsp,\oss
in contradiction\dss with\qss 
Theorem\qss \ref{antipodal-maps}.\oss  \eproof

\myuppar{Borsuk\qss Theorem.}
\emph{If\qss a\dss continuous map\dss
$f\dff \colon\dff
\mathbb{S}^n\qff \ttoo\qff \mathbb{S}^{n}$\dss
takes pairs of\dss antipodal\dss points\sss to pairs of\dss antipodal\dss points,\pss
then $f${\nsp} cannot\dss be extended\dss to\dss a map\qss
$\mathbb{B}^{\fff n\dff +\dff 1}\qff \ttoo\qff \mathbb{S}^{n}$\dnsp.\oss}

\proof
Suppose\sss that\dss 
$g\dff \colon\dff
\mathbb{B}^{\fff n\dff +\dff 1}\qff \ttoo\qff \mathbb{S}^n$\dss
is\dss such an extension.\pss
Let\dss
$p\dff \colon\dff
\mathbb{S}^{n\dff +\dff 1}\qff \ttoo\qff \mathbb{B}^{\fff n\dff +\dff 1}$\dss
be\sss the projection\dss forgetting\dss the last\sss coordinate.\pss
Let\dss us\sss define a map\qss
$h\dff \colon\dff
\mathbb{S}^{n\dff +\dff 1}\qff \ttoo\qff \mathbb{S}^n$\qss
as\dss follows\fff:\qss 
$h\dff(\dff x\trf)\off =\off g\dff \circ\dff p\trf(\dff x\trf)$\qss
if\dss $x$\dss is\dss in\dss 
the northern\dss hemisphere of\trs $\mathbb{S}^{n\dff +\dff 1}$\dss
and\qss
$h\dff(\dff x\trf)\off =\off {}-\qff g\dff \circ\dff p\trf(\dff {}-\qff x\trf)$\qss
if\trs in\sss the southern.\oss
Both\dss formulas agree on\dss $\mathbb{S}^n$\dss and\dss hence\dss $h$\dss is well\sss
defined.\oss
Clearly,\pss $h$\dss takes antipodal\sss points\sss to antipodal\dss points,\oss
in contradiction\dss with\qss Theorem\qss \ref{antipodal-maps}.\oss  \eproof

\mypar{Theorem.}{ah-lefschetz}
\emph{Let\pss
$F_{\dff 1}\fff,\pff F_{\dff 2}\fff,\pff \ldots\fff,\pff F_{\dff n\dff +\dff 1}$\qss
be closed\sss subsets of\qss $\mathbb{S}^n$\nnsp.\oss
Suppose\dss that\dss none of\qss them contains a pair of\trs antipodal\dss points\dss
and\dss the sets\qss
$F_{\dff i}\qff \cup\qff \overline{F}_{\dff i}$\qss
cover\dss $\mathbb{S}^n$\dnsp.\oss
Then\dss the intersection of\dss the sets\qss
$F_{\dff 1}\fff,\pff F_{\dff 2}\fff,\pff \ldots\fff,\pff F_{\dff n\dff +\dff 1}$\qss
is non-empty.\oss}

\proof
Since\sss  
the sets\dss $F_{\dff i}$\dss and\dss $\overline{F}_{\dff i}$\dss
are closed\sss and\sss disjoint\halfff,\oss
there are continuous functions\vspace{4.5pt}
\[
\quad
f_{\dff i}\qff \colon\dff
\mathbb{S}^n\qff \ttoo\off [\trf -\qff 1/2\dff,\pff 1/2\trf]
\]

\vspace{-7.5pt}
such\dss that\pss
$f_{\dff i}^{\dff -\dff 1}\dff(\qff 1/2\qff)
\off =\off
F_{\dff i}$\pss
and\pss
$f_{\dff i}^{\dff -\dff 1}\dff(\qff -\qff 1/2\qff)
\off =\off
\overline{F}_{\dff i}$\nsp.\qff\oss
Let\qss
$g_{\dff i}\dff(\dff x\trf)\off =\off
f_{\dff i}\dff(\dff x\trf)
\qff -\qff
f_{\dff i}\dff(\dff -\qff x\trf)$\nnsp.\oss
Then\dss $g_{\dff i}$\dss is a continuous function\dss 
$\mathbb{S}^n \ttoo [\trf -\qff 1\fff,\pff 1\trf]$\dss
such\dss that\qss
$g_{\dff i}^{\dff -\dff 1}\dff(\qff 1\qff)
\off =\off
F_{\dff i}$\qss
and\qss
$g_{\dff i}^{\dff -\dff 1}\dff(\qff -\qff 1\qff)
\off =\off
\overline{F}_{\dff i}$\nnsp.\oss
In addition,\oss
$g_{\dff i}\dff(\dff -\qff x\trf)\off =\off -\qff g_{\dff i}\dff(\dff x\trf)$\qss
for all\qss $x\qff \in\qff \mathbb{S}^n$\dnsp.\oss
Together\dss the functions\dss $g_{\dff i}$\dss
define a map\vspace{4.5pt}
\[
\quad
g\dff \colon\dff
\mathbb{S}^n
\off \ttoo\off
[\trf -\qff 1\fff,\pff 1\trf]^{\dff n\dff +\dff 1}
\qff
\]

\vspace{-7.5pt}
such\dss that\qss
$g\dff(\dff -\qff x\trf)\off =\off -\qff g\dff(\dff x\trf)$\qss
for all\qss $x\qff \in\qff \mathbb{S}^n$\dnsp.\qff\oss
Since\sss the sets\qss
$F_{\dff i}\qff \cup\qff \overline{F}_{\dff i}$\qss
cover\dss $\mathbb{S}^n$\dnsp,\oss
the image of\dss $g$\dss is contained\dss in\dss 
the boundary\dss $\mathbf{S}^{\fff n}$\dss
of\dss the cube\dss
$[\trf -\qff 1\fff,\pff 1\trf]^{\dff n\dff +\dff 1}$\nnsp.\oss

Suppose\sss that\qss
$F_{\dff 1}\qff \cap\qff F_{\dff 2}\qff \cap\qff \ldots\qff \cap\qff F_{\dff n\dff +\dff 1}
\off =\off
\varnothing$\nnsp.\qss\oss
Then\dss also\qss
$\overline{F}_{\dff 1}\qff \cap\qff \overline{F}_{\dff 2}\qff \cap\qff 
\ldots\qff \cap\qff \overline{F}_{\dff n\dff +\dff 1}
\off =\off
\varnothing$\nnsp.\oss
It\dss follows\dss that
the points\qss 
$\mathbb{1}
\off =\off
(\dff 1\fff,\pff 1\fff,\pff \ldots\fff,\pff 1\dff)$\qss 
and\qss
${}-\qff \mathbb{1}
\off =\off
(\dff -\qff 1\fff,\pff -\qff 1\fff,\pff \ldots\fff,\pff -\qff 1\dff)$\qss 
are not\dss in\dss the image of\dss $g$\nnsp.\oss

As we will see in a moment\halfff,\oss there is a continuous map\qss\vspace{4.5pt}
\[
\quad
q\dff \colon\dff
\mathbf{S}^{\fff n}
\qff \smallsetminus\qff
\{\qff \mathbb{1}\dff,\qff -\qff \mathbb{1}\pff\}
\off \ttoo\off
\mathbb{S}^{n\dff -\dff 1}
\]

\vspace{-7.5pt}
such\dss that\qss
$q\trf(\dff -\qff x\trf)\off =\off -\qff q\dff(\dff x\trf)$\qss
for all\qss $x\qff \in\qff \mathbf{S}^{\fff n}
\qff \smallsetminus\qff
\{\qff \mathbb{1}\dff,\qff -\qff \mathbb{1}\pff\}$\dnsp.\qff\oss
If\dss $q$\dss is such a map,\oss then\vspace{4.5pt}
\[
\quad
h\dff \colon\dff
x
\off \longmapsto\off
q\trf\bigl(\trf g\dff(\dff x\trf)\trf\bigr)
\]

\vspace{-7.5pt}
is a continuous map\qss
$\mathbb{S}^{\fff n}
\off \ttoo\off
\mathbb{S}^{n\dff -\dff 1}$\qss
such\dss that\qss
$h\dff(\dff -\qff x\trf)\off =\off -\qff h\dff(\dff x\trf)$\qss
for all\qss $x\qff \in\qff \mathbb{S}^n$\dnsp.\qff\oss
The contradiction\dss with\dss Corollary\qss \ref{antipodal-maps}\qss
completes\sss the proof\dss modulo\dss the existence of\trs the map\dss
$q$\nnsp.\oss

Let\dss $\mathbb{1}^{\trf \perp}$\dss be\sss the $n$\dnsp-dimensional\dss
vector subspace of\dss $\rrr^{\fff n\dff +\dff 1}$\dss
consisting of\dss vectors orthogonal\dss to\dss $\mathbb{1}$\nnsp,\oss
and\dss let\qss
$p\dff \colon\dff
\rrr^{\fff n\dff +\dff 1}
\qff \ttoo\qff
\mathbb{1}^{\trf \perp}$\qss
be\sss the orthogonal\dss projection\dss to\dss $\mathbb{1}^{\trf \perp}$\dnsp.\oss
Then\qss
$p\trf(\dff -\qff x\trf)\off =\off -\qff p\dff(\dff x\trf)$\qss
for all\qss $x$\dss and\qss
$p\dff(\dff x\trf)\off =\off 0$\qss
if\trs and\sss only\trs if\dss
$x$\dss is proportional\dss to\dss $\mathbb{1}$\nnsp.\oss
In\dss particular\halfff,\pss
$p\dff(\dff x\trf)\off \neq\off 0$\qss
if\qss 
$x\qff \in\qff \mathbf{S}^{\fff n}
\qff \smallsetminus\qff
\{\qff \mathbb{1}\dff,\qff -\qff \mathbb{1}\pff\}$\dnsp.\qff\oss
It\trs follows\dss that\dss the formula\vspace{4.5pt}
\[
\quad
q_{\dff 1}\dff(\dff x\trf)
\off =\qff\off
\frac{p\dff(\dff x\trf)}{\norm{p\dff(\dff x\trf)}}
\off,
\]

\vspace{-7.5pt}
where\dss $\norm{\bullet}$\dss is\sss the usual euclidean\dss norm,\oss
defines a map\sss from\qss
$\mathbf{S}^{\fff n}
\qff \smallsetminus\qff
\{\qff \mathbb{1}\dff,\qff -\qff \mathbb{1}\pff\}$\qss
to\sss the unit\sss sphere of\dss $\mathbb{1}^{\trf \perp}$\dss
such\dss that\qss
$q_{\dff 1}\trf(\dff -\qff x\trf)\off =\off -\qff q_{\dff 1}\dff(\dff x\trf)$\qss
for\dss all\dss $x$\nnsp.\oss
In order\dss to get\dss the map\dss $q$\nnsp,\oss
it\dss remains\sss to identify\dss this unit\sss sphere with\dss
$\mathbb{S}^{n\dff -\dff 1}$\dss
by\sss a rotation.\oss  \eproof

\mypar{Theorem.}{polya}
\emph{Let\qss
$f\dff \colon\dff
\mathbb{S}^n\qff \ttoo\qff \rrr^{\fff n}$\qss
be a continuous map.\oss
If\pss
$f\dff(\dff -\qff x\trf)\off =\off -\qff f\dff(\dff x\trf)$\qss
for every\qss $x\qff \in\qff \mathbb{S}^n$\nnsp,\oss
then\qss $f\dff(\dff x\trf)\off =\off 0$\qss for\dss some\qss
$x\qff \in\qff \mathbb{S}^n$\nnsp.}

\proof
If\qss $f\dff(\dff x\trf)\off \neq\off 0$\qss for every\qss
$x\qff \in\qff \mathbb{S}^n$\dnsp,\oss
then\qss
$x\qff \longmapsto\qff f\dff(\dff x\trf)\bigl/\fff\norm{f\dff(\dff x\trf)}$\qss
is\dss a\sss continuous map\qss
$\mathbb{S}^n\qff \ttoo\qff \mathbb{S}^{n\dff -\dff 1}$\qss
contradicting\qss Theorem\qss \ref{antipodal-maps}.\oss  \eproof

\myuppar{Historical\dss remarks.}
The applications of\qss Lusternik-Schnirelmann\qss theorems discussed\dss
in\dss this section are\dss the simplest\sss ones and are\dss the first\sss 
only\dss in\dss a\dss logical\qss (or\sss an expository\fff)\qss sense.\oss
Lusternik\dss and\dss Schnirelmann\dss proved\dss their\dss theorems\qss \cite{ls}\qss
with applications\sss to variational\dss problems in\dss mind.\oss
Their\dss main\dss result\dss is\dss a\sss solution of\dss a\sss problem
of\qss Poincar\'{e}\qss about\dss geodesics.\oss
Namely\halfff,\oss
Lusternik\dss and\dss Schnirelmann\dss proved\qss \cite{ls}\qss
that\sss on every\sss closed\sss convex surface in\dss $\rrr^{\fff 3}$\dss
there are at\dss least\dss three
closed\dss geodesics without\dss self-intersections.\oss
They did\dss not\sss stated\qss the\sss theorem\sss which\dss we\sss called\dss the\dss
first\dss Lusternik--Schnirelmann\dss theorem,\oss 
but\dss its\sss proof\trs is\sss the main part\sss of\trs their\dss proof\dss
of\trs the\sss theorem\sss which\dss we\sss called\dss the\dss
second\dss Lusternik--Schnirelmann\dss theorem.\oss
The\dss latter\dss is\qss Lemma\qss 1\qss in\dss Section\qss II.5\qss of\pss \cite{ls}.\oss

The paper\qss \cite{bo}\qss of\qss
K.\dss Borsuk,\oss published\dss in\qss 1933,\oss
three years after\dss Lusternik--Schnirelmann\dss book\qss \cite{ls},\oss
is\dss devoted\dss to\sss three\sss theorems about\dss the topology of\dss spheres.\oss
His main\dss theorem,\oss from\dss which\dss he deduced\dss the\sss two other\halfff,\oss
is\dss the\sss theorem called\qss \emph{Borsuk\dss theorem}\qss above.\oss
Strange\-ly\sss enough,\oss his second\dss theorem\dss
became known as\qss \emph{Borsuk--Ulam\dss theorem}.\oss
Borsuk\dss mentioned\dss in\sss a footnote in\qss \cite{bo}\qss
that\dss this result\dss was suggested\dss by\qss S.\dss Ulam\dss
as a conjecture.\oss
Ulam\dss wasn't\dss involved\dss in\dss the proof\halfff.\oss
G.-C.\dss Rota\qss [R]\qss tells\sss this story\dss in\dss his unique poetic style mixed\dss with\sss
bitterness.\oss
The\dss third\dss Borsuk's\dss theorem\dss is\dss the second\qss
Lusternik--Schinirel\-mann\qss theorem\dss from\qss Section\qss \ref{ls}.\oss
Apparently,\oss Borsuk\dss was unaware of\trs the work of\pss Lusternik\dss
and\dss Schnirelmann,\oss
and\dss his methods are completely\sss different\halfff.\oss
In a footnote\dss Borsuk\dss wrote \vspace{-10.5pt}

\begin{quoting}
Herr\dss H.\dss Hopf\halfff,\oss to whom\dss I\dss communicated\dss the first\dss theorem,\oss
briefly\dss indicated\dss to me\sss three other short\dss proofs of\trs these\sss theorems.\oss
But\dss these proofs are based on deep results of\dss the theory of\dss degrees of\dss maps
and\dss my\dss proof\trs basically\dss is\dss quite elementary\halfff,\oss
and\dss hence\dss I\dss believe\sss that\dss it\sss is\dss not\sss superfluous.\oss
\end{quoting}

\vspace{-10.5pt}
Borsuk's\dss main\dss tools\dss are\dss homotopies and,\oss 
this\dss footnote notwithstanding\halfff,\pss
the degrees of\dss maps.\oss
In\qss \cite{bo}\qss
he build\dss by\dss bare hands a\sss fragment\sss of\trs the homotopy\dss theory\sss
of\dss maps\sss to spheres.\oss
One may\dss speculate\sss that\sss at\dss least\dss one\sss of\pss Hopf's\qss proofs\dss
was a\sss high-brow\dss version of\qss Borsuk's\dss proof\halfff.\oss 
For a\sss modern\sss elementary\dss proof\qss in\dss this spirit\dss
see\qss \cite{ma},\oss Section\qss 2.2.\oss

Strangely\dss enough,\oss
Theorem\qss \ref{ah-lefschetz}\qss
is\dss not\dss so\dss well\dss known as\sss the other\dss versions of\pss
Borsuk--Ulam\qss theorem.\oss
J.\qss Matou\v{s}ek\qss included\dss in\qss \cite{ma},\oss Section\qss 2.1\qss 
about\sss a dozen\dss versions,\oss
but\qss not\qss  Theorem\qss \ref{ah-lefschetz}.\oss 
The author\dss learned\qss Theorem\qss \ref{ah-lefschetz}
from\qss Lefschetz's\dss book\qss \cite{lef}.\oss
See\qss \cite{lef},\oss Section\qss IV.{\dff}7,\oss Theorem\qss (21.1).\oss
Apparently\halfff,\oss this\sss theorem\dss is\dss due\sss to\qss A.W.\qss Tucker\halfff.\oss
He proved\dss the case\qss $n\off =\off 2$\qss in\qss \cite{t1}.\oss
In\dss full\dss generality\dss his results were published only\dss in\qss \cite{lef}.\oss 
A similar result\dss was proved\dss by\qss Alexandroff\qss and\dss Hopf\qss \cite{ah}.\oss
See\qss \cite{ah},\oss Section\qss XII.3,\oss  Theorem\qss X.\oss
Alexandroff\qss and\dss Hopf\qss (ibid.,\oss Theorem\qss VIII\fff)\qss 
write\sss that\qss Theorem\qss \ref{polya}\qss
was suggested\dss to\sss them\dss by\qss G.\dss P\'{o}lya.\oss

\newpage
\mysection{The\qss discrete\qss cube\qss and\qss products\qss of\qss cubical\qss cochains}{cubical-cochains}

\myuppar{The discrete cube.}
As usual\halfff,\oss let\dss $n$\dss be a natural\dss number and\qss
$I\off =\off \{\qff 1\fff,\pff 2\fff,\pff \ldots\fff,\pff n\qff \}$\nnsp.\oss
Let\dss $k$\dss be another natural\dss number\halfff.\oss
The\qss \emph{discrete\dss cube}\pss of\dss size\sss $k$\trs is\dss\vspace{3pt} 
\[
\quad
K
\off =\off
\{\qff 0\dff,\pff 1\dff,\pff \ldots\dff,\pff k \qff\}^{\dff n}
\qff.
\] 

\vspace{-9pt}
For\qss
$i\off =\off 1\fff,\pff 2\fff,\pff \ldots\fff,\pff n$\qss
let\dss $\mathcal{A}_{\dff i}$\dss and\dss $\mathcal{B}_{\dff i}$\dss 
be\sss the sets\sss of\dss all\dss points\qss\vspace{3pt}
\[
\quad
(\dff
a_{\dff 1}\fff,\pff
a_{\dff 2}\trf,\pff
\ldots\fff,\pff
a_{\dff n}
\dff)
\qff \in\qff K
\]

\vspace{-9pt}
such\dss that\qss 
$a_{\dff i}\off =\off 0$\qss
and\qss
$a_{\dff i}\off =\off k$\qss 
respectively.\oss
An\dss \emph{$(\fff n\dff -\dff 1\fff)$\dnsp-face}\qss of\trs $K$\dss is\dss  
a subset\sss of\trs $K$\dss
equal\dss to\dss $\mathcal{A}_{\dff i}$\dss or\dss $\mathcal{B}_{\dff i}$\dss for some\dss $i$\nnsp.\oss
An\dss interested\dss reader\dss may define $m$\dnsp-faces of\trs $K$\dss with\qss 
$m\off \neq\off n\qff -\qff 1$\nnsp.\oss
The faces\qss $\mathcal{A}_{\dff i}\fff,\off \mathcal{B}_{\dff i}$\nnsp,\oss
where\qss $i\qff \in\qff I$\nnsp,\oss
are said\dss to be\qss \emph{opposite}.\oss
The\qss \emph{boundary}\dss $\bd\fff K$\dss is\dss  
the union of\dss
the $(\fff n\dff -\dff 1\fff)$\dnsp-faces.\oss

\myuppar{The cubes of\trs $K$\nnsp.}
An\dss \emph{$n$\dnsp-cube}\qss of\trs $K$\dss is  
a product of\dss the form\qss\vspace{3pt}
\[
\quad
\prod_{\dff i\qff =\qff 1}^{\dff n}\qff
\{\trf a_{\dff i}\fff,\pff a_{\dff i}\qff +\qff 1\qff \}
\qff,
\]

\vspace{-9pt}
where $a_{\dff 1}\fff,\pff  
\ldots\fff,\off a_{\dff n}$
are integers between  $0$ and\dss $k\qff -\qff 1$\nnsp.\oss
A\dss \emph{cube}\qss of\trs $K$\dss 
is a product\sss of\dss the form\vspace{3pt}
\begin{equation}
\label{discrete-cubes}
\quad
\sigma
\off =\off
\prod_{\dff i\qff =\qff 1}^{\dff n}\qff
\rho_{\fff i}
\qff,
\end{equation}

\vspace{-9pt}
where\dss for every\qss $i\qff \in\qff I$\qss
either\qss
$\rho_{\fff i}
\off =\off
\{\trf a_{\dff i}\fff,\pff a_{\dff i}\qff +\qff 1\qff\}$\qss  
for a non-negative integer\qss $a_{\dff i}\qff \leq\qff k\qff -\qff 1$\nnsp,\oss
or\dss 
$\rho_{\fff i}
\off =\off
\{\dff a_{\dff i}\trf\}$\qss
for a non-negative integer\qss $a_{\dff i}\qff \leq\qff k$\nnsp.\qff\oss
The\qss \emph{dimension}\qss of\dss $\sigma$\dss is\dss  
the num\-ber of\dss 
pairs\qss $\{\trf a_{\dff i}\fff,\pff a_{\dff i}\qff +\qff 1\trf\}$\qss in\dss
the\sss product\qss (\ref{discrete-cubes}).\oss

An \emph{$m$\dnsp-cube}\qss is a cube of\dss dimension\dss $m$\nnsp.\oss 
A\qss \emph{$m$\dnsp-face}\qss of\dss a cube $\sigma$ is 
an\dss $m$\dnsp-cube contained\dss in\dss $\sigma$\nnsp.\oss
Every\dss $m$\dnsp-face of\dss the cube\qss (\ref{discrete-cubes})\qss
is obtained\dss by\dss keeping\dss
$m$\dss pairs\qss
$\{\trf a_{\dff i}\fff,\pff a_{\dff i}\qff +\qff 1\trf\}$\qss
in\dss the product\qss (\ref{discrete-cubes})\qss intact\sss
and\dss replacing\sss other such\dss pairs by\sss one of\dss their elements.\oss

\myuppar{Chains\sss and\dss cochains of\trs $K$\nnsp.}
Now one can define $m$\dnsp-chains and $m$\dnsp-cochains of\qss $K$\dss in\dss
the same manner as before.\oss
The\dss \emph{$m$\dnsp-chains of}\qss $K$ are defined as formal\sss
sums of\dss $m$\dnsp-cubes of\trs $K$\dss with coefficients in\dss $\ftwo$\nnsp.\oss
Clearly\halfff,\oss they\dss form a vector space over\dss $\ftwo$\dss
(of\dss course,\oss the $m$\dnsp-chains of\dss $Q$\sss also form such a vector
space).\oss 
Let\dss us\sss denote\sss this vector space by\dss $C_{\dff m}\dff(\dff K\fff)$\nnsp.\oss
The\dss \emph{cubical\dss $m$\dnsp-cochains}\pss of\trs $K$\nnsp,\oss
or\sss simply\dss the\sss \emph{$m$\dnsp-cochains},\oss
are defined as\sss the elements of\dss the dual\sss vector space\qss
$C^{\dff m}\dff(\dff K\fff)
\off =\off
C_{\dff m}\dff(\dff K\fff)^*$\nnsp.\oss
Since\sss the vector space\dss $C_{\dff m}\dff(\dff K\fff)$\dss 
has a canonical\dss basis,\oss
namely,\oss the basis consisting of\dss $m$\dnsp-cubes of\trs $K$\nnsp,\oss
one can\dss identify\sss cochains with\dss the formal\sss sums of\dss $m$\dnsp-cubes,\oss
similarly\dss to\sss the usual\sss cochains.\oss
But\dss viewing $m$\dnsp-cochains
as $\ftwo$\dnsp-valued\dss
functions on\dss the set\sss of\dss $m$\dnsp-cubes of\dss $K$\dss
turns out\dss to be more convenient\halfff.\oss

The boundary\dss $\partial\dff \sigma$\dss of\dss an $m$\dnsp-cube\dss $\sigma$\dss
is defined as\sss the formal sum of\dss $(\fff m\dff -\dff 1\fff)$\dnsp-faces
of\dss $\sigma$\nnsp.\oss
The map\qss
$\sigma\qff \longmapsto\qff \partial\dff \sigma$\qss
extends by\dss linearity\dss to\sss the\qss \emph{boundary\dss operator}\vspace{1.5pt}
\[
\quad
\partial\dff \colon\dff
C_{\dff m}\dff(\dff K\dff)
\off \ttoo\off
C_{\dff m\dff -\dff 1}\dff(\dff K\dff)
\qff.
\]

\vspace{-10.5pt}
By\dss the same reasons as before,\oss 
$\partial\dff \circ\dff \partial
\off =\off
0$\nnsp.\oss
The\qss \emph{coboundary\dss operator}\vspace{1.5pt}
\[
\quad
\partial^*\dff \colon\dff
C^{\dff m\dff -\dff 1}\dff(\dff K\dff)
\off \ttoo\off
C^{\dff m}\dff(\dff K\dff)
\]

\vspace{-10.5pt}
is\sss the linear operator dual\dss to\dss $\partial$\nnsp.\oss
Clearly,\pss 
$\partial^* \circ\dff \partial^*
\pff =\off
0$\nnsp.\qff\oss
If\dss $\tau$\dss is\dss an\sss $(\fff m\dff -\dff 1\fff)$\dnsp-cube,\oss
then\vspace{1.5pt}
\[
\quad
\partial^*\fff(\dff \tau\dff)
\off =\qff\off
\sum\qff \sigma
\qff,
\]

\vspace{-10.5pt}
where\sss the sum\dss is\dss taken over all $m$\dnsp-cubes\dss $\sigma$\dss
having\dss $\tau$\dss as a face.\oss

\myuppar{Products of\dss cochains.}
Let\dss $\sigma$\dss be\sss the cube given\sss by\qss the product\qss
(\ref{discrete-cubes})\qss with\dss the sets\dss $\rho_{\fff i}$\dss as above.\oss
The\qss \emph{direction}\qss of\dss $\sigma$ is\dss defined as\sss the set\dss 
$A\dff(\dff \sigma\dff)$\dss
of\trs subscripts\qss $i\qff \in\qff I$\qss such\dss that\sss $\rho_{\fff i}$\sss
has\dss the form\dss $\{\trf a_{\dff i}\fff,\pff a_{\dff i}\qff +\qff 1\trf\}$\nnsp.\oss
If\dss $\sigma$\dss is\dss an $m$\dnsp-cube,\oss
then\dss $A\dff(\dff \sigma\dff)$\dss is\dss an $m$\dnsp-element\sss set\halfff.\oss

Suppose\sss that\qss $H\qff \subset\qff A\dff(\dff \sigma\dff)$\qss
and\qss $\varepsilon\off =\off 0$\qss or\dss $1$\nnsp.\oss
Let\dss us replace in\dss the product\qss 
(\ref{discrete-cubes})\qss
every\dss factor\qss
$\rho_{\fff i}
\off =\off
\{\trf a_{\dff i}\fff,\pff a_{\dff i}\qff +\qff 1\qff\}$\qss
with\qss $i\qff \in\qff H$\pss
by\qss
$\{\dff a_{\dff i}\qff +\qff \varepsilon\qff\}$\nnsp.\oss
The resulting\dss product\dss is denoted\dss by\vspace{-1.5pt}
\[
\quad 
\lambda_{\dff H}^{\dff \varepsilon}\trf \sigma
\qff.
\]

\vspace{-13.5pt}
If\dss  $\sigma$ is an $m$\dnsp-cube and 
$H$ is\dss an \nsp$h$\dnsp-element\sss set\halfff,\pss 
then\sss $\lambda_{\dff H}^{\dff \varepsilon}\trf \sigma$\sss is an $(\fff m\dff -\dff h\fff)$\dnsp-cube.\oss 
In\dss these\sss terms\vspace{3pt}
\[
\quad
\partial\dff \sigma
\off =\off
\sum_{\dff P\fff,\pff \varepsilon\vphantom{H^4}}\pff
\lambda_{\dff P}^{\dff \varepsilon}\trf \sigma
\qff,
\]

\vspace{-15pt}
where\dss $P$\dss runs over $1$\dnsp-element\sss subsets of\trs
$A\dff(\dff \tau\dff)$\sss
and\dss $\varepsilon$\dss runs over\qss $\{\trf 0\fff,\pff 1\trf\}$\nnsp.\oss

The following definition is an adaptation of\trs the definition of\dss
products of\dss singular cubical\sss co\-chains due\sss to\pss J.-P.\qss Serre\qss \cite{se}.\oss
Let\dss us consider $m$\dnsp-cochains as\dss $\ftwo$\dnsp-valued\dss
functions on\dss the set\sss of\dss $m$\dnsp-cubes.\oss
Suppose\sss that\dss $f$\dss and\dss $g$\dss
are\dss  $p$\dnsp-cochain\dss and\dss $q$\dnsp-cochain\dss respectively.\oss
The\qss \emph{product}\dss $f\fff \cdot\dff g$\qss
is\dss the\dss $(\fff p\dff +\dff q\fff)$\dnsp-cochain 
defined\dss by\trs the\sss formula\vspace{5pt}
\begin{equation}
\label{product-cochains}
\quad
f\fff \cdot\dff g\qff(\dff \sigma\dff)
\off =\off
\sum_{\dff F\fff,\pff G\vphantom{H^4}}\off
f\dff\left(\qff \lambda_{\trf G}^{\dff 0}\qff \sigma \qff\right)
\trf \cdot\qff 
g\dff\left(\qff \lambda_{\trf F}^{\fff 1}\qff \sigma \qff\right)
\qff,
\end{equation}

\vspace{-10pt}
where\dss $\sigma$\dss is\dss an $(\fff p\dff +\dff q\fff)$\dnsp-cube,\oss
and\qss $F\fff,\pff G$\qss
runs over all\dss pairs of\dss disjoint\dss subsets of\dss
$A\dff(\dff \sigma\dff)$\dss
consisting\dss of\dss $p$\dss and\dss $q$\dss elements respectively.\oss
Equivalently,\pss $F$\dss runs over all\sss $p$\dnsp-element\sss 
subsets of\dss $A\dff(\dff \sigma\dff)$\nnsp,\oss
and\dss
$G$\dss is\sss the complement\sss of\dss $F$\dss in\dss $A\dff(\dff \sigma\dff)$\nnsp,\oss
i.e.\qss
$G
\off =\off
A\dff(\dff \sigma\dff)\qff \smallsetminus\pff F$\nnsp.\oss

\mypar{Lemma.}{asso}
\emph{The product\sss of\dss cochains is associative.\oss}

\proof
Suppose\sss that\qss
$f\nsp,\off g$\qss and\dss $h$\dss
are $p$\dnsp-cochain,\pss $q$\dnsp-cochain,\oss and $r$\dnsp-cochain\dss respectively,\oss
and\dss that\dss $\sigma$\dss is\dss an $(\fff p\dff +\dff q\dff +\dff r\fff)$\dnsp-cube.\qff\oss
Then\vspace{4.5pt}
\[
\quad
(\dff f\fff \cdot\dff g\dff)\fff \cdot\dff h\qff(\dff \sigma\dff)
\off\qff =\off
\sum_{\dff F\fff,\pff G\fff,\pff H\vphantom{H^4}}\off
f\dff\left(\qff \lambda_{\trf G}^{\dff 0}\qff \lambda_{\dff H}^{\dff 0}\sigma \qff\right)
\trf \cdot\qff 
g\dff\left(\qff \lambda_{\trf F}^{\fff 1}\qff \lambda_{\dff H}^{\dff 0}\qff \sigma \qff\right)
\trf \cdot\qff 
h\dff\left(\qff \lambda_{\trf F\dff \cup\dff G}^{\fff 1}\trf \sigma \qff\right)
\qff,
\]

\vspace{-12pt}
where\qss 
$F\fff,\off G\fff,\off H$\qss
runs\sss over\dss triples of\dss disjoint\dss subsets of\dss
$A\dff(\dff \sigma\dff)$\dss
consisting\dss of\qss $p\fff,\pff q$\qss and\dss $r$\dss elements respectively.\oss
Similarly,\vspace{4.5pt}
\[
\quad
f\fff \cdot\dff (\dff g\fff \cdot\dff h\dff)\qff(\dff \sigma\dff)
\off\qff =\off
\sum_{\dff F\fff,\pff G\fff,\pff H\vphantom{H^4}}\off
f\dff\left(\dff \lambda_{\trf G\dff \cup\dff H}^{\dff 0}\qff \sigma \dff\right)
\dff \cdot\trf 
g\dff\left(\dff \lambda_{\dff H}^{\dff 0}\qff \lambda_{\trf F}^{\fff 1}\qff \sigma \dff\right)
\dff \cdot\trf 
h\dff\left(\dff \lambda_{\trf G}^{\fff 1}\trf \lambda_{\trf F}^{\fff 1}\trf \sigma \dff\right)
\qff,
\]

\vspace{-12pt}
where\qss
$F\fff,\off G\fff,\off H$\qss
runs\sss over\dss the same\sss triples.\oss
The\sss two sums are equal\dss because,\oss
obviously,\vspace{4.5pt}
\[
\quad
\lambda_{\trf G}^{\dff 0}\qff \lambda_{\dff H}^{\dff 0}\sigma
\off =\off
\lambda_{\trf G\dff \cup\dff H}^{\dff 0}\qff \sigma\qff,
\hspace{1.5em}
\lambda_{\trf F\dff \cup\dff G}^{\fff 1}\trf \sigma
\off =\off
\lambda_{\trf G}^{\fff 1}\trf \lambda_{\trf F}^{\fff 1}\trf \sigma
\qff,
\]

\vspace{-12pt}
and\vspace{0pt}
\[
\quad
\lambda_{\trf F}^{\fff 1}\qff \lambda_{\dff H}^{\dff 0}\qff \sigma
\off =\off
\lambda_{\dff H}^{\dff 0}\qff \lambda_{\trf F}^{\fff 1}\qff \sigma
\qff.
\]

\vspace{-9pt}
The lemma follows.\oss  \eproof

\mypar{Lemma\qss ({\fff}Leibniz\dss formula\halfff).}{leibniz}
$\dis
\partial^*\dff\left(\dff f\fff \cdot\dff g \dff\right)
\off =\off
\left(\dff \partial^* f \dff\right)\fff \cdot\dff g
\off +\off
f\fff \cdot\dff \left(\dff \partial^* g \dff\right)$\oss
\emph{for every\dss cochains}\qss $f\nsp,\off g$\nnsp.

\proof
Suppose\sss that\dss $f$\dss and\dss $g$\dss
are\dss  $p$\dnsp-cochain\dss and\dss $q$\dnsp-cochain\dss respectively,\oss
and\dss that\dss $\sigma$\dss is\dss an $(\fff p\dff +\dff q\dff +\dff 1\fff)$\dnsp-cube.\oss
Then\oss\vspace{4.5pt}
\[
\quad
\partial^*\dff\left(\dff f\fff \cdot\dff g \dff\right)\dff (\dff \sigma\dff)
\off\off =\off
\sum_{\dff F\fff,\pff G\fff,\pff P\fff,\pff \varepsilon\vphantom{H^4}}\off
f\dff\left(\qff \lambda_{\trf G}^{\dff 0}\qff \lambda_{\dff P}^{\dff \varepsilon}\trf \sigma \qff\right)
\trf \cdot\qff 
g\dff\left(\qff \lambda_{\trf F}^{\fff 1}\qff \lambda_{\dff P}^{\dff \varepsilon}\trf \sigma \qff\right)
\qff,
\]

\vspace{-7.5pt}
where\qss 
$F\fff,\off G\fff,\off P$\qss
runs\sss over\dss triples of\dss disjoint\dss subsets of\dss
$A\dff(\dff \sigma\dff)$\dss
consisting\dss of\qss $p\fff,\pff q$\qss and\dss $1$\dss elements respectively,\oss
and\dss $\varepsilon$\dss runs over\qss $\{\trf 0\fff,\pff 1\trf\}$\nnsp.\oss
Similarly,\vspace{7.5pt}
\[
\quad
\left(\dff \partial^* f \dff\right)\fff \cdot\dff g\dff (\dff \sigma\dff)
\off\off =\off
\sum_{\dff F\fff,\pff G\fff,\pff P\fff,\pff \varepsilon\vphantom{H^4}}\off
f\dff\left(\qff \lambda_{\dff P}^{\dff \varepsilon}\trf \lambda_{\trf G}^{\dff 0}\qff \sigma \qff\right)
\trf \cdot\qff 
g\dff\left(\qff \lambda_{\trf F\dff \cup\dff P}^{\fff 1}\qff \trf \sigma \qff\right)
\qff,
\]

\vspace{-17pt}
and\vspace{0pt}
\[
\quad
f\fff \cdot\dff \left(\dff \partial^* g \dff\right)\dff (\dff \sigma\dff)
\off\off =\off
\sum_{\dff F\fff,\pff G\fff,\pff P\fff,\pff \varepsilon\vphantom{H^4}}\off
f\dff\left(\qff \lambda_{\trf G\dff \cup\dff P}^{\fff 0}\qff \sigma \qff\right)
\trf \cdot\qff 
g\dff\left(\qff \lambda_{\dff P}^{\dff \varepsilon}\trf \lambda_{\trf F}^{\dff 1}\qff \sigma \qff\right)
\qff,
\]

\vspace{-9pt}
where\qss 
$F\fff,\off G\fff,\off P$\qss and\dss $\varepsilon$\dss 
run over\dss the same sets.\oss\vspace{0pt}

It\dss follows\dss that\qss
$\dis
\left(\dff \partial^* f \dff\right)\fff \cdot\dff g\dff (\dff \sigma\dff)$\qss
is equal\dss to\vspace{10.5pt}
\begin{equation}
\label{df-g}
\quad
\sum_{\dff F\fff,\pff G\fff,\pff P\vphantom{H^4}}\off
f\dff\left(\qff \lambda_{\dff P}^{\dff 0}\qff \lambda_{\trf G}^{\dff 0}\qff \sigma \qff\right)
\trf \cdot\qff 
g\dff\left(\qff \lambda_{\trf F\dff \cup\dff P}^{\fff 1}\qff \sigma \qff\right)
\off +\off
\sum_{\dff F\fff,\pff G\fff,\pff P\vphantom{H^4}}\off
f\dff\left(\qff \lambda_{\dff P}^{\dff 1}\qff \lambda_{\trf G}^{\dff 0}\qff \sigma \qff\right)
\trf \cdot\qff 
g\dff\left(\qff \lambda_{\trf F\dff \cup\dff P}^{\fff 1}\qff \sigma \qff\right)
\end{equation}

\vspace{-4.5pt}
and\qss
$\dis
f\fff \cdot\dff \left(\dff \partial^* g \dff\right)\dff (\dff \sigma\dff)$\qss
is equal\dss to\vspace{10.5pt}
\begin{equation}
\label{f-dg}
\quad
\sum_{\dff F\fff,\pff G\fff,\pff P\vphantom{H^4}}\off
f\dff\left(\qff \lambda_{\trf G\dff \cup\dff P}^{\dff 0}\qff \sigma \qff\right)
\trf \cdot\qff 
g\dff\left(\qff \lambda_{\dff P}^{\dff 1}\trf \lambda_{\trf F}^{\dff 1}\qff \sigma \qff\right)
\off +\off
\sum_{\dff F\fff,\pff G\fff,\pff P\vphantom{H^4}}\off
f\dff\left(\qff \lambda_{\trf G\dff \cup\dff P}^{\dff 0}\qff \sigma \qff\right)
\trf \cdot\qff 
g\dff\left(\qff \lambda_{\dff P}^{\dff 0}\trf \lambda_{\trf F}^{\dff 1}\qff \sigma \qff\right)
\end{equation}

\vspace{-4.5pt}
Since\sss the sets\qss $F\fff,\pff G$\dss and\dss $P$\dss are disjoint\halfff,\oss\vspace{6pt}
\[
\quad
\lambda_{\dff P}^{\dff 0}\qff \lambda_{\trf G}^{\dff 0}\qff \sigma
\off =\off
\lambda_{\trf G\dff \cup\dff P}^{\dff 0}\qff \sigma
\]

\vspace{-15.25pt}
and
\vspace{-3.25pt}
\[
\quad
\lambda_{\trf F\dff \cup\dff P}^{\fff 1}\qff \sigma
\off =\off
\lambda_{\dff P}^{\dff 1}\trf \lambda_{\trf F}^{\dff 1}\qff \sigma
\qff,
\]

\vspace{-6pt}
and\dss hence\dss the first\sss sums in\qss (\ref{df-g})\qss
and\qss (\ref{f-dg})\qss 
are equal\sss and\dss cancel\dss under addition.\oss
Therefore\vspace{7.5pt}
\[
\quad
\left(\dff \partial^* f \dff\right)\fff \cdot\dff g\dff (\dff \sigma\dff)
\off +\off
f\fff \cdot\dff \left(\dff \partial^* g \dff\right)\dff (\dff \sigma\dff)
\]

\vspace{-24pt}
\[
\quad
=\off
\sum_{\dff F\fff,\pff G\fff,\pff P\vphantom{H^4}}\off
f\dff\left(\qff \lambda_{\dff P}^{\dff 1}\qff \lambda_{\trf G}^{\dff 0}\qff \sigma \qff\right)
\trf \cdot\qff 
g\dff\left(\qff \lambda_{\trf F\dff \cup\dff P}^{\fff 1}\qff \sigma \qff\right)
\off +\off
\sum_{\dff F\fff,\pff G\fff,\pff P\vphantom{H^4}}\off
f\dff\left(\qff \lambda_{\trf G\dff \cup\dff P}^{\dff 0}\qff \sigma \qff\right)
\trf \cdot\qff 
g\dff\left(\qff \lambda_{\dff P}^{\dff 0}\trf \lambda_{\trf F}^{\dff 1}\qff \sigma \qff\right)
\qff
\]

\vspace{-21pt}
\[
\quad
=\off
\sum_{\dff F\fff,\pff G\fff,\pff P\vphantom{H^4}}\off
f\dff\left(\qff \lambda_{\dff P}^{\dff 1}\qff \lambda_{\trf G}^{\dff 0}\qff \sigma \qff\right)
\trf \cdot\qff 
g\dff\left(\qff \lambda_{\dff P}^{\dff 1}\trf \lambda_{\trf F}^{\dff 1}\qff \sigma \qff\right)
\off +\off
\sum_{\dff F\fff,\pff G\fff,\pff P\vphantom{H^4}}\off
f\dff\left(\qff \lambda_{\dff P}^{\dff 0}\qff \lambda_{\trf G}^{\dff 0}\qff \sigma \qff\right)
\trf \cdot\qff 
g\dff\left(\qff \lambda_{\dff P}^{\dff 0}\trf \lambda_{\trf F}^{\dff 1}\qff \sigma \qff\right)
\]

\vspace{-21pt}
\[
\quad
=\off
\sum_{\dff F\fff,\pff G\fff,\pff P\fff,\pff \varepsilon\vphantom{H^4}}\off
f\dff\left(\qff \lambda_{\trf G}^{\dff 0}\qff \lambda_{\dff P}^{\dff \varepsilon}\qff \sigma \qff\right)
\trf \cdot\qff 
g\dff\left(\qff \lambda_{\trf F}^{\dff 1}\qff \lambda_{\dff P}^{\dff \varepsilon}\trf \sigma \qff\right)
\off.
\]

\vspace{-7.5pt}
This is nothing else but\dss the above expression for\qss
$\dis
\partial^*\dff\left(\dff f\fff \cdot\dff g \dff\right)\dff (\dff \sigma\dff)$\nnsp.\oss
The lemma follows.\oss  \eproof

\myuppar{Another\dss definition of\dss products.}
The associativity\dss implies\sss that\dss products\qss
$f\fff \cdot\dff g\fff \cdot\dff \ldots\dff \cdot\dff h$\qss
of\dss several\sss cochains are well-defined.\oss
Now we will\dss rephrase\sss the definition\qss (\ref{product-cochains})\qss
in a form\dss more suitable\sss for\dss working\dss with\dss multiple products.\oss 
Given a cube\dss $\sigma$\dnsp,\oss 
let\dss 
$A\off =\off A\dff(\dff \sigma\dff)$\nnsp,\oss 
and\dss let\vspace{4.5pt}
\[
\quad
\lambda^0\trf \sigma
\off =\off
\lambda_{\dff A}^{\dff 0}\trf \sigma
\hspace{1.5em}\mbox{and}\hspace{1.5em}
\lambda^1\trf \sigma
\off =\off
\lambda_{\dff A}^{\dff 1}\trf \sigma
\qff.
\]

\vspace{-7.5pt}
Both\sss $\lambda^0\trf \sigma$\sss and\sss $\lambda^1\trf \sigma$\sss
are $0$\dnsp-cubes,\pss which we will\sss call\dss respectively\sss the\dss
\emph{root}\qss and\sss the\dss \emph{peak}\qss of\dss $\sigma$\dnsp.\oss
In\dss these terms\sss the definition of\dss the product\dss takes\sss the following\dss form.\oss
Suppose\sss that\dss $f$\dss and\dss $g$\dss
are\dss $p$\dnsp-cochain\dss and\dss $q$\dnsp-cochain\dss respectively,\oss
and\dss let\dss $\sigma$\dss be a $(\fff p\dff +\dff q\fff)$\dnsp-cube.\oss
Then\vspace{6pt}
\begin{equation}
\label{product-cochains-other}
\quad
f\fff \cdot\dff g\qff(\dff \sigma\dff)
\off =\off
\sum_{\dff \tau\fff,\pff \rho\vphantom{\tau^a}}\off
f\dff\left(\qff \tau \qff\right)
\trf \cdot\qff 
g\dff\left(\qff \rho \qff\right)
\qff,
\end{equation}

\vspace{-9pt}
where\sss the sum\sss is\sss taken over all\dss pairs\qss $\tau\fff,\pff \rho$\qss
such\dss that $\tau$ is\dss a $p$\dnsp-face of\dss $\sigma$ 
and $\rho$ is\dss a $q$\dnsp-face of\dss $\sigma$\dnsp,\oss
the directions of\dss $\tau$ and $\rho$
are\sss disjoint\halfff,\oss 
the root\sss of\dss $\tau$\dss is equal\dss to\sss
the root\sss of\dss $\sigma$\nnsp,\oss
and\dss the peak of\dss $\tau$\dss is equal\dss to\sss the root\sss of\dss $\rho$\nnsp,\oss
i.e.\qss\vspace{4.5pt}
\[
\quad
A\dff(\dff \tau\dff)
\qff \cap\qff
A\dff(\dff \rho\dff)
\off =\off
\varnothing
\qff,
\]

\vspace{-33pt}
\[
\quad
\lambda^0\trf \tau
\off =\off
\lambda^0\trf \sigma\qff,
\hspace{1.5em}\mbox{and}\hspace{1.5em}
\lambda^1\trf \tau
\off =\off
\lambda^0\trf \rho
\]

\vspace{-7.5pt}
Indeed,\oss suppose\sss that\sss $\sigma$\sss is given\dss by\qss
(\ref{discrete-cubes}).\oss 
If\pss $F\fff,\pff G$\qss are\sss disjoint\sss subsets of\dss $A\dff(\dff \sigma\dff)$\dss
as\dss in\qss (\ref{product-cochains}),\oss
then\qss
$\tau
\off =\off 
\lambda_{\trf G}^{\dff 0}\qff \sigma$\qss
and\qss
$\rho
\off =\off 
\lambda_{\trf F}^{\qff 1}\qff \sigma$\qss
are faces of\dss $\sigma$ with\sss directions\dss $F$\dss and\dss
$G$\dss respectively,\oss and\vspace{4.5pt}
\[
\quad
\lambda^1\trf \tau
\off =\off
\lambda^0\trf \rho
\off =\off
(\qff 
a_{\dff 1}\qff +\qff \varepsilon_{\dff 1}\dff,\pff 
a_{\dff 2}\qff +\qff \varepsilon_{\dff 2}\dff,\pff 
\ldots\dff,\pff 
a_{\dff n}\qff +\qff \varepsilon_{\dff n} 
\qff)
\qff,
\]

\vspace{-7.5pt}
where\qss $\varepsilon_{\dff i}\off =\off 0$\qss
if\pss $i\qff \not\in\qff A\dff(\dff \sigma\dff)$\qss
or\qss $i\qff \in\qff G$\nnsp,\oss
and\qss $\varepsilon_{\dff i}\off =\off 1$\qss
if\pss $i\qff \in\qff F$\nnsp.\oss
Conversely,\oss if\pss $\tau\fff,\pff \rho$\qss
have\sss the above properties,\oss
then\qss
$\tau
\off =\off 
\lambda_{\trf G}^{\dff 0}\qff \sigma$\qss
and\qss
$\rho
\off =\off 
\lambda_{\trf F}^{\qff 1}\qff \sigma$\nnsp,\oss
where\qss
$F\off =\off A\dff(\dff \tau\dff)$\qss
and\qss
$G\off =\off A\dff(\dff \rho\dff)$\nnsp.\oss
It\dss follows\dss that\dss the definitions\qss
(\ref{product-cochains})\qss
and\qss (\ref{product-cochains-other})\qss
are equivalent\halfff.\oss

\myuppar{Products\sss of\dss several\dss cochains.}
The second\sss formula for\dss the product\halfff,\oss
i.e.\qss the formula\qss (\ref{product-cochains-other}),\oss 
admits a\sss transparent\sss extension\dss to\dss multiple products.\oss
Suppose\sss that\dss $f_{\dff i}$\trs is\dss a $p_{\fff i}$\nnsp-cochain\dss
for\qss
$i\off =\off 1\fff,\pff 2\fff,\pff \ldots\fff,\pff r$\nnsp.\qff\oss
If\dss $\sigma$\dss is\dss a\dss
$(\dff p_{\fff 1}\qff +\qff p_{\fff 2}\qff +\qff \ldots\qff +\qff p_{\fff r} \dff)$\dnsp-cube,\oss
then\vspace{6pt}
\begin{equation}
\label{product-multiple}
\quad
f_{\dff 1}\dff \cdot\dff f_{\dff 2}\trf \cdot\qff \ldots\dff \cdot\qff f_{\dff r}\qff(\dff \sigma\dff)
\qff\off =\off\off
\sum
\off
f_{\dff 1}\dff\left(\qff \tau_{\dff 1} \qff\right)
\trf \cdot\qff 
f_{\dff 2}\dff\left(\qff \tau_{\dff 2} \qff\right)
\trf \cdot\qff
\ldots
\trf \cdot\qff
f_{\dff r}\dff\left(\qff \tau_{\dff r} \qff\right)
\qff,
\end{equation}

\vspace{-6pt}
where\sss the sum\sss is\dss taken over all\sss sequences\qss
$\tau_{\dff 1}\fff,\pff
\tau_{\dff 2}\trf,\pff
\ldots\fff,\pff
\tau_{\dff r}$\qss
such\dss that\dss $\tau_{\dff i}$\dss is a $p_{\fff i}$\nnsp-face of\dss $\sigma$\nnsp,\oss
the directions of\dss these faces are pair-wise disjoint\halfff,\oss
the root\sss of\dss $\tau_{\dff 1}$\dss is equal\dss to\sss
the root\sss of\dss $\sigma$\nnsp,\oss
and\dss the peak of\dss $\tau_{\dff i}$\dss is equal\dss to\sss 
the root\sss of\dss $\tau_{\dff i\dff +\dff 1}$\dss
for\dss all\qss $i\qff \leq\qff r\qff -\qff 1$\nnsp.\oss

A sequence\qss
$\tau_{\dff 1}\fff,\pff
\tau_{\dff 2}\trf,\pff
\ldots\fff,\pff
\tau_{\dff r}$\qss
of\trs faces of\dss $\sigma$
subject\dss to\sss the above conditions
is uniquely\sss determined\dss by\dss the directions of\dss
the cubes\dss $\tau_{\dff i}$\nnsp.\oss
Indeed,\oss a cube is uniquely\sss determined\dss by\dss its root\sss and\sss direction.\oss
Therefore,\oss the above conditions allows\dss to find consecutively\dss 
the cubes\qss
$\tau_{\dff 1}\fff,\pff
\tau_{\dff 2}\trf,\pff
\ldots\fff,\pff
\tau_{\dff r}$\qss
if\trs their directions are known.\oss

\myuppar{Products of\dss $1$\dnsp-cochains.}
We are especially\dss interested\dss in\dss the case when\qss
$r\off =\off n$\qss
and every\dss $f_{\dff i}$\trs is\dss a $1$\dnsp-cochain.\oss
Then every\dss $\tau_{\dff i}$\dss
in\qss (\ref{product-multiple})\qss 
is a $1$\dnsp-cube
and\dss its\sss direction\sss is\sss a $1$\dnsp-element\sss subset\sss of\trs $I$\nnsp.\oss
In\dss this case\sss the sequences of\dss directions of\dss
cubes\dss $\tau_{\dff i}$\nsp,\oss
and\dss hence\sss the sequences\qss
$\tau_{\dff 1}\fff,\pff
\tau_{\dff 2}\trf,\pff
\ldots\fff,\pff
\tau_{\dff n}$\qss
themselves are in an obvious one-to-one correspondence
with\dss the permutations of\dss the set\qss
$I\off =\off \{\qff 1\fff,\pff 2\fff,\pff \ldots\fff,\pff n\qff \}$\nnsp.
Also,\oss
every such sequence\qss
$\tau_{\dff 1}\fff,\pff
\tau_{\dff 2}\trf,\pff
\ldots\fff,\pff
\tau_{\dff n}$\qss
has\sss the form\vspace{0pt}
\[
\quad
\tau_{\dff i}
\off =\off
\{\qff v_{\fff i\dff -\dff 1}\dff,\pff v_{\fff i} \pff\}
\]

\vspace{-15pt}
for a sequence\qss \vspace{-3pt} 
\[
\quad
v_{\dff 0}\fff,\pff
v_{\dff 1}\trf,\pff
\ldots\fff,\pff
v_{\dff n}
\qff \in\qff
\sigma 
\]

\vspace{-12pt}
such\dss that\dss the $0$\dnsp-cube\dss 
$\{\trf v_{\dff 0}\qff\}$\dss 
is\sss the root\sss of\dss $\sigma$\nnsp,\oss
the $0$\dnsp-cube\dss 
$\{\trf v_{\dff n}\qff\}$\dss 
is\sss the peak of\dss $\sigma$\nnsp,\oss
every\dss pair\qss
$\{\trf v_{\fff i\dff -\dff 1}\dff,\pff v_{\fff i} \qff\}$\qss
is a $1$\dnsp-cube,\oss
and\dss the directions of\dss these $1$\dnsp-cubes are all\sss different\halfff.\oss
We already\sss encountered\dss such sequences\qss
$v_{\dff 0}\fff,\pff
v_{\dff 1}\trf,\pff
\ldots\fff,\pff
v_{\dff n}$\nsp.\oss
As one can easily\sss check,\oss
they are nothing else but\dss the\sss pivot\sss sequences\sss
from\trs Section\qss \ref{tilings}.\oss
But\dss the reader\dss may\dss ignore\sss this fact\dss
and\dss treat\dss the above description as\sss the definition of\qss \emph{pivot\dss sequences}.

If\qss
$v_{\dff 0}\fff,\pff
v_{\dff 1}\trf,\pff
\ldots\fff,\pff
v_{\dff n}$\qss
is a pivot\sss sequence,\oss
passing\dss from\dss $v_{\dff i\dff -\dff 1}$\dss to\dss $v_{\dff i}$\dss
increases one of\dss the coordinates by\dss $1$\dss and\dss leaves other
coordinates unchanged.\oss
In\dss particular\halfff,\oss
no coordinate decreases along a pivot\sss sequence.\oss
Also,\oss
no pair\qss
$\{\trf v_{\fff i}\dff,\pff v_{\fff j} \qff\}$\qss
with\qss 
$\num{i\qff -\qff j}\qff \geq\qff 2$\qss
is\sss a $1$\dnsp-cube.\oss
It\dss follows\dss that\sss a pivot\sss sequence\qss
$v_{\dff 0}\fff,\pff
v_{\dff 1}\trf,\pff
\ldots\fff,\pff
v_{\dff n}$\qss
is\sss uniquely\sss determined\dss by\dss the set\qss
$\{\trf
v_{\dff 0}\fff,\pff
v_{\dff 1}\trf,\pff
\ldots\fff,\pff
v_{\dff n}
\qff\}$\nnsp.

\mypar{Lemma.}{cobound-codim-one}
\emph{Let\pss $f$\dss be\dss an\dss $(\fff n\dff -\dff 1\fff)$\dnsp-cochain.\qff\oss
Then}\vspace{1.25pt}
\[
\quad
\sum_{\sigma\vphantom{\sigma^2}}\qff\off
\partial^*\nsp f\trf(\dff \sigma\dff)
\off =\dff\off
\sum_{\tau\vphantom{\sigma^2}}\qff\off
f\trf(\dff \tau\dff)
\qff,
\]

\vspace{-12pt}
\emph{where\dss the first\sss sum\dss is\dss taken\sss 
over all\sss $n$\dnsp-cubes\dss $\sigma$\dss of\pss $K$\dss
and\dss the second\sss sum\dss is\dss taken\sss
over all\sss $(\fff n\dff -\dff 1\fff)$\dnsp-cubes\dss $\tau$\dss
contained\dss in\dss the boundary\dss $\bd\fff K$\nnsp.\oss}

\proof
Let\dss us consider cochains  
as formal\sss sums of\dss cubes.\oss 
Clearly,\oss it\dss is\dss sufficient\dss to prove\sss the\sss lemma
in\dss the case when\trs $f$\trs is equal\dss to an
$(\fff n\dff -\dff 1\fff)$\dnsp-cube\dss $\tau$\nnsp.\oss 
If\dss $\tau$\dss is contained\dss in\dss 
$\bd\fff K$\nnsp,\oss
then\dss $\tau$\dss is a face of\dss exactly\sss one $n$\dnsp-cube,\oss
and\dss both sums are equal\dss to\dss $1$\nnsp.\oss
If\dss $\tau$\dss is\dss not\sss contained\dss in\dss $\bd\fff K$\nnsp,\oss
then\dss $\tau$\dss
is a face of\dss exactly\dss two $n$\dnsp-cubes.\oss
In\dss this case\sss the left\sss sum is equal\dss to\dss $2$\dss
and\dss the right\sss sum is empty\sss and\dss hence is equal\dss to\dss $0$\nnsp.\oss
Since\qss $2\off =\off 0$\qss in\dss $\ftwo$\nsp,\oss
the\sss lemma\sss follows.\oss  \eproof

\mypar{Theorem.}{products-induction}
\emph{Let\pss
$h_{\dff 1}\fff,\pff h_{\dff 2}\dff,\pff \ldots\fff,\pff h_{\dff n}$\qss
be\dss $0$\dnsp-cochains of\oss $K$\dss and\qss let\oss
$f_{\dff i}
\off =\off
\partial^* h_{\dff i}$\pss
for every\dss
$i\qff \in\qff I$\nnsp.\oss
Let\dss $s$\dss be\dss the number of\qss $n$\dnsp-cubes\dss $\sigma$\dss of\pss $K$\dss
such\dss that}\qss\vspace{3pt}
\[
\quad
f_{\dff 1}\dff \cdot\dff f_{\dff 2}\trf \cdot\qff \ldots\dff \cdot\qff f_{\dff n}\qff(\dff \sigma\dff)
\off =\off 
1
\qff.
\]

\vspace{-9pt}
\emph{Let\dss $t$\dss be\dss the number of\qss $(\fff n\dff -\dff 1\fff)$\dnsp-cubes\qss 
$\tau\off \subset\off \bd\fff K$\qss
such\dss that\qss 
$h_{\dff 1}\dff(\qff \lambda^0\trf \tau\qff)
\off =\off 
1$\qss
and}\qss\vspace{3pt}
\[
\quad
f_{\dff 2}\trf \cdot\qff \ldots\dff \cdot\qff f_{\dff n}\qff(\dff \tau\dff)
\off =\off\dff 
1
\qff.
\]

\vspace{-9pt}
\emph{Then\qss $s\off \equiv\off t$\qss modulo\dss $2$\nnsp.\oss}

\proof
By\dss the definition of\trs the product\sss
of\dss a $0$\dnsp-cochain\sss with\sss 
an $(\fff n\dff -\dff 1\fff)$\dnsp-cochain\dss\vspace{2.25pt}
\[
\quad
f_{\dff 2}\trf \cdot\qff \ldots\dff \cdot\qff f_{\dff n}\qff(\dff \tau\dff)
\off =\off\dff 
1
\]

\vspace{-9.75pt}
if\trs and\dss only\trs if\oss
$h_{\dff 1}\dff(\qff \lambda^0\trf \tau\qff)
\off =\off 
1$\pss
and\pss
$f_{\dff 2}\trf \cdot\qff \ldots\dff \cdot\qff f_{\dff n}\qff(\dff \tau\dff)
\off =\off\dff 
1$\nnsp.\oss
It\dss follows\dss that\dss\vspace{3pt}
\[
\quad
\sum_{\tau\vphantom{\sigma^2}}\qff\off
h_{\dff 1}\dff \cdot\dff f_{\dff 2}\trf \cdot\qff \ldots\dff \cdot\qff f_{\dff n}\qff(\dff \tau\dff)
\off\qff =\off\qff
t
\hspace*{0.5em}\mbox{modulo}\hspace*{0.5em}
2
\qff,
\]

\vspace{-15pt}
where\sss the sum\dss is\dss taken over all\sss $(\fff n\dff -\dff 1\fff)$\dnsp-cubes\qss 
$\tau\off \subset\off \bd\fff K$\nnsp.\oss
Also,\oss obviously\halfff,\oss\vspace{3pt}
\[
\quad
\sum_{\sigma\vphantom{\sigma^2}}\qff\off
f_{\dff 1}\dff \cdot\dff f_{\dff 2}\trf \cdot\qff \ldots\dff \cdot\qff f_{\dff n}\qff(\dff \tau\dff)
\off\qff =\off\qff
s
\hspace*{0.5em}\mbox{modulo}\hspace*{0.5em}
2
\qff,
\]

\vspace{-15pt}
where\sss the sum\dss is\dss taken over all\sss $n$\dnsp-cubes\dss $\sigma$\dss of\trs $K$\nnsp.\oss
Therefore,\oss it\dss is\dss sufficient\dss to prove\sss that\dss the\sss two sums above
are equal.\oss
In order\dss to prove\sss this,\oss note\sss that\dss
the identity\qss
$\partial^* \circ\qff \partial^*
\pff =\off
0$\qss
implies\sss that\qss
$\partial^*\nsp f_{\dff i}
\off =\off
0$\qss
for every\dss $i$\nnsp.\oss
Together\dss with\qss
Lemma\qss \ref{leibniz}\qss
this\sss implies\dss that\vspace{3pt}
\[
\quad
\partial^*\left(\qff
h_{\dff 1}\dff \cdot\dff f_{\dff 2}\dff \cdot\dff f_{\dff 3}\trf \cdot\qff \ldots\dff \cdot\qff f_{\dff n}
\qff\right)
\off =\off
f_{\dff 1}\dff \cdot\dff f_{\dff 2}\dff \cdot\dff f_{\dff 3}\trf \cdot\qff \ldots\dff \cdot\qff f_{\dff n}
\qff.
\]

\vspace{-9pt}
Now\qss Lemma\qss \ref{cobound-codim-one}\qss
implies\sss that\dss these\sss two sums are indeed equal.\oss  \eproof

\mypar{Corollary\halfff.}{products-induction-face}
\emph{Let\pss
$h_{\dff 1}\fff,\pff h_{\dff 2}\dff,\pff \ldots\fff,\pff h_{\dff n}$\qss
be\dss $0$\dnsp-cochains of\oss $K$\dss and\qss let\oss
$f_{\dff i}
\off =\off
\partial^* h_{\dff i}$\pss
for every\dss
$i\qff \in\qff I$\nnsp.\oss
Let\dss $s$\dss 
be\dss the same number\dss as\dss in\qss Theorem\pss \ref{products-induction}.\qff\oss
Suppose\sss that}\vspace{3pt}
\[
\quad
h_{\dff i}\qff(\dff \rho\dff)
\off =\off 
1 
\hspace*{1em}\mbox{\emph{if}}\hspace*{1.2em}
\rho\off \subset\off \mathcal{A}_{\dff i}
\hspace{1em}\mbox{\emph{and}}\hspace{1.2em}
\]

\vspace{-37.5pt}
\[
\quad
h_{\dff i}\qff(\dff \rho\dff)
\off =\off 
0 
\hspace*{1em}\mbox{\emph{if}}\hspace*{1.2em}
\rho\off \subset\off \mathcal{B}_{\dff i}
\qff
\]

\vspace{-9pt}
\emph{for\dss every\qss $i\qff \in\qff I$\qss
and\qss let\dss $t_{\dff 1}$\dss be\dss the\dss number\dss of\pss 
$(\fff n\dff -\dff 1\fff)$\dnsp-cubes\oss 
$\tau\off \subset\off \mathcal{A}_{\dff 1}$\oss
such\dss that\qss}\qss\vspace{2.75pt}
\[
\quad
f_{\dff 2}\trf \cdot\qff \ldots\dff \cdot\qff f_{\dff n}\qff(\dff \tau\dff)
\off =\off\dff 
1
\qff.
\]

\vspace{-9.25pt}
\emph{Then\qss $s\off \equiv\off t_{\dff 1}$\qss modulo\dss $2$\nnsp.\oss}

\proof
Let\dss $t$\dss be\sss the number defined\dss in\qss
Theorem\pss \ref{products-induction}.\oss
It\dss is\dss sufficient\dss to prove\sss that\qss
$t_{\dff 1}\off =\off t$\nnsp.\oss
Let\dss us\dss consider\sss an $(\fff n\dff -\dff 1\fff)$\dnsp-cube\qss 
$\tau\off \subset\off \bd\fff K$\nnsp.\qff\oss
If\pss
$\tau\off \subset\off \mathcal{B}_{\dff 1}$\nnsp,\dff\oss
then\qss 
$h_{\dff 1}\dff(\qff \lambda^0\trf \tau\qff)
\off =\off 
0$\qss
and\dss hence\dss $\tau$\dss is\dss not\sss among\dss the cubes counted\trs
in\qss Theorem\pss \ref{products-induction}.\oss
If\vspace{1.5pt}
\[
\quad
\tau\off \subset\off \mathcal{A}_{\dff i}
\hspace{1em}\mbox{or}\hspace{1.2em}
\tau\off \subset\off \mathcal{B}_{\dff i}
\off,
\]

\vspace{-10.5pt}
then\qss
$f_{\dff i}\trf(\dff \varepsilon\trf)
\off =\off
0$\pss
for\sss every\sss $1$\dnsp-face\dss $\varepsilon$\dss of\dss $\tau$\nnsp.\oss
If\pss $i\qff \geq\qff 2$\nnsp,\oss
this implies\sss that\qss 
$f_{\dff 2}\trf \cdot\qff \ldots\dff \cdot\qff f_{\dff n}\qff(\dff \tau\dff)
\off =\off
0$\qss
and\dss hence\dss $\tau$\dss is\dss again\dss not\sss among\dss the cubes counted\trs
in\qss Theorem\pss \ref{products-induction}.\oss
It\dss follows\dss that\qss
$t_{\dff 1}\off =\off t$\qss
and\dss hence\sss the\sss corollary\dss follows\sss from\qss
Theorem\pss \ref{products-induction}.\oss  \eproof

\mypar{Theorem.}{products-faces}
\emph{Under\dss the assumptions of\pss Corollary\qss \ref{products-induction-face}\qss
the number\dss $s$\dss of\qss $n$\dnsp-cubes\dss $\sigma$\dss such\dss that\qss
$f_{\dff 1}\dff \cdot\dff f_{\dff 2}\trf \cdot\qff \ldots\dff \cdot\qff f_{\dff n}\qff(\dff \sigma\dff)
\off =\off 
1$\qss
is\dss odd.\oss}

\proof
Suppose\sss that\dss $n\off =\off 1$\nnsp.\oss
Then\dss we can\dss identify\dss $h_{\dff 1}$\dss
with a sequence of\dss $0$\nnsp's\qss and $1$\nnsp's\qss
starting\dss with $1$ and\sss ending\dss with $0$\nnsp.\oss
The\sss total\dss number of\dss changes from $0$\sss to $1$ and\qss from
$1$\sss to $0$ in\sss such a sequence\dss is\dss odd.\oss
But\dss the number of\dss changes\dss is\dss equal\dss to\dss 
$s$\nnsp.\oss
This proves\sss the\sss theorem\dss for\qss $n\off =\off 1$\nnsp.\oss
The assumptions of\qss Corollary\qss \ref{products-induction-face}\qss
imply\dss that\dss the same assumptions hold\dss for\qss $n\qff -\qff 1$\qss
in\dss the role of\dss $n$\nnsp,\oss
the face\dss $\mathcal{A}_{\dff 1}$\dss in\dss the role of\qss $K$\nnsp,\oss
and\dss the restrictions of\qss 
$h_{\dff 2}\fff,\pff h_{\dff 3}\dff,\pff \ldots\fff,\pff h_{\dff n}$\qss
to\dss $\mathcal{A}_{\dff 1}$\dss
in\dss the role of\qss 
$h_{\dff 1}\fff,\pff h_{\dff 2}\dff,\pff \ldots\fff,\pff h_{\dff n}$\nsp.\oss
Therefore\dss an\dss induction\dss by\dss $n$\dss
completes\sss the proof\halfff.\oss  \eproof

\prooftitle{Another\dss proof\halfff}
For every\qss
$i\qff \in\qff I$\qss
let\dss us define a $0$\dnsp-cochain\dss $u_{\dff i}$\dss as\sss follows\fff:\pss
$u_{\dff i}\trf(\dff \rho\dff)
\off =\off 
1$\qss
if\qss the $0$\dnsp-cube\dss $\rho$\dss
is contained\dss in\qss $\mathcal{A}_{\dff i}$\qss
and\qss
$u_{\dff i}\trf(\dff \rho\dff)
\off =\off 
0$\qss
otherwise.\oss
Clearly\halfff,\oss
the assumptions of\trs the\dss Corollary\qss \ref{products-induction-face}\qss
hold\dss for\qss
$h_{\dff i}\off =\off u_{\dff i}$\nsp.\oss
For every\dss
$i\qff \in\qff I$\nnsp,\oss
let\oss
$g_{\dff i}
\off =\off
\partial^* u_{\dff i}$\nsp.\oss
Then\qss $g_{\dff i}\dff(\dff \varepsilon\dff)
\off =\off 
1$\qss
if\dss and\dss only\dss if\trs the root\sss of\trs the $1$\dnsp-cube\dss $\varepsilon$\dss
belongs\dss to\dss $\mathcal{A}_{\dff i}$\dss and\dss its direction\dss is\dss $\{\trf i\qff\}$\nnsp.\oss
Hence\dss the summand\vspace{1.5pt}
\[
\quad
g_{\trf 1}\dff\left(\qff \tau_{\dff 1} \qff\right)
\trf \cdot\qff 
g_{\trf 2}\dff\left(\qff \tau_{\dff 2} \qff\right)
\trf \cdot\qff
\ldots
\trf \cdot\qff
g_{\dff n}\dff\left(\qff \tau_{\dff n} \qff\right)
\]

\vspace{-10.5pt}
of\qss
$g_{\dff 1}\dff \cdot\dff g_{\dff 2}\trf \cdot\qff \ldots\dff \cdot\qff g_{\dff n}\qff(\dff \sigma\dff)$\qss
for some $n$\dnsp-cube\dss $\sigma$\dss
is equal\dss to\dss $1$\dss
if\trs and\dss only\trs if\trs for every\qss
$i\qff \in\qff I$\qss the root\sss of\trs the $1$\dnsp-cube\dss $\tau_{\dff i}$\dss
belongs\sss to\dss $\mathcal{A}_{\dff i}$\dss and\dss the direction 
of\dss $\tau_{\dff i}$\dss is\dss $\{\trf i\qff\}$\nnsp.\oss

We claim\dss that\dss there is only\sss one $n$\dnsp-cube\dss $\sigma$\dss
and only\sss one sequence\qss
$\tau_{\dff 1}\fff,\pff
\tau_{\dff 2}\trf,\pff
\ldots\fff,\pff
\tau_{\dff n}$\qss
with\dss these properties.\oss
Let\qss
$v_{\dff 0}\fff,\pff v_{\dff 1}\dff,\pff \ldots\fff,\pff v_{\dff n}$\qss
be\sss the corresponding\sss pivot\sss sequence.\oss
Then\dss the $i${\nnsp}th coordinate of\dss $v_{\dff i\dff -\dff 1}$\dss
is equal\dss to\dss $0$\dss because\dss $v_{\dff i\dff -\dff 1}$\dss
is\dss the root\sss of\dss $\tau_{\dff i}$\dss and\dss hence\qss
$v_{\dff i\dff -\dff 1}\qff \in\qff \mathcal{A}_{\dff i}$\nsp.\oss
Similarly,\oss the $i${\nnsp}th coordinate of\dss $v_{\dff i}$\dss
is equal\dss to\dss $1$\dss because\dss $v_{\dff i}$\dss
is\dss the peak\sss of\dss $\tau_{\dff i}$\nsp.\oss
Since no coordinate decreases along a pivot\sss sequence,\oss
this implies\sss that\dss the $i${\nnsp}th coordinate of\dss $v_{ j}$\dss
is equal\dss to\dss $0$\dss if\qss
$j\qff \leq\qff i\qff -\qff 1$\qss
and\dss is equal\dss to\dss $1$\dss if\qss
$j\qff \geq\qff i$\nnsp.\oss
It\dss follows\dss that\qss
$v_{ j}
\off =\off 
(\qff 1\fff,\pff \ldots\fff,\pff 1\fff,\pff 0\fff,\pff \ldots\fff,\pff 0 \qff)$\nnsp,\oss
where\sss $j$\sss coordinates equal\dss to $1$ are followed\dss by\qss $n\qff -\qff j$\qss 
coordinates equal\dss to $0$\nnsp.\oss
This proves our\sss claim.\oss
Clearly,\oss this claim\dss implies\sss that\vspace{1.5pt}
\[
\quad
\sum_{\sigma\vphantom{\sigma^2}}\qff\off
g_{\dff 1}\dff \cdot\dff g_{\dff 2}\trf \cdot\qff \ldots\dff \cdot\qff g_{\dff n}\qff(\dff \sigma\dff)
\off =\dff\off
1
\qff,
\]

\vspace{-13.5pt}
where\dss $\sigma$\dss runs over\sss all\sss $n$\dnsp-cubes.\oss
Therefore\sss the\sss theorem\dss is\dss true\dss
for\qss 
$h_{\dff i}
\off =\off
u_{\dff i}$\nnsp.\oss
Let\dss us compare\sss the general\sss case with\dss this one.\oss
Let\qss
$d_{\dff i}
\off =\off
h_{\dff i}\qff -\qff u_{\dff i}$\nsp.\oss
Then\qss 
$\partial^* d_{\dff i}
\off =\off
f_{\dff i}\qff -\qff g_{\dff i}$\qss 
for every\qss
$i\qff \in\qff I$\nnsp.\oss
Together\dss with\dss
Lemma\qss \ref{leibniz}\qss
and\dss the identity\qss
$\partial^* \circ\qff \partial^*
\pff =\off
0$\qss
this implies\dss that\vspace{7.5pt}
\[
\quad
\partial^*\left(\qff
d_{\dff 1}\dff \cdot\dff f_{\dff 2}\dff \cdot\dff f_{\dff 3}\trf \cdot\qff \ldots\qff \cdot\qff f_{\dff n}
\qff\right)
\off\qff =\off\qff
f_{\dff 1}\dff \cdot\dff f_{\dff 2}\dff \cdot\dff f_{\dff 3}\trf \cdot\qff \ldots\qff \cdot\qff f_{\dff n}
\off -\off
g_{\dff 1}\dff \cdot\dff f_{\dff 2}\dff \cdot\dff f_{\dff 3}\trf \cdot\qff \ldots\qff \cdot\qff f_{\dff n}
\qff,
\]

\vspace{-30.0pt}
\[
\quad
\partial^*\left(\qff
g_{\dff 1}\dff \cdot\dff d_{\dff 2}\dff \cdot\dff f_{\dff 3}\trf \cdot\qff \ldots\qff \cdot\qff f_{\dff n}
\qff\right)
\off\qff =\off\qff
g_{\dff 1}\dff \cdot\dff f_{\dff 2}\dff \cdot\dff f_{\dff 3}\trf \cdot\qff \ldots\qff \cdot\qff f_{\dff n}
\off -\off
g_{\dff 1}\dff \cdot\dff g_{\dff 2}\dff \cdot\dff f_{\dff 3}\trf \cdot\qff \ldots\qff \cdot\qff f_{\dff n}
\qff,
\]

\vspace{-41.25pt}
\[
\quad
\hspace*{10em}
\ldots\ldots\ldots\ldots
\]

\vspace{-38.75pt}
\[
\quad
\partial^*\left(\qff
g_{\dff 1} \dff \cdot\dff g_{\dff 2}\trf \cdot
\qff \ldots\qff \cdot\qff 
g_{\dff n\dff -\dff 1}\dff \cdot\dff d_{\dff n}
\qff\right)
\off\qff =\off\qff
g_{\dff 1}\dff \cdot\dff g_{\dff 2}\trf \cdot\qff \ldots\qff \cdot\qff 
g_{\dff n\dff -\dff 1}\dff \cdot\dff f_{\dff n}
\off -\off
g_{\dff 1}\dff \cdot\dff g_{\dff 2}\trf \cdot\qff \ldots\qff \cdot\qff 
g_{\dff n\dff -\dff 1}\dff \cdot\dff g_{\dff n}
\qff.
\]

\vspace{-3pt}
By\sss summing\dss these equalities we see\sss that\vspace{6pt}
\begin{equation}
\label{comparing-fg}
\quad
\partial^* \Sigma
\off\qff =\off\qff
f_{\dff 1}\dff \cdot\dff f_{\dff 2}\trf \cdot\qff \ldots\qff \cdot\qff f_{\dff n}
\off -\off
g_{\dff 1}\dff \cdot\dff g_{\dff 2}\trf \cdot\qff \ldots\qff \cdot\qff g_{\dff n}
\qff,
\end{equation}

\vspace{-7.5pt}
where\vspace{-1.5pt}
\[
\quad
\Sigma
\off\qff =\off\qff
d_{\dff 1}\dff \cdot\dff f_{\dff 2}\dff \cdot\dff f_{\dff 3}\trf \cdot\qff \ldots\qff \cdot\qff f_{\dff n}
\qff\off +\qff\off
g_{\dff 1}\dff \cdot\dff d_{\dff 2}\dff \cdot\dff f_{\dff 3}\trf \cdot\qff \ldots\qff \cdot\qff f_{\dff n}
\qff\off +\qff\off
\ldots
\qff\off +\qff\off
g_{\dff 1}\dff \cdot\dff g_{\dff 2}\trf \cdot\qff 
\ldots
\qff\cdot\dff g_{\dff n\dff -\dff 1}\dff \cdot\dff d_{\dff n}
\qff.
\]

\vspace{-7.5pt}
We claim\dss that\dss every\sss summand\sss of\dss $\Sigma$\dss is equal\dss to\dss $0$\dss
on $(\fff n\dff -\dff 1\fff)$\dnsp-cubes contained\dss in\dss $\bd\fff K$\nnsp.\oss
Let\dss us consider\dss the\dss $j${\dnsp}th\dss summand,\oss
where\qss $j\qff \in\qff I$\nnsp.\oss
It\dss is equal\dss to\qss $g\fff \cdot\dff d\dff \cdot\fff f$\nnsp,\oss
where\vspace{4.5pt}
\[
\quad
g
\off =\off
g_{\dff 1}\dff \cdot\qff 
\ldots\qff \cdot\qff 
g_{\dff j\dff -\dff 1}
\qff,
\hspace{1.2em}
d\off =\off d_{\dff j}
\qff,
\hspace{1.2em}\mbox{and}\hspace{1.2em}
f
\off =\off
f_{\dff j\dff +\dff 1}\dff \cdot\qff 
\ldots\dff \cdot\qff 
f_{\dff n}
\off.
\]

\vspace{-7.5pt}
Let\dss $\tau$\dss be\dss an $(\fff n\dff -\dff 1\fff)$\dnsp-cube.\oss
Then\vspace{6pt}
\begin{equation}
\label{ghf-sum}
\quad
g\fff \cdot\dff d\dff \cdot\fff f\qff(\dff \tau\dff)
\qff\off =\off\off
\sum
\off
g\dff\left(\qff \tau_{\dff 1} \qff\right)
\trf \cdot\qff 
d\dff\left(\qff \tau_{\dff 2} \qff\right)
\trf \cdot\qff
f\dff\left(\qff \tau_{\dff 3} \qff\right)
\end{equation}

\vspace{-6pt}
where\sss the sum\dss is\dss taken over\dss the\sss triples\qss
$\tau_{\dff 1}\dff,\pff
\tau_{\dff 2}\dff,\pff
\tau_{\dff 3}$\qss
subject\halfff,\oss up\sss to\sss the notations,\oss 
to\sss the same conditions as\sss the sequences\qss
$\tau_{\dff 1}\fff,\pff
\tau_{\dff 2}\trf,\pff
\ldots\fff,\pff
\tau_{\dff r}$\qss
in\qss (\ref{product-multiple}).\oss

Arguing as in\dss the case of\dss the product\qss
$g_{\dff 1}\dff \cdot\dff g_{\dff 2}\trf \cdot\qff \ldots\dff \cdot\qff g_{\dff n}$\nsp,\oss
we see\sss that\dss if\qss
$g\dff\left(\qff \tau_{\dff 1} \qff\right)
\off \neq\off
0$\nnsp,\oss
then\dss the direction of\dss $\tau_{\dff 1}$\dss is\qss
$\{\trf 1\fff,\pff 2\fff,\pff \ldots\fff,\pff j\qff -\qff 1\trf\}$\qss
and\dss the peak of\dss $\tau_{\dff 1}$\dss
has\sss the form\dss $\{\dff v\qff\}$\nnsp,\oss where\qss\vspace{3pt}
\[
\quad
v
\off =\off
(\qff 1\fff,\pff \ldots\fff,\pff 1\fff,\pff a_{\dff j}\fff,\pff \ldots\fff,\pff a_{\dff n} \qff)
\qff
\]

\vspace{-9pt}
with\dss the first\qss $j\qff -\qff 1$\qss coordinates equal\dss to $1$\nnsp.\oss
Since\qss $v\qff \in\qff \tau$\nnsp,\oss
it\dss follows\dss that\dss
in\dss this case\dss $\tau$\dss is\dss not\sss contained\dss in\dss any\dss face\qss
$\mathcal{A}_{\dff i}\dff,\off \mathcal{B}_{\dff i}$\qss
with\qss $i\qff \leq\qff j\qff -\qff 1$\nnsp.\oss

The peak\dss $\{\dff v\qff\}$\dss of\dss $\tau_{\dff 1}$\dss is also\sss the root\sss of\trs
the $0$\dnsp-cube\dss $\tau_{\dff 2}$\nnsp.\oss 
If\qss
$d\dff\left(\qff \tau_{\dff 2} \qff\right)
\off \neq\off
0$\nnsp,\oss
then\qss
$0\qff <\qff a_{\dff j}\qff <\qff k$\nnsp.\oss
Similarly\dss to\sss the previous paragraph,\oss
this implies\sss that $\tau$
is not\sss contained\dss in\dss the faces\qss $\mathcal{A}_{\fff j}\dff,\off \mathcal{B}_{\fff j}$\nsp.

Finally,\oss arguing as in\dss the\sss proof\dss of\qss
Corollary\qss \ref{products-induction-face},\oss
we see\sss that\dss if\qss
$f\dff\left(\qff \tau_{\dff 3} \qff\right)
\off \neq\off
0$\nnsp,\oss
then\dss $\tau_{\dff 3}$\dss and\dss hence\dss $\tau$\dss
is\dss not\sss contained\dss in\dss any\dss face\qss
$\mathcal{A}_{\dff i}\dff,\off \mathcal{B}_{\dff i}$\qss
with\qss $i\qff \geq\qff j\qff +\qff 1$\nnsp.\oss

By\sss collecting all\dss this information\dss together\sss
we see\sss that\dss if\dss at\dss least\sss one of\dss terms
of\trs sum\qss (\ref{ghf-sum})\qss is\qss $\neq\off 0$\nnsp,\oss
then\dss $\tau$\dss is not\sss contained\dss in any
$(\fff n\dff -\dff 1\fff)$\dnsp-face of\trs $K$\nnsp,\oss
i.e.\qss not\sss contained\dss in\dss the boundary\dss $\bd\fff K$\nnsp.\oss
Since\sss this is\sss true for\sss every\qss $j\qff \in\qff I$\nnsp,\oss
it\dss follows\dss that\qss
$\Sigma\dff(\dff \tau\dff)
\off =\off
0$\qss
for\sss every $(\fff n\dff -\dff 1\fff)$\dnsp-cube\dss $\tau$\dss
contained\dss in\dss $\bd\fff K$\nnsp.\oss
Now\qss Lemma\qss \ref{cobound-codim-one}\qss implies\sss that\vspace{4.5pt}
\[
\quad
\sum_{\sigma\vphantom{\sigma^2}}\qff\off
\partial^*\dff \Sigma\qff(\dff \sigma\dff)
\off =\dff\off
0
\qff,
\]

\vspace{-13.5pt}
where\sss the sum is\sss taken over all\sss $n$\dnsp-cubes\dss $\sigma$\dss of\trs $K$\nnsp.\oss

Together\dss with\qss (\ref{comparing-fg})\qss this implies\sss that\vspace{1.5pt}
\[
\quad
\sum_{\sigma\vphantom{\sigma^2}}\qff\off
f_{\dff 1}\dff \cdot\dff f_{\dff 2}\trf \cdot\qff \ldots\dff \cdot\qff f_{\dff n}\qff(\dff \sigma\dff)
\off\qff =\pff\off
\sum_{\sigma\vphantom{\sigma^2}}\qff\off
g_{\dff 1}\dff \cdot\dff g_{\dff 2}\trf \cdot\qff \ldots\dff \cdot\qff g_{\dff n}\qff(\dff \sigma\dff)
\off\qff =\pff\off
1
\qff.
\]

\vspace{-15pt}
The\sss theorem\dss follows.\oss  \eproof

\myuppar{Remark.}
The second\dss proof\trs is\dss longer\halfff,\oss
but\dss more straightforward.\oss  
While\sss the first\dss proof\dss compares\dss
the situation\sss with a similar one in\dss lower dimension,\oss
the second proof compares\dss the general\sss situation with\dss
the simplest\sss one in\dss the same dimension and avoids\sss induction\dss by\dss $n$\nnsp.\oss
Only\dss the second\dss method\dss is\dss available in\dss
Lusternik--Schnirelmann\dss situation.\oss
See\dss Section\qss \ref{sphere-cochains}.

\mypar{Theorem.}{separation-products}
\emph{Suppose\dss that\pss
$c_{\dff 1}\fff,\pff c_{\dff 2}\dff,\pff \ldots\fff,\pff c_{\dff n}$\qss
are subsets of\oss $K$\dss and\qss let\oss
$\overline{c}_{\dff i}
\qff\off =\qff\off
K\off \smallsetminus\off c_{\dff i}$\nsp.\qff\oss
If}\vspace{3pt}
\[
\quad
\mathcal{A}_{\dff i}\off \subset\off c_{\dff i}
\hspace{1em}\mbox{\emph{and}}\hspace{1.2em}
\mathcal{B}_{\dff i}\off \subset\off \overline{c}_{\dff i}
\]

\vspace{-9pt}
\emph{for\dss every\qss $i\qff \in\qff I$\nnsp,\oss
then\dss there is a pivot\dss sequence\pss
$v_{\dff 0}\fff,\pff v_{\dff 1}\dff,\pff \ldots\fff,\pff v_{\dff n}$\qss
such\dss that}\vspace{3pt}
\[
\quad
\{\qff v_{\fff i\dff -\dff 1}\dff,\pff v_{\fff i} \pff\}
\hspace{0.8em}\mbox{\emph{intersects\qss both}}\hspace{0.8em}
c_{\dff i}
\hspace{0.8em}\mbox{\emph{and}}\hspace{1em}
\overline{c}_{\dff i}
\]

\vspace{-9pt}
\emph{for\dss every\qss $i\qff \in\qff I$\nnsp.\oss
In\dss particular\halfff,\oss
there\dss is\dss an\dss $n$\dnsp-cube\dss intersecting\trs
all\trs sets\qss $c_{\dff i}$\qss
and\oss $\overline{c}_{\dff i}$\nsp,\oss
$i\qff \in\qff I$\nnsp.\oss
The\dss number\dss of\trs such\dss pivot\dss sequences\dss is\dss odd.\oss}

\proof
For every\qss
$i\qff \in\qff I$\qss
let\dss us define a $0$\dnsp-cochain\dss $h_{\dff i}$\dss as\sss follows\fff:\pss
$h_{\dff i}\trf(\dff \rho\dff)
\off =\off 
1$\qss
if\qss the $0$\dnsp-cube\dss $\rho$\dss
is contained\dss in\qss $c_{\dff i}$\qss
and\qss
$h_{\dff i}\trf(\dff \rho\dff)
\off =\off 
0$\qss
otherwise.\oss
Let\oss
$f_{\dff i}
\off =\off
\partial^* h_{\dff i}$\pss
for every\dss
$i\qff \in\qff I$\nnsp.\oss
Then\qss
$f_{\dff i}\trf(\dff \varepsilon\dff)
\off =\off 
1$\qss
if\qss the $1$\dnsp-cube\dss $\varepsilon$\dss intersects\sss both\dss
$c_{\dff i}$\dss
and\qss 
$\overline{c}_{\dff i}$\nsp,\oss
and\pss
$f_{\dff i}\trf(\dff \varepsilon\dff)
\off =\off 
0$\qss
otherwise.\oss
It\dss follows\dss that\dss for\sss every\sss 
$n$\dnsp-cube\dss $\sigma$\vspace{1.5pt}
\[
\quad
f_{\dff 1}\dff \cdot\dff f_{\dff 2}\trf \cdot\qff \ldots\dff \cdot\qff f_{\dff n}\qff(\dff \sigma\dff)
\qff
\]

\vspace{-10.5pt}
is\dss the image\sss in\dss $\ftwo$\dss of\qss 
the number of\qss pivot\sss sequences\qss 
$v_{\dff 0}\fff,\pff
v_{\dff 1}\trf,\pff
\ldots\fff,\pff
v_{\dff n}
\qff \in\qff
\sigma$\qss
such\dss that\dss the $1$\dnsp-cube\qss
$\{\qff v_{\fff i\dff -\dff 1}\dff,\pff v_{\fff i} \pff\}$\qss
intersects\sss both\dss $c_{\dff i}$\dss
and\qss $\overline{c}_{\dff i}$\dss
for every\qss $i\qff \in\qff I$\nnsp.\oss
It\dss follows\dss that\dss 
the sum\vspace{1.5pt}
\begin{equation*}
\quad
\sum_{\sigma\vphantom{\sigma^2}}\qff\off
f_{\dff 1}\dff \cdot\dff f_{\dff 2}\trf \cdot\qff \ldots\dff \cdot\qff f_{\dff n}\qff(\dff \sigma\dff)
\qff,
\end{equation*}

\vspace{-15.0pt}
where\dss $\sigma$\dss runs over\sss all\sss $n$\dnsp-cubes,\oss
is equal\dss to\sss the\sss 
number of\trs pivot\sss sequences with\dss required\dss properties\sss
taken\dss modulo\dss $2$\nnsp.\oss
Therefore\dss
it\dss is\sss sufficient\dss to prove\sss that\dss this sum\dss
is\sss equal\dss to\dss $1$\nnsp.\oss
Clearly\halfff,\oss this sum\dss is\dss equal\dss to\sss the number\dss $s$\dss 
of\dss $n$\dnsp-cubes\dss $\sigma$\dss such\dss that\vspace{1.5pt}
\[
\quad
f_{\dff 1}\dff \cdot\dff f_{\dff 2}\trf \cdot\qff \ldots\dff \cdot\qff f_{\dff n}\qff(\dff \sigma\dff)
\off =\off 
1
\qff.
\]

\vspace{-10.5pt}
taken\dss modulo\dss $2$\nnsp.\oss
The assumptions about\dss the subsets\dss $c_{\dff i}$\dss
are equivalent\dss to\sss the assumptions of\qss Theorem\qss \ref{products-faces}\qss
about\dss the cochains\dss $h_{\dff i}$\nsp.\oss
Therefore\qss Theorem\qss \ref{products-faces}\qss implies\sss 
that\dss $s$\dss is\dss odd.\oss  \eproof

\mypar{Theorem.}{separation-products-sets}
\emph{Let\pss
$D_{\dff 1}\fff,\pff D_{\dff 2}\fff,\pff \ldots\fff,\pff D_{\dff n}$\qss
be closed\sss subsets of\oss $[\trf 0\fff,\qff 1\trf]^{n}$\dss such\dss that\qss 
for every\dss $i$\dss the set\qss $D_{\dff i}$\dss
contains\dss $A_{\dff i}$\dss and\dss is\dss disjoint\trs from\dss $B_{\dff i}$\nsp.\oss
Let\oss $\overline{D}_{\dff i}$\dss be\sss the closure of\pss
$[\trf 0\fff,\qff 1\trf]^{n}\qff \smallsetminus\qff D_{\dff i}$\nsp.\qff\oss
Then}\oss\vspace{3pt}
\[
\quad
\bigl(\dff
D_{\dff 1}\qff \cap\qff 
D_{\dff 2}\qff \cap\qff
\ldots
\qff \cap\qff
D_{\dff n}
\qff\bigr)
\off \cap\off
\bigl(\qff
\overline{D}_{\dff 1}\qff \cap\qff 
\overline{D}_{\dff 2}\qff \cap\qff
\ldots
\qff \cap\qff
\overline{D}_{\dff n}
\qff\bigr)
\off\qff \neq\off\qff
\varnothing 
\qff.
\]

\vspace{-12.0pt}
\proof
The sets\qss
$D_{\dff 1}\fff,\pff D_{\dff 2}\fff,\pff \ldots\fff,\pff D_{\dff n}$\qss
together\dss with\dss the sets\qss
$\overline{D}_{\dff 1}\fff,\pff \overline{D}_{\dff 2}\fff,\pff \ldots\fff,\pff \overline{D}_{\dff n}$\qss
form a closed covering of\qss $[\trf 0\fff,\qff 1\trf]^{n}$\qss
(of\dss course,\oss every\dss pair\qss 
$D_{\dff i}\fff,\pff \overline{D}_{\dff i}$\qss
already\dss forms a covering).\oss
Let\qss $\varepsilon\qff >\qff 0$\qss be a\dss Lebesgue\dss number of\dss this covering.\oss
Suppose\sss that\dss the number\dss $k$\dss is chosen\dss to be so large\sss that\dss
an $n$\dnsp-dimensional\sss cube in\dss $\rrr^{\fff n}$\dss with\dss
the sides of\dss the length\dss $1/k$\dss has diameter\qss
$<\qff \varepsilon$\nnsp.\oss
Let\qss 
$p\dff \colon\dff
K\qff \ttoo\qff [\trf 0\fff,\qff 1\trf]^{n}$\qss
be\sss the map defined\dss by\vspace{3pt}
\[
\quad
p\dff(\dff a_{\dff 1}\fff,\pff a_{\dff 2}\fff,\pff \ldots\fff,\pff a_{\dff n}\dff)
\off =\off
(\dff a_{\dff 1}/k\fff,\off a_{\dff 2}/k\fff,\off \ldots\fff,\off a_{\dff n}/k\qff)
\qff.
\]

\vspace{-9pt}
For each\qss $i\qff \in\qff I$\qss let\qss 
$c_{\dff i}
\off =\off
p^{\dff -\dff 1}\dff (\trf D_{\dff i}\trf)$\qss
and\qss
$\overline{c}_{\dff i}
\qff\off =\qff\off
K\off \smallsetminus\off c_{\dff i}$\nsp.\oss
Then\qss\vspace{3pt}
\[
\quad
\mathcal{A}_{\dff i}\off \subset\off c_{\dff i}
\hspace{1em}\mbox{and}\hspace{1.2em}
\mathcal{B}_{\dff i}\off \subset\off \overline{c}_{\dff i}
\]

\vspace{-9pt}
for every\qss $i\qff \in\qff I$\nnsp.\oss
Therefore,\oss the assumptions of\qss
Theorem\qss \ref{separation-products}\qss
hold\sss and\dss this\sss theorem\dss
implies\sss that\dss there exists an $n$\dnsp-cube\sss $\sigma$\sss of\trs $K$\dss 
intersecting\trs all\trs sets\qss $c_{\dff i}$\nsp,\pss
$\overline{c}_{\dff i}$\nsp,\oss
$i\qff \in\qff I$\nnsp.\oss
Since,\oss obviously,\vspace{3pt}
\[
\quad
p\dff(\dff c_{\dff i}\dff)\qss \subset\qss D_{\dff i}
\hspace{1em}\mbox{and}\hspace{1.2em}
p\dff(\qff \overline{c}_{\dff i}\dff)\qss \subset\qss \overline{D}_{\dff i}
\]

\vspace{-9pt}
for every\qss $i\qff \in\qff I$\nnsp,\oss
the image\dss $p\dff(\dff \sigma\dff)$\dss intersects\sss all\dss sets\qss $D_{\dff i}$\qss
and\oss $\overline{D}_{\dff i}$\nsp,\oss
$i\qff \in\qff I$\nnsp.\oss
The image\dss $p\dff(\dff \sigma\dff)$\dss is equal\dss to\sss the set\sss of\dss vertices of\dss an
$n$\dnsp-dimensional\sss cube with\dss
the sides of\dss the length\dss $1/k$\nnsp.\oss
By\dss the choice of\trs $k$\dss this implies\sss that\dss the diameter of\dss
$p\dff(\dff \sigma\dff)$\dss is\qss
$<\qff \varepsilon$\nnsp.\oss
It\dss follows\dss that\dss every\dss point\sss of\dss
$p\dff(\dff \sigma\dff)$\dss is\dss at\dss the distance\qss $<\qff \varepsilon$\qss
from each of\trs the\sss 
sets\qss $D_{\dff i}$\qss
and\oss $\overline{D}_{\dff i}$\nsp,\oss
$i\qff \in\qff I$\nnsp.\oss
By\dss the choice of\dss $\varepsilon$\dss this implies\dss
that\dss the intersection of\dss all\dss these sets\dss is\dss non-empty.\oss  \eproof

\myuppar{Theorem about\dss partitions.}
\emph{Let\pss
$C_{\dff 1}\fff,\pff C_{\dff 2}\fff,\pff \ldots\fff,\pff C_{\dff n}$\qss
be closed\sss subsets of\qss $[\trf 0\fff,\qff 1\trf]^{n}$\dss such\dss that\dss 
for every\dss $i$\dss the set\dss $C_{\dff i}$\dss
is\dss a\dss partition\dss between\dss $A_{\dff i}$\dss and\dss $B_{\dff i}$\dss 
in\dss the sense\sss that\dss $A_{\dff i}$\dss and\dss $B_{\dff i}$\dss
are contained\dss in\sss different\sss components of\dss the complement\qss
$[\trf 0\fff,\qff 1\trf]^{n}\pff \smallsetminus\qss C_{\dff i}$\nsp.\qff\oss
Then\oss
$C_{\dff 1}\qff \cap\qff 
C_{\dff 2}\qff \cap\qff
\ldots
\qff \cap\qff
C_{\dff n}
\off \neq\off
\varnothing$\nnsp.\oss}

\proof
This\dss theorem was already\dss proved\dss in\dss Section\qss \ref{lebesgue}.\oss
Now\dss we\dss would\dss like\sss to deduce it\dss from\qss Theorem\qss \ref{separation-products}\qss
using\trs Theorem\qss \ref{separation-products-sets}\qss as an\dss intermediary.\oss
Let\dss $U_{\fff i}$\dss be\sss the component\sss of\qss
$[\trf 0\fff,\qff 1\trf]^{n}\pff \smallsetminus\qss C_{\dff i}$\qss
containing\dss $A_{\dff i}$\nsp,\oss
and\dss let\dss $D_{\dff i}$\dss be its closure.\oss
Let\dss $\overline{D}_{\dff i}$\dss be\sss the closure of\qss
$[\trf 0\fff,\qff 1\trf]^{n}\pff \smallsetminus\qss D_{\dff i}$\nsp.\oss
Then\dss $B_{\dff i}$\dss is contained\dss in\dss $\overline{D}_{\dff i}$\nsp.\oss
If\dss a\sss point\qss
$x\qff \in\qff
[\trf 0\fff,\qff 1\trf]^{n}$\qss  
is\sss not\sss contained\dss in\dss $C_{\dff i}$\nsp,\oss
then\sss either\qss $x\qff \in\qff U_{\fff i}$\nsp,\oss
or\dss $x$\dss belongs\sss to some other component\sss of\dss the complement\qss
$[\trf 0\fff,\qff 1\trf]^{n}\pff \smallsetminus\qss C_{\dff i}$\nsp.\oss
In\dss the first\sss case\qss
$x\qff \not\in\qff \overline{D}_{\dff i}$\nsp,\oss
in\dss the second case\qss
$x\qff \not\in\qff D_{\dff i}$\nsp.\oss
It\dss follows\dss that\qss\vspace{2.25pt}
\[
\quad
D_{\dff i}\qff \cap\qff \overline{D}_{\dff i}
\off \subset\off
C_{\dff i} 
\qff.
\]

\vspace{-9.75pt}
It\dss remains\sss to apply\qss Theorem\qss \ref{separation-products-sets}.\oss  \eproof

\myuppar{Lebesgue\dss first\sss covering\dss theorem\dss revisited.}
Theorem\dss about\dss partitions
seems\sss to be weaker\dss than\trs Lebesgue\dss first\sss covering\dss theorem.\oss
But\dss the latter can\dss be deduced\dss from\trs
Theorem\dss about\dss partitions
by\dss using\qss Lebesgue\dss fusion of\dss sets\qss
(see\dss Section\qss \ref{lebesgue})\qss 
and\sss some elementary\dss topology.

\mypar{Lemma.}{extensions}
\emph{Let\pss $Z\off \subset\off [\trf 0\fff,\qff 1\trf]^{n}$\qss 
be\dss a\sss closed set\dss 
and\pss $E$\nnsp,\dss $F\off \subset\off Z$\qss be\sss its\sss closed subsets\dss 
such\dss that\pss $E\qff \cup\qff F\off =\off Z$\nnsp.\oss 
If\pss the\dss sets\pss
$E$\nnsp,\dss $F$\dss are\dss disjoint\qss from\pss $B_{\dff i}$\nsp,\dss
$A_{\dff i}$\qss respectively\halfff,\oss
then\dss there exists\sss a\sss closed\sss set\qss 
$C_{\dff i}\off \subset\off [\trf 0\fff,\qff 1\trf]^{n}$\qss
such\dss that\pss 
$C_{\dff i}\qff \cap\qff Z
\off =\off
E\qff \cap\qff F$\qss
and\qss $C_{\dff i}$\dss is\dss
separating\qss $A_{\dff i}$\qss from\pss $B_{\dff i}$\nsp.\oss}

\proof
Clearly,\oss $E\qff \cap\qff A_{\dff i}$\qss and\qss $E\qff \cap\qff F$\qss
are disjoint\sss closed subsets of\trs $E$\nnsp.\oss
It\dss follows\dss that\dss there exists a continuous function\qss
$\varphi\dff \colon\dff
E\qff \ttoo\qff [\dff {}-\qff 1\fff,\pff 0\qff]$\qss
such\dss that\vspace{4.5pt}
\[
\quad
\varphi^{\dff -\dff 1}\dff(\dff -\qff 1\dff)
\off =\off
E\qff \cap\qff A_{\dff i}
\hspace{1em}\mbox{and}\hspace{1.2em}
\varphi^{\dff -\dff 1}\dff(\dff 0\dff)
\off =\off
E\qff \cap\qff F
\qff.
\]

\vspace{-7.5pt}
Similarly,\oss there exists a continuous function\qss
$\psi\dff \colon\dff
F\qff \ttoo\qff [\trf 0\fff,\pff 1\trf]$\qss
such\dss that\vspace{4.5pt}
\[
\quad
\psi^{\dff -\dff 1}\dff(\dff 0\dff)
\off =\off
E\qff \cap\qff F
\hspace{1em}\mbox{and}\hspace{1.2em}
\varphi^{\dff -\dff 1}\dff(\dff 1\dff)
\off =\off
F\qff \cap\qff B_{\dff i}
\qff.
\]

\vspace{-7.5pt}
Clearly,\oss there is a continuous function\qss
$Z\qff \ttoo\qff [\dff {}-\qff 1\fff,\pff 1\qff]$\qss
equal\dss to\dss $\varphi$\dss on\dss $E$\dss and\dss to\dss $\psi$\dss on\dss $F$\nnsp.\oss
Moreover\halfff,\oss this function can be extended\dss to a function\qss
$A_{\dff i}\qff \cup\qff Z\qff \cup\qff B_{\dff i}
\qff \ttoo\qff
[\dff {}-\qff 1\fff,\pff 1\qff]$\qss
equal\dss to\dss $-\qff 1$\dss on\dss $A_{\dff i}$\dss
and\dss to\dss $1$\dss on\dss $B_{\dff i}$\nsp.\oss
Finally,\oss the\sss latter\dss function can be extended\dss to a continuous function\qss
$\rho\dff \colon\dff
[\trf 0\fff,\qff 1\trf]^{n}\qff \ttoo\qff [\dff {}-\qff 1\fff,\pff 1\qff]$\nnsp.\oss
A\sss trivial\sss verification\dss shows\sss that\dss
one can\dss take\qss
$C_{\dff i}\off =\off
\rho^{\dff -\dff 1}\dff(\dff 0\dff)$\nnsp.\oss  \eproof

\mypar{Lemma.}{sequences-of-sets}
\emph{Let\pss
$Z_{\dff 0}\off \supset\off 
Z_{\dff 1}\off \supset\off
\ldots\off \supset\off
Z_{\dff n}$\qss
be a decreasing\sss sequence of\dss sets.\oss
Suppose\dss that\pss
$C_{\dff 1}\dff,\pff C_{\dff 2}\dff,\pff \ldots\dff,\pff C_{\dff n}$\qss
are subsets of\pss $Z_{\dff 0}$\dss
such\dss that\qss
$C_{\dff m}\qff \cap\qff Z_{\dff m\dff -\dff 1}
\off =\off
Z_{\dff m}$\qss
for all\qss
$m\off =\off 1\fff,\pff 2\fff,\pff \ldots\fff,\pff n$\nnsp.\oss
Then\oss
$C_{\dff 1}\qff \cap\qff C_{\dff 2}\qff \cap\qff \ldots\qff \cap\qff C_{\dff m}
\off =\off
Z_{\dff m}$\oss
for all\qss
$m\off =\off 1\fff,\pff 2\fff,\pff \ldots\fff,\pff n$\nnsp.\oss} \eproof

\myuppar{Lebesgue\dss first\dss covering\dss theorem.} 
\emph{Let\qss
$D_{\dff 1}\fff,\pff D_{\dff 2}\dff,\pff \ldots\fff,\pff D_{\dff r}$\qss
be a covering\dss of\pss the unit\sss cube\qss
$[\trf 0\fff,\qff 1\trf]^{n}$\qss
by closed\dss sets.\oss
Suppose\dss that\dss none of\qss the sets\dss $D_{\dff i}$sss intersects\sss
two opposite\sss $(\fff n\dff -\dff 1\fff)$\dnsp-faces of\pss $[\trf 0\fff,\qff 1\trf]^{n}$\dnsp.\oss
Then\dss among\dss the sets\qss $D_{\dff i}$\qss there are\qss 
$n\qff +\qff 1$\qss sets with non-empty\dss intersection.\oss}

\proof
Now\dss we will\sss deduce\dss this\dss theorem\dss from\qss 
Theorem\dss about\dss partitions.\oss
Lebesgue\dss fusion of\dss sets construction\dss
from\dss Section\qss \ref{lebesgue}\qss leads\sss to\sss a covering\sss of\qss
$[\trf 0\fff,\qff 1\trf]^{n}$\qss
by\dss sets\pss
$E_{\dff 1}\fff,\pff E_{\dff 2}\dff,\pff \ldots\fff,\pff E_{\dff n\dff +\dff 1}$\qss
such\dss that\dss these sets are unions of\dss disjoint\sss collections of\dss
the sets\dss $D_{\dff i}$\dss
and\dss the conditions\qss ({\fff}i{\fff})\qss and\qss ({\fff}i{\fff}i{\fff})\qss
of\qss Theorem\qss \ref{e-coverings}\qss  
hold\qss
({\dff}for\dss $E_{\dff i}$\dss in\dss the role of\dss $e_{\dff i}$\nsp).\oss
In\dss particular\halfff,\oss
$E_{\dff i}$\dss is\dss disjoint\dss from\dss
$A_{\fff j}$\dss if\pss $i\qff >\qff j$\qss
and\dss is\dss disjoint\dss from\dss $B_{\dff i}$\dss if\pss $i\qff \leq\qff n$\nnsp.\oss
For\qss $m\qff \in\qff I$\pss let\vspace{4.5pt}
\[
\quad
X_{\dff m}
\off =\off
E_{\dff 1}\qff \cap\qff
E_{\dff 2}\qff \cap\qff
\ldots\qff \cap\qff
E_{\dff m}
\qff,
\]

\vspace{-37.5pt}
\[
\quad
Y_{\dff m}
\off =\off
E_{\dff m\dff +\dff 1}\qff \cup\qff
E_{\dff m\dff +\dff 2}\qff \cup\qff
\ldots\qff \cup\qff
E_{\dff n\dff +\dff 1}
\qff.
\]

\vspace{-7.5pt}
Then\dss 
$X_{\dff m}$\dss is\dss disjoint\dss from\dss $B_{\dff m}$\dss
and\dss $Y_{\dff m}$\dss is\dss disjoint\dss from\dss $A_{\dff m}$\dss
for every\qss $m\qff \in\qff I$\nnsp.\oss

Let\qss
$Z_{\dff 0}
\off =\off
[\trf 0\fff,\qff 1\trf]^{n}$\nnsp.\oss
Obviously,\oss $Z_{\dff 0}
\off =\off
X_{\dff 1}\qff \cup\qff Y_{\dff 1}$\nsp.\oss
The conditions\qss ({\fff}i{\fff})\qss and\qss ({\fff}i{\fff}i{\fff})\qss
imply\dss that\dss $X_{\dff 1}$\dss contains $A_{\dff 1}$
and\dss is\dss disjoint\dss from\dss $B_{\dff 1}$\nsp,\oss
and\dss $Y_{\dff 1}$\dss contains $B_{\dff 1}$
and\dss is\dss disjoint\dss from\dss $A_{\dff 1}$\nsp.\oss
Hence\vspace{4.5pt}
\[
\quad
A_{\dff 1}
\off \subset\off 
X_{\dff 1}\qff \smallsetminus\qff \bigl(\qff X_{\dff 1}\qff \cap\qff Y_{\dff 1} \qff\bigr)
\hspace{1em}\mbox{and}\hspace{1.2em}
B_{\dff 1}
\off \subset\off 
Y_{\dff 1}\qff \smallsetminus\qff \bigl(\qff X_{\dff 1}\qff \cap\qff Y_{\dff 1} \qff\bigr)
\qff.
\]

\vspace{-9pt}
Obviously\halfff,\oss the sets\qss\vspace{3pt}
\[
\quad
X_{\dff 1}\qff \smallsetminus\qff \bigl(\qff X_{\dff 1}\qff \cap\qff Y_{\dff 1} \qff\bigr)
\hspace{1em}\mbox{and}\hspace{1.2em}
Y_{\dff 1}\qff \smallsetminus\qff \bigl(\qff X_{\dff 1}\qff \cap\qff Y_{\dff 1} \qff\bigr)
\qff.
\]

\vspace{-7.5pt}
are closed\dss 
in\qss
$Z_{\dff 0}\qff \smallsetminus\qff \bigl(\qff X_{\dff 1}\qff \cap\qff Y_{\dff 1} \qff\bigr)$\nnsp.\oss
It\dss follows\dss that\qss
$X_{\dff 1}\qff \cap\qff Y_{\dff 1}$\qss
separates $A_{\dff 1}$\sss from\dss $B_{\dff 1}$\sss in\dss $Z_{\dff 0}$\nsp.

Suppose\dss that\qss $1\qff \leq\qff m\qff \leq\qff n$\qss and\dss let\pss
$Z_{\dff m}
\off =\off
X_{\dff m}\qff \cap\qff Y_{\dff m}$\nsp.\oss
In\dss particular\halfff,\oss
$Z_{\dff 1}
\off =\off
X_{\dff 1}\qff \cap\qff Y_{\dff 1}$\nsp.\oss
The sets\dss $Z_{\dff m}$\dss are analogues of\trs the cubical\sss sets\dss
$\num{\gamma_{\fff m}}$\dss
from\dss the\sss proof\dss of\pss Theorem\qss \ref{e-coverings}.\oss 
Similarly\dss to\sss that\dss proof\halfff,\oss
we would\dss like\sss to prove\sss that\qss 
$Z_{\dff n}\off \neq\off \varnothing$\nnsp.\oss
In order\dss to do\sss this,\oss 
we will\dss ``extend\fff''\qss each\sss set\dss $Z_{\dff m}$\dss
to\sss a\dss set\dss $C_{\dff m}$\dss 
separating\dss $A_{\dff m}$\dss and\dss $B_{\dff m}$\dss
in\dss $[\trf 0\fff,\qff 1\trf]^{n}$\dnsp.\oss
Since\sss the set\dss $Z_{\dff 1}$\dss itself\qss has\sss this property\halfff,\oss
we can\sss set\qss
$C_{\dff 1}\off =\dff\off Z_{\dff 1}$\nsp.\oss

To begin\dss with,\oss let\dss us\dss note\sss that\qss 
$Y_{\dff m\dff -\dff 1}
\off =\off
E_{\dff m}\qff \cup\qff Y_{\dff m}$\qss
and\dss hence\vspace{4.5pt}
\[
\quad
Z_{\dff m\dff -\dff 1}
\off =\off
X_{\dff m\dff -\dff 1}
\off \cap\off 
\bigl(\qff E_{\dff m}\qff \cup\qff Y_{\dff m} \qff\bigr)
\]

\vspace{-36pt}
\[
\quad
\phantom{Z_{\dff m\dff -\dff 1}
\off }
=\off
\bigl(\qff X_{\dff m\dff -\dff 1}\qff \cap\qff E_{\dff m}\qff\bigr)
\off \cup\off
\bigl(\qff X_{\dff m\dff -\dff 1}\qff \cap\qff Y_{\dff m}\qff\bigr)
\qff
\]

\vspace{-36pt}
\[
\quad
\phantom{Z_{\dff m\dff -\dff 1}
\off }
=\off
X_{\dff m}
\off \cup\off
\bigl(\qff X_{\dff m\dff -\dff 1}\qff \cap\qff Y_{\dff m}\qff\bigr)
\qff.
\]

\vspace{-7.5pt}
Clearly,\oss 
$X_{\dff m\dff -\dff 1}\qff \cap\qff Y_{\dff m}$\qss
is disjoint\dss from\dss $A_{\dff m}$\qss
together with\dss $Y_{\dff m}$\nsp.\oss
By\dss applying\qss Lemma\qss \ref{extensions}\qss to\vspace{4.5pt}
\[
\quad
Z\off =\off Z_{\dff m\dff -\dff 1}\qff,
\hspace{1.2em}
E\off =\off X_{\dff m}\qff,
\hspace{1em}\mbox{and}\hspace{1.2em}
F\off =\off X_{\dff m\dff -\dff 1}\qff \cap\qff Y_{\dff m}
\qff
\]

\vspace{-7.5pt}
we see\sss that\dss 
there exists\sss a\sss closed set\dss $C_{\dff m}$\dss
separating\dss $A_{\dff m}$\dss from\dss $B_{\dff m}$\dss
and such\dss that\vspace{4.5pt}
\[
\quad
C_{\dff m}\qff \cap\qff Z_{\dff m\dff -\dff 1}
\off =\off
X_{\dff m}\qff \cap\qff X_{\dff m\dff -\dff 1}\qff \cap\qff Y_{\dff m}
\off =\off
Z_{\dff m}
\qff.
\]

\vspace{-7.5pt}
Therefore\qss
$C_{\dff m}\qff \cap\qff Z_{\dff m\dff -\dff 1}
\off =\off
Z_{\dff m}$\qss
for every\qss $m\qff \in\qff I$\nnsp.\oss
Since,\oss obviously\halfff,\oss
$Z_{\dff 0}\off \supset\off 
Z_{\dff 1}\off \supset\off
\ldots\off \supset\off
Z_{\dff n}$\nsp,\oss
Lemma\qss \ref{sequences-of-sets}\qss implies\sss that\oss
$C_{\dff 1}\qff \cap\qff C_{\dff 2}\qff \cap\qff \ldots\qff \cap\qff C_{\dff n}
\off =\off
Z_{\dff n}$\nsp.\oss
Theorem\dss about\dss partitions\sss
implies\sss that\dss the intersection\qss
$C_{\dff 1}\qff \cap\qff C_{\dff 2}\qff \cap\qff \ldots\qff \cap\qff C_{\dff n}$\oss
is non-empty\halfff.\oss
It\dss follows\dss that\qss
$Z_{\dff n}$\dss is\dss non-empty\halfff.\oss 
But\vspace{3pt} 
\[
\quad
Z_{\dff n}
\off =\off
E_{\dff 1}\qff \cap\qff
E_{\dff 2}\qff \cap\qff
\ldots\qff \cap\qff
E_{\dff n\dff +\dff 1}
\]

\vspace{-9pt}
and\dss hence\sss the intersection of\trs the sets\qss 
$E_{\dff 1}\fff,\pff E_{\dff 2}\dff,\pff \ldots\fff,\pff E_{\dff n\dff +\dff 1}$\qss
is non-empty.\oss
Since\sss these sets are unions of\dss disjoint\sss collections of\dss
the sets\qss 
$D_{\dff 1}\fff,\pff D_{\dff 2}\dff,\pff \ldots\fff,\pff D_{\dff r}$\nsp,\oss
it\dss follows\dss that\sss among\dss the sets\qss 
$D_{\dff 1}\fff,\pff D_{\dff 2}\dff,\pff \ldots\fff,\pff D_{\dff r}$\qss there are\qss 
$n\qff +\qff 1$\qss sets\sss with\dss non-empty\dss intersection.\oss  \eproof

\newpage
\mysection{The\qss discrete\qss sphere\qss and\qss products\qss of\qss cubical\qss cochains}{sphere-cochains}

\myuppar{The discrete sphere.}
Since we are going\dss to explore\sss the properties of\dss the central\sss
symmetry of\dss an $n$\dnsp-dimensional\sss sphere,\oss it\dss is\dss convenient\dss to replace\sss
the cube\dss $K$\dss from\dss Section\qss \ref{cubical-cochains}\qss
by a centrally\sss symmetric $(\dff n\dff +\dff 1\fff)$\dnsp-dimensional\sss 
cube with\dss the center at\qss $0\qff \in\qff \rrr^{\fff n\dff +\dff 1}$\dnsp.\oss
Let\dss $I$\dss be\sss the set\qss 
$\{\trf 1\fff,\pff 2\fff,\pff \ldots\fff,\pff n\qff +\qff 1 \trf\}$\qss
and\dss let\dss $k$\dss be a natural\dss number\halfff.\oss
Our symmetric cube is\vspace{1.5pt} 
\[
\quad
C
\off =\off
\{\pff -\qff k\dff,\pff \ldots\dff,\pff -\qff 1\dff,\pff 0\dff,\pff 1\dff,\pff 
\ldots\dff,\pff k \pff\}^{\dff n\dff +\dff 1}
\qff.
\] 

\vspace{-10.5pt}
The\qss \emph{discrete sphere}\qss $S$\dss 
of\trs dimension\qss $n\qff -\qff 1$\qss and\dss size\dss $2\fff k$\qss is\dss
the set\sss of\dss all\dss points\qss\vspace{0pt}
\[
\quad
(\dff
a_{\dff 1}\fff,\pff
a_{\dff 2}\trf,\pff
\ldots\fff,\pff
a_{\dff n\dff +\dff 1}
\dff)
\qff \in\qff C
\]

\vspace{-12pt}
such\dss that\qss 
$a_{\dff i}\off =\off -\qff k$\pss
or\pss
$a_{\dff i}\off =\off k$\qss 
for some\qss $i\qff \in\qff I$\nnsp.\oss
All\dss notions and\dss results from\dss Section\qss \ref{cubical-cochains},\oss
except\sss of\trs products,\oss
have obvious analogues for\dss the symmetric cube\dss $C$\nnsp.\oss
In\dss particular\halfff,\oss the boundary\dss $\bd\fff C$\dss is defined,\oss
and,\dss obviously,\qss 
$S\off =\off \bd\fff C$\nnsp.\oss
The\dss \emph{$m$\dnsp-cubes}\qss of\dss $S$\dss are defined as
$m$\dnsp-cubes of\dss $C$\dss contained\dss in\dss $S$\nnsp.\oss
The $m$\dnsp-chains of\dss $S$\dss are\sss the formal\sss sums of\dss
$m$\dnsp-cubes of\dss $S$\dss with coefficients in\dss $\ftwo$\nsp,\oss
and\dss the $m$\dnsp-cochains are\sss the $\ftwo$\dnsp-valued\dss
functions on\dss the\sss set\sss of\sss $m$\dnsp-cubes\sss of\dss $S$\nnsp.\oss
The boundary operator\dss $\partial$\dss on chains and\dss the coboundary operator\dss
$\partial^*$ on cochains are defined as before,\oss
but\dss taking\sss only\dss the cubes of\dss $S$\dss into account\halfff.\oss

\myuppar{Symmetric and\sss asymmetric\sss cochains.}
Recall\dss that\qss
$\iota\dff \colon\dff
\rrr^{\dff n\dff +\dff 1}\qff \ttoo\qff \rrr^{\dff n\dff +\dff 1}$\qss 
is\dss the antipodal\dss involution,\oss
$\iota\dff(\dff x\trf)\off =\off -\qff x$\nnsp.\oss
Obviously,\oss $\iota$\dss leaves\sss the sphere\dss $S$\dss invariant\sss
and\dss maps every $m$\dnsp-cube of\dss $S$\dss into an $m$\dnsp-cube of\dss $S$\nnsp.\oss
This leads\sss to a self-map\dss $\iota_{\dff *}$\dss of\dss the $\ftwo$\dnsp-vector
space of\dss $m$\dnsp-chains of\dss $S$\nnsp.\oss
The dual\sss self-map of\dss the $\ftwo$\dnsp-vector
space of\dss $m$\dnsp-cochains of\dss $S$\dss is denoted\dss by\dss $\iota^*$\nnsp.\oss
If\dss $f$\dss is an $m$\dnsp-cochain considered as an $\ftwo$\dnsp-valued\dss
functions on\dss the\sss set\sss of\sss $m$\dnsp-cubes,\oss then\qss\vspace{2pt}
\[
\quad
\iota^*\left(\qff f \qff\right)\dff (\dff \sigma\dff)
\off\dff =\off\qff
f\qff\bigl(\qff \iota_{\dff *}\dff(\dff \sigma\dff)\qff\bigr)
\]

\vspace{-10pt}
for every $m$\dnsp-cube $\sigma$ of\qss $S$\nnsp.\pss
Obviously,\oss a cube\dss $\tau$\dss is a face of\dss another cube\dss $\sigma$\dss
if\dss and\dss only\dss if\qss 
$\iota\qff(\dff \tau\trf)$\dss is a face of\trs  
$\iota\qff(\dff \sigma\trf)$\nnsp.\oss
It\dss follows\dss that\dss $\partial$\dss and\trs $\partial^*$\sss 
commute\sss with\sss $\iota_{\dff *}$ and\trs $\iota^*$\sss respectively,\oss
i.e.\vspace{3pt}
\[
\quad
\partial\trf \circ\qff \iota_{\dff *}
\off\dff =\off\qff
\iota_{\dff *}\dff \circ\qff \partial
\hspace{1.2em}\mbox{and}\hspace{1.2em}
\partial^*\fff \circ\pff \iota^*
\off\dff =\off\qff
\iota^*\fff \circ\pff \partial^*
\qff.
\]

\vspace{-9pt}
An $m$\dnsp-cochain $f$\dss is\sss said\dss to be\qss \emph{symmetric}\oss
if\oss 
$\dis
\iota^*\left(\qff f \qff\right)
\off =\off
f$\nnsp,\oss
i.e.\oss if\vspace{2.25pt}
\[
\quad
f\qff\bigl(\qff \iota_{\dff *}\dff(\dff \sigma\dff)\qff\bigr)
\off =\dff\off
f\qff(\dff \sigma\dff)
\]

\vspace{-9.75pt}
for\dss every $m$\dnsp-cube $\sigma$ of\qss $S$\nnsp.\oss
If\dss $f$\dss is symmetric,\oss
then\dss $\partial^* (\qff f\qff)$\dss is also symmetric.\oss

Let\dss us denote by\dss $\iota$\dss also\sss the map\qss
$\ftwo\qff \ttoo\qff \ftwo$\qss defined\dss as follows\fff:\qss
$\iota\qff(\dff 0\dff)\off =\off 1$\qss
and\qss
$\iota\qff(\dff 1\dff)\off =\off 0$\nnsp.\oss
Clearly,\oss $\iota\dff \circ\iota$\qss is\dss the identity\dss map,\oss
i.e.\qss
$\iota\dff \colon\dff
\ftwo\qff \ttoo\qff \ftwo$\qss
is an\dss involution.\oss
But\dss $\iota$\dss is\dss not\sss an au\-to\-morphism of\qss $\ftwo$\nsp.\oss
An $m$\dnsp-cochain $f$\dss is\sss said\dss to be\qss \emph{asymmetric}\oss
if\oss 
$\dis
\iota^*\left(\qff f \qff\right)
\off =\off
\iota\dff \circ\dff f$\nnsp,\oss
i.e.\oss if\vspace{2.25pt}
\[
\quad
f\qff\bigl(\qff \iota_{\dff *}\dff(\dff \sigma\dff)\qff\bigr)
\off =\dff\off
\iota\pff \bigl(\qff f\qff(\dff \sigma\dff)\qff\bigr)
\]

\vspace{-9.75pt}
for\dss every $m$\dnsp-simplex $\sigma$ of\qss $S$\nnsp.\oss
This property\dss is\dss not\dss preserved\dss by\dss 
the coboundary\sss operator\dss $\partial^*$\nsp\dnsp.\oss

\mypar{Lemma.}{asymmetric-coboundary}
\emph{If\pss $h$\dss is\dss an\dss asymmetric\qss $0$\dnsp-cochain,\oss
then\dss the cochain\qss $\partial^* h$\dss
is\trs symmetric.\oss}

\proof
Let\qss 
$\sigma
\off =\off
\{\trf u\fff,\pff v\trf\}$\qss 
be\dss a $1$\dnsp-cube of\dss $S$\nnsp.\oss
Then\vspace{3pt}
\[
\quad
\partial^* h\qff(\dff \sigma\dff)
\off =\off
h\qff(\dff u\trf)\qff +\qff h\qff(\dff v\trf)
\qff
\]

\vspace{-9pt}
and\dss hence\qss
$\partial^* h\qff(\dff \sigma\dff)
\off =\off 1$\qss
if\dss and\dss only\dss if\qss
one of\dss the values\qss 
$h\qff(\dff u\trf)$\nnsp,\pss $h\qff(\dff v\trf)$\qss
is equal\dss to\dss $1$\dss and\dss the other\dss to\dss $0$\nnsp.\oss
Equivalently,\pss\vspace{3.625pt}
\[
\quad
\partial^* h\qff(\dff \sigma\dff)
\off =\off 1
\hspace{1.2em}\mbox{if\qss and\trs only\qss if}\hspace{1.2em}
h\qff(\dff u\trf)\off \neq\off h\qff(\dff v\trf)
\qff.
\]

\vspace{-8.375pt}
Since\dss $h$\dss is asymmetric,\oss
the latter condition is equivalent\dss to\qss\vspace{2.5pt}
\[
\quad
h\qff\bigl(\pff \iota\qff(\dff u\trf)\qff\halfff\bigr)
\off \neq\off\fff 
h\qff\bigl(\pff \iota\qff(\dff v\trf)\qff\halfff\bigr)
\]

\vspace{-10pt}
and\dss hence\sss to\oss
$\dis
\partial^* h\pff\bigl(\pff \iota_{\dff *}\dff(\dff \sigma\dff) \qff\bigr)
\off =\off 1$\nnsp.\oss 
It\dss follows\dss that\oss
$\dis
\iota^*\dff\left(\qff \partial^* h \qff\right) 
\off\dff =\off\qff
\partial^* h$\nnsp.\oss \eproof

\mypar{Lemma.}{sums-asymmetric}
\emph{If\pss $f\nsp,\pff g$\qss are\dss asymmetric\dss $m$\dnsp-cochains,\oss
then\dss the\dss $m$\dnsp-cochain\qss $f\qff +\qff g$\qss is symmetric.}

\proof
Let\dss $\sigma$\dss be an $m$\dnsp-cube of\trs $S$\nnsp.\oss
Then\vspace{2.25pt}
\[
\quad
f\qff(\dff \sigma\trf)\qff +\qff g\qff(\dff \sigma\trf)
\off =\dff\off
1
\]

\vspace{-9.75pt}
if\dss and\dss only\dss if\qss
one of\dss the values\qss 
$f\qff(\dff \sigma\trf)$\nnsp,\pss $g\qff(\dff \sigma\trf)$\qss
is equal\dss to\dss $1$\dss and\dss the other\dss to\dss $0$\nnsp,\oss
i.e.\qss if\dss and\dss only\dss if\pss
$f\qff(\dff \sigma\trf)\off \neq\off g\qff(\dff \sigma\trf)$\nnsp.\oss
Since\qss $f\fff,\pff g$\dss are asymmetric,\oss
the latter condition is equivalent\dss to\qss\vspace{2.25pt}
\[
\quad
f\qff\bigl(\pff \iota_{\dff *}\dff(\dff \sigma\trf)\qff\halfff\bigr)
\off \neq\off\fff 
g\qff\bigl(\pff \iota_{\dff *}\dff(\dff \sigma\trf)\qff\halfff\bigr)
\]

\vspace{-9.75pt}
and\dss hence\sss to\oss
$\dis
\iota^*\left(\qff f \qff\right)\qff (\dff \sigma\dff)
\pff +\pff\dff
\iota^*\bigl(\trf g \qff\bigr)\qff (\dff \sigma\dff)
\off =\off\dff 
1$\nnsp.\oss 
I\dss follows\dss that\oss
$\dis
\iota^*\left(\qff f \qff\right)
\pff +\pff\dff
\iota^*\dff\bigl(\trf g \qff\bigr)
\off =\off
f\pff +\pff g$\nnsp.\oss  \eproof

\myuppar{Summing\dss values of\dss symmetric $n$\dnsp-cochains.}
For a symmetric $n$\dnsp-cochain\dss $f$\qss 
we\sss will\dss interpret\vspace{0pt}
\[
\quad
\sum\nolimits_{\qff \sigma}\off
f\trf(\dff \sigma\dff)
\]

\vspace{-12pt}
as\sss the sum over any set\sss of\dss $n$\dnsp-cubes $\sigma$
containing one element\sss of\dss every\dss pair\sss of\qss
the\dss form\qss $\rho\fff,\dff\off \iota_{\dff *}\dff(\dff \rho\trf)$\nnsp.\oss
Clearly,\oss this sum does not\sss depend on\dss the choice of\dss such set\halfff.\oss

\mypar{Lemma.}{coboundary-sphere-sum}
\emph{If\pss $f$\dss is\dss a\sss symmetric\dss $(\fff n\dff -\dff 1\fff)$\dnsp-cochain,\oss
then}\oss 
$\dis
\sum\nolimits_{\qff \sigma}\off
\partial^* f\qff(\dff \sigma\dff)
\off =\dff\off
0$\nnsp.\oss

\proof
It\dss is\dss similar\dss to\sss the proof\dss of\qss Lemma\qss 
\ref{cobound-codim-one}.\oss
If\dss one considers $(\fff n\dff -\dff 1\fff)$\dnsp-cochains 
as formal\sss sums of\dss $(\fff n\dff -\dff 1\fff)$\dnsp-cubes,\oss
it\dss is\dss sufficient\dss to prove\sss the\sss lemma
in\dss the case when\dss 
$f\off =\off \tau\qff +\qff \iota_{\dff *}\dff(\dff \tau\trf)$\trs 
for\dss an
$(\fff n\dff -\dff 1\fff)$\dnsp-cube\dss $\tau$\nnsp.\oss
The $(\fff n\dff -\dff 1\fff)$\dnsp-cube\dss $\tau$\dss is\dss a\sss face of\dss
exactly\dss two $n$\dnsp-cubes,\oss
say,\oss the $n$\dnsp-cubes\dss $\sigma_{\dff 1}$\dss and\dss $\sigma_{\dff 2}$\nsp.\oss
Then\sss $\partial^*\nsp f\trf(\dff \sigma\dff)$ is\dss $1$\dss
if\pss
$\sigma\off =\off \sigma_{\dff 1}$\nsp,\qss
$\sigma_{\dff 2}$\nsp,\qss $\iota_{\dff *}\dff(\dff \sigma_{\dff 1}\trf)$\nnsp,\qss
$\iota_{\dff *}\dff(\dff \sigma_{\dff 2}\trf)$\nsp,\pss
and\trs is\dss $0$\dss otherwise.\oss
Obviously,\oss
$\iota_{\dff *}\dff(\dff \sigma_{\dff 1}\trf)
\off \neq\dff\off
\sigma_{\dff 2}$\qss
and\dss hence we may assume\sss that\sss $\sigma$ in\dss the sum\dss runs over
a set\sss containing\qss $\sigma_{\dff 1}\dff,\off \sigma_{\dff 2}$\nsp.\oss
Clearly,\oss in\dss this case\sss the sum\dss is\sss equal\dss to\qss
$1\qff +\qff 1\off =\off 0\qff \in\qff \ftwo$\nsp.\oss  \eproof

\myuppar{Products of\dss cochains of\dss $C$\dss and\dss $S$\nnsp.}
The most\sss obvious analogue of\dss the products from\dss Section\qss \ref{cubical-cochains}\qss
results from\dss identifying\dss $C$\dss with\dss
$\{\qff 0\fff,\pff 1\fff,\pff 
\ldots\dff,\pff 2\fff k \qff\}^{\dff n\dff +\dff 1}$\dss
using\dss the\sss translation\dss by\dss 
$(\dff k\fff,\pff k\fff,\pff \ldots\fff,\pff k\trf)$\nnsp.\oss
Unfortunately,\oss this\sss leads\sss to a\sss product\sss 
which\dss is\dss
not\dss $\iota$\dnsp-invariant\halfff.\oss
More precisely,\oss if\trs the product\dss is defined\dss in\dss this way,\oss
then\qss
$\iota^*\dff(\dff g\dff \cdot\dff h\dff)
\off =\off
\iota^*\dff(\dff h\dff)\dff \cdot\qff \iota^*\dff(\dff g\dff)$\nnsp,\dff\oss
i.e.\qss $\iota^*$\dss is not\sss a homomorphism,\oss but\dss an\sss
anti-homomorphism\qss ({\fff}the product\sss of\dss cochains\sss is\dss not\sss
commutative).\oss

Let\dss $\norm{\cdot}_{\dff 1}$\dss be\sss the 
$l_{\dff 1}$\nsp\dnsp-norm on\dss $\rrr^{\fff n\dff +\dff 1}$\nsp\dnsp.\oss
Recall\dss that\dss
if\qss
$x\off =\off
(\dff
x_{\dff 1}\fff,\pff
x_{\dff 2}\trf,\pff
\ldots\fff,\pff
x_{\dff n\dff +\dff 1}
\dff)$\nnsp,\oss
then\vspace{0pt}
\[
\quad
\norm{x}_{\dff 1}
\off =\off
\num{x_{\dff 1}}\qff +\qff
\num{x_{\dff 2}}\qff +\qff
\ldots\qff +\qff
\num{x_{\dff n\dff +\dff 1}}
\qff.
\]

\vspace{-12pt}
Let\dss $\sigma$\dss be\sss a cube of\trs $S$\dss or\dss $C$\nnsp.\oss
Let\dss us define\sss
the\qss \emph{root}\pss $\lambda^0\trf \sigma$\sss of\dss $\sigma$\dss
as\sss the $0$\dnsp-face\qss
$\lambda^0\trf \sigma
\off =\off
\{\qff r\qff\}$\qss
of\dss $\sigma$\dss such\dss that\dss
the $l_{\dff 1}$\nsp\dnsp-norm\dss $\norm{r\fff}_{\dff 1}$\dss
is\sss minimal\sss among $l_{\dff 1}$\nsp\dnsp-norms\sss of\dss elements of\dss $\sigma$\nnsp.\oss
Similarly,\oss the\qss \emph{peak}\pss $\lambda^1\trf \sigma$\sss of\dss $\sigma$\dss
is defined as\sss the $0$\dnsp-face\qss
$\lambda^1\trf \sigma
\off =\off
\{\qff p\qff\}$\qss
of\dss $\sigma$\dss such\dss that\dss $\norm{p\fff}_{\dff 1}$\dss
is\sss maximal\sss among $l_{\dff 1}$\nsp\dnsp-norms\sss of\dss elements of\dss $\sigma$\nnsp.\oss
More generally,\oss every cube\dss $\sigma$\dss of\dss $C$\dss has\sss the form\dss\vspace{-1.5pt}
\[
\quad
\sigma
\off =\off
\prod_{i\qff =\qff 1}^{n\dff +\dff 1}\qff
\rho_{\dff i}
\qff,
\]

\vspace{-13.5pt}
where for every\dss $i$\dss either\qss
$\rho_{\dff i}\off =\off \{\trf a_{\dff i}\fff,\pff b_{\dff i} \trf\}$\qss 
for some integers\qss $a_{\dff i}\fff,\pff b_{\dff i}$\qss
between\dss $-\qff k$\dss and\dss $k$\dss
such\dss that\qss
$\num{a_{\dff i}\qff -\qff b_{\dff i}}\off =\off 1$\nnsp,\oss
or\qss
$\rho_{\dff i}\off =\off \{\trf a_{\dff i} \trf\}$\qss 
for some integer\dss $a_{\dff i}$\dss
between\dss $-\qff k$\dss and\dss $k$\nnsp.\oss
We may assume\sss that\dss in\dss the first\sss case\qss
$\num{a_{\dff i}}\qff <\qff \num{b_{\dff i}}$\nnsp.\oss
The\qss \emph{direction}\dss $A\dff(\dff \sigma\dff)$\dss of\dss $\sigma$\dss
is\dss the set\sss of\dss subscripts\sss $i$\sss such\dss that\dss $\rho_{\dff i}$\dss
consists of\trs two elements.\oss
Suppose\sss that\qss $H\qff \subset\qff A\dff(\dff \sigma\dff)$\qss
and\qss $\varepsilon\off =\off 0$\qss or\dss $1$\nnsp.\oss
Let\dss us replace in\dss the above product\sss 
every\dss factor\qss
$\rho_{\fff i}
\off =\off
\{\trf a_{\dff i}\fff,\pff b_{\dff i} \qff\}$\qss
with\qss $i\qff \in\qff H$\pss
by\dss
$\{\dff a_{\dff i}\qff\}$\dss if\pss $\varepsilon\off =\off 0$\qss
and\dss by\dss
$\{\dff b_{\dff i}\qff\}$\dss if\pss $\varepsilon\off =\off 1$\nnsp.\oss
The resulting\dss product\dss is\dss is denoted\dss by\vspace{-0.75pt}
\[
\quad 
\lambda_{\dff H}^{\dff \varepsilon}\trf \sigma
\qff.
\]

\vspace{-12.75pt}
Clearly,\oss if\pss 
$A\off =\off A\dff(\dff \sigma\dff)$\nnsp,\oss
then\pss
$\lambda^0\trf \sigma
\off =\off
\lambda_{\dff A}^{\dff 0}\trf \sigma$\pss
and\pss
$\lambda^1\trf \sigma
\off =\off
\lambda_{\dff A}^{\dff 1}\trf \sigma$\nnsp.\oss 

The key advantage of\dss these definitions  
is\dss their\dss
$\iota$\dnsp-invariance.\oss
In\dss particular\halfff,\oss if\dss $\tau$\dss is\dss the root\sss or\dss peak of\dss $\sigma$\nnsp,\oss
then,\oss obviously,\oss
$\iota_{\dff *}\dff(\dff \tau\dff)$\dss is\dss the\sss root\sss or\dss peak of\trs
$\iota_{\dff *}\dff(\dff \sigma\dff)$\dss respectively\qss
({\fff}in contrast\dss with\dss the definition\dss based on\dss the above\sss translation,\oss
which\dss leads\sss to interchanging\dss roots and\dss peaks\sss by\dss $\iota_{\dff *}$\nsp).\oss
More generally,\oss if\pss
$H\qff \subset\qff A\dff(\dff \sigma\dff)$\nnsp,\oss
then\qss\vspace{1.5pt}
\[
\quad
\lambda_{\dff H}^{\dff \varepsilon}\trf \bigl(\qff \iota_{\dff *}\dff(\dff \sigma\dff)\qff\bigr)
\off =\dff\off
\iota_{\dff *}\dff\bigl(\qff \lambda_{\dff H}^{\dff \varepsilon}\trf\sigma\qff\bigr)
\qff.
\]

\vspace{-10.5pt}
Now we can define\sss the product\sss of\dss cochains by\dss 
one of\trs the formulas\qss (\ref{product-cochains})\qss or\qss 
(\ref{product-cochains-other})\qss
from\qss Section\qss \ref{cubical-cochains}.\oss 
The resulting definition\dss is
$\iota$\dnsp-invariant\halfff,\oss i.e.\qss
if\pss $f\nsp,\pff g$\qss are cochains
of\qss $S$\dss or\dss $C$\nnsp,\oss
then\vspace{1.5pt}
\[
\quad
\iota^*\left(\qff f\dff \cdot\qff g \qff\right)
\off =\off\qff
\iota^*\left(\qff g \qff\right)
\qff \cdot\pff
\iota^*\left(\qff f \qff\right)
\qff.
\]

\vspace{-10.5pt}
In\dss particular\halfff,\oss
the product\sss of\trs two or\dss more\sss symmetric cochains\sss
of\trs $S$\dss is\sss symmetric.\oss
As\sss in\qss Section\qss \ref{cubical-cochains},\oss
the products\sss of\dss cochains on\dss $C$\dss and\dss $S$\dss are associative
and satisfy\qss Leibniz\dss formula\qss
({\fff}i.e.\qss analogues of\qss
Lemmas\qss \ref{asso}\qss and\qss \ref{leibniz}\qss
hold\halfff).\oss
The proofs are\sss the same.\oss

\myuppar{Products of\dss $1$\dnsp-cochains.}
If\pss
$f_{\dff 1}\dff,\off f_{\dff 2}\trf,\off \ldots\dff,\off f_{\dff n}$\dss
are $1$\dnsp-cochains and\dss $\sigma$\dss is\dss an
$n$\dnsp-cube of\trs $S$\nnsp,\oss
then\vspace{2.25pt}
\[
\quad
f_{\dff 1}\dff \cdot\dff f_{\dff 2}\trf \cdot\qff \ldots\dff \cdot\qff f_{\dff n}\qff(\dff \sigma\dff)
\qff\off =\off\off
\sum
\off
f_{\dff 1}\dff\left(\qff \tau_{\dff 1} \qff\right)
\trf \cdot\qff 
f_{\dff 2}\dff\left(\qff \tau_{\dff 2} \qff\right)
\trf \cdot\qff
\ldots
\trf \cdot\qff
f_{\dff n}\dff\left(\qff \tau_{\dff n} \qff\right)
\qff,
\]

\vspace{-9.75pt}
where\sss the sum\sss is\dss taken over all\sss sequences\qss
$\tau_{\dff 1}\fff,\pff
\tau_{\dff 2}\trf,\pff
\ldots\fff,\pff
\tau_{\dff n}$\qss
such\dss that\dss $\tau_{\dff i}$\dss is a $1$\nnsp-face of\dss $\sigma$\nnsp,\oss
the directions of\dss these faces are pair-wise disjoint\halfff,\oss
the root\sss of\dss $\tau_{\dff 1}$\dss is equal\dss to\sss
the root\sss of\dss $\sigma$\nnsp,\oss
and\dss the peak of\dss $\tau_{\dff i}$\dss is equal\dss to\sss 
the root\sss of\dss $\tau_{\dff i\dff +\dff 1}$\dss
for\dss all\qss $i\qff \leq\qff n\qff -\qff 1$\nnsp.\oss
These conditions imply\dss that\dss the\sss peak\dss of\dss $\tau_{\dff n}$\dss
is equal\dss to\sss the peak\sss of\dss $\sigma$\nnsp.\oss
Every such sequence\qss
$\tau_{\dff 1}\fff,\pff
\tau_{\dff 2}\trf,\pff
\ldots\fff,\pff
\tau_{\dff n}$\qss
has\sss the form\vspace{0pt}
\[
\quad
\tau_{\dff i}
\off =\off
\{\qff v_{\fff i\dff -\dff 1}\dff,\pff v_{\fff i} \pff\}
\]

\vspace{-12pt}
for a sequence\qss 
$v_{\dff 0}\fff,\pff
v_{\dff 1}\trf,\pff
\ldots\fff,\pff
v_{\dff n}
\qff \in\qff
\sigma$\qss
such\dss that\dss the $0$\dnsp-cube\dss 
$\{\trf v_{\dff 0}\qff\}$\dss 
is\sss the root\sss of\dss $\sigma$\nnsp,\oss
the $0$\dnsp-cube\dss 
$\{\trf v_{\dff n}\qff\}$\dss 
is\sss the peak of\dss $\sigma$\nnsp,\oss
and\dss 
passing\dss from\dss $v_{\dff i\dff -\dff 1}$\dss to\dss $v_{\dff i}$\dss
increases\sss the\qss \emph{absolute\dss value}\qss 
of\dss one of\dss the coordinates by\dss $1$\dss
({\fff}the coordinate\sss itself\dss may decrease by\dss $1$\nnsp)\qss 
and\dss leaves other
coordinates unchanged.\oss
We\sss will\sss call\sss such sequences\qss
$v_{\dff 0}\fff,\pff
v_{\dff 1}\trf,\pff
\ldots\fff,\pff
v_{\dff n}$\qss
\emph{pivot\dss sequences}.\oss

\mypar{Theorem.}{power-of-generator}
\emph{Let\qss
$h_{\dff 1}\dff,\off h_{\dff 2}\trf,\off \ldots\dff,\off h_{\dff n}$\qss
be\dss asymmetric\dss $0$\dnsp-cochains\dss of\oss $S$\dss and\qss let\oss
$f_{\dff i}\off =\dff\off \partial^* h_{\dff i}$\oss
for every\qss
$i\off =\off 1\fff,\pff 2\fff,\pff \ldots\fff,\pff n$\nnsp.\oss
Then\dss there\sss exists\dss an\dss $n$\dnsp-cube\dss $\sigma$\dss 
of\oss $S$\dss such\dss that}\vspace{1.5pt}
\begin{equation}
\label{product-evaluation-1}
\quad
f_{\dff 1}\dff \cdot\dff f_{\dff 2}\trf \cdot\qff \ldots\dff \cdot\qff f_{\dff n}\qff(\dff \sigma\dff)
\off =\off
1
\qff.
\end{equation}

\vspace{-10.5pt}
\emph{Moreover\halfff,\oss 
the set\dss of\qss such\dss $n$\dnsp-cubes\dss is\dss invariant\dss under\dss
the\sss involution\qss $\iota$\dss
and\dss consists\sss of\qss an\dss odd\dss number\dss of\qss pairs\dss of\dss
$n$\dnsp-cubes\sss of\qss
the\dss form\qss $\rho\fff,\dff\off \iota_{\dff *}\dff(\dff \rho\trf)$\nnsp.\oss}

\proof
By\qss Lemma\qss \ref{asymmetric-coboundary}\qss every\sss cochain\dss $f_{\dff i}$\dss
is\dss symmetric.\oss
This\sss implies\sss that\dss
that\dss the product\sss cochain\qss
$f_{\dff 1}\dff \cdot\dff f_{\dff 2}\trf \cdot\qff \ldots\dff \cdot\qff f_{\dff n}$\qss
is\dss symmetric.\oss
In\dss turn,\oss this implies\sss that\dss the set\sss of\dss $n$\dnsp-cubes\sss $\sigma$\sss
satisfying\qss (\ref{product-evaluation-1})\qss
is\sss equal\dss to\sss the union\sss of\qss several\dss disjoint\dss pairs\sss of\qss
the\dss form\qss $\rho\fff,\dff\off \iota_{\dff *}\dff(\dff \rho\trf)$\nnsp.\oss
The number of\dss such\dss pairs\dss is\dss equal\dss
modulo\dss $2$\dss to\vspace{1.5pt}
\begin{equation}
\label{product-sum}
\quad
\sum\nolimits_{\qff \sigma}\off
f_{\dff 1}\dff \cdot\dff f_{\dff 2}\trf \cdot\qff \ldots\dff \cdot\qff f_{\dff n}\qff(\dff \sigma\dff)
\qff,
\end{equation}

\vspace{-10.5pt}
where\sss the sum\sss is\sss understood\dss as was
explained\dss before\dss Lemma\qss \ref{coboundary-sphere-sum}.\oss 
Our\dss plan\dss is\dss to prove\sss that\qss (\ref{product-sum})\qss
does not\sss depend on\dss the choice of\dss
asymmetric $0$\dnsp-cochains\sss
$h_{\dff i}$\dss
and\dss then compute\sss this sum\dss for a particular choice of\trs these cochains,\oss
as\sss in\dss the second\dss proof\dss of\qss Theorem\qss \ref{separation-products}.\oss

Let\dss us\sss look\dss what\dss happens\dss if\qss $h_{\dff 1}$\dss
is\dss replaced\dss by\sss
some asymmetric $0$\dnsp-chain\sss $g$\nsp.\oss 
Let\qss
$f\off =\dff\off \partial^* g$\oss
and\oss
$d_{\dff 1}
\off =\off
h_{\dff 1}\qff -\pff g
\off =\off
h_{\dff 1}\qff +\pff g$\nnsp.\oss
Lemma\qss \ref{sums-asymmetric}\qss implies\dss that\sss $d_{\dff 1}$\sss 
is\dss a\sss symmetric $0$\dnsp-chain.\oss
Clearly,\vspace{3pt}
\[
\quad
f_{\dff 1}\qff -\qff f
\off =\off
\partial^* d_{\dff 1}
\qff.
\]

\vspace{-9pt}
Since\qss
$f_{\dff i}
\off =\off
\partial^* h_{\dff i}$\pss
for every\dss
$i$\nnsp,\oss
the identity\qss
$\partial^* \circ\qff\fff \partial^*
\pff =\off
0$\qss
implies\sss that\qss
$\partial^* f_{\dff i}
\off =\off
0$\qss
for every\dss $i$\nnsp.\oss
Together\dss with\qss Leibniz\dss formula\dss this\sss implies\sss that\vspace{1.5pt}
\[
\quad
\partial^*\left(\qff
d_{\dff 1}\dff \cdot\dff f_{\dff 2}\trf \cdot\qff \ldots\qff \cdot\qff f_{\dff n}
\qff\right)
\off\qff =\off\qff
f_{\dff 1}\dff \cdot\dff f_{\dff 2}\trf \cdot\qff \ldots\qff \cdot\qff f_{\dff n}
\off -\off
f\dff \cdot\dff f_{\dff 2}\trf \cdot\qff \ldots\qff \cdot\qff f_{\dff n}
\qff.
\]

\vspace{-10.5pt}
Since\sss the $(\fff n\dff -\dff 1\fff)$\dnsp-cochain\qss
$d_{\dff 1}\dff \cdot\dff f_{\dff 2}\trf \cdot\qff \ldots\qff \cdot\qff f_{\dff n}$\qss
is\sss symmetric,\oss
Lemma\qss \ref{coboundary-sphere-sum}\qss implies\sss that\vspace{3pt}
\[
\quad
\sum\nolimits_{\qff \sigma}\off
f_{\dff 1}\dff \cdot\dff f_{\dff 2}\trf \cdot\qff \ldots\dff \cdot\qff f_{\dff n}\qff(\dff \sigma\dff)
\qff\off =\off\off
\sum\nolimits_{\qff \sigma}\off
f\dff \cdot\dff f_{\dff 2}\trf \cdot\qff \ldots\dff \cdot\qff f_{\dff n}\qff(\dff \sigma\dff)
\qff,
\]

\vspace{-9pt}
i.e.\qss replacing\dss $h_{\dff 1}$\dss by another asymmetric
$0$\dnsp-chain does not\sss change\dss the sum\qss (\ref{product-sum}).\oss
Similar arguments\sss apply\dss to other\dss  
$h_{\dff i}$\nsp,\oss
and\dss hence\qss (\ref{product-sum})\qss is\dss  
independent\sss
on\dss the choice of\dss cochains\dss $h_{\dff i}$\nsp.\oss
It\dss remains\sss to compute\sss the sum\qss (\ref{product-sum})\qss
for a particular\dss choice of\dss cochains\dss $h_{\dff i}$\nsp.\oss
We will\sss do\sss this in\dss the case when\qss
$h_{\dff i}\off =\off h$\qss for all\dss $i$\nnsp,\oss
where\dss $h$\dss is\dss defined as follows.\oss
Let\dss $H$\dss be\sss the set\sss of\trs points\vspace{1.5pt}
\[
\quad
a
\off =\off
(\dff
a_{\dff 1}\fff,\pff
a_{\dff 2}\trf,\pff
\ldots\fff,\pff
a_{\dff n\dff +\dff 1}
\dff)
\qff \in\qff S
\]

\vspace{-10.5pt}
such\dss that\dss the first\dss non-zero coordinate\sss $a_{\dff i}$\dss
of\dss $a$\qss
({\dff}i.e.\qss the first\dss non-zero number\dss in\dss the sequence\dss
$a_{\dff 1}\fff,\pff
a_{\dff 2}\trf,\pff
\ldots\fff,\pff
a_{\dff n\dff +\dff 1}$\nsp)\qss
is\dss negative,\oss
let\qss 
$\overline{H}\off =\off S\qff \smallsetminus\qff H$\nnsp,\oss
and\dss let\vspace{1.5pt}
\[
\quad
h\dff(\dff a\trf)\off =\off 1
\hspace{1.2em}\mbox{if}\hspace{1.2em}
a\qff \in\qff
H
\qff,
\]

\vspace{-40.5pt}
\[
\quad
h\dff(\dff a\trf)\off =\off 0
\hspace{1.2em}\mbox{if}\hspace{1.2em}
a\qff \in\qff
\overline{H}
\qff.
\]

\vspace{-10.5pt}
Clearly,\oss $\overline{H}\off =\off \iota\qff(\qff H\qff)$\qss
and\dss hence\dss the cochain\dss $h$\dss is antisymmetric.\oss
Let\qss
$f\off =\dff\off \partial^*\nsp h$\nnsp.\oss
Now\dss we are going\dss to\sss 
compute\sss the $n$\dnsp-fold\dss product\qss
$f\dff \cdot\dff f\dff \cdot\dff \ldots\dff \cdot\dff f$\nnsp.\oss

Let\dss $\sigma$\dss be\dss an $n$\dnsp-cube of\trs $S$\dss
and\qss
$v_{\dff 0}\fff,\pff
v_{\dff 1}\trf,\pff
\ldots\fff,\pff
v_{\dff n}
\qff \in\qff
\sigma$\qss
be\sss a pivot\sss sequence.\oss
For\qss $1\qff \leq\qff i\qff \leq\qff n$\qss let\qss
$\tau_{\dff i}
\off =\off
\{\qff v_{\fff i\dff -\dff 1}\dff,\pff v_{\fff i} \pff\}$\nnsp.\oss
Suppose\sss that\vspace{1.5pt}
\[
\quad
f\dff\left(\qff \tau_{\dff 1} \qff\right)
\trf \cdot\qff 
f\dff\left(\qff \tau_{\dff 2} \qff\right)
\trf \cdot\qff
\ldots
\trf \cdot\qff
f\dff\left(\qff \tau_{\dff n} \qff\right)
\off =\dff\off
1
\qff.
\]

\vspace{-10.5pt}
We will\sss show\sss that\dss the $n$\dnsp-cube\dss $\sigma$\dss and\dss the sequence\qss
$v_{\dff 0}\fff,\pff
v_{\dff 1}\trf,\pff
\ldots\fff,\pff
v_{\dff n}$\qss
are determined\dss by\dss this condition\sss up\sss to
simultaneously\dss replacing\dss them\dss by\dss 
$\iota_{\dff *}\dff(\dff \sigma\dff)$\dss and\qss
$\iota\dff(\dff v_{\dff 0}\trf)\fff,\pff
\iota\dff(\dff v_{\dff 1}\trf)\fff,\pff
\ldots\fff,\pff
\iota\dff(\dff v_{\dff n}\trf)$\nnsp.\oss

By\dss the definition,\oss
$f\dff(\dff \tau_{\dff i}\trf)\off =\off 1$\qss
if\dss and\dss only\dss if\dss
one of\trs the points\qss $v_{\fff i\dff -\dff 1}\dff,\pff v_{\fff i}$\qss
belongs\sss to\dss $H$\dss and\dss the other\dss to\dss
$\overline{H}$\nnsp.\oss
It\dss follows\dss that\dss
the points\qss 
$v_{\dff 0}\fff,\pff
v_{\dff 1}\trf,\pff
\ldots\fff,\pff
v_{\dff n}$\qss
alternate\dss between\dss $H$\dss and\dss $\overline{H}$\nnsp.\oss
It\dss will\dss be convenient\dss to say\dss that\sss a point\dss $v$\dss
is\dss a point\sss of\dss the\sss type\dss $H$\dss if\qss $v\qff \in\qff H$\nnsp,\oss
and a point\sss of\dss the\sss type\dss $\overline{H}$\dss if\qss $v\qff \in\qff \overline{H}$\nnsp.\oss
In\sss such\dss terms,\oss
when one moves from\dss $v_{\dff 0}$\dss to\dss $v_{\dff n}$\dss
along\dss the sequence\qss
$v_{\dff 0}\fff,\pff
v_{\dff 1}\trf,\pff
\ldots\fff,\pff
v_{\dff n}$\nsp,\oss
the\sss type has\sss to change at\sss each\sss step
and\dss hence has\sss to change\sss $n$\dss times.\oss

Let\qss
$v_{\fff t}
\off =\off
(\qff
a_{\dff 1}\fff,\pff
a_{\dff 2}\trf,\pff
\ldots\fff,\pff
a_{\dff n\dff +\dff 1}
\qff)$\nnsp,\oss
and\dss let\dss us\dss treat\dss 
$a_{\dff 1}\fff,\pff
a_{\dff 2}\trf,\pff
\ldots\fff,\pff
a_{\dff n\dff +\dff 1}$\qss
as integer-valued variables depending\dss on\dss
$t\off =\off 0\fff,\pff 1\fff,\pff \ldots\fff,\pff n$\nnsp.\oss
Since\qss 
$v_{\dff 0}\fff,\pff
v_{\dff 1}\trf,\pff
\ldots\fff,\pff
v_{\dff n}$\qss
is a pivot\sss sequence,\oss
each of\dss the coordinates\dss $a_{\dff i}$\dss may\sss change at\dss no more\sss than
one step.\oss
Since\dss $\sigma$\dss is\sss a cube of\trs $S$\nnsp,\oss
one of\dss the coordinates 
does not\sss depend\sss on\dss $t$\dss and\dss 
is\sss equal\dss to either\dss ${}-\qff k$\dss or\qss to\dss $k$\trs
for\dss all\dss $t$\nnsp.\oss
If\dss $a_{\dff i}$\dss is\dss the unchanging coordinate,\oss
then\dss the type of\dss $v_{\fff t}$\dss may change no more than\qss
$i\qff -\qff 1$\qss times.\oss
Since\sss the\sss type changes $n$\sss times,\oss
the unchanging coordinate is\dss $a_{\dff n\dff +\dff 1}$\dss
and\dss it\dss is equal\dss to either\dss ${}-\qff k$\dss or\dss $k$\nnsp.\oss
Suppose\sss that\qss
$a_{\dff n\dff +\dff 1}
\off =\off
-\qff k$\nnsp,\oss
the case\qss
$a_{\dff n\dff +\dff 1}
\off =\off
k$\qss
being completely\sss similar\halfff.\oss

Since\qss
$v_{\dff 0}\fff,\pff
v_{\dff 1}\trf,\pff
\ldots\fff,\pff
v_{\dff n}$\qss
is\dss a pivot sequence,\oss
for every\dss $i$\dss
the absolute value\dss $\num{a_{\dff i}}$\dss
of\trs the coordinate\dss $a_{\dff i}$\dss of\dss $v_{\fff t}$\dss does not\dss decrease
when\sss $t$\trs increases from\dss $0$\dss to\dss $n$\nnsp.\oss
It\dss follows\dss that\dss 
if\qss $a_{\dff i}\off \neq\off 0$\qss for\qss 
$t\off =\off u$\nsp,\oss
then\dss the sign of\dss $a_{\dff i}$\dss remains\sss the same for\sss
all\qss $t\qff \geq\qff u$\nsp.\oss
In\dss turn,\oss
this\sss implies\dss that\dss the\sss type of\dss $v_{\fff t}$\dss
may change no more\sss than\qss $i\qff -\qff 1$\qss times\sss
in\dss the sequence\qss
$v_{\dff u}\fff,\pff
v_{\dff u\dff +\dff 1}\trf,\pff
\ldots\fff,\pff
v_{\dff n}$\dss
({\fff}because\sss there are only\qss $i\qff -\qff 1$\qss coordinates
available\sss to change\sss the\sss type).\oss
On\dss the other\dss hand,\oss the\sss type should change at\sss
each step,\oss i.e.\qss $n\qff -\qff u$\qss times.\oss
Therefore,\oss if\pss
$a_{\dff i}\off \neq\off 0$\qss for\qss 
$t\off =\off u$\nsp,\oss
then\qss
$i\qff -\qff 1\qff \geq\qff n\qff -\qff u$\nnsp,\oss
i.e.\qss
$i\qff \geq\qff n\qff +\qff 1\qff -\qff u$\nnsp.\oss

In\dss particular\halfff,\oss
if\trs the coordinate\dss $a_{\dff i}$\dss of\dss $v_{\dff 0}$\dss
is\dss non-zero,\oss
then\qss $i\qff \geq\qff n\qff +\qff 1$\nnsp.\oss
Therefore\sss the only\dss non-zero coordinate of\dss $v_{\dff 0}$\dss
is\qss $a_{\dff n\dff +\dff 1}
\off =\off
-\qff k$\qss
and\dss hence\qss
$v_{\dff 0}
\off =\off
(\dff 0\fff,\pff \ldots\fff,\pff 0\fff,\pff {}-\qff k\dff)
\qff \in\qff H$\nnsp.\oss

If\trs the coordinate\dss $a_{\dff i}$\dss for\qss $t\off =\off 1$\qss
is\dss non-zero,\oss
then\qss $i\qff \geq\qff n$\nnsp,\oss
and\trs if\pss
$a_{\dff n}\off =\off 0$\qss for\qss $t\off =\off 1$\nnsp,\oss
then\dss $v_{\dff 1}$\dss and\dss $v_{\dff 0}$\dss have\sss
the same\sss type,\oss contrary\dss to\sss the assumption.\oss
Therefore\qss $a_{\dff n}\off \neq\off 0$\qss for\qss $t\off =\off 1$\nnsp.\oss
Since\dss $v_{\dff 0}$\dss has\dss the\dss type\dss $H$\nnsp,\oss
the\dss type\dss of\dss $v_{\dff 1}$\dss is\qss $\overline{H}$\nnsp.\oss
It\trs follows\dss that\qss $a_{\dff n}\off =\off 1$\qss for\qss $t\off =\off 1$\nnsp,\oss
and\dss hence\qss
$v_{\dff 1}
\off =\off
(\dff 0\fff,\pff \ldots\fff,\pff 0\fff,\pff 1\fff,\pff {}-\qff k\dff)
\pff \in\pff \overline{H}$\nnsp.\oss

If\trs the coordinate\dss $a_{\dff i}$\dss for\qss $t\off =\off 2$\qss
is\dss non-zero,\oss
then\qss $i\qff \geq\qff n\qff -\qff 1$\nnsp,\oss
and\trs if\pss
$a_{\dff n\dff -\dff 1}\off =\off 0$\qss for\qss $t\off =\off 2$\nnsp,\oss
then\dss $v_{\dff 2}$\dss and\dss $v_{\dff 2}$\dss have\sss
the same\sss type,\oss contrary\dss to\sss the assumption.\oss
Therefore\qss $a_{\dff n\dff -\dff 1}\off \neq\off 0$\qss for\qss $t\off =\off 2$\nnsp.\oss
Since\dss $v_{\dff 1}$\dss has\dss the\dss type\dss $\overline{H}$\nnsp,\oss
the\dss type\dss of\dss $v_{\dff 2}$\dss is\qss $H$\nnsp.\oss
It\trs follows\dss that\qss $a_{\dff n\dff -\dff 1}\off =\off {}-\qff 1$\qss for\qss $t\off =\off 2$\nnsp,\oss
and\dss hence\qss
$v_{\dff 2}
\off =\off
(\dff 0\fff,\pff \ldots\fff,\pff 0\fff,\pff {}-\qff 1\fff,\pff 1\fff,\pff {}-\qff k\dff)
\pff \in\pff H$\nnsp.\oss

By\sss continuing\dss to argue\dss in\dss this way,\oss
we see\sss that\dss the sequence\qss
$v_{\dff 0}\fff,\pff
v_{\dff 1}\trf,\pff
\ldots\fff,\pff
v_{\dff n}$\qss
with\dss the unchanging\sss coordinate\dss
$a_{\dff n\dff +\dff 1}
\off =\off
-\qff k$\pss
is\dss indeed\dss u\-nique\-ly\sss determined.\oss
Clearly,\oss there\dss is\sss unique $n$\dnsp-cube $\sigma$
containing\dss this sequence.\oss
Moreover\halfff,\oss the sequence\qss
$v_{\dff 0}\fff,\pff
v_{\dff 1}\trf,\pff
\ldots\fff,\pff
v_{\dff n}$\qss
constructed\dss in\dss the course of\trs this argument\dss
is\sss a pivot\sss sequence\sss such\dss that\dss the\sss type of\dss $v_{\dff i}$\dss
changes at\sss each step.\oss
Since\sss there is only one such sequence contained\dss in\dss $\sigma$\nnsp,\oss
it\dss follows\dss that\qss\vspace{1.5pt}
\[
\quad
f\dff \cdot\dff f\trf \cdot\qff \ldots\dff \cdot\qff f\qff(\dff \sigma\dff)
\off =\off
1
\qff.
\]

\vspace{-10.5pt}
The same is\sss true under\dss the assumption\qss
$a_{\dff n\dff +\dff 1}
\off =\off
k$\nnsp.\oss
Hence\sss there are exactly\dss two $n$\dnsp-cubes\dss $\sigma$\dss
such\dss that\qss 
$f\dff \cdot\dff f\trf \cdot\qff \ldots\dff \cdot\qff f\qff(\dff \sigma\dff)
\off =\off
1$\qss
and\dss $\iota$\dss takes each of\trs them\dss
to\sss the other\halfff,\oss
as one can either\dss verify\sss directly,\oss
or\dss deduce from\dss the $\iota$\dnsp-invariance of\dss this condition.\oss
It\dss follows\dss that\vspace{1.75pt}
\[
\quad
\sum\nolimits_{\qff \sigma}\off
f\dff \cdot\dff f\trf \cdot\qff \ldots\dff \cdot\qff f\qff(\dff \sigma\dff)
\qff\off =\off\off
1
\qff.
\]

\vspace{-10.25pt}
As we saw,\oss this implies\sss that\qss (\ref{product-sum})\qss is equal\dss to\dss $1$\dss
for\sss every choice of\qss
$h_{\dff 1}\dff,\off h_{\dff 2}\trf,\off \ldots\dff,\off h_{\dff n}$\nsp.\oss  \eproof

\mypar{Theorem.}{subsets-of-sphere-products}
\emph{Suppose\dss that\qss
$c_{\dff 1}\fff,\pff c_{\dff 2}\fff,\pff \ldots\fff,\pff c_{\dff n}$\qss
are subsets of\oss $S$\dss such\dss that\pss
$\iota\qff(\dff c_{\dff i}\trf)\off =\off S\qff \smallsetminus\qff c_{\dff i}$\qss
for every\qss 
$i\off =\off 1\fff,\pff 2\fff,\pff \ldots\fff,\pff n$\nnsp.\oss
Then\dss there exits an\dss $n$\dnsp-cube\dss $\sigma$\dss of\oss $S$\dss
such\dss that\qss for\dss every\dss $i$}\vspace{2.625pt}
\[
\quad
\sigma\qff \cap\qff c_{\dff i}
\off \neq\off
\varnothing
\hspace{1em}\mbox{\emph{and}}\hspace{1.25em}
\sigma\qff \cap\pff \bigl(\qff S\qff \smallsetminus\qff c_{\dff i} \qff\bigr)
\off \neq\off
\varnothing
\qff.
\]

\vspace{-9.375pt}
\proof
This\dss is\dss just\sss a\sss restatement\sss of\qss Theorem\qss \ref{power-of-generator}\qss
in\dss terms of\dss subsets of\trs $S$\nnsp.\oss
In\dss more details,\oss
let\qss $1\qff \leq\qff i\qff \leq\qff n$\nnsp,\oss
and\dss let\dss us define a $0$\dnsp-cochain\dss 
$h_{\dff i}$\dss as\dss follows\fff:\vspace{1.5pt}
\[
\quad
h_{\dff i}\trf(\dff a\trf)\off =\off 1
\hspace{1.2em}\mbox{if}\hspace{1.2em}
a\qff \in\qff
c_{\dff i}
\hspace{1.2em}\mbox{and}\hspace{1.2em}
h_{\dff i}\trf(\dff a\trf)\off =\off 0
\hspace{1.2em}\mbox{if}\hspace{1.2em}
a\qff \in\qff
S\qff \smallsetminus\qff c_{\dff i}
\qff.
\]

\vspace{-10.5pt}
Since\qss
$\iota\qff(\dff c_{\dff i}\trf)\off =\off S\qff \smallsetminus\qff c_{\dff i}$\nsp,\oss
every\sss cochain\dss $h_{\dff i}$\dss is\dss asymmetric.\oss
Let\qss
$f_{\dff i}\off =\dff\off \partial^* h_{\dff i}$\nnsp.\oss
If\qss $\tau$\dss is\dss a $1$\dnsp-cube,\oss
then\qss $f_{\dff i}\dff(\dff \tau\dff)\off =\off 1$\qss
if\qss and\dss only\qss if\qss
$\tau$\dss intersects both\dss $c_{\dff i}$\dss and\qss
$S\qff \smallsetminus\qff c_{\dff i}$\nsp.\oss
Therefore,\oss if\dss\vspace{1.5pt} 
\[
\quad
f_{\dff 1}\dff \cdot\dff f_{\dff 2}\trf \cdot\qff \ldots\dff \cdot\qff f_{\dff n}\qff(\dff \sigma\dff)
\off =\off
1 
\]

\vspace{-10.5pt}
for an $n$\dnsp-cube $\sigma$\nnsp,\oss
then $\sigma$ intersects\sss both\dss $c_{\dff i}$\dss and\qss
$S\qff \smallsetminus\qff c_{\dff i}$\qss
for\dss every\qss 
$i\off =\off 1\fff,\pff 2\fff,\pff \ldots\fff,\pff n$\nnsp.\oss
It\dss remains\dss to apply\qss Theorem\qss \ref{power-of-generator}.\oss  \eproof

\mypar{Theorem.}{ls-product}
\emph{Let\qss $\mathbb{S}^n$\dss be\sss the standard\dss unit\dss sphere in\qss 
$\rrr^{\fff n\dff +\dff 1}$\dnsp.\oss
Let\pss
$F_{\dff 1}\fff,\pff F_{\dff 2}\fff,\pff \ldots\fff,\pff F_{\dff n}$\qss
be closed\sss subsets of\pss $\mathbb{S}^n$\nnsp.\oss
Suppose\dss that\qss none of\qss them contains a pair of\trs antipodal\dss points.\oss
Then\dss the sets\qss
$F_{\dff i}\qff \cup\pff \iota\trf(\trf F_{\dff i}\trf)$\qss
do not\dss cover\qss $\mathbb{S}^n$\dnsp.\oss}

\proof
This\dss 
is\dss the\sss first\qss Lusternik--Schnirelmann\dss theorem\sss
and\dss was proved\dss in\dss Section\qss \ref{ls}.\oss
Now\dss we\sss will\dss deduce it\dss from\qss Theorem\qss \ref{subsets-of-sphere-products}.\oss
Let\sss use\sss the\dss $l_{\dff \infty}$\dnsp-distance
on\dss $\rrr^{\fff n\dff +\dff 1}$\dss and\dss $S$\dss
(see\dss Section\qss \ref{ls},\pss for example).\oss
Since\sss for\sss every\dss $i$\dss the sets\dss $F_{\dff i}$\dss 
and\dss $\iota\trf(\trf F_{\dff i}\trf)$\dss
are compact\sss and\sss disjoint\halfff,\oss
there is a real\dss number\qss
$\varepsilon\qff >\qff 0$\qss
such\dss that\dss the distance between\dss 
$F_{\dff i}$\dss and\qss $\iota\trf(\trf F_{\dff i}\trf)$\dss
is\qss $>\qff \varepsilon$\qss
for every\dss $i$\nnsp.\oss
Let\qss
$r\dff \colon\dff
S\qff \ttoo\qff \mathbb{S}^n$\qss
be\sss the restriction of\trs the radial\dss projection\qss
$\rrr^{\fff n\dff +\dff 1}\qff \smallsetminus\qff \{\trf 0\trf\}
\qff \ttoo\qff
\mathbb{S}^n$\qss
to\dss $S$\dss and\dss let\vspace{0pt}
\[
\quad
e_{\dff i}
\off =\off
r^{\dff -\dff 1}\qff\bigl(\qff F_{\dff i}\trf\bigr)
\qff,
\]

\vspace{-12pt}
where\qss
$i\off =\off 1\fff,\pff 2\fff,\pff \ldots\fff,\pff n$\nnsp.\oss
Obviously\halfff,\oss $r\dff \circ\qff \iota\off =\off \iota\qff \circ\dff r$\qss 
and\dss hence\qss
$\iota\qff(\trf e_{\dff i}\trf)
\off =\off
r^{\dff -\dff 1}\qff\bigl(\qff\fff \iota\trf(\qff F_{\dff i}\trf)\qff\bigr)$\nnsp.\oss
Clearly,\oss $e_{\dff i}$ is disjoint\dss from
$\iota\qff(\trf e_{\dff i}\trf)$\nnsp.\oss
Let\sss $E_{\dff i}$\sss be\sss the result\sss of\dss adding\dss to $e_{\dff i}$
all\dss points of\trs $S$\dss at\dss the\sss $l_{\dff \infty}$\dnsp-dis\-tance 
$1$ from $e_{\dff i}$\nnsp.\oss 
Equivalently\halfff,\qss 
$E_{\dff i}$\dss is\dss  
the union of\dss $e_{\dff i}$
and\sss all\sss $n$\dnsp-cubes of\dss $S$\sss intersecting $e_{\dff i}$\nsp.

If\qss $k$\dss is large enough\qss
(say,\oss if\pss $k\off >\off 3/\varepsilon$\nsp),\oss
then\dss the distance between\dss the sets\dss
$e_{\dff i}$\dss and\dss $\iota\qff(\trf e_{\dff i}\trf)$\qss 
is\qss $>\qff 3$\qss 
and\dss hence\dss 
the sets\dss $E_{\dff i}$\dss
and\dss $\iota\trf(\qff E_{\dff i}\trf)$\dss are disjoint\dss
for\dss every\dss $i$\nnsp.\oss
The complement\sss of\dss the union\qss
$E_{\dff i}\qff \cup\qff \iota\trf(\qff E_{\dff i}\trf)$\qss
in\dss $S$\dss consists of\dss several\sss disjoint\dss pairs of\trs the form\qss
$a\fff,\off \iota\trf(\dff a\dff)$\nnsp.\oss
Let\dss $c_{\dff i}$\dss be\sss the result\sss of\dss adding\dss to\dss $E_{\dff i}$\dss
one point\dss from each such\sss pair\halfff.\oss
Then\qss $\iota\qff(\dff c_{\dff i}\trf)\off =\off S\qff \smallsetminus\qff c_{\dff i}$\nsp.\oss

By\qss Theorem\qss \ref{subsets-of-sphere-products}\qss 
there\sss is\sss an $n$\dnsp-cube\dss $\sigma$\dss
of\trs $S$\dss intersecting\dss $c_{\dff i}$\dss and\qss
$S\qff \smallsetminus\qff c_{\dff i}$\qss for every\dss $i$\nnsp.\oss
If\dss $\sigma$\dss intersects\dss $e_{\dff i}$\nsp,\oss
then $\sigma$ is\dss contained\dss in\dss $E_{\dff i}\qff \subset\qff c_{\dff i}$\dss
and\dss hence is\sss disjoint\dss from\qss
$S\qff \smallsetminus\qff c_{\dff i}$\nsp.\oss
Similarly,\oss if\sss $\sigma$ intersects $\iota\qff(\trf e_{\dff i}\trf)$\nsp,\oss
then $\sigma$ is\sss disjoint\dss from $e_{\dff i}$\nsp.\oss
Therefore $\sigma$ is\sss disjoint\dss from 
$e_{\dff i}
\qff \cup\qff \iota\qff(\trf e_{\dff i}\trf)$\dss
and\dss hence $r\dff(\dff \sigma\dff)$
is\dss disjoint\dss from\dss 
$F_{\dff i}
\qff \cup\qff \iota\trf(\trf F_{\dff i}\trf)$\dss for\sss every\dss $i$\nnsp.\oss
The\dss theorem\dss follows.\oss  \eproof

\newpage
\mysection{Kuhn's\qss cubical\qss Sperner\qss lemmas}{kuhn}

\myuppar{In search of\trs the lost\dss cubes\fff:\qss I.}
In\dss the introduction\dss to his paper\qss \cite{ku}\qss
H.W.\dss Kuhn wrote\vspace{-10.5pt}

\begin{quoting}
The formal description of\trs the subdivision of\dss a\sss triangle\qss
(or\halfff,\oss more generally,\oss a simplex)\qss
is cumbersome.\oss
Does an analogue of\dss the\dss Sperner\dss Lemma\dss holds\sss for\dss the
cube,\oss for which subdivision is a\sss trivial\dss formal\sss operation\fff?\oss\vspace{-1.5pt}

The central\sss result\sss of\trs this paper is a combinatorial\dss lemma
which is a cubical\sss analogue of\dss Sperner's\dss Lemma\oss \ldots 
\end{quoting}

\vspace{-10.5pt}
Apparently,\oss Kuhn\dss wasn't\dss aware of\qss Lebesgue's\dss work\pss \cite{l-one},\pss
\cite{l-two}.\oss
Lebesgue's\dss paper\qss \cite{l-two}\qss includes\dss results which could\dss be called\qss
``cubical\sss analogues of\trs Sperner's\dss lemma''\qss if\trs not\dss for\dss the fact\dss that\dss
they\sss preceded\dss Sperner's\dss work\dss by\dss many\dss years
and\dss served as a\sss major\dss motivation for\dss it\halfff.\oss
Kuhn's\qss ``cubical\sss analogue''\qss turned out\dss to be very close\sss to\dss
Lebesgue\dss results.\oss 
It\dss is\dss stated\dss in\dss terms of\dss
a partial\sss order\dss $\leq$\dss on\dss $\rrr^{\fff n}$\dss
from\qss Section\qss \ref{tilings}.\oss
Recall\dss that\dss for\dss two points\vspace{1.5pt}
\[
\quad
a
\off =\off
(\dff
a_{\dff 1}\fff,\pff
a_{\dff 2}\trf,\pff
\ldots\fff,\pff
a_{\dff n}
\dff)
\off \in\off \rrr^{\fff n}
\hspace{1.2em}\mbox{and}\hspace{1.5em}
b
\off =\off
(\dff
b_{\dff 1}\fff,\pff
b_{\dff 2}\trf,\pff
\ldots\fff,\pff
b_{\dff n}
\dff)
\off \in\off \rrr^{\fff n}
\qff
\]

\vspace{-10.5pt}
$a\qff \leq\qff b$\qss means\sss that\qss
$a_{\dff i}\qff \leq\qff b_{\dff i}$\qss
for every\qss $i\qff \in\qff I$\nnsp.\oss
As usual,\qss
$a\qff <\qff b$\qss
means\dss that\qss
$a\qff \leq\qff b$\qss
and\qss
$a\qff \neq\qff b$\nnsp.\oss
Let\oss
$\mathbb{1}
\off =\off
(\dff
1\fff,\pff
1\trf,\pff
\ldots\fff,\pff
1
\dff)$\pss
and\oss
$\mathbb{K} 
\off =\off
\{\qff 0\dff,\pff 1\qff\}^{\dff n}$\dnsp.\oss

\myuppar{Kuhn's\qss lemma.} 
\emph{Let\oss
$L\dff \colon\dff
K\off \ttoo\off \mathbb{K}$\oss
be\dss a\dss map\dss with components\oss
$l_{\dff i}\qff \colon\qff
K\off \ttoo\off \{\qff 0\dff,\pff 1\qff\}$\nnsp,\oss
where\qss
$i\off =\off 1\fff,\pff 2\fff,\pff \ldots\fff,\pff n$\nnsp.\oss
Suppose\sss that}\dss\vspace{0pt}\vspace{-1.25pt}
\[
\quad
l_{\dff i}\dff
(\dff
a_{\dff 1}\fff,\pff
a_{\dff 2}\trf,\pff
\ldots\fff,\pff
a_{\dff n}
\trf)
\off =\off
0
\hspace{1.5em}\mbox{\emph{if}}\hspace{1.5em}
a_{\dff i}
\off =\off
0
\hspace{1.5em}\mbox{\emph{and}}\hspace{1.5em}
\]

\vspace{-40.75pt}
\[
\quad
l_{\dff i}\dff
(\dff
a_{\dff 1}\fff,\pff
a_{\dff 2}\trf,\pff
\ldots\fff,\pff
a_{\dff n}
\trf)
\off =\off
1
\hspace{1.5em}\mbox{\emph{if}}\hspace{1.5em}
a_{\dff i}
\off =\off
k
\]

\vspace{-10.75pt}
\emph{for every\qss $i\qff \in\qff I$\nnsp.\oss
Then\dss there\dss is\dss a\dss sequence\qss
$v_{\dff 0}\fff,\pff
v_{\dff 1}\trf,\pff
\ldots\fff,\pff
v_{\dff n}
\qff \in\qff
K$\qss
such\dss that}\vspace{1.25pt}
\begin{equation}
\label{kuhn-simplex}
\quad
v_{\dff 0}\off <\off
v_{\dff 1}\off <\off
\ldots\off <\off
v_{\dff n}
\off \leq\off
v_{\dff 0}\qff +\qff \mathbb{1}
\end{equation}

\vspace{-10.75pt}
\emph{and\dss each\dss map\dss $l_{\dff i}$\dss
takes\dss both\dss values\dss $0$\dss and\dss $1$\dss
on\dss the\dss set\pss
$\{\trf
v_{\dff 0}\fff,\pff
v_{\dff 1}\trf,\pff
\ldots\fff,\pff
v_{\dff n}
\qff\}$\nnsp.\oss}

\myuppar{The assumptions and\dss the conclusion of\qss Kuhn's\qss lemma.}
It\dss is only\dss natural\dss to introduce\sss 
the sets\dss $c_{\dff i}$\dss of\trs points\qss
$v\qff \in\qff K$\qss
such\dss that\qss
$l_{\dff i}\dff(\dff v\trf)\off =\off 0$\nnsp.\oss
Then\dss the complement\qss $\overline{c}_{\dff i}
\off =\off
K\qff \smallsetminus\qff c_{\dff i}$\qss
is\dss the set\sss of\trs points\qss
$v\qff \in\qff K$\qss
such\dss that\qss
$l_{\dff i}\dff(\dff v\trf)\off =\off 1$\qss
and\dss hence\dss $L$\dss is determined\dss by\dss the sets\qss
$c_{\dff 1}\fff,\pff
c_{\dff 2}\trf,\pff
\ldots\fff,\pff
c_{\dff n}$\nsp.\oss
In\dss terms\sss of\trs these sets\sss the assumptions of\qss Kuhn's\dss lemma\qss
mean\dss that\qss
$\mathcal{A}_{\dff i}\off \subset\off c_{\dff i}$\qss
and\qss
$\mathcal{B}_{\dff i}\off \subset\off \overline{c}_{\dff i}$\qss
for every\dss $i$\nnsp.\oss
The conclusion means\sss that\dss the set\qss
$\{\trf
v_{\dff 0}\fff,\pff
v_{\dff 1}\trf,\pff
\ldots\fff,\pff
v_{\dff n}
\qff\}$\qss
intersects\dss both\dss $c_{\dff i}$\dss and\qss $\overline{c}_{\dff i}$\dss
for every\dss $i$\nnsp.\oss
The inequalities\dss 
$v_{\dff 0}\qff \leq\qff v_{\dff i}\qff \leq\qff v_{\dff 0}\qff +\qff \mathbb{1}$\qss
mean\dss that\dss 
$\{\trf
v_{\dff 0}\fff,\pff
v_{\dff 1}\trf,\pff
\ldots\fff,\pff
v_{\dff n}
\qff\}$\qss
is contained\dss in an $n$\dnsp-cube of\trs $K$\dss
and\dss hence\qss (\ref{kuhn-simplex})\qss means\sss that\qss
$v_{\dff 0}\fff,\pff
v_{\dff 1}\trf,\pff
\ldots\fff,\pff
v_{\dff n}$\qss
is\sss a pivot\sss sequence in\dss the sense of\qss Sections\qss
\ref{tilings}\qss and\qss \ref{cubical-cochains}.\oss 
Now\qss Kuhn's\qss lemma\sss takes\sss the following\dss form.\oss

\myuppar{An\dss equivalent\qss form\dss of\pss Kuhn's\qss lemma.} 
\emph{Suppose\dss that\pss
$c_{\dff 1}\fff,\pff c_{\dff 2}\dff,\pff \ldots\fff,\pff c_{\dff n}$\qss
are subsets of\oss $K$\dss and\qss let\oss
$\overline{c}_{\dff i}
\qff\off =\qff\off
K\off \smallsetminus\off c_{\dff i}$\nsp.\qff\oss
Suppose also\dss that}\vspace{0pt}
\[
\quad
\mathcal{A}_{\dff i}\off \subset\off c_{\dff i}
\hspace{1em}\mbox{\emph{and}}\hspace{1.2em}
\mathcal{B}_{\dff i}\off \subset\off \overline{c}_{\dff i}
\]

\vspace{-12pt}
\emph{for\dss every\qss $i\qff \in\qff I$\nnsp.\qff\oss
Then\dss there is a pivot\dss sequence\pss
$v_{\dff 0}\fff,\pff v_{\dff 1}\dff,\pff \ldots\fff,\pff v_{\dff n}$\qss
such\dss that\qss 
$\{\qff
v_{\dff 0}\fff,\pff v_{\dff 1}\dff,\pff \ldots\fff,\pff v_{\dff n}
\qff\}$\qss
intersects\dss both\dss
$c_{\dff i}$\dss and\qss
$\overline{c}_{\dff i}$\qss
for\dss every\qss $i\qff \in\qff I$\nnsp.\oss}

\myuppar{Kuhn's\qss lemma and\dss products of\dss cochains.}
The above form of\pss Kuhn's\qss lemma\dss immediately\dss follows from\qss
Theorem\qss \ref{separation-products}\qss
and\dss is\dss a\sss somewhat\dss weakened\dss form of\trs the\sss latter\halfff.\oss
In\dss fact\halfff,\oss
the\sss theory\dss presented\dss in\dss Section\qss
\ref{cubical-cochains}\qss was partially\dss motivated\dss by\dss the desire\sss 
to
elucidate\qss Kuhn's\qss lemma.\oss

\myuppar{Kuhn's\dss reduced\dss labelings.}
Kuhn\dss proves his lemma using an\dss induction\dss by\dss $n$\nnsp.\oss
As\dss it\dss is\dss often\dss the case,\pss
using\dss induction\dss forces\dss to strengthen\dss the statement\halfff.\oss
Kuhn\qss thinks about\dss the map\dss $L$\dss as\sss a\sss labeling of\trs $K$\nnsp,\oss
and\dss his strengthening
is\dss based on associating\dss
to\dss $L$\dss a\qss \emph{reduced\dss labeling}\qss 
$r\dff \colon\dff
K\qff \ttoo\qff
\{\trf
1\fff,\pff
2\trf,\pff
\ldots\fff,\pff
n\qff +\qff 1
\qff\}$\qss
defined as follows.\oss
Let\qss
\[
\quad
c_{\dff n\dff +\dff 1}
\off =\off\dff
K\off \smallsetminus\off
(\qff
c_{\dff 1}\qff \cup\qff c_{\dff 2}\qff \cup\qff\ldots\qff \cup\qff c_{\dff n}
\qff)
\qff.
\]
Then\qss
$r\dff(\dff v\trf)$\dss
is\dss equal\dss to\sss the
minimal\dss integer\dss
$i$\dss such\dss that\qss
$v\qff \in\qff c_{\dff i}$\nsp.\oss

\myuppar{Kuhn's\qss strong\dss lemma.} 
\emph{Let\pss
$L\dff \colon\dff
K\off \ttoo\off \mathbb{K}$\oss
be\dss a\dss map satisfying\dss
the\sss assumptions of\pss Kuhn's\qss lemma,\oss
and\dss let\sss
$r$\sss be\sss the\dss reduced\dss labeling\dss
associated\dss with\qss $L$\nnsp.\oss
Then\dss there\dss exist\dss a\sss pivot\sss sequence\qss
$v_{\dff 0}\fff,\pff
v_{\dff 1}\trf,\pff
\ldots\fff,\pff
v_{\dff n}
\qff \in\qff
K$\qss
such\dss that\qss 
the reduced\dss labeling\qss $r$\dss
maps\qss
$\{\trf
v_{\dff 0}\fff,\qff
v_{\dff 1}\trf,\qff
\ldots\fff,\qff
v_{\dff n}
\trf\}$\qss
onto\dss the\dss set\qss
$\{\qff
1\dff,\off
2\dff,\off
\ldots\dff,\off
n\qff +\qff 1
\qff\}$\nnsp.\oss
Moreover\halfff,\oss
the\dss number\dss of\qss such\sss pivot\sss sequences\dss is\dss odd.\oss}\vspace{0.25pt}

\myuppar{Triangulating\sss cubes.}
Regardless of\dss his interest\dss in working\sss with cubes
as opposed\dss to simplices,\oss
Kuhn's\qss proofs\sss of\trs his lemmas starts with\dss triangulating\dss cubes,\pss
albeit\dss in a canonical\dss way.\oss
The\sss canonical\dss triangulations of\dss cubes used\dss by\qss Kuhn\dss were discovered\dss by\dss
H.\dss Freudenthal\qss \cite{f},\pss
who worked\dss with solid cubes.\oss
The geometry\sss of\dss these\sss triangulations\sss is\sss
more\sss transparent\sss in\sss the case of\dss solid cubes,\oss
but\dss in\dss the discrete context\dss it\sss
is more natural\dss to work\dss with\dss the corresponding\dss 
abstract\dss simplicial\sss complexes.\oss
We refer\dss to\qss Appendix\qss \ref{freudenthal}\qss
for a discussion of\dss geometry\halfff.\oss

Let\dss us\dss turn\dss $K$\dss into an abstract\sss simplicial\sss complex.\oss
Naturally,\oss its set\sss of\dss vertices\dss is\dss $K$\dss
itself\halfff.\oss
For every\dss pivot\sss sequence\qss
$v_{\dff 0}\fff,\pff
v_{\dff 1}\trf,\pff
\ldots\fff,\pff
v_{\dff n}
\qff \in\qff
K$\qss
the set\qss
$\{\trf
v_{\dff 0}\fff,\pff
v_{\dff 1}\trf,\pff
\ldots\fff,\pff
v_{\dff n}
\qff\}$\qss
is\dss a\sss simplex,\oss as also every\dss its subset\halfff.\oss
There are no other simplices.\oss
By an abuse of\dss notation,\oss we denote\sss this simplicial\sss
complex also by\dss $K$\nnsp.\oss
Also,\qss we denote\sss by\sss $\bd\fff K$\sss the simplicial\sss complex\sss
having\dss the boundary\sss $\bd\fff K$\sss as\sss the set\sss of\dss vertices
and\dss the simplices of\trs $K$\dss contained\dss in\sss $\bd\fff K$\sss
as\dss its\dss simplices.\oss

As one may expect\halfff,\oss the simplicial\sss complex\dss $K$\dss
is a pseudo-manifold\sss and\dss its boundary\dss is\dss the simplicial\sss
complex\sss $\bd\fff K$\nnsp.\oss
By\dss the very definition,\oss it\dss is\sss dimensionally\dss homogenous\fff:\oss
every simplex is a face of\dss an $n$\dnsp-simplex.\oss
We will\dss not\dss use\sss the fact\dss that\dss $K$\dss is strongly connected
and\dss leave its proof\trs to\sss the interested\dss readers.\oss
Let\dss us prove\sss the non-branching\dss property.\oss

\mypar{Lemma.}{cube-non-branching}
\emph{An\sss $(\fff n\dff -\dff 1\fff)$\dnsp-simplex of\oss $K$\dss
is a\dss face of\dss exactly\dss two $n$\dnsp-simplices of\oss $K$\dss
if\pss it\dss is\sss not\dss contained\dss in\pss $\bd\fff K$\nnsp,\oss
and\dss of\qss exactly\sss one $n$\dnsp-simplex if\pss it\dss is\dss contained\dss
in\pss $\bd\fff K$\nnsp.\oss}

\proof
Let\qss
$\varepsilon_{\dff 1}\fff,\pff
\varepsilon_{\dff 2}\trf,\pff
\ldots\fff,\pff
\varepsilon_{\dff n}$\qss
be\sss the standard\dss basis of\trs the vector space\dss $\rrr^{\fff n}$\dnsp.\oss

Let\qss
$v_{\dff 0}\fff,\pff
v_{\dff 1}\trf,\pff
\ldots\fff,\pff
v_{\dff n}$\qss
be a\sss pivot\sss sequence of\trs vertices of\trs $K$\nnsp.\oss
Removing a point\dss $v_{\dff i}$\dss from\dss the $n$\dnsp-simplex\qss
$\sigma
\off =\off
\{\trf
v_{\dff 0}\fff,\pff
v_{\dff 1}\trf,\pff
\ldots\fff,\pff
v_{\dff n}
\qff\}$\qss
results\dss in\dss an $(\fff n\dff -\dff 1\fff)$\dnsp-simplex\halfff,\oss
which we denote by\dss $\tau_{\dff i}$\nnsp.\oss
By\dss the definition,\oss every\sss $(\fff n\dff -\dff 1\fff)$\dnsp-simplex
can\dss be obtained\dss in\dss this way.\oss

Suppose first\dss that\qss 
$i\off \neq\off 0\fff,\pff n$\nnsp.\oss
Then\dss $\tau_{\dff i}$\dss contains both\dss $v_{\dff 0}$\dss
and\qss $v_{\dff n}\off =\off v_{\dff 0}\qff +\qff \mathbb{1}$\qss
and\dss hence is not\sss contained\dss in\dss $\bd\fff K$\nnsp.\qff\oss
In a\sss pivot\sss sequence 
different\sss coordinates are increasing\dss by $1$
at\sss different\sss steps and\dss
therefore\vspace{3.5pt}
\[
\quad
v_{\dff i}\off =\off v_{\dff i\dff -\dff 1}\qff +\pff \varepsilon_{\dff p}
\hspace{1.2em}\mbox{and}\hspace{1.2em}
v_{\dff i\dff +\dff 1}\off =\off v_{\dff i}\qff +\qff \varepsilon_{\dff q}
\]

\vspace{-8.5pt}
for some\qss $p\off \neq\off q$\nnsp.\oss
Let\qss
$w\off =\off v_{\dff i\dff -\dff 1}\qff +\pff \varepsilon_{\dff q}$\nsp.\oss
Then\qss
$v_{\dff i\dff +\dff 1}\off =\off w\qff +\qff \varepsilon_{\dff p}$\qss
and\dss hence replacing\dss $v_{\dff i}$\dss by\dss $w$\dss
in\dss the sequence\qss
$v_{\dff 0}\fff,\pff
v_{\dff 1}\trf,\pff
\ldots\fff,\pff
v_{\dff n}$\qss
results\sss in\sss a\sss pivot\sss sequence satisfying.\oss
This,\oss in\dss turn,\oss implies\sss that\sss
replacing\dss $v_{\dff i}$\dss by\dss $w$\dss in\qss
$\sigma$\dss results\dss in\dss an $n$\dnsp-simplex\dss having\sss $\tau_{\dff i}$\dss
as\dss a\dss face.\oss
Since in a\sss pivot\sss sequence 
different\sss coordinates are increasing\dss 
at\sss different\sss steps,\oss
$v_{\dff i}$\dss can\dss be replaced\sss only\dss by\dss $w$\nsp.\oss
Therefore\dss there are exactly\dss two $n$\dnsp-simplices\sss
having\dss $\tau_{\dff i}$\dss as a face.\oss
This completes\sss the proof\qss in\dss the case when\qss
$i\off \neq\off 0\fff,\pff n$\nnsp.\oss

Suppose now\dss that\dss $i\off =\off 0$\nnsp.\oss
If\qss $j\qff >\qff 0$\nnsp,\oss
then\qss 
$v_{\dff j}\off =\off v_{\dff j\dff -\dff 1}\qff +\pff \varepsilon_{\dff p}$\qss
for some\dss $p$\nnsp,\oss
and\dss hence\dss there are no vertices\dss $w$\dss such\dss that\qss
$v_{\dff j\dff -\dff 1}\qff <\qff w\qff <\qff v_{\dff j}$\nsp.\oss
It\dss follows\dss that\trs
if\qss replacing\dss the vertex\dss $v_{\dff 0}$\dss 
in\dss $\sigma$\dss 
by\qss $w\off \neq\off v_{\dff 0}$\qss 
results in an $n$\dnsp-simplex\halfff,\oss
then either\qss
$w\qff <\qff v_{\dff 1}$\nsp,\dff\oss
or\pss
$v_{\dff n}\qff <\qff w$\nnsp.\oss
Let\dss $q$\dss be such\dss that\vspace{3.25pt}
\[
\quad
v_{\dff 1}\off =\off v_{\dff 0}\qff +\pff \varepsilon_{\dff q}
\qff.
\]

\vspace{-8.5pt}
Since different\sss coordinates are increasing at\sss different\sss steps,\oss
if\qss
$w\qff <\qff v_{\dff 1}$\nsp,\dff\oss
then\qss
$v_{\dff 1}\off =\off w\qff +\pff \varepsilon_{\dff q}$\qss
and\dss hence\qss
$w\off =\off v_{\dff 0}$\nsp,\oss
contrary\dss to\sss the assumption.\oss
By\dss the same reason,\oss
if\qss
$v_{\dff n}\qff <\qff w$\nnsp,\dff\oss
then\qss
$w\off =\off v_{\dff n}\qff +\pff \varepsilon_{\dff q}$\nsp.\oss
Therefore,\oss 
if\pss
$v_{\dff n}\qff +\pff \varepsilon_{\dff q}
\pff \not\in\pff K$\nnsp,\oss
then\dss $\sigma$\dss is\dss the\sss only\sss $n$\dnsp-simplex
containing\dss $\tau_{\dff 0}$\nsp.\oss
But\qss
$v_{\dff n}\qff +\pff \varepsilon_{\dff q}
\pff \not\in\pff K$\qss
only\qss if\qss the $q${\nnsp}th\dss coordinate of\dss $v_{\dff n}$\dss is equal\dss to\dss $k$\nnsp.\oss
Since\dss the $q${\nnsp}th\dss coordinate does not\dss change along\dss the sequence\qss
$v_{\dff 1}\fff,\pff
v_{\dff 2}\trf,\pff
\ldots\fff,\pff
v_{\dff n}$\nsp,\oss
in\dss this case\qss 
$\tau_{\dff 0}\off \subset\off
\mathcal{B}_{\dff q}\off \subset\off
\bd\fff K$\nnsp.\oss

Suppose now\dss that\pss
$v_{\dff n}\qff +\pff \varepsilon_{\dff q}
\qff \in\qff K$\nnsp.\oss
Then\qss
$\tau_{\dff 0}\off \not\subset\off \mathcal{B}_{\dff q}$\nsp.\oss
Also,\pss
$\tau_{\dff 0}\off \not\subset\off \mathcal{A}_{\dff q}$\qss
because\qss
$v_{\dff 1}\qff \in\qff \tau_{\dff 0}$\qss
and\qss
$v_{\dff 1}\off =\off v_{\dff 0}\qff +\pff \varepsilon_{\dff r}$\nnsp,\oss
and\dss $\tau_{\dff 0}$\dss
cannot\dss be contained\dss in\dss $\mathcal{A}_{\dff p}$\dss or\dss $\mathcal{B}_{\dff p}$\dss
for\qss $p\off \neq\off q$\qss
because only\dss the $q${\nnsp}th\dss coordinate does not\dss change\sss
along\dss the sequence\qss
$v_{\dff 1}\fff,\pff
v_{\dff 2}\trf,\pff
\ldots\fff,\pff
v_{\dff n}$\nsp.\oss
At\dss the same\sss time\vspace{3.5pt}
\[
\quad
\{\trf
v_{\dff 1}\fff,\pff
v_{\dff 2}\trf,\pff
\ldots\fff,\pff
v_{\dff n}\fff,\pff
v_{\dff n}\qff +\pff \varepsilon_{\dff q}
\pff\}
\qff
\]

\vspace{-8.5pt}
is\dss the second $n$\dnsp-simplex\dss having\dss $\tau_{\dff 0}$\dss
as a face and\dss there are no other such simplices.\oss
Therefore in\dss this case\dss $\tau_{\dff 0}$\dss
is not\sss contained\dss in\dss $\bd\fff K$\dss
and\dss is\dss a face of\dss exactly\dss two $n$\dnsp-simplices.\oss
This completes\sss the proof\qss in\dss the case of\qss
$i\off =\off 0$\nnsp.\oss
The case of\qss $i\off =\off n$\qss is completely\sss similar\halfff.\oss  \eproof

\prooftitle{Kuhn's\qss proof\qss of\pss Kuhn's\qss strong\dss lemma}
Let\dss $\triangle$\dss be\sss the abstract\sss simplicial\sss complex\sss having\qss
$1\dff,\off
2\dff,\off
\ldots\dff,\off
n\qff +\qff 1$\qss
as its vertices and all\sss subsets of\qss
$\{\qff
1\dff,\off
2\dff,\off
\ldots\dff,\off
n\qff +\qff 1
\qff\}$\qss
as its simplices.\oss
The reduced\dss labeling\sss $r$\sss can\dss be considered as a simplicial\dss map\qss
$K\qff \ttoo \triangle$\nnsp.\oss
Let\dss $e$\dss be\sss the number of\dss $n$\dnsp-simplices\dss $\sigma$\dss of\trs $K$\dss
such\dss that\qss
$r\dff(\dff \sigma\dff)\off =\off \triangle$\nnsp.\oss
We need\dss to prove\sss that\dss $e$\dss is odd.\oss

Let\dss 
$\triangle_{\dff 1}$\dss
be\sss the $(\fff n\dff -\dff 1\fff)$\dnsp-face\qss
$\{\qff
2\dff,\off
\ldots\dff,\off
n\qff +\qff 1
\qff\}$\qss
of\qss $\triangle$\nnsp,\oss
and\dss let\dss $h$\dss be\dss the number of\dss 
$(\fff n\qff -\qff 1\fff)$\dnsp-simplices\dss $\tau$\dss of\trs $K$\dss
such\dss that\qss $\tau\qff \subset\qff \mathcal{B}_{\dff 1}$\qss and\qss
$r\dff(\dff \tau\dff)\off =\off \triangle_{\dff 1}$\nsp.\oss
We claim\dss that\qss
$e\qff \equiv\qff h$\qss modulo $2$\nnsp.\oss
Clearly\halfff,\oss
$r\dff(\dff v\trf)\qff \leq\qff i$\qss
if\qss $v\qff \in\qff \mathcal{A}_{\dff i}$\nsp,\oss
and\qss
$r\dff(\dff v\trf)\qff \neq\qff i$\qss
if\qss $v\qff \in\qff \mathcal{B}_{\dff i}$\nsp.\oss
It\dss follows\dss that\qss
$r\dff(\dff \tau\dff)\off \neq\off \triangle_{\dff 1}$\qss
if\dss $\tau$\dss is\dss an
$(\fff n\qff -\qff 1\fff)$\dnsp-simplex of\trs $K$\dss
contained\dss in\dss $\bd\fff K$\dss
but\dss not\dss in\dss $\mathcal{B}_{\dff 1}$\nsp.\oss
Using\qss Lemma\qss \ref{cube-non-branching},\oss
we can\dss follow\dss now\sss either\dss the combinatorial\sss or\dss
a cochain-based\dss proof\dss of\qss Sperner's\dss lemma\sss
and conclude\sss that\qss
$e\qff \equiv\qff h$\qss modulo $2$\nnsp.\oss

By\dss forgetting\dss the first\sss coordinate we can\dss
identify\dss $\mathcal{B}_{\dff 1}$\dss with\dss the cube defined\dss
in\dss the same way as\dss $K$\nnsp,\oss but\dss with\qss $n\qff -\qff 1$\qss
in\dss the role of\dss $n$\nnsp.\oss
Let\dss us consider\dss the map\qss
$\mathcal{B}_{\dff 1}\qff \ttoo\qff \{\qff 0\dff,\pff 1\qff\}^{\dff n\dff -\dff 1}$\qss
having\pss
$l_{\dff 2}\dff,\pff \ldots\dff,\pff l_{\dff n}$\qss
as its components.\oss
Then\qss
$r'\dff \colon\dff
v\qff \longmapsto\qff r\dff(\dff v\trf)\qff -\qff 1$\qss
is\dss the corresponding\dss reduced\dss labeling\halfff.\oss
Obviously,\pss $h$\dss is equal\dss to\sss
the number of\dss 
$(\fff n\qff -\qff 1\fff)$\dnsp-simplices\dss 
$\tau$\dss of\trs $\mathcal{B}_{\dff 1}$\dss such\dss that\qss
$r'\dff(\dff \tau\dff)\off =\off 
\{\qff
1\dff,\pff
2\dff,\off
\ldots\dff,\off
n
\qff\}$\nsp.\oss
By\dss using an\dss induction\dss by\dss $n$\dss we can assume\sss that\dss $h$\dss is odd\qss
(the case of\qss $n\off =\off 0$\qss being\dss trivial\dff).\oss
By\dss the previous paragraph\dss this implies\sss that\dss $e$\dss is odd.\oss  \eproof

\myuppar{Kuhn's\qss strong\dss lemma,\oss Hurewicz's\trs theorems,\oss
and\trs Lebesgue\trs tilings.}
Let\qss
$e_{\dff i}
\off =\off
r^{\dff -\dff 1}\dff(\trf i\trf)$\nnsp,\oss
where\qss
$i\off =\off 1\fff,\pff 2\fff,\pff \ldots\fff,\pff n\qff +\qff 1$\nnsp.\oss
Clearly\halfff,\oss the sets\qss
$e_{\dff 1}\fff,\pff e_{\dff 2}\fff,\pff \ldots\fff,\pff e_{\dff n\dff +\dff 1}$\qss
are pairwise disjoint\sss and\dss their\dss union\sss is\sss equal\dss to\dss $K$\nnsp.\oss
The assumptions of\qss Kuhn's\qss lemma\dss imply\dss that\qss\vspace{2.25pt}\vspace{-2.25pt}
\[
\quad
\mathcal{A}_{\dff i}
\off \subset\off
e_{\dff 1}\off \cup\off \ldots\off \cup\off e_{\dff i}
\]

\vspace{-9.75pt}
and\dss $e_{\dff i}$\dss
is\dss disjoint\qss from\dss $\mathcal{B}_{\dff i}$\dss
for every\qss $i\qff \leq\qff n$\nnsp.\oss
If\qss $i\qff >\qff j\qff \geq\qff 1$\nnsp,\oss
then\dss $e_{\dff i}$\dss is disjoint\dss from\dss $\mathcal{A}_{\fff j}$\qss
because\qss
$\mathcal{A}_{\fff j}\off \subset\off
e_{\dff 1}\qff \cup\qff \ldots\qff \cup\qff e_{\fff j}$\qss
and\dss the sets\qss 
$e_{\dff 1}\dff,\pff e_{\dff 2}\dff,\pff \ldots\fff,\pff e_{\dff n\dff +\dff 1}$\qss
are\dss pairwise\dss disjoint\halfff.\qff\oss
Therefore\qss
$e_{\dff 1}\fff,\pff e_{\dff 2}\dff,\pff \ldots\fff,\pff e_{\dff n\dff +\dff 1}$\qss
satisfy\dss the assumptions of\qss Theorem\qss \ref{wh-lebesgue-tilings}\qss 
and\dss this\sss theorem\dss implies\qss Kuhn's\qss strong\dss lemma.\oss
In\sss other\dss words,\oss Kuhn's\qss strong\sss lemma\sss is\qss
an\dss immediate\sss corollary\sss of\qss Hurewicz's\dss theorems\sss as applied\dss
to\dss Lebesgue\dss tilings.\oss
Apparently\halfff,\oss this connection\dss wasn't\dss noticed\dss before.\oss

\myuppar{Kuhn's\qss lemmas\sss and\qss Lebesgue\dss results.}
In order\dss to compare\dss Kuhn's\dss results with\dss Lebesgue\dss ones,\pss
one needs\sss to express\sss them\dss in\dss the same language.\oss
This\dss is\dss done\sss in\dss Appendix\qss \ref{two-cubes},\oss
where\dss Lebesgue\dss results are\sss
transported\dss from\dss $Q$\dss to\dss $K$\nnsp.\oss
It\dss seems\sss that\dss the only\dss essential\sss conclusion 
of\qss Kuhn's\qss lemma\sss and\qss 
Theorem\qss \ref{separation-products}\qss 
is\dss the existence of\dss an $n$\dnsp-cube intersecting\dss
$c_{\dff i}$\dss and\qss $\overline{c}_{\dff i}$\qss
for every\qss $i\qff \in\qff I$\nnsp.\oss
This part\dss is nothing\dss else\sss but\qss Theorem\qss \ref{cubes-partitions}\qss
transported\dss from\dss $Q$\dss to\dss $K$\nnsp.\oss
See\qss Theorem\qss \ref{separation-products-weak}.\oss
Similarly\halfff,\oss
it\dss seems\dss that\dss the only\sss essential\sss conclusion 
of\qss Kuhn's\qss strong\dss lemma\dss is\dss the
existence of\trs an $n$\dnsp-cube\dss $\sigma$\dss of\qss $K$\dss
such\dss that\qss 
$r\dff(\dff \sigma\dff)
\off =\off
\{\qff
1\dff,\off
2\dff,\off
\ldots\dff,\off
n\qff +\qff 1
\qff\}$\nnsp.\oss
This\dss is\dss equivalent\dss to\sss the existence of\dss an $n$\dnsp-cube\dss $\sigma$\dss
intersecting\sss every\sss set\qss
$e_{\dff i}
\off =\off
r^{\dff -\dff 1}\dff(\trf i\trf)$\nnsp.\oss
By\dss the previous subsection,\oss the sets\qss
$e_{\dff 1}\fff,\pff e_{\dff 2}\dff,\pff \ldots\fff,\pff e_{\dff n\dff +\dff 1}$\qss
satisfy\dss the assumptions of\qss Theorem\qss 
\ref{discrete-coverings}\qss
and\dss hence\sss this\sss theorem\dss implies\sss the existence of\dss such an $n$\dnsp-cube.\oss
On\dss the other hand\qss Theorem\qss 
\ref{discrete-coverings}\qss
is\dss the\sss result\sss of\trs transporting\qss Theorem\qss \ref{e-coverings}\qss
from\dss $Q$\dss to\dss $K$\nnsp.\oss 
We see\sss that\dss the essential\dss part\sss of\qss Kuhn's\qss lemmas\sss
was\sss contained\sss already\dss in\dss Lebesgue\dss results.

\myuppar{In search of\trs the lost\dss cubes\fff:\qss II.}
The above proof\dss is\qss Kuhn's\qss own\dss proof\dss up\sss to notations and\dss
terminology.\oss
With\dss the simplicial\sss complex\dss $K$\dss at\dss hand,\oss
it\dss is a fairly\sss straightforward adaptation of\trs standard\dss proofs\sss
of\trs Sperner's\dss lemma.\oss
Given\qss Kuhn's\qss goal\dss to replace simplices by\sss cubes,\oss
it\dss is\dss hardly\sss satisfactory,\oss
and\qss Kuhn\qss \cite{ku}\qss
posed,\oss among other\halfff,\oss the following question.\vspace{-9.25pt}

\begin{quoting}
Is\dss it\dss possible\sss to prove\sss the\dss 
Cubical\dss Sperner\dss Lemma\dss
without\dss resorting\dss to a simplicial\sss decomposition
of\trs the $n$\dnsp-cube\fff?\oss
\end{quoting}

\vspace{-9.5pt}
Since\dss Freudenthal's\dss simplicial\sss decomposition of\dss the cube is\dss
included\sss already\dss in\dss the statement\sss of\qss Kuhn's\qss lemma\dss 
in\dss the form of\trs
the condition\qss (\ref{kuhn-simplex}),\oss
this question\dss hardly can\dss have a positive answer\halfff,\oss
if\qss understood\dss literally.\oss 
But\sss one can\dss hide\sss the simplices of\dss this decomposition fairly\dss well.\oss
A positive answer\dss was suggested\dss by\qss
L.A.\trs Wolsey\qss \cite{w}.\oss
In\qss \cite{w}\qss the simplices are hidden\dss in\dss the recursive definition of\pss
\emph{completely\dss labeled}\dss $m$\dnsp-cubes.\oss 
See\qss \cite{w},\oss Definition\qss 1.\oss

Theorem\qss \ref{separation-products}\qss 
also may\sss qualify\sss as a positive answer\halfff.\oss
In\dss its\dss proof\trs the simplices are hidden\dss in\dss the formula\sss for\dss the
multiplication\sss of\dss several\sss $1$\dnsp-cochains.\oss
But\dss this formula\sss is\sss just\sss a special\sss case of\dss the general\dss formula\qss
(\ref{product-multiple})\qss for\dss the multiplication\sss of\dss several\sss cochains,\oss
which\dss immediately\dss follows\dss from\dss the definition of\dss 
products\sss of\qss two cubical\sss cochains.\oss
Neither\qss (\ref{product-multiple}),\oss nor\dss this definition,\pss
adapted\dss from\dss the celebrated\dss work of\qss J.-P.\dss Serre\qss \cite{se},\pss
involve any\dss triangulations.\oss

Instead of\trs hiding simplices one can weaken\sss a\sss little\sss the conclusions of\qss
Kuhn's\qss lemmas.\oss
Namely,\oss one can\dss replace\sss 
the sets\sss of\dss terms of\dss pivot\sss sequences\qss
by $n$\dnsp-cubes.\oss 
As we saw,\oss 
these weakened\dss results 
were essentially\dss proved\dss by\qss
Lebesgue\dss 40\dss years earlier\dss than\dss Kuhn's\dss asked\dss his question.\oss

\myuppar{Missed\sss opportunities.}
At\dss the end of\dss his paper\trs Kuhn\qss wrote\sss that\vspace{-9.25pt}

\begin{quoting}
\ldots\oss cubical\sss complexes have\sss the crucial\sss advantage of\trs their
ease of\dss description\dss in\dss binary\dss notation\fff;\oss
for digital\sss computation,\oss therefore,\oss they are a natural\sss object\dss
to study\halfff.\oss
This advantage,\oss however\halfff,\oss is offset\dss by\dss the absence of\dss
an appropriate homology\dss theory\dss with a boundary\sss operator suited\dss
to our purposes.
\end{quoting}

\vspace{-9.25pt}
This was written about\dss 10\dss years after\dss the cubical\sss homology
and cohomology\dss theories appeared\dss in\dss the celebrated\trs Serre's\trs paper\qss
\cite{se}.\oss
At\dss the\sss time\qss H.W.\qss Kuhn\dss was an associate professor
of\dss mathematics and economics
at\dss Princeton and a colleague of\pss J.W.\qss Milnor\dss and\qss J.C.\qss Moore\qss
(see\qss \cite{go}\fff).\oss
Serre's\dss theory\dss and other\dss
methods of\qss French school\dss was\qss ``a\dss Princeton\dss speciality'',\oss
and\qss J.C.\qss Moore\dss was a leading\sss proponent\sss of\trs these methods\qss
(see\qss \cite{bgk}\fff).\oss
A.W.\qss Tucker\halfff,\oss a pupil\sss of\qss S.\qss Lefschetz,\oss
was\dss the chairman of\trs Princeton\dss Mathematics\dss Department\dss 
from\dss 1952\dss to\dss 1963.\oss

In\dss 1960\qss H.W.\qss Kuhn\qss and\qss A.W.\qss Tucker\dss
were experts in\dss game\sss theory\sss and\dss mathematical\dss pro\-gram\-ming\halfff.\oss
Still,\pss even\dss the\sss title of\qss Kuhn's\qss paper\qss \cite{ku}\qss
indicates\sss that\dss it\dss deals with\dss topology\halfff.\oss
Did\dss he ever\dss discuss\dss it\dss with a\sss topologist\dff?\oss
The situation\dss resembles\sss a story\dss related\dss by\trs F.\dss Dyson\dss in\dss his
famous and\dss influential\sss essay\qss \cite{dy}.\oss
Being a physicist\halfff,\oss Dyson\dss did\sss not\dss discuss\dss his ideas in\dss
number\dss theory\dss with\qss I.\dss MacDonald,\oss a mathematician,\oss
who happened\dss to work at\dss the same\sss time and\dss place on almost\dss
the same questions.\oss

\mysection{Ky\qss Fan's\qss cubical\qss Sperner\qss lemma}{kyfan}

\myuppar{Adjacency-preserving\dss maps.}
Let\dss us\dss keep\sss the notations of\qss Section\qss \ref{kuhn}.\oss
Two points\qss $v\fff,\pff w$\qss of\trs $K$\dss or\dss $\mathbb{K}$\dss  
are said\dss to be\qss
\emph{adjacent}\pss if\qss $\{\dff v\fff,\qff w\trf\}$\qss is\dss a
$1$\dnsp-cube.\oss
A map\qss
$\varphi\dff \colon\dff
K\qff \ttoo\qff \mathbb{K}$\qss
is\dss said\dss to be\qss \emph{adjacency-preserving}\pss if\qss it\dss takes
adjacent\sss vertices of\trs $K$\dss to adjacent\sss or equal\sss vertices of\trs $\mathbb{K}$\nnsp.\oss
Ky\dss Fan\qss \cite{kf}\qss called\sss such\dss maps\qss
\emph{cubical\dss vertex\dss map},\oss
but\dss in\dss the opinion of\trs the present\sss author\dss these map do
not\dss really\sss deserve\sss such\dss name.\oss
The notion of\dss adjacency-preserving\dss maps\dss
is\qss \emph{not}\qss
a\sss genuine cubical\sss analogue of\trs the notion of\dss simplicial\dss maps.\oss
For example,\oss the image of\dss a cube of\trs $K$\dss under an adjacency-preserving\dss map
is\dss not\sss a cube\sss of\trs $\mathbb{K}$\dss in\dss general.\oss

\myuppar{Symmetries.} 
A\trs \emph{reflection}\pss of\trs the unit\sss discrete cube\dss $\mathbb{K}$\dss is\dss a\dss bijection\qss
$\mathbb{K}\qff \ttoo\qff \mathbb{K}$\qss
acting\dss independently\sss on each of\trs $n$\dss coordinates by\dss maps
equal\dss either\dss to\qss
$a\qff \longmapsto\qff 1\qff -\qff a$\pss or\dss to\sss the identity\pss
$a\qff \longmapsto\qff a$\nnsp.\oss
A\trs \emph{symmetry}\pss of\trs $\mathbb{K}$\dss is\dss a\dss
composition of\dss a\sss reflection and\sss a permutation of\dss coordinates.\oss 
A\qss \emph{reflection}\qss or\sss a\qss \emph{symmetry}\qss of\dss 
an $n$\dnsp-cube\dss $\sigma$\dss of\trs $K$\dss
is\dss a\dss bijection\qss
$\sigma\qff \ttoo\qff \sigma$\qss
such\dss that\dss identifying\dss $\sigma$\dss with\dss $\mathbb{K}$\dss
by\sss a\sss translation\dss turns it\sss into a\sss reflection\sss or\sss 
a symmetry\sss of\trs $\mathbb{K}$\dss respectively\halfff.\oss

\myuppar{Opposite vertices and\dss free\dss pivot\dss sequences.}
Let\dss $\tau$\dss be an $m$\dnsp-cube of\trs $K$\nnsp.\oss
Two vertices of\dss $\tau$\dss are said\dss to be\qss \emph{opposite}\pss
if\trs they differ\dss in $m$ coordinates.\oss
A\qss \emph{free\dss pivot\dss sequence}\qss of\dss vertices of\dss $\tau$\dss
is\dss a\sss sequence\qss
$v_{\dff 0}\fff,\pff
v_{\dff 1}\trf,\pff
\ldots\fff,\pff
v_{\dff m}
\qff \in\qff
\tau$\qss
such\dss that\dss only\sss one coordinate changes when one passes from
$v_{\dff i\dff -\dff 1}$\sss
to $v_{\dff i}$ 
and\dss different\sss coordinates change at\sss different\sss steps.\oss
Clearly\halfff,\qss $v_{\dff 0}$ and $v_{\dff m}$
are opposite vertices of\dss $\tau$\nnsp,\qss
and\dss two
opposite vertices can\dss be connected\dss by a\sss free pivot\sss sequence.

Suppose\sss now\dss that\qss $m\off =\off n$\nnsp.\oss
In\dss this case,\oss
if\dss $s\dff(\dff i\trf)$\dss is\dss the number of\trs the coordinate
changing\dss from\dss
$v_{\dff i\dff -\dff 1}$\dss
to\dss $v_{\dff i}$\nsp,\oss
then\dss $s$\dss is\dss a\sss permutation of\qss
$1\fff,\pff 2\fff,\pff \ldots\fff,\pff n$\nnsp,\oss
and\qss 
$v_{\dff 0}\fff,\pff
v_{\dff 1}\trf,\pff
\ldots\fff,\pff
v_{\dff n}$\qss
is\dss unique\-ly\dss determined\dss by\dss $v_{\dff 0}$\dss
and\dss $s$\nnsp.\oss
Also,\oss if\qss $m\off =\off n$\nnsp,\oss
then\sss a\sss sequence\qss
$v_{\dff 0}\fff,\pff
v_{\dff 1}\trf,\pff
\ldots\fff,\pff
v_{\dff n}
\qff \in\qff
\tau$\qss
is\dss a\sss free\dss pivot\sss sequence\dss if\trs and\dss only\trs if\trs
some\sss reflection\sss of\dss $\tau$ turns\sss it\dss into a pivot\sss sequence.\oss

\myuppar{Adjacency-preserving\dss maps of\dss cubes.}
Let\dss us\dss fix\sss an $n$\dnsp-cube\dss $\sigma$\dss of\trs $K$\nnsp.\oss
Let\qss
$\lambda^0\trf \sigma\off =\off \{\trf r\qff\}$\qss
and\qss
$\lambda^1\trf \sigma\off =\off \{\trf p\qff\}$\qss
be\sss the root\dss and\trs the peak\dss of\dss $\sigma$\dss respectively\halfff.\oss
Let\dss us\dss fix also an adjacency-pre\-serv\-ing\dss map\qss
$\varphi\dff \colon\dff
\sigma\qff \ttoo\qff \mathbb{K}$\nnsp.\oss
We are interested\dss when\dss $\varphi$\dss is\dss a\dss bijection.\oss

\mypar{Lemma.}{opposite-pivot}
\emph{Suppose\dss that\qss
$\varphi\dff(\dff r\trf)
\off =\off 
0$\qss
and\qss
$\varphi\dff(\dff p\trf)
\off =\off 
\mathbb{1}$\nnsp.\oss
Then\dss $\varphi$\dss takes every\dss pivot\dss 
sequence\sss to a pivot\sss sequence.\oss}\vspace{-0.1pt}

\proof
After\dss identifying\dss $\sigma$\dss with\dss $\mathbb{K}$\dss
by\sss a\sss translation\dss we may assume\sss that\qss
$\sigma\off =\off \mathbb{K}$\nnsp.\oss
Then\qss $r\off =\off 0$\qss and\qss $p\off =\off \mathbb{1}$\nnsp.\oss
Every\dss pivot\dss sequence in\dss $\mathbb{K}$\dss starts with\dss $0$\dss
and\dss ends with\dss $\mathbb{1}$\nnsp.\oss
Let\qss\vspace{0pt}
\[
\quad
0
\off =\off
v_{\dff 0}\dff,\off\off
v_{\dff 1}\trf,\off\off
\ldots\dff,\off\off
v_{\dff n}
\off =\off
\mathbb{1}
\]

\vspace{-12pt}
be\dss a\sss pivot\sss sequence.\oss
Then\qss
$\{\trf v_{\dff i\dff -\dff 1}\dff,\pff v_{\dff i}\trf\}$\qss
is\dss a $1$\dnsp-cube for every\qss $i\qff \geq\qff 1$\nnsp.\oss
Since\dss $\varphi$\dss is\dss adjacency-preserving,\oss
this implies\sss that\dss no more\sss than one coordinate changes
when one passes\sss from \dss $\varphi\dff(\dff v_{\dff i\dff -\dff 1}\dff)$\dss
to\dss $\varphi\dff(\dff v_{\dff i}\dff)$\nnsp,\oss
and\dss the changing coordinate,\pss if\qss any\halfff,\pss changes by\dss $1$\nnsp.\oss
Since we have\sss to get\dss from\sss $0$\sss to\sss $\mathbb{1}$\sss
in $n$ such steps,\oss
some coordinate should\dss increase\sss by\dss $1$\dss at\sss each step
and\dss different\sss coordinates should\dss increase at\sss different\sss steps.\oss
Therefore\qss
$\varphi\dff(\dff v_{\dff 0}\dff)\dff,\off\off
\varphi\dff(\dff v_{\dff 1}\dff)\trf,\off\off
\ldots\dff,\off\off
\varphi\dff(\dff v_{\dff n}\dff)$\qss
is\dss a\sss pivot\sss sequence.\oss  \eproof

\mypar{Corollary\halfff.}{one-all}
\emph{If\qss $\varphi$\dss takes some pivot\dss sequence\sss to a pivot\sss sequence,\oss
then\dss $\varphi$\dss takes every\dss pivot\dss sequence\sss to a pivot\sss sequence.\oss}

\proof
Every\dss pivot\sss sequence in $\sigma$ starts\sss at\sss $r$ and ends\dss at\sss $p$\nnsp,\oss
and every\dss pivot\sss sequence in\dss $\mathbb{K}$\dss starts\sss at\sss $0$ 
and ends\dss at\sss $\mathbb{1}$\nnsp.\oss  \eproof

\mypar{Lemma.}{bijective-pivots}
\emph{If\qss $\varphi$\dss is\dss bijective,\oss 
then\dss $\varphi$\dss takes some\dss free\sss pivot\dss sequence\sss
to\sss a\dss pivot\sss sequence.\oss}

\proof
Let\qss
$v\off =\off \varphi^{\dff -\dff 1}\dff(\dff 0\dff)$\qss
and\qss
$w\off =\off \varphi^{\dff -\dff 1}\dff(\dff \mathbb{1}\trf)$\nnsp.\oss
Suppose\sss that\dss $w$\dss differs from\dss $v$\dss
in\dss $m$\dss coordinates.\oss
Then\dss there is a sequence\qss\vspace{2.5pt}
\[
\quad
v
\off =\off
v_{\dff 0}\dff,\off\off
v_{\dff 1}\trf,\off\off
\ldots\dff,\off\off
v_{\dff m}
\off =\off
w
\]

\vspace{-9.5pt}
of\dss points in $\sigma$ such\dss that\sss $v_{\dff i}$
differs from $v_{\dff i\dff -\dff 1}$in only one coordinate for every\qss $i\qff \leq\qff m$\nnsp.\oss
Equiv\-a\-lent\-ly\halfff,\pss
$\{\trf v_{\dff i\dff -\dff 1}\dff,\pff v_{\dff i} \qff\}$\dss
is\dss a $1$\dnsp-cube for every\qss $i\qff \leq\qff m$\nnsp.\oss
Since\dss $\varphi$\dss is\dss adjacency-preserving,\oss
this implies\sss that\dss $\varphi\dff(\dff v_{\dff m}\dff)$\dss
differs from\dss $\varphi\dff(\dff v_{\dff 0}\dff)$\dss
in\qss $\leq\qff m$\qss coordinates.\oss 
But\sss actually\qss
$\varphi\dff(\dff v_{\dff m}\dff)\off =\off \mathbb{1}$\qss
and\qss
$\varphi\dff(\dff v_{\dff 0}\dff)\off =\off 0$\qss
differ\dss in $n$ coordinates.\oss
Therefore\qss $m\off =\off n$\nnsp,\oss
i.e.\dss $w\off =\off v_{\dff n}$\dss
differs\dss from\dss  $v\off =\off v_{\dff 0}$\dss
in $n$ coordinates.\oss
In\dss turn,\pss this implies\sss that\sss different\sss coordinates
are changed\sss at\sss different\sss steps\sss and\dss hence\dss 
$v_{\dff 0}\fff,\pff
v_{\dff 1}\trf,\pff
\ldots\fff,\pff
v_{\dff n}$\qss
is\dss a\sss free pivot\sss sequence.\oss
By\dss the same reason\qss\vspace{2.5pt}
\[
\quad
\varphi\dff(\dff v_{\dff 0}\dff)\dff,\off\off
\varphi\dff(\dff v_{\dff 1}\dff)\dff,\off\off
\ldots\dff,\off\off
\varphi\dff(\dff v_{\dff n}\dff)
\]

\vspace{-9.5pt}
is\dss a\sss free\dss pivot\sss sequence.\oss
Since all\sss $n$\sss coordinates are increasing\dss between\qss
$\varphi\dff(\dff v_{\dff 0}\dff)\off =\off 0$\qss
and\qss
$\varphi\dff(\dff v_{\dff n}\dff)\off =\off \mathbb{1}$\nnsp,\oss
it\dss is\dss actually\sss a\sss pivot\sss sequence.\oss  \eproof

\mypar{Lemma.}{bijective-free-pivots}
\emph{If\qss $\varphi$\dss is\dss bijective,\oss 
then\dss $\varphi$\dss takes some\sss pivot\dss sequence\sss
to\sss a\sss free\dss pivot\sss sequence.}

\proof
As above,\oss we may assume\sss that\qss
$\sigma\off =\off \mathbb{K}$\nnsp.\oss
Since\dss $\varphi$\dss is\dss bijective and adjacency-pre\-serv\-ing\halfff,\oss
$\varphi$\dss induces a\sss bijective self-map of\trs the set\sss of\dss
$1$\dnsp-cubes of\trs $\mathbb{K}$\nnsp.\oss
It\dss follows\sss that\dss $\varphi^{\dff -\dff 1}$\dss is\dss
a\sss bijective adjacency-preserving\dss map.\oss
It\dss remains\dss to apply\qss Lemma\qss \ref{bijective-pivots}\qss
to\dss $\varphi^{\dff -\dff 1}$\dnsp.\oss  \eproof

\mypar{Corollary\halfff.}{pivot-correspondence}
\emph{Let\qss $v\qff \in\qff \sigma$\nnsp.\oss
If\qss $\varphi$\dss is\dss bijective,\oss
then\dss $\varphi$\dss induces a one-to-one
correspondence between\dss the\dss free pivot\dss sequences starting\dss at\sss $v$ 
and\dss the\dss free\sss pivot\sss
sequences starting\dss at\dss $\varphi\dff(\dff v\dff)$\nnsp.\oss}

\proof
After composing\dss $\varphi$\dss with a reflection\dss we may assume\sss that\dss $v$\dss
is\dss the $<$\dnsp-minimal\dss point\sss of\dss $\sigma$\nnsp.\oss
Then\dss the free pivot\sss sequences starting\dss at $v$ are simply\dss pivot\sss sequences.\oss 
Lemma\qss \ref{bijective-free-pivots}\qss implies\sss that\dss there\dss is\dss a\sss reflection $\rho$
such\dss that\dss $\rho\dff \circ\dff \varphi$\dss
takes some pivot\dss sequence\sss to a pivot\sss sequence.\oss
Now\dss Corollary\qss \ref{one-all}\qss implies\sss that\qss
$\rho\dff \circ\dff \varphi$\dss
takes every\dss pivot\dss sequence\sss to a pivot\sss sequence,\oss
and\dss hence $\varphi$ takes every\dss pivot\sss sequence\sss
to\sss the image of\dss a\sss pivot\sss sequence under\dss
$\rho^{\dff -\dff 1}\pff =\off \rho$\nnsp.\oss
Such an image\dss is a free\sss pivot\sss sequence starting\sss at\dss
$\varphi\dff(\dff v\dff)$\nnsp.\oss
This proves\sss the first\sss statement\sss of\trs the\sss lemma.\oss
The second\dss follows\sss from\dss the first\sss and\dss the fact\dss
that\dss $\varphi$\dss is\dss a\sss bijection.\oss  \eproof

\myuppar{Cochains associated\dss with $\varphi$\nnsp.}
Let\qss 
$h_{\dff i}\qff \colon\qff
\sigma\off \ttoo\off \{\qff 0\dff,\pff 1\qff\}$\nnsp,\oss
where\qss
$i\off =\off 1\fff,\pff 2\fff,\pff \ldots\fff,\pff n$\nnsp,\oss
be\sss the components of\dss $\varphi$\dss as\sss a\sss map\sss to\qss
$\mathbb{K}\off =\off \{\qff 0\dff,\pff 1\qff\}^{\fff n}$\nnsp.\oss
Let\dss us\dss identify\qss $\{\qff 0\dff,\pff 1\qff\}$\qss with\qss $\ftwo$\qss
and consider\dss the maps\dss $h_{\dff i}$\dss also as $0$\dnsp-cochains of\trs
$\sigma$\nnsp.\oss
Let\qss\vspace{2.625pt}
\[
\quad
f_{\dff i}
\off =\off
\partial^* h_{\dff i}
\qff.
\]

\vspace{-9.375pt}
If\qss
$\tau\off =\off \{\qff v\fff,\qff w \qff\}
\off \subset\off
\sigma$\qss is\dss a $1$\dnsp-cube,\oss
then\qss
$f_{\dff i}\trf(\dff \tau\dff)\off =\off 1$\qss
if\qss and\dss only\trs if\pss 
$h_{\dff i}\dff (\dff v\dff)\off \neq\off h_{\dff i}\dff (\dff w\dff)$\nnsp.\oss

Two points of\trs $\mathbb{K}$\dss are adjacent\qss if\qss and\dss only\trs if\qss
they differ\dss in\dss no more\sss than one coordinate.\oss
Therefore,\oss the adjacency-preserving\dss property\sss of\dss $\varphi$ means\sss that\dss
for every $1$\dnsp-cube\qss $\tau\qff \subset\qff \sigma$\qss
there\sss is\dss no more\sss than one $i$ such\dss that\qss
$f_{\dff i}\trf(\dff \tau\dff)\pff =\off 1$\nnsp.\oss
Equivalently\halfff,\pss all\sss components\dss $h_{\dff i}$\dss of\qss $\varphi$\nnsp,\oss
except\halfff,\oss perhaps,\oss one of\trs them,\oss
are constant\dss on\dss $\tau$\nnsp.\oss

\mypar{Theorem.}{bijection-to-product}
\emph{If\qss $\varphi$\dss is\dss bijective,\oss 
then\oss
$f_{\dff 1}\dff \cdot\dff f_{\dff 2}\trf \cdot\qff \ldots\dff \cdot\qff 
f_{\dff n}\qff(\dff \sigma\dff)
\off =\off 
1$\nnsp.\oss}

\proof
Let\qss
$v_{\dff 0}\fff,\pff
v_{\dff 1}\trf,\pff
\ldots\fff,\pff
v_{\dff n}
\qff \in\qff
\sigma$\qss
be a pivot\sss sequence.\oss
It\dss defines a sequence\qss
$\tau_{\dff 1}\fff,\pff
\tau_{\dff 2}\trf,\pff
\ldots\fff,\pff
\tau_{\dff n}$\qss
of\dss $1$\dnsp-cubes,\oss
where\qss
$\tau_{\dff i}
\off =\off
\{\qff v_{\fff i\dff -\dff 1}\dff,\pff v_{\fff i} \pff\}$\qss
for each\qss $i\qff \geq\qff 1$\nnsp.\oss
Suppose\sss that\vspace{4.5pt}
\begin{equation}
\label{product-is-1}
\quad
f_{\dff 1}\dff\left(\qff \tau_{\dff 1} \qff\right)
\trf \cdot\qff 
f_{\dff 2}\dff\left(\qff \tau_{\dff 2} \qff\right)
\trf \cdot\qff
\ldots
\trf \cdot\qff
f_{\dff n}\dff\left(\qff \tau_{\dff n} \qff\right)
\off =\off
1
\qff.
\end{equation}

\vspace{-7.5pt}
Then\qss 
$f_{\dff i}\dff\left(\qff \tau_{\dff i} \qff\right)\off =\off 1$\qss
for every\dss $i$\nnsp.\oss
This means\dss that\dss
$\varphi\qff(\dff v_{\dff i}\trf)$\dss
differs\dss from\dss
$\varphi\qff(\dff v_{\dff i\dff -\dff 1}\trf)$\dss
in $i${\dnsp}th\dss coordinate.\oss
Since\dss $\varphi$\dss is\dss adjacency-preserving\halfff,\oss
$\{\qff \varphi\qff(\dff v_{\dff i\dff -\dff 1}\trf)\fff,\off
\varphi\qff(\dff v_{\dff i}\trf) \trf\}$\qss
is\dss a $1$\dnsp-cube
and\dss hence\dss
$\varphi\qff(\dff v_{\dff i}\trf)$\dss
differs\dss from\dss
$\varphi\qff(\dff v_{\dff i\dff -\dff 1}\trf)$\dss
\emph{only}\pss in $i${\dnsp}th\dss coordinate.\oss
It\dss follows\dss that\dss the sequence\vspace{4.5pt}
\[
\quad
\varphi\qff(\dff v_{\dff 0}\trf)\fff,\off\off
\varphi\qff(\dff v_{\dff 1}\trf)\fff,\off\off
\ldots\fff,\off\off
\varphi\qff(\dff v_{\dff n}\trf)
\]

\vspace{-7.5pt}
is\dss an\dss free\dss pivot\sss sequence and\dss the corresponding\dss
permutation\dss $s$\dss is\dss the identity\dss permutation.\oss
Therefore\sss this sequence\dss
is\dss uniquely\sss determined\dss by\qss
$\varphi\qff(\dff v_{\dff 0}\trf)$\qss
and\dss the property\qss (\ref{product-is-1}).\oss
In view of\qss Corollary\qss \ref{pivot-correspondence},\oss
this implies\sss that\dss the pivot\sss sequence\qss
$v_{\dff 0}\fff,\pff
v_{\dff 1}\trf,\pff
\ldots\fff,\pff
v_{\dff n}
\qff \in\qff
\sigma$\qss
is\dss uniquely\dss determined\dss by\dss
the property\qss (\ref{product-is-1}),\oss i.e.\qss
there\dss is\dss exactly\sss one\sss pivot\sss sequence\qss
$v_{\dff 0}\fff,\pff
v_{\dff 1}\trf,\pff
\ldots\fff,\pff
v_{\dff n}
\qff \in\qff
\sigma$\qss
such\dss that\qss (\ref{product-is-1})\qss holds.\oss

As we saw\sss in\dss
Section\qss \ref{cubical-cochains},\oss
the value\qss
$f_{\dff 1}\dff \cdot\dff f_{\dff 2}\trf \cdot\qff \ldots\dff \cdot\qff 
f_{\dff n}\qff(\dff \sigma\dff)$\qss
is\dss equal\dss to\sss the sum of\dss all\dss products\qss
$f_{\dff 1}\dff\left(\qff \tau_{\dff 1} \qff\right)
\trf \cdot\qff 
f_{\dff 2}\dff\left(\qff \tau_{\dff 2} \qff\right)
\trf \cdot\qff
\ldots
\trf \cdot\qff
f_{\dff n}\dff\left(\qff \tau_{\dff n} \qff\right)$\qss
corresponding\dss to\sss pivot\sss sequences.\oss
Each such\dss product\dss is\dss equal\sss either\dss to\dss $1$\dss or\dss to\dss $0$\nnsp.\oss
By\dss the previous paragraph,\oss 
this product\dss is\dss equal\dss to $1$ for exactly one pivot\sss sequence.\oss
It\dss follows\dss that\qss
$f_{\dff 1}\dff \cdot\dff f_{\dff 2}\trf \cdot\qff \ldots\dff \cdot\qff 
f_{\dff n}\qff(\dff \sigma\dff)
\off =\off 
1$\nnsp.\oss  \eproof

\mypar{Lemma.}{face-product}
\emph{Let\dss $\tau$\dss be an $(\fff n\dff -\dff 1\fff)$\dnsp-face 
of\trs $\sigma$\dss such\dss that}\qss\vspace{3pt}
\begin{equation}
\label{2-n-product}
\quad
f_{\dff 2}\dff \cdot\dff f_{\dff 3}\trf \cdot\qff \ldots\dff \cdot\qff f_{\dff n}\qff(\dff \tau\dff)
\off =\off
1
\qff.
\end{equation}

\vspace{-9pt}
\emph{Then\dss $\varphi$\dss maps\dss $\tau$\dss into\qss
$\mathbb{A}_{\dff 1}
\off =\off
0\dff \times\dff 
\{\qff 0\dff,\pff 1\qff\}^{\fff n\dff -\dff 1}$\qss
or\dff\oss
$\mathbb{B}_{\dff 1}
\off =\off
0\dff \times\dff 
\{\qff 0\dff,\pff 1\qff\}^{\fff n\dff -\dff 1}$\nnsp.\oss
Moreover\halfff,\oss}\vspace{3pt}
\[
\quad
\varphi\dff(\dff \tau\dff)
\off \subset\dff\off 
\mathbb{B}_{\dff 1}
\]

\vspace{-9pt}
\emph{if\qss and\dss only\qss if\oss
$h_{\dff 1}\dff(\qff \lambda^0\trf \tau\qff)
\off =\off 
1$\nnsp,\oss
where\qss $\lambda^0\trf \tau$\qss is\dss the root\dss of\pss $\tau$\nnsp.\oss}

\proof
Let\qss
$\lambda^0\trf \tau\off =\off \{\trf v\qff\}$\qss
and\qss
$\lambda^1\trf \tau\off =\off \{\trf w\qff\}$\qss 
be\sss the root\sss and\dss the peak of\dss
$\tau$\nnsp.\oss
The equality\qss (\ref{2-n-product})\qss implies\sss that\dss there\dss
is\dss a\sss free
pivot\sss sequence\qss 
$v
\off =\off
v_{\dff 1}\dff,\off
v_{\dff 2}\dff,\off
\ldots\dff,\off
v_{\dff n}
\off =\off
w$\pss
in\dss $\tau$\dss
such\dss that\qss\vspace{4.5pt}
\[
\quad
f_{\dff 2}\dff\left(\qff \tau_{\dff 2} \qff\right)
\trf \cdot\qff 
f_{\dff 3}\dff\left(\qff \tau_{\dff 3} \qff\right)
\trf \cdot\qff
\ldots
\trf \cdot\qff
f_{\dff n}\dff\left(\qff \tau_{\dff n} \qff\right)
\off =\off
1
\qff,
\]

\vspace{-7.5pt}
where\qss
$\tau_{\dff i}
\off =\off
\{\qff v_{\fff i\dff -\dff 1}\dff,\pff v_{\fff i} \pff\}$\qss
for each\qss $i\qff \geq\qff 2$\nnsp.\oss
It\dss follows\dss that\qss 
$f_{\dff i}\dff\left(\qff \tau_{\dff i} \qff\right)\off =\off 1$\qss
and\dss hence\dss
$\varphi\qff(\dff v_{\dff i}\trf)$\dss
differs\dss from\dss
$\varphi\qff(\dff v_{\dff i\dff -\dff 1}\trf)$\dss
in\dss the $i${\dnsp}th\dss coordinate
for every\dss $i$\nnsp.\oss
Since\dss $\varphi$\dss is\dss adjacency-preserving\halfff,\oss
$\varphi\qff(\dff v_{\dff i}\trf)$\dss
cannot\dss differ\dss from\dss
$\varphi\qff(\dff v_{\dff i\dff -\dff 1}\trf)$\dss
in\sss any\sss other coordinate,\oss
in\dss particular\halfff,\oss in\dss the first\sss one.\oss
It\dss follows\dss that\dss the first\sss coordinate
of\dss
$\varphi\qff(\dff v_{\dff i}\trf)$\dss
is\dss independent\sss of\dss $i$\dss
and\dss that\dss
$\varphi\qff(\dff w\trf)
\off =\off 
\varphi\qff(\dff v_{\dff n}\trf)$\dss
differs\sss from\dss
$\varphi\qff(\dff v\trf)
\off =\off 
\varphi\qff(\dff v_{\dff 0}\trf)$\dss
in all\sss other coordinates.\oss
In\sss other\dss words,\pss
$\varphi\qff(\dff v\trf)$
and\dss
$\varphi\qff(\dff w\trf)$\dss
both\dss belong\dss either\dss to\dss
$\mathbb{A}_{\dff 1}$\dss
or\qss
$\mathbb{B}_{\dff 1}$\dss
and are opposite vertices of\qss
$\mathbb{A}_{\dff 1}$\dss
or\qss
$\mathbb{B}_{\dff 1}$\dss
respectively\halfff.\oss

Suppose now\dss that\qss
$v
\off =\off
w_{\dff 1}\dff,\off
w_{\dff 2}\dff,\off
\ldots\dff,\off
w_{\dff n}
\off =\off
w$\pss
is\dss some other\dss free pivot\sss sequence in\dss $\tau$\dss
with\dss the same starting and ending\dss points.\oss
Since\dss $\varphi$\dss is\dss adjacency-preserving\halfff,\oss
$\varphi\qff(\dff w_{\dff i}\trf)$\dss
cannot\dss differ\dss from\dss
$\varphi\qff(\dff w_{\dff i\dff -\dff 1}\trf)$\dss
in\dss more\sss than one coordinate.\oss
Since\dss
$\varphi\qff(\dff w_{\dff n}\trf)
\off =\off 
\varphi\qff(\dff w\trf)$\dss
differs\sss from\dss
$\varphi\qff(\dff w_{\dff 0}\trf)
\off =\off 
\varphi\qff(\dff v\trf)$\dss
in\dss all\sss coordinates except\dss the first\sss one,\oss
it\dss follows\sss that\dss the fist\sss coordinate of\dss
$\varphi\qff(\dff w_{\dff i}\trf)$\dss
is\dss independent\sss of\dss $i$\dss and\dss hence\dss is\dss
equal\dss to\sss the first\sss coordinate of\dss
$\varphi\qff(\dff v\trf)$\nnsp.\oss
Since every\dss point\sss of\dss $\tau$\dss is\dss a\sss term\sss of\dss a\sss
pivot\sss sequence starting\sss at $v$\sss and ending\sss at $w$\nnsp,\oss
every\dss point\dss of\dss $\varphi\dff(\dff \tau\dff)$\dss has\sss
the same\sss first\sss coordinate as\dss $\varphi\qff(\dff v\trf)$\nnsp.\oss
Hence\dss $\varphi\dff(\dff \tau\dff)$\dss is\dss contained\dss in\sss either\dss  
$\mathbb{A}_{\dff 1}$\dss or\qss
$\mathbb{B}_{\dff 1}$\nsp.\oss

This proves\sss the first\sss statement\sss of\trs the lemma.\oss
The second\sss statement\dss follows\dss from\dss the first\sss one and\sss the fact\dss 
that\dss the root\dss $\lambda^0\trf \tau$\dss is\dss a\sss $0$\dnsp-face\sss of\dss
$\tau$\nnsp.\oss  \eproof\vspace{2.365pt}

\mypar{Theorem.}{product-to-bijection}
\emph{If\oss
$f_{\dff 1}\dff \cdot\dff f_{\dff 2}\trf \cdot\qff \ldots\dff \cdot\qff 
f_{\dff n}\qff(\dff \sigma\dff)
\off =\off 
1$\nnsp,\oss
then\dss $\varphi$\dss is\dss a\sss bijection.\oss}\vspace{2.365pt}

\proof
After\dss identifying\dss $\sigma$\dss with\dss $\mathbb{K}$\dss
by\sss a\sss translation\dss 
we may\sss apply\qss Theorem\qss \ref{products-induction}\qss
to\dss $\sigma$\dss in\dss the role of\qss $K$\dss and\dss the cochains\qss
$h_{\dff 1}\fff,\pff h_{\dff 2}\dff,\pff \ldots\fff,\pff h_{\dff n}$\nsp.\oss
Clearly\halfff,\oss in\dss this situation\qss
$s\off =\off 1$\nnsp.\oss
Also,\pss every $(\fff n\dff -\dff 1\fff)$\dnsp-cube contained\dss in\dss $\sigma$\dss
is\trs automatically\sss contained\dss in\dss $\bd\fff \sigma$
and\dss is\dss an $(\fff n\dff -\dff 1\fff)$\dnsp-face of\dss $\sigma$\nnsp.\oss

Therefore\dss Theorem\qss \ref{products-induction}\qss implies\sss that\dss there\dss is\dss
and\dss odd\dss
number\sss of\sss $(\fff n\dff -\dff 1\fff)$\dnsp-faces\dss $\tau$\dss
of\qss $\sigma$\dss 
such\dss that\qss
$h_{\dff 1}\dff(\qff \lambda^0\trf \tau\qff)
\off =\off 
1$\qss
and\pss (\ref{2-n-product})\pss holds.\oss 
In view of\qss Lemma\qss \ref{face-product}\qss this means\sss that\dss
the number of\sss $(\fff n\dff -\dff 1\fff)$\dnsp-faces\dss $\tau$\dss
of\qss $\sigma$\dss 
such\dss that\qss
$\varphi\dff(\dff \tau\dff)\off \subset\dff\off \mathbb{B}_{\dff 1}$\qss
and\pss (\ref{2-n-product})\pss holds\dss
is\dss odd.\oss

We claim\dss that\dss
the number\dss $N$\dss of\dss $(\fff n\dff -\dff 1\fff)$\dnsp-faces\dss $\tau$\dss
such\dss that\qss
(\ref{2-n-product})\qss holds\sss
is\dss even.\oss
By\qss Lemma\qss \ref{leibniz}\vspace{4.5pt}
\[
\quad
\partial^*\left(\qff
h_{\dff 2}\dff \cdot\dff f_{\dff 3}\trf \cdot\qff \ldots\dff \cdot\qff f_{\dff n}
\qff\right)
\off =\off
f_{\dff 2}\dff \cdot\dff f_{\dff 3}\trf \cdot\qff \ldots\dff \cdot\qff f_{\dff n}
\qff.
\]

\vspace{-7.5pt}
The identity\qss $\partial^* \circ\qff \partial^*
\off =\off
0$\qss implies\sss that\qss
$\partial^*\left(\qff
f_{\dff 2}\dff \cdot\dff f_{\dff 3}\trf \cdot\qff \ldots\dff \cdot\qff f_{\dff n}
\qff\right)
\off =\off 0$\qss
and\dss hence\vspace{4.5pt}
\[
\quad
0
\off =\off
\partial^*\left(\qff
f_{\dff 2}\dff \cdot\dff f_{\dff 3}\trf \cdot\qff \ldots\dff \cdot\qff f_{\dff n}
\qff\right)\dff
(\qff \sigma\qff)
\]

\vspace{-31.75pt}
\[
\quad
\phantom{0
\off }
=\off\dff
f_{\dff 2}\dff \cdot\dff f_{\dff 3}\trf \cdot\qff \ldots\dff \cdot\qff f_{\dff n}\qff
(\qff \partial\dff\sigma\qff)
\off\dff =\off\qff
\sum_{\tau\vphantom{\sigma^2}}\qff\off
f_{\dff 2}\dff \cdot\dff f_{\dff 3}\trf \cdot\qff \ldots\dff \cdot\qff f_{\dff n}\qff
(\qff \tau\qff)
\qff,
\]

\vspace{-12pt}
where\dss the sum\dss is\dss taken over all\sss $(\fff n\dff -\dff 1\fff)$\dnsp-faces $\tau$
of\trs $\sigma$\nnsp.\oss
It\dss follows\dss that\dss $N$\dss is\dss indeed\sss even.\oss
So,\oss the number of\dss $(\fff n\dff -\dff 1\fff)$\dnsp-faces\dss $\tau$\dss
such\dss that\qss
(\ref{2-n-product})\qss holds\sss
is\dss even,\oss
and\dss there\dss is\dss an\sss odd\dss number of\trs faces\dss $\tau$\dss such\dss that\qss
$\varphi\dff(\dff \tau\dff)\off \subset\dff\off \mathbb{B}_{\dff 1}$\qss
among\dss them.\oss
Therefore,\oss there\dss is\dss also an\sss odd\dss number of\trs faces\dss $\tau$\dss such\dss that\qss
$\varphi\dff(\dff \tau\dff)\off \subset\dff\off \mathbb{A}_{\dff 1}$\qss
among\dss them.\oss
In\dss particular\halfff,\oss
there\dss exist\dss  $(\fff n\dff -\dff 1\fff)$\dnsp-faces\dss $\tau_{\dff \mathbb{A}}$\dss
and\dss $\tau_{\dff \mathbb{B}}$\dss
such\dss that\qss
$\varphi\dff(\dff \tau_{\dff \mathbb{A}}\dff)\off \subset\dff\off \mathbb{A}_{\dff 1}$\qss
and\qss
$\varphi\dff(\dff \tau_{\dff \mathbb{B}}\dff)\off \subset\dff\off \mathbb{B}_{\dff 1}$\nnsp.\oss
Clearly\halfff,\pss
$\tau_{\dff \mathbb{A}}$\qss
and\qss $\tau_{\dff \mathbb{B}}$\qss
are disjoint\halfff.\oss

Now\dss it\dss is\dss time\sss to apply\sss an\dss induction\dss by\dss $n$\nnsp.\oss
The\sss theorem\dss is\dss trivially\dss true for\qss $n\off =\off 1$\nnsp.\oss
Arguing\dss by\dss induction,\oss
we may assume\sss that\dss it\dss is\dss true with\qss $n\qff -\qff 1$\qss 
in\dss the role of\dss $n$\nnsp.\oss
Let\dss $\tau$\dss be an $(\fff n\dff -\dff 1\fff)$\dnsp-face of\dss $\sigma$\dss
such\dss that\qss (\ref{2-n-product})\qss holds.\oss
Consider\dss the map\qss
$\psi
\dff \colon\dff
\tau\qff \ttoo\qss
\{\qff 0\dff,\pff 1\qff\}^{\fff n\dff -\dff 1}$\qss
having\qss
$h_{\dff 2}\dff,\off h_{\dff 3}\dff,\off \ldots\dff,\off h_{\dff n}$\qss
as\dss its\dss components.\oss
The inductive assumption\dss implies\sss that\dss $\psi$\dss is\dss a\sss bijection.\oss
Together\dss with\dss Lemma\qss \ref{face-product}\qss
this\sss implies\sss that\dss $\varphi$\dss maps\dss $\tau$\dss bijectively\sss onto
either\dss $\mathbb{A}_{\dff 1}$\dss or\qss $\mathbb{B}_{\dff 1}$\nsp.\oss
By\sss applying\dss this result\dss to\qss
$\tau\off =\off \tau_{\dff \mathbb{A}}$\qss
and\qss
$\tau\off =\off \tau_{\dff \mathbb{B}}$\qss
we see\sss that\dss $\varphi$\dss maps\qss 
$\tau_{\dff \mathbb{A}}$\qss
and\qss $\tau_{\dff \mathbb{B}}$\qss
bijectively\sss onto\dss
\dss $\mathbb{A}_{\dff 1}$\dss and\pss $\mathbb{B}_{\dff 1}$\dss
respectively\halfff.\oss
Since\qss 
$\tau_{\dff \mathbb{A}}$\qss
and\qss $\tau_{\dff \mathbb{B}}$\qss
are\sss two disjoint\sss
$(\fff n\dff -\dff 1\fff)$\dnsp-face of\dss $\sigma$\nnsp,\oss
it\dss follows\dss that\dss $\varphi$\dss is\dss a\sss bijection.\oss  \eproof

\mypar{Theorem.}{equivalence}
\emph{The\dss map\dss $\varphi$\dss is\dss a\dss bijection\qss
if\pss and\trs only\qss if\qff\oss
$f_{\dff 1}\dff \cdot\dff f_{\dff 2}\trf \cdot\qff \ldots\dff \cdot\qff 
f_{\dff n}\qff(\dff \sigma\dff)
\off =\off\dff 
1$\nnsp.\oss}  \eproof

\myuppar{Ky\dss Fan's\dss lemma.}
\emph{Let\qss
$\Phi\dff \colon\dff
K\qff \ttoo\qff \mathbb{K}$\qss
be an adjacency-preserving\dss map,\oss
and\dss let\dss $\mathbb{B}$\dss be some\dss $(\fff n\qff -\qff 1\fff)$\dnsp-face of\oss $\mathbb{K}$\nnsp.\oss
Let\dss $s$\dss be\dss the number of\qss $n$\dnsp-cubes of\oss $K$\dss
which are mapped\dss by\dss $\Phi$\dss bijectively\dss onto\qss $\mathbb{K}$\nnsp,\oss
and\trs let\qss $t$\dss be\trs the number of\pss $(\fff n\qff -\qff 1\fff)$\dnsp-cubes
of\pss $K$\dss contained\dss in\dss $\bd\fff K$\dss and\dss bijectively\dss mapped\dss by\dss $\Phi$\dss 
onto\qss $\mathbb{B}$\nnsp.\oss
Then\qss $s\off \equiv\off t$\qss modulo\dss $2$\nnsp.\oss}

\proof
We may assume\sss that\qss
$\mathbb{B}
\off =\off
\mathbb{B}_{\dff 1}
\off =\off
1\dff \times\dff 
\{\qff 0\dff,\pff 1\qff\}^{\fff n\dff -\dff 1}$\dnsp.\oss
Let\qss 
$h_{\dff 1}\dff,\off
h_{\dff 2}\dff,\off
\ldots\dff,\off
h_{\dff n}$\qss 
be\sss the components of\dss $\Phi$\dss considered
as $0$\dnsp-cochains of\trs $K$\nnsp,\pss
and\dss let\qss
$f_{\dff i}
\off =\off
\partial^* h_{\dff i}$\nsp.\oss
Theorem\qss \ref{equivalence}\qss implies\sss that\dss $s$\dss
is\dss equal\dss to\sss the number of\dss $n$\dnsp-cubes\dss $\sigma$\dss of\trs $K$\dss
such\dss that\qss
$f_{\dff 1}\dff \cdot\dff f_{\dff 2}\trf \cdot\qff \ldots\dff \cdot\qff 
f_{\dff n}\qff(\dff \sigma\dff)
\off =\off\dff 
1$\nnsp.\oss
By\dss combining\qss Theorem\qss \ref{equivalence}\qss
with\qss $n\qff -\qff 1$\qss in\dss the role of\dss $n$\dss
with\qss Lemma\qss \ref{face-product}\qss we see\sss that\dss
$t$\dss is\dss equal\dss to\sss the number of\dss 
$(\fff n\qff -\qff 1\fff)$\dnsp-cubes\qss 
$\tau 
\off \subset\off \bd\fff K$\qss such\dss that\qss
$h_{\dff 1}\dff(\qff \lambda^0\trf \tau\qff)
\off =\off 
1$\qss and\qss (\ref{2-n-product})\pss holds.\oss
These observations show\dss that\dss the\sss theorem\dss is\dss
a\sss special\sss case of\qss
Theorem\qss \ref{products-induction}.\oss  \eproof

\myuppar{Corollary\halfff.}
\emph{Let\oss
$L\dff \colon\dff
K\off \ttoo\off \mathbb{K}$\oss
be an adjacency-preserving\dss map satisfying\dss the assumptions of\pss
Kuhn's\qss lemma.\oss
Then\dss the number of\qss $n$\dnsp-cubes of\pss $K$\dss mapped\dss by\dss $L$\dss
bijectively\dss onto\dss $\mathbb{K}$\dss is\dss odd.\oss}

\proof
It\dss is\dss sufficient\dss to use\qss Ky\dss Fan\qss theorem\dss together\dss with 
an\dss induction\dss by\dss $n$\nnsp.\oss  \eproof

\myuppar{Adjacency-preserving\trs and\trs transversality\halfff.}
Let\qss
$\Phi\dff \colon\dff
K\off \ttoo\off \mathbb{K}$\qss
be an arbitrary\dss map and\dss let\qss
$h_{\dff 1}\dff,\off
h_{\dff 2}\dff,\off
\ldots\dff,\off
h_{\dff n}$\qss
be its components.\oss
As\dss in\dss Section\qss \ref{kuhn},\oss
let\dss 
$c_{\dff i}$\dss be\sss the set\sss of\trs points\qss
$v\qff \in\qff K$\qss
such\dss that\qss
$h_{\dff i}\dff(\dff v\trf)\off =\off 0$\nnsp.\oss
Then\dss $\Phi$\dss is\dss determined\dss by\dss the sets\qss
$c_{\dff 1}\fff,\pff
c_{\dff 2}\trf,\pff
\ldots\fff,\pff
c_{\dff n}$\nsp.\oss
Let\qss
$\overline{c}_{\dff i}
\off =\off
K\qff \smallsetminus\qff c_{\dff i}$\nsp.\oss
Considering\dss the maps\dss $h_{\dff i}$\dss as $0$\dnsp-chains\dss of\trs $K$\dss
allows\sss us\sss define $1$\dnsp-cochains\dss
$f_{\dff i}
\off =\off
\partial^* h_{\dff i}$\nsp.\oss

Let\dss $Q$\dss be\sss the solid cube of\dss size\qss
$l\off =\off k\qff +\qff 1$\nnsp,\oss 
where\dss $k$\dss is\dss the size of\trs $K$\nnsp.\qff\oss
Let\vspace{3pt}
\[
\quad
s_{\dff i}
\off\dff =\off
*\dff c_{\dff i}
\hspace*{1.2em}\mbox{and}\hspace*{1.5em}
\overline{s}_{\dff i}
\off\dff =\off
*\dff \overline{c}_{\dff i}
\off
\]

\vspace{-9pt}
be\sss the cubic subsets of\dss $Q$\dss dual\dss in\dss the sense of\qss
Appendix\qss \ref{two-cubes}\pss to\dss
$c_{\dff i}$\dss and\qss $\overline{c}_{\dff i}$\dss
respectively\halfff.\oss
For every\qss
$i\off =\off 1\fff,\pff 2\fff,\pff \ldots\fff,\pff n$\qss
let\dss us\sss consider\dss the\sss formal\sss sum\dss $\gamma_{\dff i}$\dss
of\dss
$(\fff n\dff -\dff 1\fff)$\dnsp-cubes\dss $\tau$\dss of\dss $Q$\dss
such\dss that\sss $\tau$\sss is\dss a\sss common\dss face\sss of\trs an $n$\dnsp-cube 
contained\dss in $s_{\dff i}$\dss
and\dss an $n$\dnsp-cube 
contained\trs in\qss $\overline{s}_{\dff i}$\hnsp.\oss
In\sss a\sss natural\sss sense\dss $\gamma_{\dff i}$\dss
is\dss the common\dss boundary\sss of\trs the cubical\sss sets\dss
$s_{\dff i}$\dss and\qss
$\overline{s}_{\dff i}$\hnsp.\oss
One may also say\dss that\dss $\gamma_{\dff i}$\dss is\dss the\sss internal\dss part\sss
of\qss the boundary\sss of\dss $s_{\dff i}$\nsp.\oss

If\qss
$\tau\off =\off \{\qff v\fff,\qff w \qff\}$\qss is\dss a $1$\dnsp-cube,\oss
then\qss
$f_{\dff i}\trf(\dff \tau\dff)\off =\off 1$\qss
if\qss and\dss only\trs if\pss
one of\trs the $n$\dnsp-cubes\dss
$*\dff v$\nnsp,\qss $*\dff w$\dss is\dss contained\dss in\dss 
$*\dff c_{\dff i}$\dss and\dss the other\sss does not\halfff.\oss
By\qss Lemma\qss \ref{two-cubes}\qss
the dual\sss $(\fff n\dff -\dff 1\fff)$\dnsp-cube\dss $*\dff \tau$\dss
is\dss the common\dss face\sss of\qss $*\dff v$\dss and\dss $*\dff w$\nnsp.\oss
It\dss follows\dss that\qss 
$*\qff \gamma_{\dff i}
\off =\off\dff
f_{\dff i}$\nsp.\oss

As we already\dss saw,\pss
$\Phi$\dss is\dss adjacency-preserving\qss if\qss and\dss only\trs if\pss
for every $1$\dnsp-cube\dss $\tau$\dss
there\sss is\dss no more\sss than one $i$ such\dss that\qss
$f_{\dff i}\trf(\dff \tau\dff)\pff =\off 1$\nnsp.\oss
In\dss the language of\dss $Q$\dss this means\sss that\dss
every\sss $(\fff n\dff -\dff 1\fff)$\dnsp-cube of\dss $Q$\dss
enters no more\dss than one\dss $\gamma_{\dff i}$\nsp.\oss
In other\dss words,\oss if\qss $i\off \neq\off j$\nnsp,\oss
then\dss $\gamma_{\dff i}$\dss and\dss $\gamma_{\fff j}$\dss
have no common $(\fff n\dff -\dff 1\fff)$\dnsp-cubes.\oss
If\dss one\sss thinks about\dss $\gamma_{\dff i}$\dss and\dss $\gamma_{\fff j}$\dss
as hypersurfaces,\oss
then\dss this means\dss that\dss they are nowhere\sss tangent\halfff.\oss
Two smooth\dss hypersurfaces are nowhere\sss tangent\dss if\qss and\dss only\trs if\qss
they are\sss transverse.\oss
We see\sss that\dss the adjacency-preserving\dss property\dss
is\dss a\sss sort\sss of\trs transversality\halfff.\oss

The pivot\dss sequences\qss
$v_{\dff 0}\fff,\pff v_{\dff 1}\dff,\pff \ldots\fff,\pff v_{\dff n}$\qss
such\dss that\qss 
$\{\qff
v_{\dff 0}\fff,\pff v_{\dff 1}\dff,\pff \ldots\fff,\pff v_{\dff n}
\qff\}$\qss
intersects\dss both\dss
$c_{\dff i}$\dss and\qss
$\overline{c}_{\dff i}$\qss
for\dss every\dss $i$\dss
correspond\dss to\sss the points of\dss intersection of\trs all\qss
hypersurfaces\dss $\gamma_{\dff i}$\nsp,\pss
and\qss Kuhn's\dss lemma ensures\sss that\dss these hypersurfaces do intersect\halfff.\oss
But\dss it\dss provides\sss no information about\qss \emph{how}\qss they\dss intersect\halfff.\oss
As\dss is\dss well\dss known,\oss
under\dss appropriate\sss transversality\dss assumptions\dss 
the intersections can\dss be understood\dss better\dss than\sss otherwise.\oss
If\dss $\Phi$\dss is\dss adjacency-preserving\halfff,\oss
then\qss Ky\dss Fan\dss theorem\dss ensures\sss the existence of\dss cubes\dss $\sigma$\dss of\trs $K$
such\dss that\dss $\Phi$\dss
maps\dss $\sigma$\dss bijectively\sss onto\dss $\mathbb{K}$\nnsp.\oss
Such cubes
also correspond\dss to\sss the points of\dss intersection of\trs all\qss
hypersurfaces\dss $\gamma_{\dff i}$\qss
(in\dss fact\halfff,\pss every\dss pivot\sss sequence from\qss Kuhn's\qss lemma\dss
is\dss contained\dss in\sss such a cube).\oss
The\sss bijectivity\dss tells us\dss that\trs
$2^{\fff n}$\sss cubes\dss $c$\dss meeting at\dss such\dss point\sss of\dss intersection\dss
realize all\trs
$2^{\fff n}$\sss compatible combinations of\dss conditions\qss
$c\qff \subset\qff s_{\dff i}$\nsp,\oss
$c\qff \subset\qff \overline{s}_{\fff j}$\qss
(where\qss $1\qff \leq\qff i\fff\halfff,\pff j\qff \leq\qff n$\nsp).\oss
Therefore at\dss least\dss in\dss this respect\dss the intersection of\trs the\dss hypersurfaces\dss
defined\dss by\sss an\sss adjacency-preserving\dss map\dss $\Phi$\dss
is\dss similar\dss to\dss the simplest\sss case,\oss
the intersection of\dss $n$\dss hyperplanes parallel\dss to\sss the coordinate
hyperplanes.\oss

\myuppar{A question of\qss Kuhn.}
Kuhn\qss \cite{ku}\qss asked\dss if\qss there\dss is\dss any\dss relation\dss
between\qss Kuhn's\trs lemma\dss and\qss Ky\dss Fan's\qss theorem.\oss
The\sss author\dss believes\dss that\dss the above discussion\dss provides\sss an\sss answer\halfff.\oss

\newpage
\myappend{Freudenthal's\pss triangulations\qss of\qss cubes\qss and\qss simplices}{freudenthal}

\myuppar{The unit\dss cube\dss $[\trf 0\fff,\pff 1\trf]^{\fff n}$\dnsp.}
As usual,\oss let\qss
$I\off =\off
\{\trf 1\fff,\pff 2\trf,\pff \ldots\fff,\pff n \qff\}$\nnsp.\oss
For every permutation\qss ({\fff}i.e\qss bijection)\dss 
$\omega\dff \colon\dff
I \qff \ttoo\qff I$\qss 
let\dss $\Delta_{\dff \omega}$\dss be\sss the set\sss of\dss points\qss
$(\dff x_{\dff 1}\fff,\pff x_{\dff 2}\trf,\pff \ldots\fff,\pff x_{\dff n} \dff)
\qff \in\fff\pff
\rrr^{\fff n}$\qss
such\dss that\vspace{2.5pt}
\begin{equation}
\label{decreasing-omega}
\quad
1\off \geq\off
x_{\trf \omega\dff(\dff 1\dff)}\off \geq\off
x_{\trf \omega\dff(\dff 2\dff)}\off \geq\off
\ldots\off \geq\off
x_{\trf \omega\dff(\dff n\dff)}
\off \geq\off
0
\qff.
\end{equation}

\vspace{-9.5pt}
Obviously,\oss the cube\dss $[\trf 0\fff,\pff 1\trf]^{\fff n}$\dss is equal\dss to\dss
the union of\dss the sets\dss $\Delta_{\dff \omega}$\dss over all\dss permutations\dss $\omega$\nnsp.\oss
If\dss $\omega$\dss is\dss the identity,\oss then\qss 
$\Delta\off =\off \Delta_{\dff \omega}$\qss
is\dss the set\sss of\dss points\qss
$(\dff x_{\dff 1}\fff,\pff x_{\dff 2}\trf,\pff \ldots\fff,\pff x_{\dff n} \dff)
\qff \in\fff\pff
\rrr^{\fff n}$\qss
such\dss that\vspace{2.5pt}
\begin{equation}
\label{decreasing}
\quad
1\off \geq\off
x_{\dff 1}\off \geq\off
x_{\dff 2}\off \geq\off
\ldots\off \geq\off
x_{\dff n}
\off \geq\off
0
\qff.
\end{equation}

\vspace{-9.5pt}
\myapar{Lemma.}{main-simplex}
\emph{\dnsp$\Delta$\qss
is\dss the geometric $n$\dnsp-simplex 
with\dss the vertices\qss
$u_{\dff 0}\fff,\pff u_{\dff 1}\trf,\pff \ldots\fff,\pff u_{\dff n}
\qff \in\qff 
\rrr^{\fff n}$\dnsp,\oss
where}\vspace{2.5pt}
\[
\quad
u_{\dff i}
\off =\off
(\dff 1\fff,\pff \ldots\fff,\pff 1\fff,\pff 0\trf,\pff \ldots\fff,\pff 0 \dff)
\qff
\]

\vspace{-9.5pt}
\emph{is\dss the point\dss with\dss the\dss first\dss $i$\sss coordinates equal\dss to\dss $1$
and\dss the last\trs $n\qff -\qff i$\trs coordinates equal\dss to\dss $0$\nnsp.\oss
Every\qss face of\qss $\Delta$\dss is\trs defined\dss by\qss
(\ref{decreasing})\qss with several\dss inequality signs\pss $\geq$\qss
replaced\dss by\pss $=$\nsp.\oss}

\proof
Indeed,\oss
every\qss
$(\dff x_{\dff 1}\fff,\pff x_{\dff 2}\trf,\pff \ldots\fff,\pff x_{\dff n} \dff)
\qff \in\qff 
\rrr^{\fff n}$\qss
has a unique presentation\dss in\dss the form\vspace{2.5pt}
\begin{equation}
\label{y-to-x}
\quad
(\dff x_{\dff 1}\fff,\pff x_{\dff 2}\trf,\pff \ldots\fff,\pff x_{\dff n} \dff)
\off =\off
y_{\dff 0}\dff u_{\dff 0}
\qff +\qff
y_{\dff 1}\dff u_{\dff 1}
\qff +\qff
\ldots
\qff +\qff
y_{\dff n\dff -\dff 1}\dff u_{\dff n\dff -\dff 1}
\qff +\qff
y_{\dff n}\dff u_{\dff n}
\qff,
\end{equation}

\vspace{-9.5pt}
where\qss
$y_{\dff 0}\dff,\pff y_{\dff 1}\dff,\pff \ldots\dff,\pff y_{\dff n}
\qff \in\qff \rrr$\qss
and\qss\vspace{2.5pt}
\[
\quad
y_{\dff 0}\qff +\qff y_{\dff 1}\qff +\qff \ldots\qff +\qff y_{\dff n}
\off =\off 
1
\qff.
\]

\vspace{-9.5pt}
In\dss fact\halfff,\oss (\ref{y-to-x})\qss holds\dss 
if\dss and\dss only\trs if\dss\vspace{2.5pt} 
\[
\quad
x_{\dff i}
\off =\off 
y_{\dff i}\qff +\qff  y_{\dff i\dff +\dff 1}\qff +\qff\ldots\qff +\qff y_{\dff n}
\]

\vspace{-9.5pt}
for all\qss
$i\off =\off 1\fff,\pff 2\fff,\pff \ldots\fff,\pff n$\qss
and\dss hence\dss if\dss and\dss only\trs if\pss\vspace{2.5pt} 
\[
\quad
y_{\dff n}\off =\off x_{\dff n}
\hspace{1em}\mbox{and}\hspace{1.2em}
y_{\dff i}\off =\off x_{\dff i}\qff -\qff x_{\dff i\dff +\dff 1}
\]

\vspace{-9.5pt}
for all\qss $i\qff \leq\qff n\qff -\qff 1$\nnsp.\oss
These equalities imply\dss that\qss
$x_{\dff 1}\off =\off y_{\dff 1}\qff +\qff y_{\dff 2}\qff +\qff \ldots\qff +\qff y_{\dff n}$\qss
and\qss
$y_{\dff 0}\off =\off 1\qff -\qff x_{\dff 1}$\nsp.\oss
It\dss follows\dss that\qss
$y_{\dff 0}\dff,\pff y_{\dff 1}\dff,\pff \ldots\dff,\pff y_{\dff n}
\qff \geq\qff 0$\qss
if\dss and\dss only\dss if\qss (\ref{decreasing})\qss holds\dss and\dss hence\dss
$\Delta$\dss is\dss
the convex\dss hull\sss of\trs the points\qss
$u_{\dff 0}\fff,\pff u_{\dff 1}\trf,\pff \ldots\fff,\pff u_{\dff n}$\nsp.\oss
Clearly,\oss these points are affinely\dss independent\sss
and\dss hence\dss $\Delta$\dss is\dss indeed\dss a geometric simplex
and\dss has\sss these points as its vertices.\oss
Every\dss face of\dss $\Delta$\dss is\dss defined\dss by\sss
requiring several\dss barycentric coordinates\dss $y_{\dff i}$\dss to be\dss $0$\nnsp.\oss 
But\dss the equality\qss $y_{\dff i}\off =\off 0$\qss
is\dss equivalent\dss to\qss 
$x_{\dff n}\off =\off 0$\qss if\qss $i\off =\off n$\nnsp,\oss
to\qss 
$x_{\dff 1}\off =\off 1$\qss if\qss $i\off =\off 0$\nnsp,\oss
and\dss to\qss 
$x_{\dff i}\off =\off x_{\dff i\dff +\dff 1}$\qss if\qss $0\qff <\qff i\qff <\off n$\nnsp.\oss  \eproof

\myapar{Corollary.}{other-simplices}
\emph{\dnsp$\Delta_{\dff \omega}$\dss is\dss a\sss geometric $n$\dnsp-simplex\dss 
for every\dss $\omega$\nnsp.\oss 
Its vertices\qss
$v_{\dff 0}\fff,\pff
v_{\dff 1}\trf,\pff
\ldots\fff,\pff
v_{\dff n}$\qss
belong\dss to\qss
$\{\qff 0\dff,\pff 1\qff\}^{\dff n}$\qss
and\sss can\dss be ordered\dss in\dss such a way\dss that}\oss\vspace{0pt} 
\[
\quad
0
\off =\off
v_{\dff 0}\off <\off
v_{\dff 1}\off <\off
\ldots\off <\off
v_{\dff n}
\off =\off
\mathbb{1}
\qff,
\]

\vspace{-12pt}
\emph{where\qss $<$\qss is\dss the order\dss from\qss Sections\qss \ref{tilings}\qss and\qss \ref{kuhn}.\oss
Every\qss face of\qss $\Delta_{\dff \omega}$\dss is\trs defined\dss by\qss
(\ref{decreasing-omega})\qss with several\dss inequality signs\pss $\geq$\qss
replaced\dss by\pss $=$\nsp.\oss}

\proof
Obviously,\oss the vertices\qss
$u_{\dff 0}\fff,\pff u_{\dff 1}\trf,\pff \ldots\fff,\pff u_{\dff n}$\qss
of\dss $\Delta$\dss belong\dss to\qss $\{\qff 0\dff,\pff 1\qff\}^{\dff n}$\nsp.\oss
Moreover\halfff,\oss\vspace{0pt} 
\[
\quad
0
\off =\off
u_{\dff 0}\off <\off
u_{\dff 1}\off <\off
\ldots\off <\off
u_{\dff n}
\off =\off
\mathbb{1}
\qff.
\]

\vspace{-12pt}
Every\sss set\dss $\Delta_{\dff \omega}$\dss is\dss the image of\dss
$\Delta$\dss under a map\qss
$\rrr^{\fff n}\qff \ttoo\qff \rrr^{\fff n}$\qss
permuting\sss coordinates.\oss
It\dss remains\sss to notice\sss that\sss such maps preserve\sss 
the order\qss $<$\qss and other\dss relevant\dss properties.\oss  \eproof

\myapar{Corollary\halfff.}{faces-intersections}
\emph{For\dss every\qss  
$\omega\fff,\off \omega'$\dss
the intersection\qss
$\Delta_{\dff \omega}\qff \cap\qff \Delta_{\dff \omega'}$\dss
is\sss a\sss common\dss face of\pss the simplices\qss 
$\Delta_{\dff \omega}$\dss and\qss $\Delta_{\dff \omega'}$\nsp.}\vspace{-0.75pt}

\proof
Clearly,\qss 
$\Delta_{\dff \omega}\qff \cap\qff \Delta_{\dff \omega'}$\qss
can\dss be defined\dss by\qss
(\ref{decreasing-omega})\qss with several\dss  
signs\pss $\geq$\qss
replaced\dss by\pss $=$\nsp.\oss  \eproof

\myuppar{The canonical\dss triangulations of\qss unit\dss cubes.}
Since\dss $[\trf 0\fff,\pff 1\trf]^{\fff n}$\dss
is equal\dss to\sss the union of\dss simplices\dss $\Delta_{\dff \omega}$\nsp,\oss
Corollary\qss \ref{faces-intersections}\qss implies\sss that\sss
simplices\dss  $\Delta_{\dff \omega}$\dss together\sss with\dss their\dss faces
form a\sss triangulation of\trs the unit\sss cube\dss $[\trf 0\fff,\pff 1\trf]^{\fff n}$\dnsp.\oss
This\dss is\dss the\qss \emph{canonical\dss triangulation}\qss of\qss 
$[\trf 0\fff,\pff 1\trf]^{\fff n}$\dnsp.\oss
It\dss is\dss invariant\dss under all\dss permutations of\trs the coordinates,\oss
but\dss not\dss under\dss most\sss of\dss other symmetries of\dss $[\trf 0\fff,\pff 1\trf]^{\fff n}$\dnsp.\oss
An $m$\dnsp-cube $c$ of\dss $Q$\dss can\dss be identified\dss with\dss
$[\trf 0\fff,\pff 1\trf]^{\fff m}$\dss by a\sss translation followed\dss by a
permutation of\dss coordinates.\oss
Clearly,\oss the involved\dss translation\dss is\dss uniquely\sss determined\dss
by\sss $c$\nnsp.\oss
One can\dss use such an\sss identification\dss to\sss transplant\dss
the canonical\dss triangulation of\qss $[\trf 0\fff,\pff 1\trf]^{\fff m}$\dss
to\dss $c$\nnsp.\oss
Since\sss the canonical\dss triangulation of\qss $[\trf 0\fff,\pff 1\trf]^{\fff m}$\dss
is\dss invariant\dss under\sss permutation of\dss coordinates,\oss
the\sss result\sss does not\sss depend on\dss the choice of\dss
identification.\oss
It\dss is\sss called\dss the\qss \emph{canonical\dss triangulation}\qss of\dss $c$\nnsp.

\myapar{Lemma.}{faces-canonical}
\emph{If\qss $c$\dss is\dss a\dss face\dss of\oss $[\trf 0\fff,\pff 1\trf]^{\fff n}$\dnsp,\dff\oss
then\dss the collection of\dss 
simplices of\trs the canonical\dss
triangulation of\pss  $[\trf 0\fff,\pff 1\trf]^{\fff n}$\sss contained\dss in\dss $c$\dss
is\dss the canonical\dss triangulation of\qss $c$\nnsp.\oss}

\proof
The\sss face $c$ is\dss 
defined\dss by\dss requiring\sss some coordinates\sss to be $0$ and\sss some other\sss
to be $1$\nnsp.\oss
After\dss permuting coordinates we may assume\dss that\dss the coordinates
required\dss to be $0$ are\sss the\sss last\sss ones,\oss
and\dss required\dss to be $1$
are\sss the first\sss ones.\oss
In\dss this case\sss the\sss lemma\dss is\dss trivial.\oss  \eproof

\myuppar{The canonical\dss triangulation of\dss $Q$\nnsp.}
We\sss see\sss that\dss the canonical\dss triangulations of\dss cubes 
agree on\dss common\dss faces.\oss 
Hence\dss the union of\dss these canonical\sss triangulations
is\sss a\sss triangulation of\dss $Q$\nnsp.\oss
It\dss is\dss called\dss the\qss \emph{canonical\dss triangulation}\pss of\dss $Q$\nnsp.\oss
The corresponding abstract\sss simplicial\sss complex\dss has as\sss its set\sss of\dss
vertices\sss the discrete cube
of\dss the same size as\dss $Q$\nnsp.\oss
Corollary\qss \ref{other-simplices}\qss implies\sss that\dss it\dss is\dss
nothing else but\dss the abstract\sss simplicial\sss complex
from\dss Section\qss \ref{kuhn}.\oss

\myuppar{Freudenthal's\dss triangulations of\dss standard\sss simplices.}
We will\sss assume\sss that\dss the size $l$\sss of\trs 
the big cube\dss $Q$\dss is\qss $\geq\qff 2$\nnsp.\oss
Let\dss $l\fff \Delta$\dss the simplex defined\dss by\dss the inequalities\qss
(\ref{decreasing})\qss with $1$ replaced\dss by $l$\nnsp.\oss
Alternatively,\pss $l\fff \Delta$\dss is\dss the image of\dss $\Delta$\dss
under\dss the map\qss
$\rrr^{\fff n} \ttoo\qff \rrr^{\fff n}$\qss
multiplying each coordinate by $l$\nnsp.\oss
As we will\sss see in a moment\qss
(see\qss Theorem\qss \ref{permutations-inclusions}\qss 
and\dss Corollary\qss \ref{triangulations-all}\qss below),\pss
$l\fff \Delta$\dss is equal\dss to\sss the union of\dss simplices of\dss
the canonical\dss triangulation of\dss $Q$\dss contained\dss in\dss $l\fff \Delta$\nnsp.\oss
Therefore\sss these simplices form a\sss triangulation of\dss $l\fff \Delta$\nnsp,\oss
which we will\sss call\qss \emph{Freudenthal's\dss triangulation}\pss of\dss  $l\fff \Delta$\nnsp.\oss 
Freudenthal\qss himself\pss \cite{f}\qss considered only\dss the case\qss $l\off =\off 2$\nnsp.\oss

It\dss will\dss be convenient\dss to use instead of\dss permutations\dss $\omega$\dss
their\dss inverses.\oss
For every\dss permutation\dss $\omega\dff \colon\dff I\qff \ttoo\qff I$\pss
let\qss $\Delta\dff(\dff \omega\dff)\off =\off \Delta_{\qff \omega^{\dff -\dff 1}}$\nsp.\oss
Then\dss $\Delta\dff(\dff \omega\dff)$\dss consists\sss of\dss points\qss
$(\dff x_{\dff 1}\fff,\pff x_{\dff 2}\fff,\pff \ldots\fff,\off x_{\dff n} \dff)
\off \in\off
[\trf 0\fff,\pff 1\trf]^{\fff n}$\qss
such\dss that\qss $x_{\dff i}\off \geq\off x_{\fff j}$\qss
if\trs and\dss only\trs if\pss 
$\omega\dff(\dff i\dff)\off \geq\off \omega\dff(\dff j\dff)$\nnsp.\oss
Let\dss $l\fff \Delta\dff(\dff \omega\dff)$\dss be\sss the subset\sss of\qss
$Q\off =\off [\trf 0\fff,\pff l\qff]^{\fff n}$\qss
defined\dss by\dss the same inequalities.\oss
Every simplex of\trs the canonical\dss triangulation of\dss $Q$\dss
has\sss the form\qss
$a\qff +\qff \Delta\dff(\dff \omega\dff)$\nnsp,\oss
where\qss
$a\off =\off
(\dff a_{\dff 1}\fff,\pff a_{\dff 2}\fff,\pff \ldots\fff,\off a_{\dff n}\dff)$\qss
and\qss
$a_{\dff 1}\fff,\pff a_{\dff 2}\fff,\pff \ldots\fff,\off a_{\dff n}$\qss
are integers between $0$ and\qss $l\qff -\qff 1$\nnsp,\oss 
and $\omega$ is\dss a permutation.\oss
Let\dss us\dss look when\qss
$a\qff +\qff \Delta\dff(\dff \omega\dff)$\qss
is\dss contained\dss in\qss
$l\fff \Delta\dff(\dff \omega'\dff)$\nnsp.\oss

\myapar{Theorem.}{permutations-inclusions}
\emph{For every\dss permutation\dss $\omega'$\dss the simplex\qss
$a\qff +\qff \Delta\dff(\dff \omega\dff)$\qss
is\dss contained\dss in\qss
$l\fff \Delta\dff(\dff \omega'\dff)$\qss
if\qss and\dss only\qss if\oss
$\omega'\off =\off \omega\nsp\mid\nsp a$\nsp,\oss
where\dss $\omega\nsp\mid\nsp a$\dss is\dss the unique permutation such\dss that\qss}\vspace{3pt}
\[
\quad
\omega\nsp\mid\nsp a\qff(\dff i\dff)\off >\off \omega\nsp\mid\nsp a\qff(\dff j\dff)
\]

\vspace{-9pt}
\emph{if\qss either\qss
$a_{\dff i}\off >\off a_{\fff j}$\nsp,\oss
or\pss
$a_{\dff i}\off =\off a_{\fff j}$\qss
and\pss
$\omega\dff(\dff i\dff)\off >\off \omega\dff(\dff j\dff)$\nnsp.\oss
Every\dss simplex\qss $l\fff \Delta\dff(\dff \omega'\dff)$\qss
is\dss equal\dss to\sss the union of\qss simplices of\qss
the\qss form\qss $a\qff +\qff \Delta\dff(\dff \omega\dff)$\dss
contained\dss in\dss it\halfff.\oss}

\proof
Every\dss point\sss of\dss $Q$\dss has\sss the form\qss\vspace{3pt}
\[
\quad
a\qff +\qff x
\off =\off
(\dff 
a_{\dff 1}\qff +\qff x_{\dff 1}\dff,\off 
a_{\dff 2}\qff +\qff x_{\dff 2}\dff,\off 
\ldots\dff,\off 
a_{\dff n}\qff +\qff x_{\dff n} 
\dff)
\qff,
\]

\vspace{-9pt}
where\qss
$a\off =\off
(\dff a_{\dff 1}\fff,\pff a_{\dff 2}\fff,\pff \ldots\fff,\off a_{\dff n}\dff)$\qss
is\dss as above and\qss
$x
\off =\off
(\dff x_{\dff 1}\fff,\pff x_{\dff 2}\fff,\pff \ldots\fff,\off x_{\dff n}\dff) 
\off \in\off
[\trf 0\fff,\pff 1\trf]^{\fff n}$\dnsp.\oss
Clearly,\oss if\pss $a_{\dff i}\off >\off a_{\fff j}$\nsp,\oss
then\qss
$a_{\dff i}\qff +\qff x_{\dff i}
\off \geq\off
a_{\fff j}\qff +\qff x_{\fff j}$\nsp.\oss
If\pss
$a_{\dff i}\off =\off a_{\fff j}$\nsp,\oss
then\pss
$a_{\dff i}\qff +\qff x_{\dff i}
\off \geq\off
a_{\fff j}\qff +\qff x_{\fff j}$\qss
is\dss trivially\sss equivalent\dss to\pss
$x_{\dff i}
\off \geq\off
x_{\fff j}$\nsp.\oss
It\dss follows\dss that\qss
$a\qff +\qff \Delta\trf(\dff \omega\trf)
\off \subset\off
l\fff \Delta\dff(\dff \omega\nsp\mid\nsp a\trf)$\nnsp.\oss
Since different\sss simplices\dss $l\fff \Delta\dff(\dff \omega'\dff)$\dss
have disjoint\dss interiors,\oss
$a\qff +\qff \Delta\dff(\dff \omega\trf)$\qss
is\dss not\sss contained\dss in any other\dss
$l\fff \Delta\dff(\dff \omega'\dff)$\nnsp.\oss
Moreover\halfff,\oss
if\pss $\omega'\off \neq\off \omega\nsp\mid\nsp a$\nsp,\oss
then\dss the intersection of\sss
$a\qff +\qff \Delta\dff(\dff \omega\trf)$\sss
and\dss $l\fff \Delta\dff(\dff \omega'\dff)$\dss
is\dss contained\dss in\dss the boundary of\qss
$l\fff \Delta\dff(\dff \omega'\dff)$\nnsp.\oss
Since\dss $Q$\dss is\dss equal\dss to\sss the union of\dss
simplices of\trs the form\qss
$a\qff +\qff \Delta\dff(\dff \omega\trf)$\nnsp,\oss
it\dss follows\dss that\dss $l\fff \Delta\dff(\dff \omega'\dff)$\dss
is\dss equal\dss to\sss the union of\dss simplices\dss
$a\qff +\qff \Delta\dff(\dff \omega\trf)$\dss
contained\dss in\dss $l\fff \Delta\dff(\dff \omega'\dff)$\nnsp.\oss  \eproof

\myapar{Corollary.}{triangulations-all}
\emph{For every\dss permutation\dss $\omega'$\dss the simplices\dss
$a\qff +\qff \Delta\dff(\dff \omega\dff)$\dss
corresponding\dss to pairs\qss $a\fff,\pff \omega$\qss
such\dss that\qss $\omega\nsp\mid\nsp a\off =\off \omega'$\nnsp,\oss
together\dss with\dss their\trs faces,\oss
form\dss a\dss triangulation of\pss $l\fff \Delta\dff(\dff \omega'\dff)$\nnsp.\oss}  \eproof

\myuppar{Remark.}
The simplices\dss $a\qff +\qff \Delta\dff(\dff \omega\dff)$\dss
corresponding\dss to pairs\qss $a\fff,\pff \omega$\qss
such\dss that\qss $\omega\nsp\mid\nsp a\off =\off \operatorname{id}$\qss 
triangulate\dss $l\fff \Delta$\nnsp.\oss
If\qss $\omega\nsp\mid\nsp a\off =\off \operatorname{id}$\nnsp,\oss
then\qss $i\qff >\qff j$\qss implies\sss that\qss
$a_{\dff i}\qff \geq\qff a_{\fff j}$\nsp,\oss
i.e.\qss
$a_{\dff 1}\qff \geq\qff
a_{\dff 2}\qff \geq\qff
\ldots\qff \geq\qff
a_{\dff n}$\qss
and\dss hence\qss $a\qff \in\qff l\fff \Delta$\nnsp.\oss
Of\dss course,\oss geometrically\dss this\dss is\sss obvious.\oss
Indeed,\pss $0\qff \in\qff \Delta\dff(\dff \omega\dff)$\qss
and\dss hence\qss
$a\qff +\qff \Delta\dff(\dff \omega\dff)
\qff \subset\qff
l\fff \Delta$\qss
implies\qss
$a\qff \in\qff l\fff \Delta$\nnsp.\oss 
Since\qss $a_{\dff i}\qff \leq\qff l\qff -\qff 1$\nnsp,\oss
even\dss more\dss is\sss true\fff:\pss
$a\qff +\qff \mathbb{1}\qff \in\qff l\fff \Delta$\nnsp.\oss

\myuppar{The case\qss $l\off =\off 2$\nnsp.}
This\dss is\sss the only\sss case considered\dss by\dss Freudenthal\qss \cite{f}.\oss
If\qss $l\off =\off 2$\nnsp,\oss
then\qss
$1\qff \geq\qff
a_{\dff 1}\qff \geq\qff
a_{\dff 2}\qff \geq\qff
\ldots\qff \geq\qff
a_{\dff n}\qff \geq\qff
0$\qss
and\dss hence\sss there\sss is\sss a\sss number\qss $q\qff \leq\qff n$\qss
such\dss that\qss $a_{\dff i}\off =\off 1$\qss for\qss $i\qff \leq\qff q$\qss
and\qss $a_{\dff j}\off =\off 0$\qss for\qss $i\qff >\qff q$\nnsp.\oss
For such\sss $a$\sss the condition\dss 
$\omega\nsp\mid\nsp a\off =\off \operatorname{id}$\dss
is\sss equivalent\dss to\sss the following\fff:\pss
if\dss either\qss $i\fff,\pff j\qff \leq\qff q$\qss
or\qss $i\fff,\pff j\qff >\qff q$\nnsp,\oss
then\qss $i\qff >\qff j$\qss implies\qss
$\omega\dff(\dff i\dff)\off >\off \omega\dff(\dff j\dff)$\nnsp.\oss
Permutations\dss $\omega$\dss satisfying\dss this condition are known as\dss
\emph{$(\fff p\fff,\qff q\dff)$\dnsp-shuffles},\oss
where\qss $p\off =\off n\qff -\qff q$\nnsp.\oss
Clearly,\oss a\dss $(\fff p\fff,\qff q\dff)$\dnsp-shuffle\dss $\omega$\dss
is\sss determined\dss by\dss the set\qss
$\{\qff\fff \omega\dff(\dff i\dff)\qff \mid\qff i\qff \leq\qff q\pff\}$\nnsp.\oss
Conversely,\oss this set\dss is\dss trivially\sss determined\dss by\dss 
$\omega$\dss and\dss $q$\nnsp,\oss
but\dss not\dss by\dss $\omega$\dss alone,\oss
as\sss the example of\dss the identity\dss permutation shows.\oss
It\dss follows\dss that\dss the simplices of\qss Freudenthal's\dss
triangulation of\dss $2\fff \Delta$\dss correspond\dss to subsets of\qss
$\{\trf 1\fff,\pff 2\fff,\pff \ldots\fff,\pff n\trf\}$\nnsp.\oss
In\dss particular\halfff,\oss the number of\dss $n$\dnsp-simplices
of\trs this\sss triangulation\dss is\dss $2^n$\dnsp.\oss

The vertices of\dss every $n$\dnsp-simplex of\qss Freudenthal's\dss
triangulation of\dss $2\fff \Delta$\dss can\dss be obtained\dss as\dss
follows.\oss
Let\dss the $0${\dnsp}th\dss vertex\dss be\dss
$(\dff 1\fff,\pff \ldots\fff,\pff 1\fff,\pff  0\fff,\pff \ldots\fff,\off 0\dff)$\dss
with $q$ coordinates $1$ followed\dss by\qss $n\qff -\qff q$\qss coordinates $0$\nnsp.\oss
The $1${\dnsp}st\sss vertex results\dss from adding $1$ either\dss to the first\sss $1$\sss
or\dss to\sss the first\sss $0$\nnsp.\oss
Continuing\dss in\dss this way\halfff,\oss one adds $1$ at\dss each step either\dss to
the leftmost\sss $1$ among\dss the first\sss $q$\sss coordinates,\oss
or\dss to\sss the leftmost\sss $0$\sss among\dss the last\qss $n\qff -\qff q$\qss
coordinates.\oss
After $n$ steps one\dss gets\sss all\dss vertices of\sss a simplex
of\qss Freudenthal's\dss triangulation of\dss $2\fff \Delta$\nnsp.\vspace{-0.0625pt}

\myuppar{Freudenthal's\qss triangulations of\trs arbitrary\dss simplices and\dss polyhedra.}
Every $n$\dnsp-simplex\sss $\Gamma$\sss in\dss a euclidean space\sss $\rrr^{\fff m}$\sss
is\dss the image of\dss $l\fff \Delta$\sss under an\qss \emph{affine map}\pss
$\rrr^{\fff n} \ttoo\qff \rrr^{\fff m}$\dnsp,\oss i.e.\qss
a map of\trs the form\qss
$x\qff \longmapsto\qff L\dff (\dff x\dff)\qff +\qff a$\nnsp,\oss
where\qss
$L\dff \colon\dff
\rrr^{\fff n} \ttoo\qff \rrr^{\fff m}$\qss
is\dss a\sss linear\dss map and\qss $a\qff \in\qff \rrr^{\fff m}$\dnsp.\oss
The image of\qss Freudenthal's\dss triangulation\sss of\dss  $l\fff \Delta$\sss
under such map is\sss a\sss triangulation of\dss $\Gamma$\dnsp.\oss
If\trs the vertices of\dss $\Gamma$\dss are ordered,\oss
then\dss there\sss is\sss a preferred\dss map\qss
$l\fff \Delta\dff \ttoo\dff \Gamma$\nnsp,\oss
namely\halfff,\oss
the map\sss taking\dss the vertex\dss $l\fff u_{\dff i}$\dss of\qss $l\fff \Delta$\trs
to\sss the $i${\nnsp}th\dss vertex of\dss $\Gamma$\nnsp.\oss
This\sss leads\sss to a preferred\dss triangulation of\dss $\Gamma$\nnsp.\oss 

Let\dss $S$\dss is a geometric simplicial\sss complex and\dss
let\dss us choose a\sss linear order on\dss the set\sss of\trs its vertices.\oss
This order\dss induces an\sss order on\dss the set\sss of\dss vertices
of\dss each simplex of\dss $S$\nnsp.\oss
The corresponding\dss preferred\dss triangulations of\dss simplices
agree on\sss common\dss faces
and\dss hence define a\sss triangulation of\trs the polyhedron\dss
$\norm{S}$\dss of\dss $S$\nnsp.\oss
For\dss large $l$\dss one\sss gets\sss triangulations
of\trs $\norm{S}$\dss into arbitrarily\sss small\sss simplices.\oss
Instead of\dss using\dss large $l$\qss Freudenthal\dss iterated\dss
the construction\dss with\qss $l\off =\off 2$\nnsp.\vspace{-0.0625pt}

\myuppar{The geometry\sss of\qss Freudenthal's\trs triangulations.}
Every $n$\dnsp-simplex of\trs Freudenthal's\dss triangulation
of\dss $l\fff \Delta$\sss is\sss the image of\dss $\Delta$\dss
under a permutation of\dss coordinates followed\dss by\sss a\sss translation.\oss
In\dss particular\halfff,\oss all $n$\dnsp-simplices of\trs this\sss triangulation
are isometric\sss to $\Delta$ and\dss $\Delta$\dss is\dss homothetic\sss to\dss $l\fff \Delta$\nnsp.\oss
Less formally,\oss one can say\dss that\dss these simplices have\sss the same shape as\dss
$l\fff \Delta$\dss and are $l$ times smaller\dss than\dss $l\fff \Delta$\nnsp.\oss
Taking\dss the images under affine maps preserves\sss these properties,\oss
and\dss hence\qss Freudenthal's\dss construction\dss leads\sss to\sss triangulations
of\dss any\dss polyhedron\dss $\norm{S}$\dss into
into arbitrarily\sss small\sss simplices,\oss
each\dss having\dss the same shape as one of\trs the simplices of\dss $S$\nnsp.\oss

Freudenthal's\dss construction\qss \cite{f}\qss
answered a question of\qss L.E.J.\dss Brouwer\dss
about\sss a\qss \emph{simple}\qss construction of\dss triangulations of\dss
polyhedra into arbitrarily\sss small\sss simplices
which do not\dss eventually\dss include nearly\dss flat\sss simplices,\oss
i.e.\dss $n$\dnsp-simplices with arbitrarily\sss small\dss ratio\dss $V/d^{\dff n}$\dnsp,\oss
where $V$ is\dss the volume and $d$ is\dss the diameter\halfff.\oss
In contrast\dss with\dss Freudenthal's\dss triangulations,\oss
the iterated\dss barycentric subdivisions eventually\sss
include nearly\dss flat\sss $n$\dnsp-simplices.\oss

\newpage
\myappend{The\qss solid\qss cube\qss  and\qss the\qss discrete\qss cube}{two-cubes}

\myuppar{A naive approach\dss to relations between\dss $Q$\dss and\trs $K$\nnsp.}
Let\dss $Q\fff,\pff l$\dss and\trs $K\dff,\pff k$\dss be as in\dss Sections\qss \ref{lebesgue}\qss
and\qss \ref{cubical-cochains}\qss respectively.\oss
Naively,\oss one would\dss think\dss that\sss one should set\qss
$k\off =\off l$\qss and\dss treat\dss $K$\dss as\dss the set\sss of\dss vertices
of\dss cubes of\trs $Q$\nnsp.\oss
An $m$\dnsp-cube of\dss $Q$\sss is determined\dss by\sss its set\sss of\dss vertices,\oss
which\sss is\sss an $m$\dnsp-cube of\trs $K$\nnsp,\oss
and each $m$\dnsp-cube of\dss $K$\dss is\dss the set\sss of\dss vertices of\dss
an $m$\dnsp-cube of\dss $Q$\nnsp.\oss
This gives as a canonical\sss one-to-one correspondence between $m$\dnsp-cubes of\dss $Q$\dss
and of\dss $K$\nnsp.\oss
In\dss turn,\oss this leads\sss to a canonical\sss 
one-to-one correspondence between $m$\dnsp-chains of\dss $Q$\dss
and of\dss $K$\nnsp.\oss
A cubical\sss subset\sss of\dss $Q$\dss can\dss be identified\dss with\dss
the formal\sss sum of\dss $n$\dnsp-cubes contained\dss in\dss it\halfff,\oss
and\dss hence with an $n$\dnsp-chain.\oss
This allows\sss to interpret\sss the cubical\sss subsets of\dss $Q$\dss
in\dss terms of\trs $K$\nnsp.\oss

The proofs of\dss all\dss results of\qss Section\qss \ref{lebesgue}\qss
up\sss to\dss Theorem\qss \ref{collecting-sets}\qss
({\fff}i.e\qss up\sss turning\dss to closed sets)\qss
can\dss be straightforwardly\dss translated\dss into\sss purely\sss
combinatorial\dss language of\trs $K$\dss
with\dss exception\dss of\trs the intersections-based\dss proof\dss 
of\qss Lemma\qss \ref{essential}.\oss
One can save\sss the day\dss with\dss
the projections-based\dss proof\halfff,\pss which
admits a straightforward\dss translation.\oss
In\qss Lusternik--Schnirelmann\qss context\dss
there is no obvious analogue of\dss projections.\oss
In\dss fact\halfff,\pss
the intersections-based\dss proofs also
can\dss be\sss translated\dss into\sss the combinatorial\dss language of\trs $K$\nnsp,\oss
but\sss at\dss the cost\sss of\dss obscuring\dss the ideas.

\myuppar{Duality\halfff.}
The proper\dss way\dss to understand\dss the relations between\dss the proofs
using\dss the solid cube\dss $Q$\dss and\dss the proofs using\dss the
discrete cube\dss $K$\dss is different\halfff.\oss
The discrete counterpart\dss to\dss $Q$\dss is not\dss the set\sss of\dss the vertices
of\dss $Q$\nnsp,\oss 
but\dss the set\sss of\dss the centers of\dss $n$\dnsp-cubes of\dss $Q$\nnsp,\oss
and\dss the discrete cube\dss $K$\dss should\dss be\sss thought\sss of\dss
as\sss the set\sss of\trs these centers,\oss 
for convenience shifted\dss in such a way\dss that\dss
the coordinates of\dss the points of\dss $K$\dss are integers.\oss
Let\dss us set\qss $k\off =\off l\qff -\qff 1$\qss
and consider\dss the\sss translation\qss
$t\dff \colon\dff
\rrr^{\fff n} \ttoo\qff \rrr^{\fff n}$\qss
by\dss the vector\qss
$(\dff 1/2\fff,\pff 1/2\fff,\pff \ldots\fff,\pff 1/2\dff)$\nnsp,\oss
i.e.\qss the map\vspace{2.25pt}
\[
\quad
t\dff \colon\dff
(\dff x_{\dff 1}\dff,\pff x_{\dff 2}\dff,\pff \ldots\dff,\pff x_{\dff n} \dff)
\off \longmapsto\off
(\qff x_{\dff 1}\qff +\qff 1/2\dff,\pff 
x_{\dff 2}\qff +\qff 1/2\dff,\pff 
\ldots\dff,\pff 
x_{\dff n}\qff +\qff 1/2 \qff)
\qff.
\]

\vspace{-9.75pt}
Then $t\trf(\trf K\trf)$ is\dss the set\sss of\dss 
centers of\dss $n$\dnsp-cubes of\dss $Q$\nnsp.\oss
For\dss $v\qff \in\qff K$\dss let\dss $*v$\dss
be\sss the $n$\dnsp-cube of\dss $Q$\dss with\dss the center\dss $t\trf(\dff v\trf)$\nnsp.\oss
In\dss more details,\oss
if\qss
$v
\off =\off
(\dff
a_{\dff 1}\fff,\pff
a_{\dff 2}\trf,\pff
\ldots\fff,\pff
a_{\dff n}
\dff)$\nnsp,\oss
then\vspace{-0.5pt}
\[
\quad
*v
\off =\off
\prod_{\dff i\qff =\qff 1\vphantom{K}}^{\dff n}\qff
[\trf a_{\dff i}\fff,\pff a_{\dff i}\qff +\qff 1\qff ]
\qff.
\]

\vspace{-12.5pt}
The map\dss
$v\qff \longmapsto\qff *v$\dss
is a bijection\dss between\dss $K$\dss and\dss the set\sss of\dss $n$\dnsp-cubes of\dss $Q$\nnsp.\oss
It\dss should\dss be considered as a sort\sss of\dss duality\fff:\oss
it\dss assigns\sss to a $0$\dnsp-cube of\trs $K$\dss an $n$\dnsp-cube of\trs $Q$\nnsp.\oss
Experts will\dss immediately\dss recognize here an\dss instance of\qss Poincar\'{e}\dss duality\halfff.\oss
For\qss $s\qff \subset\qff K$\pss let\vspace{0pt}
\[
\quad
*s
\off =\off
\bigcup_{v\qff \in\qff s\vphantom{K}}\qff *v
\qff.
\]

\vspace{-12pt}
Then\qss
$s\qff \longmapsto\qss *s$\qss
is a bijection\dss between subsets of\trs $K$\dss
and\sss cubical\sss subsets\sss of\dss $Q$\nnsp.\oss
The maps\qss
$\bullet\qff \longmapsto\qff *\bullet$\qss
allow\sss to\sss translate\sss results about\dss $K$\dss
into results about\dss $Q$\dss and\sss vice versa.\oss
Let\dss us\dss begin\dss with simplest\dss examples.\oss

\myapar{Lemma.}{dual-main-faces}
\emph{The sets\qss
$*\mathcal{A}_{\dff i}$\trs and\qss $*\mathcal{B}_{\dff i}$\trs
are\dss the unions of\dss all\sss $n$\dnsp-cubes of\trs $Q$\dss
intersecting\dss $A_{\dff i}$\dss and\dss $B_{\dff i}$\dss respectively\halfff.\qff\oss
A\dss subset\qss  
$s\off \subset\off K$\qss
contains\qss
$\mathcal{A}_{\dff i}$\trs or\pss $\mathcal{B}_{\dff i}$\qss
if\qss and\dss only\qss if\pss the subset\oss
$*s\pff \subset\off Q$\sss contains\qss $A_{\dff i}$\trs or\pss $B_{\dff i}$\dss
respectively.\oss
Similarly,\oss a\dss subset\qss $s\off \subset\off K$\qss intersects\qss
$\mathcal{A}_{\dff i}$\trs or\pss $\mathcal{B}_{\dff i}$\qss
if\qss and\dss only\qss if\pss
the subset\oss
$*s\pff \subset\off Q$\sss intersects\qss $A_{\dff i}$\trs or\pss $B_{\dff i}$\dss
respectively.\oss}

\proof
The first\sss statement\dss is obvious,\pss
and\dss the rest\dss immediately\dss follows\sss from\dss it\halfff.\oss  \eproof

\myapar{Lemma.}{cubes-intersections}
\emph{Let\qss
$s_{\dff 1}\fff,\pff
s_{\dff 2}\trf,\pff
\ldots\fff,\pff
s_{\fff r}$\qss
be subsets of\pss
$K$\nnsp.\oss
There exists an $n$\dnsp-cube\sss of\oss $K$\dss 
intersecting\trs every\dss set\pss 
$s_{\dff 1}\fff,\pff
s_{\dff 2}\trf,\pff
\ldots\fff,\pff
s_{\fff r}$\pss
if\qss and\dss only\qss if\dff\oss 
$*s_{\dff 1}\qff \cap\qff
*s_{\dff 2}\qff \cap\qff
\ldots\qff \cap\qff
*s_{\fff r}
\qff\off \neq\qff\off
\varnothing$\nsp.\oss}

\proof
Let\dss us consider an $n$\dnsp-cube\vspace{3pt}
\[
\quad
\sigma
\off =\off
\prod\nolimits_{\dff i\qff =\qff 1}^{\dff n}\off
\{\trf a_{\dff i}\fff,\pff a_{\dff i}\qff +\qff 1\qff \}
\qff
\]

\vspace{-9pt}
of\trs $K$\nnsp.\qss\off
If\qss $v\qff \in\qff \sigma$\nnsp,\oss
then,\oss obviously,\oss 
$(\qff
a_{\dff 1}\qff +\qff 1\fff,\pff
a_{\dff 2}\qff +\qff 1\trf,\pff
\ldots\fff,\pff
a_{\dff n}\qff +\qff 1
\qff)
\qff \in\qff
v^*$\nnsp.\qff\oss
Therefore\trs\vspace{3pt}  
\[
\quad
*v_{\dff 1}\qff \cap\qff
*v_{\dff 2}\qff \cap\qff
\ldots\qff \cap\qff
*v_{\fff r}
\off \neq\off
\varnothing
\qff
\]

\vspace{-9pt}
if\qss all\qss  
$v_{\dff 1}\fff,\pff
v_{\dff 2}\trf,\pff
\ldots\fff,\pff
v_{\fff r}$\qss
are contained\dss in\sss an $n$\dnsp-cube of\trs $K$\nnsp.\oss
It\dss follows\dss that\dss if\qss
$s_{\dff 1}\fff,\pff
s_{\dff 2}\trf,\pff
\ldots\fff,\pff
s_{\fff r}$\qss
intersect\dss $\sigma$\nnsp,\oss
then\qss
$*s_{\dff 1}\qff \cap\qff
*s_{\dff 2}\qff \cap\qff
\ldots\qff \cap\qff
*s_{\fff r}
\qff\off \neq\qff\off
\varnothing$\nsp.\oss
This proves\sss the\qss ``only\dss if{\trf}''\qss part\sss of\trs lemma.\oss

Suppose\sss now\dss that\pss
$*s_{\dff 1}\qff \cap\qff
*s_{\dff 2}\qff \cap\qff
\ldots\qff \cap\qff
*s_{\fff r}
\off \neq\off
\varnothing$\nsp.\oss
By\dss the definition of\dss the sets\dss $*s_{\dff i}$ \dss
in\dss this case\pss
$*v_{\dff 1}\qff \cap\qff
*v_{\dff 2}\qff \cap\qff
\ldots\qff \cap\qff
*v_{\fff r}
\off \neq\off
\varnothing$\pss 
for\sss some\sss points\qss $v_{\fff i}\qff \in\qff s_{\dff i}$\nsp.\oss
Suppose\sss that\vspace{3pt}
\[
\quad
(\dff
a_{\dff 1}\fff,\pff
a_{\dff 2}\trf,\pff
\ldots\fff,\pff
a_{\dff n}
\dff)
\off \in\off 
*v_{\dff 1}\qff \cap\qff
*v_{\dff 2}\qff \cap\qff
\ldots\qff \cap\qff
*v_{\fff r}
\qff
\]

\vspace{-9pt}
and\sss consider\dss the $n$\dnsp-cube\vspace{3pt}
\[
\quad
\sigma
\off =\off
\prod\nolimits_{\dff i\qff =\qff 1}^{\dff n}\off
\{\trf a_{\dff i}\qff -\qff 1\fff,\pff a_{\dff i}\qff \}
\qff
\]

\vspace{-9pt}
of\trs $K$\nnsp.\oss
Clearly,\oss
$v_{\dff i}\qff \in\qff \sigma$\qss for every\qss
$i\qff \in\qff I$\nnsp.\oss  
This proves\sss the\qss ``{\halfff}if{\trf}''\qss part\sss of\trs lemma.\oss  \eproof

\myapar{Theorem.}{discrete-coverings}
\emph{Let\qss
$s_{\dff 1}\fff,\pff s_{\dff 2}\dff,\pff \ldots\fff,\pff s_{\dff n\dff +\dff 1}
\off \subset\oss K$\nnsp.\oss 
Suppose\sss that\qss
$s_{\dff 1}\qff \cup\qff s_{\dff 2}\qff \cup\qff \ldots\qff \cup\qff s_{\dff n\dff +\dff 1}
\qff =\off
K$\dss and}\vspace{1pt}

\hspace*{1em}({\fff}i{\fff})\phantom{i{\dff}}\oss 
\emph{$s_{\dff 1}\qff \cup\qff s_{\dff 2}\qff \cup\qff \ldots\qff \cup\qff s_{\dff i}$\qss
contains\qss 
$\mathcal{A}_{\dff i}$\qss
and\qss $s_{\dff i}$\dss 
is\dss disjoint\dss from\qss
$\mathcal{B}_{\dff i}$\oss for\qss $i\qff \leq\qff n$\nnsp,\oss and}
\vspace{0pt}

\hspace*{1em}({\fff}i{\fff}i{\fff})\oss
\emph{$s_{\dff i}$\dss is\dss disjoint\qss from\qss
$\mathcal{A}_{\fff j}$\pss if\oss $i\qff >\qff j$\nnsp.\oss}
\vspace{1pt}

\emph{Then\dss there exists an $n$\dnsp-cube\sss of\oss $K$\dss 
intersecting\trs every\dss set\pss
$s_{\dff 1}\fff,\pff s_{\dff 2}\dff,\pff \ldots\fff,\pff s_{\dff n\dff +\dff 1}$\nsp.\oss}

\proof\dnsp
This\sss theorem is\sss a\sss duality-based\dss translation of\qss 
Theorem\qss \ref{e-coverings}\qss 
into\sss the language of\trs $K$\nsp.
Let\qss
$e_{\dff i}
\off =\off
*s_{\dff i}$\pss
for every\qss
$i\off =\off 1\fff,\pff 2\fff,\pff \ldots\fff,\pff n\qff +\qff 1$\nnsp.\oss
Since\qss
$s_{\dff 1}\qff \cup\qff s_{\dff 2}\qff \cup\qff \ldots\qff \cup\qff s_{\dff n\dff +\dff 1}
\qff =\off
K$\nnsp,\oss
the sets\qss
$e_{\dff 1}\fff,\pff e_{\dff 2}\dff,\pff \ldots\fff,\pff e_{\dff n\dff +\dff 1}$\qss
form a covering of\dss $Q$\nnsp.\oss
The assumption\qss ({\fff}i{\fff})\qss 
together\dss with\dss Lemma\qss \ref{dual-main-faces}\qss
imply\dss that\qss
$e_{\dff 1}\qff \cup\qff e_{\dff 2}\qff \cup\qff \ldots\qff \cup\qff e_{\dff i}$\qss
contains\dss $A_{\dff i}$\dss and\dss is\sss disjoint\dss from\dss $B_{\dff i}$\dss
for\qss $i\qff \leq\qff n$\nnsp.\oss
Similarly,\oss the assumption\qss ({\fff}i{\fff}i{\fff})\qss
together\dss with\dss Lemma\qss \ref{dual-main-faces}\qss
imply\dss that\dss $e_{\dff i}$\dss is\sss disjoint\dss from\dss $A_{\dff j}$\dss
for\qss $i\qff >\qff j$\nnsp.\oss
Therefore\sss the sets\qss
$e_{\dff 1}\fff,\pff e_{\dff 2}\dff,\pff \ldots\fff,\pff e_{\dff n\dff +\dff 1}$\qss
satisfy\dss the assumptions of\qss
Theorem\qss \ref{e-coverings}.\oss 
By\dss this\dss theorem\qss
$e_{\dff 1}\qff \cap\qff e_{\dff 2}\qff \cap\qff \ldots\qff \cap\qff e_{\dff n\dff +\dff 1}
\qff \neq\off
\varnothing$\nnsp.\oss
It\dss remains\dss to\sss apply\qss Lemma\qss \ref{cubes-intersections}.\oss  \eproof

\myapar{Theorem.}{separation-products-weak}
\emph{Let\pss
$s_{\dff 1}\fff,\pff s_{\dff 2}\dff,\pff \ldots\fff,\pff s_{\dff n}
\off \subset\off K$\qss and\oss
$\overline{s}_{\dff i}
\qff\off =\qff\off
K\off \smallsetminus\off s_{\dff i}$\nsp.\qff\oss
If\oss 
$\mathcal{A}_{\dff i}
\off \subset\off 
s_{\dff 1}\qff \cup\qff 
s_{\dff 2}\qff \cup\qff
\ldots
\qff \cup\qff 
s_{\dff i}$\qss 
and\qss
$\mathcal{B}_{\dff i}\off \subset\off \overline{s}_{\dff i}$\oss
for\dss every\qss $i\qff \in\qff I$\nnsp,\oss
then\dss there\dss is\dss an\dss $n$\dnsp-cube\dss of\oss $K$\dss  
intersecting\trs
all\trs sets\pss $s_{\dff i}\dff,\off\qff
\overline{s}_{\dff i}$\nsp.\oss}

\proof\dnsp
This\sss theorem\dss is\sss a\sss duality-based\dss translation of\qss Theorem\qss \ref{cubes-partitions}\qss 
into\sss the language of\trs $K$\nsp.
Let\qss
$s_{\dff n\dff +\dff 1}
\off =\dff\off
\overline{s}_{\dff 1}\qff \cap\qff \ldots\qff \cap\pff \overline{s}_{\dff n}$\qss
and\qss
$d_{\dff i}
\off =\off
*s_{\dff i}$\pss
for\dss every\qss
$i\off =\off 1\fff,\pff 2\fff,\pff \ldots\fff,\pff n\qff +\qff 1$\nnsp.\oss
Then\qss the\sss sets\qss
$d_{\dff 1}\fff,\pff d_{\dff 2}\dff,\pff \ldots\fff,\pff d_{\dff n\dff +\dff 1}$\qss
form a covering\dss of\dss $Q$\nnsp.\oss
Lemma\qss \ref{dual-main-faces}\qss implies\sss that\qss
$A_{\dff i}\qff \subset\qff d_{\dff 1}\qff \cup\qff d_{\dff 2}\qff \cup\qff
\ldots\qff \cup\qff d_{\dff i}$\qss
and\dss $d_{\dff i}$\dss is\dss disjoint\dss from\dss $B_{\dff i}$\dss
for every\qss $i\qff \leq\qff n$\nnsp.\oss
By\qss Theorem\qss \ref{cubes-partitions}\qss this implies\sss that\dss the intersection\qss
$d_{\dff 1}\qff \cap\qff d_{\dff 2}\qff \cap\qff \ldots\qff \cap\qff d_{\dff n\dff +\dff 1}$\qss
is\dss non-empty\halfff.\oss
It\dss remains\dss to\sss apply\qss Lemma\qss \ref{cubes-intersections}.\oss  \eproof

\myuppar{A\sss duality-based\sss approach\dss to\dss Lebesgue\dss theorems.}
In\qss Section\qss \ref{cubical-cochains}\qss
we presented a proof\dss of\pss
Lebesgue\dss first\dss covering\dss theorem\qss
(which\dss immediately\dss implies\dss the\dss second one)\qss
based on\dss the multiplication of\dss cochains.\oss
A\sss duality-based\dss version of\dss this approach
allows\dss to move\sss to\sss  
$Q$\dss 
a\sss little bit\sss earlier and\dss prove a combinatorial\sss analogue of\qss
Theorem\qss \ref{separation-products-sets}.\oss 
Here\dss it\dss is.

\myapar{Theorem.}{separation-lebesgue}
\emph{Suppose\dss that\pss
$d_{\dff 1}\fff,\pff d_{\dff 2}\dff,\pff \ldots\fff,\pff d_{\dff n}$\qss
are cubical\sss subsets of\oss $Q$\nsp.\qff\oss 
If}\vspace{3pt}
\[
\quad
A_{\dff i}\off \subset\off d_{\dff i}
\hspace{1em}\mbox{\emph{and}}\hspace{1.2em}
B_{\dff i}\off \subset\off \overline{d}_{\dff i}
\]

\vspace{-9pt}
\emph{\textup{(}where\qss
$\overline{d}_{\dff i}$\dss
is\dss the\sss union of\qss all $n$\dnsp-cubes of\oss $Q$\dss 
not\sss contained\dss in\dss $d_{\dff i}$\hnsp\textup{)}\oss
for\dss every\qss $i\qff \in\qff I$\nnsp,\oss
then}\vspace{4.5pt}
\[
\quad
\bigl(\dff
d_{\dff 1}\qff \cap\qff 
d_{\dff 2}\qff \cap\qff
\ldots
\qff \cap\qff
d_{\dff n}
\qff\bigr)
\off \cap\off
\bigl(\qff
\overline{d}_{\dff 1}\qff \cap\qff 
\overline{d}_{\dff 2}\qff \cap\qff
\ldots
\qff \cap\qff
\overline{d}_{\dff n}
\qff\bigr)
\off\qff \neq\off\qff
\varnothing
\qff.
\]

\vspace{-9.0pt}
\proof
For each\qss $i\qff \in\qff I$\pss
let\qss $c_{\dff i}\qff \subset\pff K$\qss
be\dss such\trs that\qss
$*c_{\dff i}\off =\dff\off d_{\dff i}$\nsp.\qff\oss
Let\oss 
$\overline{c}_{\dff i}\off =\dff\off K\qff \smallsetminus\qff c_{\dff i}$\nsp.\qff\oss
Then\qss\vspace{3pt}
\[
\quad
*\dff\left(\qff \overline{c}_{\dff i}\trf\right)
\off =\qff\off 
\overline{d}_{\dff i} 
\qff. 
\]

\vspace{-9pt}
Lemma\qss \ref{dual-main-faces}\qss implies\sss that\qss 
$\mathcal{A}_{\dff i}\off \subset\off c_{\dff i}$\qss
and\qss
$\mathcal{B}_{\dff i}\off \subset\off \overline{c}_{\dff i}$\qss
for every\qss $i\qff \in\qff I$\nnsp.\oss
Therefore\sss the sets\qss
$c_{\dff i}\dff,\pff \overline{c}_{\dff i}$\qss 
satisfy\dss the assumptions of\qss Theorem\qss \ref{separation-products}.\oss
By\qss this\dss theorem\dss  
there exists an $n$\dnsp-cube\sss
intersecting\dss all\dss these sets.\oss
It\trs remains\dss to apply\qss
Lemma\qss \ref{cubes-intersections}.\oss  \eproof

\myuppar{Deducing\dss Lebesgue\dss theorems\dss from\qss Theorem\qss \ref{separation-lebesgue}.}
First\sss of\dss all,\oss
one can deduce\qss
Theorem\qss \ref{separation-products-sets}\qss 
from\dss its\sss combinatorial\sss analogue\qss 
({\fff}i.e.\qss from\dss the\sss last\dss theorem)\qss 
in\sss exactly\dss the same manner as\dss Lebesgue\dss first\dss covering\dss theorem\dss
was deduced\dss from\dss its\sss combinatorial\sss analogue\qss
({\fff}i.e.\qss from\qss Theorem\qss \ref{collecting-sets}).\oss 
We\dss leave details\sss to\sss interested\dss readers.\oss
After\dss this one can\dss prove\dss Lebesgue\dss fist\sss covering\dss theorem\dss in
exactly\dss the same way\sss as in\dss Section\qss \ref{cubical-cochains}.\oss

\myuppar{Duality\dss in\sss arbitrary\dss dimensions.}
Suppose\sss that\vspace{3pt}
\begin{equation}
\label{cube-for-duality}
\quad
\sigma
\off =\off
\prod_{\dff i\qff =\qff 1}^{\dff n}\off
\rho_{\fff i}
\qff,
\end{equation}

\vspace{-9pt}
is\dss a cube of\trs $K$\nnsp.\oss
Thus\qss $i\qff \in\qff I$\qss
either\qss
$\rho_{\fff i}
\off =\off
\{\trf a_{\dff i}\fff,\pff a_{\dff i}\qff +\qff 1\qff\}$\qss  
for a non-negative integer\qss $a_{\dff i}\qff \leq\qff k\qff -\qff 1$\nnsp,\oss
or\dss 
$\rho_{\fff i}
\off =\off
\{\dff a_{\dff i}\trf\}$\qss
for a non-negative integer\qss $a_{\dff i}\qff \leq\qff k$\nnsp.\qff\oss
Let\vspace{4.5pt}
\[
\quad
*\dff\rho_{\fff i}
\off =\off
\{\dff a_{\dff i}\qff +\qff 1\trf\}
\hspace*{1.5em}\mbox{if}\hspace*{1.7em}
\rho_{\fff i}
\off =\off
\{\trf a_{\dff i}\fff,\pff a_{\dff i}\qff +\qff 1\trf\}
\qff,
\]

\vspace{-31.5pt}
\[
\quad
*\dff\rho_{\fff i}
\off =\off
[\trf a_{\dff i}\fff,\pff a_{\dff i}\qff +\qff 1 \trf]
\hspace*{1.5em}\mbox{if}\hspace*{1.7em}
\rho_{\fff i}
\off =\off
\{\trf a_{\dff i}\trf\}
\qff,
\]

\vspace{-13.5pt}
and\vspace*{-9pt}
\[
\quad
*\dff \sigma
\off =\off
\prod_{\dff i\qff =\qff 1}^{\dff n}\qff
*\dff \rho_{\fff i}
\off.
\]

\vspace{-10.5pt}
Clearly\halfff,\pss $*\dff \sigma$\dss is\dss a\sss cube of\trs $Q$\nnsp,\oss
called\dss the\qss \emph{dual\dss cube\sss of}\dss $\sigma$\nnsp.\oss
If\dss $\sigma$\dss is\dss an $m$\dnsp-cube,\oss
then\dss $*\dff \sigma$\dss is\dss an $(\fff n\dff -\dff m\fff)$\dnsp-cube.\oss

A\sss cube\sss $\sigma$\sss of\trs $K$\dss is uniquely\sss determined\dss by\dss
its\sss dual\sss cube\sss $*\dff \sigma$\nnsp,\oss
and a\sss cube\dss $c$\dss of\trs $Q$ has\sss the form\dss $*\dff \sigma$\qss
if\trs and\dss only\trs if\dss $c$\dss
is\dss not\dss contained\dss in\sss $\bd\fff Q$\nnsp.\oss
We\sss will\sss also say\dss that\dss the cube\sss $\sigma$\dss is\dss the 
\emph{dual\dss cube}\qss of\dss $*\dff \sigma$\nnsp.\oss
If\qss $c\off =\off *\dff \sigma$\nsp\dnsp,\dff\oss we\sss will\dss also write\qss
$\sigma\off =\off *\fff c$\nnsp.\oss
So,\oss 
$\sigma\off =\off **\fff \sigma$\qss
if\dss $\sigma$\dss is\dss a cube of\trs $K$\nnsp,\oss
and\qss
$c\off =\off **\fff c$\qss
if\dss $c$\dss is\dss a\sss cube of\trs $Q$\dss not\dss contained\dss in\sss $\bd\fff Q$\nnsp.\oss

\myapar{Lemma.}{face-duality}
\emph{Suppose\dss that\qss $\sigma\fff,\pff \tau$\dss are\dss two cubes of\pss $K$\nnsp.\oss
Then\sss $\sigma$\dss is\dss a\qss face\dss of\qss $\tau$\dss
if\qss and\dss only\qss if\qss
$*\dff \tau$\dss is\dss a\qss face\dss of\qss $*\dff \sigma$\nnsp.\oss}

\proof
Let\dss $\sigma$\dss be\sss an $m$\dnsp-cube of\trs $K$\dss given\dss by\dss the 
product\qss (\ref{cube-for-duality})\qss
and\dss let\dss $\tau$\sss be\sss an $(\fff m\dff -\dff r\dff)$\dnsp-face\dss of\dss $\sigma$\dss
for some\qss $r\qff \leq\qff m$\nnsp.\oss
Then\sss $\tau$\sss can\dss be obtained\dss by\dss 
replacing\dss $r$\dss two-point\sss sets\qss\vspace{4pt}
\[
\quad
\rho_{\fff i}
\off =\off
\{\trf a_{\dff i}\fff,\pff a_{\dff i}\qff +\qff 1\trf\}
\]

\vspace{-8pt}
in\qss
(\ref{cube-for-duality})\qss
by\dss either\dss $\{\dff a_{\dff i}\trf\}$\nnsp,\oss
or\dss $\{\dff a_{\dff i}\qff +\qff 1\trf\}$\nnsp,\oss
and\dss $*\dff \tau$\dss can\dss be obtained\dss by\dss replacing\dss
$r$\dss one-point\sss sets\vspace{3pt}
\[
\quad
*\dff \rho_{\fff i}
\off =\off
\{\trf a_{\dff i}\qff +\qff 1 \trf\}
\]

\vspace{-9pt}
in\dss the product\sss defining\dss $*\dff \sigma$\dss
by\dss the\sss intervals\vspace{4pt}
\[
\quad
[\trf a_{\dff i}\dff,\pff a_{\dff i}\qff +\qff 1 \qff]
\hspace{1.5em}\mbox{or}\hspace{1.5em}
[\trf a_{\dff i}\qff +\qff 1\fff,\pff a_{\dff i}\qff +\qff 2 \qff]
\]

\vspace{-8pt}
respectively.\oss
Clearly,\pss $*\dff \sigma$\dss is\dss a\dss face\sss of\qss $*\dff \tau$\nnsp.\oss
A similar\sss argument\dss shows\sss that\dss $\tau$\dss is\dss a\sss face\sss of\dss $\sigma$\dss
if\qss $*\dff \sigma$\dss is\dss a\sss face\sss of\dss $*\dff \tau$\nnsp.\oss
The lemma follows.\oss  \eproof

\myuppar{Corollary\halfff.}
\emph{Suppose\sss that\qss
$v_{\dff 1}\fff,\pff
v_{\dff 2}\trf,\pff
\ldots\fff,\pff
v_{\fff r}$\qss
are\qss vertices of\qss an $m$\dnsp-cube\dss $\sigma$\dss of\oss $K$\nnsp.\oss
Then\dss the\dss $n$\dnsp-cubes\qss
$*v_{\dff 1}\fff,\pff
*v_{\dff 2}\trf,\pff
\ldots\fff,\pff
*v_{\fff r}$\qss
of\oss $Q$\dss
have a common $(\fff n\dff -\dff m\fff)$\dnsp-face.\oss}

\proof
Indeed,\oss $*\dff \sigma$\dss is\dss a\sss face of\dss $*v_{\dff i}$\dss
for every\qss $i\qff \leq\qff r$\nnsp.\oss  \eproof

\myuppar{Poincar\'{e}\dss duality.}
Let\dss us\dss
indicate\sss some deeper aspects of\dss
relations\sss between\dss $K$\dss and\dss $Q$\nnsp.\oss

Let\dss us\sss consider\dss the cubes of\trs $K$\dss 
as chains and\dss the cubes of\trs $Q$\dss
as cochains.\oss
A cochain of\trs $Q$\dss is called\dss a\qss \emph{relative\dss cochain}\qss 
of\dss
$(\dff Q\fff,\pff \bd\fff Q\dff)$\dss
if\trs it\dss is\dss equal\dss to a formal\sss sum of\dss cubes of\trs $Q$\dss
not\sss contained\dss in\dss $\bd\fff Q$\nnsp.\oss
As\sss linear\dss functionals\qss
$C_{\dff n\dff -\dff m}\trf(\dff Q\dff)\qff \ttoo\qff \ftwo$\nsp,\oss
relative cochains are functionalss vanishing on\dss the subspace\dss
$C_{\dff n\dff -\dff m}\trf(\trf \bd\fff Q\dff)$\nnsp.\oss
The duality\dss map\qss $\sigma\qff \longmapsto\qff *\dff \sigma$\qss
is\dss a bijection\dss between $m$\dnsp-cubes of\trs $K$\dss
and $(\fff n\dff -\dff m\fff)$\dnsp-cubes of\trs $Q$\dss not\dss in\dss $\bd\dff Q$\dss
and\dss hence\sss
leads\sss to an\sss isomorphism\vspace{2.25pt}
\[
\quad
C_{\dff m}\dff(\dff K\dff)
\off \ttoo\off
C^{\dff n\dff -\dff m}\dff(\dff Q\fff,\pff \bd\fff Q\dff)
\qff.
\]

\vspace{-9.75pt}
If\dss a cube of\dss $Q$\dss is\sss contained\dss in\dss $\bd\fff Q$\nnsp,\oss
then all\dss its faces are contained\dss in\dss $\bd\fff Q$\nnsp.\oss
It\dss follows\sss that\dss the coboundary operators\dss $\partial^*$\dss
define  maps\vspace{2.25pt}
\[
\quad
\partial^*\dff \colon\dff
C^{\dff k}\dff(\dff Q\fff,\pff \bd\fff Q\dff)
\off \ttoo\off
C^{\dff k\dff +\dff 1}\dff(\dff Q\fff,\pff \bd\fff Q\dff)
\qff,
\]

\vspace{-9.75pt}
also called\qss \emph{coboundary operators}.\oss
As it\dss turns out\halfff,\oss the duality\dss map\sss
turns\dss the coboundary operators\dss $\partial^*$\dss
into\sss the boundary operators\dss $\partial$\nnsp.\oss
In\dss more details,\oss
if\dss $\gamma$\dss is\dss a chain\dss of\trs $K$\nnsp,\oss 
then\vspace{2.25pt}
\begin{equation}
\label{dual-boundaries}
\quad
*\dff\left(\dff \partial\dff(\dff \gamma\dff)\dff\right)
\off =\qff\off
\partial^*\fff(\trf *\dff \gamma\trf)
\qff.
\end{equation}

\vspace{-9.75pt}
This immediately\dss 
follows\dss from\dss Lemma\qss \ref{face-duality}.\oss

One can also proceed\dss in other direction and start\dss with\dss
the duality\dss map\qss
$c\qff \longmapsto\qff *\fff c$\qss
from cubes of\trs $Q$\dss to cubes of\trs $K$\nnsp.\oss
The only difficulty\sss is\dss the fact\dss that\dss $*\fff c$
is not\sss defined\dss if\dss $c$ is\sss a\sss cube of\trs $Q$\sss
contained\dss in\dss $\bd\fff Q$\nnsp.\oss
The solution is\dss to set\dss in\dss this case\qss
$*\fff c\off =\off 0$\nnsp.\oss
Let\dss us\sss introduce\sss the spaces\vspace{2.25pt}
\[
\quad
C_{\dff m}\dff(\dff Q\fff,\pff \bd\fff Q\dff)
\off =\off
C_{\dff m}\dff(\dff Q\dff)\left/\dff C_{\dff m}\dff(\trf \bd\fff Q\dff)\right.
\]

\vspace{-9.75pt}
of\pss \emph{relative chains}\qss of\trs
$(\dff Q\fff,\pff \bd\fff Q\dff)$\nnsp.\oss
The duality\dss map\qss
$c\qff \longmapsto\qff *\fff c$\qss
leads\sss to an\sss isomorphism\vspace{2.25pt}
\[
\quad
C_{\dff m}\dff(\dff Q\fff,\pff \bd\fff Q\dff)
\off \ttoo\off
C^{\dff n\dff -\dff m}\dff(\dff K\dff)
\qff.
\]

\vspace{-9.75pt}
The boundary\sss operators define maps\vspace{2.25pt}
\[
\quad
\partial\dff \colon\dff
C_{\dff m}\dff(\dff Q\fff,\pff \bd\fff Q\dff)
\off \ttoo\off
C_{\dff m\dff -\dff 1}\dff(\dff Q\fff,\pff \bd\fff Q\dff)
\qff,
\]

\vspace{-9.75pt}
also called\qss \emph{boundary operators}.\oss
Lemma\qss \ref{face-duality}\qss implies\sss that\dss the identity\qss
(\ref{dual-boundaries})\qss holds\sss also when\dss $\gamma$\dss is\dss
a\dss relative chain\dss of\dss $(\dff Q\fff,\pff \bd\fff Q\dff)$\nnsp.\oss

\vspace*{12pt}
\begin{flushright}
July\qss 17,\oss 2019
 
https\halfff:/\!/\hspace*{-0.06em}nikolaivivanov.com

E-mail\halfff:\oss nikolai.v.ivanov{\fff}@{\dff}icloud.com
\end{flushright}

\end{document}